\documentclass[a4paper, 12pt]{scrreport}  

\addtokomafont{chapter}{\fontsize{20.74}{24}\selectfont}

\usepackage[utf8]{inputenc}  
\usepackage[english]{babel}  

\usepackage{amsthm}  
\usepackage{amsmath}
\usepackage{amssymb}  
\usepackage{mathtools}  
\usepackage{mathrsfs}  

\usepackage{tikz-cd}	
\tikzcdset{scale cd/.style={
	every label/.append style={scale=#1},
	cells={nodes={scale=#1}}
}}

\usepackage{mleftright}  

\usepackage{enumerate}  

\usepackage[most, skins]{tcolorbox}  

\usepackage{hyperref}
\hypersetup{
	colorlinks = true,
	citecolor  = blue,
	linkcolor  = blue, 
	filecolor  = blue,      
	urlcolor   = blue,
}

\newcommand{\ie}{i.\,e.\ }
\newcommand{\Ie}{I.\,e.\ }
\newcommand{\eg}{e.\,g.\ }
\newcommand{\Eg}{E.\,g.\ }
\newcommand{\confer}{cf.\ }
\renewcommand{\quote}[1]{``#1''}
\newcommand{\emphquote}[1]{\emph{\quote{#1}}}

\newcommand{\assumptionPreWord}{\textcolor{red}{assumption(!)}\,}

\newcommand{\ldef}{\coloneqq}
\newcommand{\rdef}{\eqqcolon}
\newcommand{\setseparator}{\,\vert\,}
\newcommand{\powerset}[1]{{\mathcal P(#1)}}

\newcommand{\N}{{\mathbb{N}}}
\newcommand{\A}{{\mathbb{A}}}
\newcommand{\Z}{{\mathbb{Z}}}
\newcommand{\Q}{{\mathbb{Q}}}
\newcommand{\R}{{\mathbb{R}}}
\newcommand{\C}{{\mathbb{C}}}
\newcommand{\adicsQ}[1]{{\Q_{#1}}}

\newcommand{\yoneda}{{\mathcal Y}}
\newcommand{\iso}{\cong}
\renewcommand{\equiv}{\simeq}
\newcommand{\grothendieckaxiom}[1]{(AB#1)}

\DeclareMathAlphabet\mathbfcal{OMS}{cmsy}{b}{n}

\newcommand{\category}[1]{{\mathbfcal{#1}}}
\newcommand{\categoryname}[1]{{\mathbf{#1}}}

\newcommand{\Hom}{{\operatorname{Hom}}}
\newcommand{\eHom}{{\underline{\operatorname{Hom}}}}
\newcommand{\intHom}{{\operatorname{\mathscr{H}\!om}}}

\newcommand{\eStructHom}[1]{\eHom_{\struct{#1}}}
\newcommand{\eRHom}{\RightD\eHom}
\newcommand{\intRHom}{\RightD\intHom}
\newcommand{\eStructRHom}[1]{\eRHom_{\struct{#1}}}

\newcommand{\id}{\mathrm{id}}
\newcommand{\shortisorightarrow}{\stackrel{\sim\;\;}{\smash{\to}\rule{0pt}{0.2ex}}}
\newcommand{\isorightarrow}{\stackrel{\sim\;\;\;\;\;}{\smash{\longrightarrow}\rule{0pt}{0.2ex}}}

\newcommand{\isoleftarrow}{\stackrel{\;\;\;\;\,\sim\,}{\smash{\longleftarrow}\rule{0pt}{0.2ex}}}

\newcommand{\op}{{\mathrm{op}}}
\newcommand{\Overcategory}[2]{{#1\categoryname{/}#2}}
\newcommand{\functorCategory}[2]{{\categoryname{Func}\mleft(#1, #2\mright)}}

\newcommand{\preSheavesSymbol}{{\categoryname{PSh}}}
\newcommand{\sheavesSymbol}{{\categoryname{Sh}}}
\newcommand{\Presheaves}[2]{{\preSheavesSymbol{\mleft(#1, #2\mright)}}}
\newcommand{\Sheaves}[2]{{\sheavesSymbol{\mleft(#1, #2\mright)}}}

\newcommand{\smallPresheaves}[2]{{\categoryname{P\underline{s}h}{\mleft(#1, #2\mright)}}}
\newcommand{\smallSheaves}[2]{{\categoryname{\underline{s}h}{\mleft(#1, #2\mright)}}}

\newcommand{\ProCategory}[1]{{\categoryname{Pro}{\mleft(#1\mright)}}}

\newcommand{\ladj}{\dashv}

\newcommand{\longxrightarrow}[1]{\xrightarrow{\;#1\;}}
\newcommand{\xleftrightarrows}[2]{\xleftrightharpoons[\;#2\;]{\;#1\;}}

\newcommand{\Set}{{\categoryname{\underline{S}et}}}
\newcommand{\boundedSet}[1]{{\Set_{#1}}}
\newcommand{\set}{{\categoryname{\underline{s}et}}}
\newcommand{\Top}{{\categoryname{Top}}}
\newcommand{\Schemes}{{\categoryname{Schemes}}}
\newcommand{\Ab}{{\categoryname{Ab}}}
\newcommand{\TopAb}{{\categoryname{TopAb}}}

\newcommand{\commutativeRings}{\categoryname{CR1ng}}

\newcommand{\limit}{\varprojlim}
\newcommand{\colimit}{\varinjlim}
\newcommand{\homotopyLimit}{\mathrm{h}\!\limit\nolimits}
\newcommand{\homotopyColimit}{\mathrm{h}\!\colimit\nolimits}

\DeclareMathOperator{\equalizer}{Eq}
\DeclareMathOperator{\kernel}{Ker}
\DeclareMathOperator{\coKernel}{CoKer}
\DeclareMathOperator{\image}{Img}

\newcommand{\site}[1]{{\categoryname{#1}}}

\newcommand{\covers}[1]{{\operatorname{Cov}(#1)}}

\newcommand{\open}[1]{{\categoryname{Ouv}(#1)}}
\newcommand{\etalesite}[1]{{{#1}_{\text{ét}}}}
\newcommand{\proetalesite}[1]{{{#1}_{\text{pro-ét}}}}

\newcommand{\CHaus}{{\categoryname{CHaus}}}
\newcommand{\boundedCHaus}[1]{{\CHaus_{#1}}}

\newcommand{\proFiniteSets}{\categoryname{ProFin}}
\newcommand{\boundedProFiniteSets}[1]{\proFiniteSets_{#1}}

\newcommand{\extremallyDisconnectedSets}{\categoryname{ExDisc}}
\newcommand{\boundedExtremallyDisconnectedSets}[1]{\extremallyDisconnectedSets_{#1}}

\newcommand{\condensedSymbol}{{\categoryname{Cond}}}

\newcommand{\boundedCondensedSets}[1]{{\condensedSymbol\Set_{#1}}}
\newcommand{\associated}[1]{\underline{#1}}
\newcommand{\condensed}[1]{{\condensedSymbol(#1)}}
\newcommand{\condensedSets}{{\condensedSymbol\Set}}
\newcommand{\condensedAb}{{\condensedSymbol\Ab}}
\newcommand{\topologizedUnderlying}[1]{{#1(\ast)_\mathrm{top}}}
\newcommand{\resTo}[1]{{\vert_{#1}}}
\newcommand{\enlargementsymbol}[2]{{#1\rightsquigarrow #2}}
\newcommand{\CGWH}{C.G.W.H. }

\newcommand{\sheaf}[1]{{\mathcal{#1}}}
\newcommand{\sheafificationsymbol}{\sharp}
\newcommand{\sheafification}[1]{{\mleft(#1\mright)^\sheafificationsymbol}}
\renewcommand{\restriction}{\vert}

\newcommand{\Ordinals}{{\categoryname{Ord}}}
\newcommand{\Cardinals}{{\categoryname{Card}}}
\newcommand{\cofinality}[1]{{\operatorname{cof}(#1)}}

\newcommand{\topology}[1]{{\mathscr{#1}}}
\newcommand{\boundedkSpaces}[1]{{\mathrm{k}\Top^{#1}}}
\newcommand{\kSpaces}{{\mathrm{k}\Top}}
\newcommand{\kificationsymbol}{{\mathrm{cg}}}
\newcommand{\boundedkificationsymbol}[1]{{#1\mathrm{-cg}}}
\newcommand{\stonecechsymbol}{{\beta}}
\newcommand{\stonecech}[1]{{\stonecechsymbol #1}}
\newcommand{\filterkernel}{{\operatorname{ker}}}
\newcommand{\tOne}{{T_1}}
\newcommand{\tOneSpaces}{{\Top_1}}

\newcommand{\structsymbol}{{\mathcal{O}}}
\newcommand{\struct}[1]{{\structsymbol_{#1}}}
\newcommand{\stalkStruct}[2]{{\structsymbol_{#1, #2}}}

\newcommand{\structTensor}[1]{\otimes_{\struct{#1}}}
\newcommand{\associatedSheaf}[1]{\widetilde{#1}}

\newcommand{\constantSheaf}[1]{{\underline{#1}}}
\newcommand{\free}[2]{{#1[#2]}}
\newcommand{\freeAbelian}[1]{{\free {\associated{\Z}} {#1}}}

\newcommand{\scalarRestriction}[1]{{#1^\ast}}
\newcommand{\scalarExtension}[1]{{#1_{(!)}}}

\newcommand{\pMod}[1]{{\categoryname{pMod}(#1)}}
\renewcommand{\mod}[1]{{\categoryname{Mod}(#1)}}
\newcommand{\maybePMod}[1]{{\categoryname{(p)Mod}(#1)}}

\newcommand{\structMod}[1]{{\mod{\struct{#1}}}}
\newcommand{\qCoh}[1]{{\categoryname{qCoh}(\struct{#1})}}
\newcommand{\Coh}[1]{{\categoryname{Coh}(\struct{#1})}}

\newcommand{\ideal}[1]{{\langle #1\rangle}}
\newcommand{\Spec}{\operatorname{Spec}}
\newcommand{\Spa}{\operatorname{Spa}}
\newcommand{\ad}{{\mathrm{ad}}}

\newcommand{\preAnaRing}[1]{{\mathfrak{#1}}}
\newcommand{\underlying}[1]{{\underline{\preAnaRing{#1}}}}
\newcommand{\ufree}[2]{{\free {\underlying{#1}} {#2}}}
\newcommand{\pfree}[2]{{\free {\preAnaRing{#1}} {#2}}}

\newcommand{\preAnalyticRings}{\categoryname{PAnaRing}}
\newcommand{\analyticRings}{\categoryname{AnaRing}}
\newcommand{\presentableMod}[1]{{\categoryname{Pres}(#1)}}
\newcommand{\pPresentableMod}[1]{{\categoryname{Pres}(\preAnaRing{#1})}}
\newcommand{\completeMod}[1]{{\categoryname{Comp}(#1)}}
\newcommand{\pCompleteMod}[1]{{\categoryname{Comp}(\preAnaRing{#1})}}

\newcommand{\completionSymbol}[1]{#1}
\newcommand{\derivedCompletionSymbol}[1]{{\LeftD#1}}
\newcommand{\completion}[2]{\mleft(#1\mright)^{\completionSymbol{#2}}}
\newcommand{\completionAsTensor}[3]{#1 \otimes_{#2} #3}
\newcommand{\derivedCompletion}[2]{\mleft(#1\mright)^{\derivedCompletionSymbol{#2}}}
\newcommand{\derivedCompletionAsTensor}[3]{#1 \Dotimes_{#2} #3}
\newcommand{\completedScalarExtension}[2]{{#1_{(!)}^{#2}}}
\newcommand{\derivedCompletedScalarExtension}[2]{{\LeftD#1_{(!)}^{#2}}}

\newcommand{\continuous}[2]{{C(#1, #2)}}
\newcommand{\continuousZero}[2]{{C_0(#1, #2)}}
\newcommand{\measures}[2]{{\mathcal{M}(#1, #2)}}
\renewcommand{\d}{\,\mathrm{d}}
\newcommand{\laurent}[2]{{#1(\!(#2)\!)}}
\newcommand{\series}[2]{{#1[\![#2]\!]}}

\newcommand{\derived}[1]{{\categoryname{D}(#1)}}
\newcommand{\derivedPlus}[1]{{\categoryname{D}^+(#1)}}
\newcommand{\derivedMinus}[1]{{\categoryname{D}^-(#1)}}

\newcommand{\RightD}{\mathcal{R}}
\newcommand{\LeftD}{\mathscr{L}}

\newcommand{\Ext}{{\operatorname{Ext}}}

\newcommand{\Dotimes}{\otimes^\LeftD }

\newcommand{\nerveSymbol}{{\mathcal N}}
\newcommand{\nerve}[1]{{\nerveSymbol(#1)}}

\newcommand{\complexes}[1]{{\categoryname{Ch}(#1)}}

\newcommand{\sTruncation}[1]{{\sigma_{#1}}}
\newcommand{\cTruncation}[1]{{\tau_{#1}}}
\newcommand{\homotopyCategoryOfComplexes}[1]{{K(#1)}}

\newcommand{\inftyD}{{\mathcal D}}

\newcommand{\dual}{\vee}
\newcommand{\exterior}{\bigwedge\nolimits}
\newcommand{\cotangent}{\Omega}
\newcommand{\canonical}{\omega}

\newtheoremstyle
{normal}        
{-2pt}          
{5pt}           
{\normalfont}   
{}              
{}     
{}              
{\newline}      
{\hspace{-1.95pt}\underline{\hspace{1.95pt}\bfseries\thmname{#1}\bfseries\thmnumber{ #2}\normalfont\thmnote{ (#3)}}\vspace{0.2em}}              

\newtheoremstyle
{proof}         
{3pt}           
{5pt}           
{\normalfont}   
{}              
{\itshape}      
{:}             
{\newline}      
{}              

\newtheoremstyle
{shortproof}    
{3pt}           
{5pt}           
{\normalfont}   
{}              
{\itshape}      
{:}             
{ }      		
{}              

\makeatletter
\newlength{\savedparskip}
\newlength{\savedparindent}
\makeatother

\newenvironment{pbox}{
	
	\setlength{\savedparskip}{\parskip}
	\setlength{\savedparindent}{\parindent}
	
	\begin{tcolorbox}[
		frame hidden,
		borderline west = {0.4pt}{0pt},
		boxrule=0pt,
		colback=white,
		sharp corners,
		overlay={
			
		},
		left=-0.1ex,
		right=-0.25ex,
		top=-0.0ex,
		bottom=-0.0ex,
		enhanced,
		breakable,
		before upper={
			\setlength{\parskip}{\savedparskip}
			\setlength{\parindent}{\savedparindent}
		}	
	]
}{
	\end{tcolorbox}
}

\newenvironment{assumbox}{
	
	\setlength{\savedparskip}{\parskip}
	\setlength{\savedparindent}{\parindent}
	
	\begin{tcolorbox}[
		colframe=black!75!white,
		colback=red!10!white,
		before skip=2pt,
		after skip=4pt,
		toprule=0.4mm,
		bottomrule=0.4mm,
		leftrule=1mm,
		rightrule=0.5mm,
		opacityback=1.0,
		sharp corners,
		enhanced,
		breakable,
		before upper={
			\setlength{\parskip}{\savedparskip}
			\setlength{\parindent}{\savedparindent}
		}
	]
}{
	\end{tcolorbox}
}

\theoremstyle{normal}
\newtheorem{customtheorem}{Theorem}[section]
\newenvironment{theorem}{\pbox\customtheorem}{\endcustomtheorem\endpbox}

\theoremstyle{normal}
\newtheorem{customcorollary}[customtheorem]{Corollary}
\newenvironment{corollary}{\pbox\customcorollary}{\endcustomcorollary\endpbox}

\theoremstyle{normal}
\newtheorem{customlemma}[customtheorem]{Lemma}
\newenvironment{lemma}{\pbox\customlemma}{\endcustomlemma\endpbox}

\theoremstyle{normal}
\newtheorem{customproposition}[customtheorem]{Proposition}
\newenvironment{proposition}{\pbox\customproposition}{\endcustomproposition\endpbox}

\theoremstyle{normal}
\newtheorem{customdefinition}[customtheorem]{Definition}
\newenvironment{definition}{\pbox\customdefinition}{\endcustomdefinition\endpbox}

\theoremstyle{proof}
\newtheorem*{customproof}{Proof}
\renewenvironment{proof}{\customproof}{\qed \endcustomproof} 

\theoremstyle{shortproof}
\newtheorem*{customshortproof}{Proof}
\newenvironment{shortproof}{\customshortproof}{\qed \endcustomshortproof} 

\theoremstyle{proof}
\newtheorem*{customproofsketch}{Proof sketch}
\newenvironment{proofsketch}{\customproofsketch}{\hfill $\blacksquare$ \endcustomproofsketch}  

\theoremstyle{normal}
\newtheorem{customnotation}[customtheorem]{Notation}
\newenvironment{notation}{\pbox\customnotation}{\endcustomnotation\endpbox}

\theoremstyle{normal}
\newtheorem{customconvention}[customtheorem]{Convention}
\newenvironment{convention}{\pbox\customconvention}{\endcustomconvention\endpbox}

\theoremstyle{normal}
\newtheorem{customexample}[customtheorem]{Example}
\newenvironment{example}{\pbox\customexample}{\endcustomexample\endpbox}

\theoremstyle{normal}
\newtheorem{customremark}[customtheorem]{Remark}
\newenvironment{remark}{\pbox\customremark}{\endcustomremark\endpbox}

\theoremstyle{normal}
\newtheorem{customconstruction}[customtheorem]{Construction}
\newenvironment{construction}{\pbox\customconstruction}{\endcustomconstruction\endpbox}

\theoremstyle{normal}
\newtheorem{customconstructionsketch}[customtheorem]{Sketch of construction}
\newenvironment{construction-sketch}{
	\pbox\customconstructionsketch
}{
	\endcustomconstructionsketch\endpbox
}

\theoremstyle{normal}
\newtheorem{customassumption}[customtheorem]{Assumption}
\newenvironment{assumption}{\assumbox\customassumption}{\endcustomassumption\endassumbox}

\title{
	A Proof of Coherent Duality\\
	via Techniques\\
	from Condensed Mathematics
}
\date{October, 2025}
\author{Cedric Brendel}

\begin{document}
	
	{\let\newpage\relax\maketitle}  
	
	\thispagestyle{empty} 

	\begin{center}
		
		\vfill

		A thesis presented for the degree of\\
		M. Sc. Mathematics\\
		under supervision of Prof. Dr. Ulrich Thiel
		
		\vspace{0.8cm}

		\includegraphics[width=0.4\textwidth]{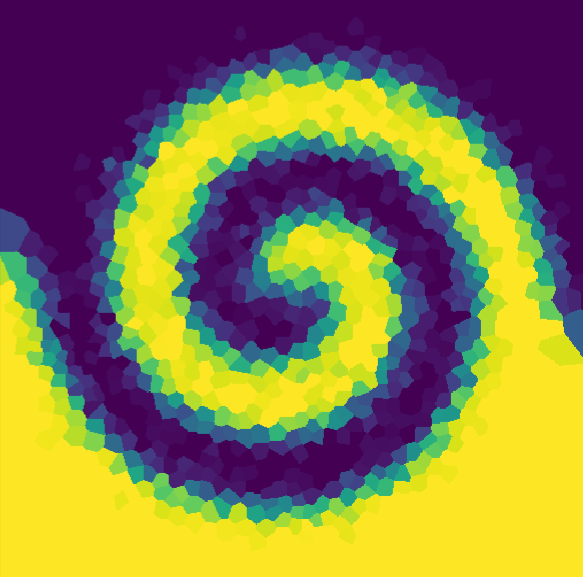}
		
		\vspace{0.8cm}
		
		Department of Mathematics\\
		Rheinland-Pfälzische Technische Universität Kaiserslautern-Landau\\ \vspace{1cm}
		
		
	\end{center}
	
	\cleardoublepage
	
	\newpage
	 
\section*{Introduction: Why work with condensed mathematics?}

The content of this thesis will be two-fold.

Initially, we will focus purely on the theory of condensed mathematics. In chapter \ref{chapter: sites, sheaves, topoi and condensed set} we will introduce the necessary background needed to define condensed sets and to state the formal properties their category possesses. In particular, we will see how condensed sets generalize (good) topological spaces. In chapter \ref{chapter: rings, modules and condensed algebra} the interplay between condensed sets and algebra is established -- one of the motivating drives that led Clausen and Scholze to introduce condensed mathematics, the main reference for which seems to be Clausen and Scholze's \cite[lecture notes]{condenseddotpdf}. In chapter \ref{chapter: analytic rings and completeness} analytic rings and their associated notion of complete modules are introduced. They will play a major part in proving coherent duality in the second part of the thesis, but have vast applications outside of \quote{traditional algebraic geometry}. For example, Clausen and Scholze introduce in their \cite[Lectures on Analytic Geometry]{analyticdotpdf} a special (class of) analytic rings, whose associated complete modules provide a well-rounded notion of topological vector space in functional analysis (\eg that their category is abelian, that certain tensor products arise naturally from the theory, ...). To cite them in the introduction of said lecture notes:

\emphquote{The goal of this course is to [turn] functional analysis into a branch of commutative algebra, and various types of analytic geometry (like manifolds) into algebraic geometry.}

We will however not cover this quite ambitious aspect of condensed mathematics and instead restrict ourselves to that part of the theory of analytic rings which is needed to prove the duality theorem later on. Overall, the first part of thesis aims to establish a well-rounded picture of the theory of condensed mathematics and tries to show what can be gained \emph{theoretically} from working with condensed mathematics.

The basic question motivating condensed mathematics is:
\begin{center}
	\emphquote{How does one do algebra when the algebraic objects carry a topology?}
\end{center}
This question often arises naturally. We list a few scenarios given by Clausen and Scholze, where the underlying topology of objects is an integral part of the theory.
\begin{enumerate}
	\item
	The representation theory of topological groups, \eg matrix groups over $\R$ or $\adicsQ p$, where we might be interested in representations given by a continuous action on a Banach space or more generally a topological vector space.
	
	\item 
	Algebraic/analytic geometry over topological fields such as again $\R$ or $\adicsQ p$.
	
	\item
	Constructing cohomological theories, explicitly incorporating and respecting the underlying topology, \eg Tate's continuous group cohomology:	For a topological group $G$ acting continuously on a topological abelian group $M$, one can construct a complex $C^\bullet(G, M)$ given by $\Hom_\Top(G^n, M)$ in degree $n\ge 0$, yielding cohomology groups $H^n(G, M) \ldef H^n(C^\bullet(G, M))$.
\end{enumerate}

While in many cases solutions have been found to work with this extra structure, in general some problems arise that make life hard. Here are a few examples and why we might care about them:
\begin{enumerate}
	\item
	A frequent problem is that categories of such objects are ill-behaved. A common (and important) example is \eg the category $\TopAb$ of topological abelian groups. While abelian groups on their own of course form an abelian category, this is no longer the case when the groups carry a topology. Consider for example the natural morphism of topological abelian groups
	$$(\R, \text{discrete topology}) \to  (\R, \text{natural topology}).$$
	This is certainly not an isomorphism, yet its kernel and cokernel in $\TopAb$ are trivial -- this is impossible in an abelian category. Hence if we wish to have an ambient abelian category, we would need to alter kernel and cokernel in such a way to explain the failure of being an isomorphism.
	
	\item
	Base changes and tensor products do not generalize well to the topological setting. If $A\to B$ is a morphism of topological $R$-algebras it is in general not clear how to construct a good base change $M\otimes_A^\text{top} B$ of a topological $A$-module $M$. The implications are manifold, \eg in the context of geometry, a theory of topological (quasi-)coherent sheaves lacks a base change. Similarly, in representation theory if one would like to consider some topological analogue of a $p$-modular system $(F, \mathcal O, k)$, the associated functor of $p$-modular reduction (sending topological modules in characteristic $0$ to topological modules in characteristic $p$) should be a base change along $\mathcal O \twoheadrightarrow \mathcal O/\mathfrak m = k$.
	
	\item 
	Very natural definitions of (co)homology fail to produce satisfactory theories, \ie topological structure does not mix well with the standard theory of derived categories. For example, in Tate's continuous group cohomology a short exact sequence of topological modules fails to produce a long exact sequence in cohomology. Hence, the theory does not produce a $\delta$-functor in the sense of Grothendieck and many important homological tools expected from such a theory fail (for example dimensional shifting arguments).
\end{enumerate}

How do we deal with these problems? Consider again the example of the category $\TopAb$ of topological abelian groups which we saw to not be abelian. This category essentially consists of \emph{abelian group objects} inside the category $\Top$ of topological spaces. This is where the root cause of the above problems lies: Algebraic objects internal to an arbitrary category are not as well-behaved as those internal to the category $\Set$ of sets. It is however well known that for a \emph{(Grothendieck) topos} -- a category of sheaves on a small site -- the situation is much better. For example abelian group objects in any topos form an abelian category! This is where the fundamental insight of condensed mathematics lies: One should replace $\Top$ by a better behaved, \emphquote{topos-like} category -- one should replace topological spaces by condensed sets. Only then can one expect that arbitrary internal algebraic objects behave well.

The question remains how to construct such a category. Roughly speaking, a \emph{condensed set} is as (special) $\Set$-valued sheaf on the \emph{pro-étale site} $\proetalesite{\ast}$ of a point. This site is a variant of the well known étale site on a scheme, where loosely speaking the étale condition on affine morphisms is relaxed to a pro-étale condition, that is we instead allow morphisms that are limits of a cofiltered system of pro-étale morphisms (see construction \ref{Construction: Pro-categories} for a formalization of pro-categories).

The pro-étale site was originally introduced by Bhatt and Scholze in \cite{bhatt2014proetaletopologyschemes} to give a new definition of the derived category of constructible $\ell$-adic sheaves.

In condensed mathematics however, only the pro-étale site of a point is of much interest, so we can restrict to that case. It turns out, this site can equivalently be given by the site of \emph{profinite spaces} with finite, jointly surjective families of continuous maps as covers. Calling this site $\proFiniteSets$ we obtain a Yoneda-like functor
\begin{align*}
	\Top &\longrightarrow \Sheaves{\proFiniteSets}{\Set},\\
	X&\longmapsto \Hom_\Top(-, X)
\end{align*}
into the category of sheaves on this site $\proFiniteSets$. Any condensed set will be such a sheaf! For technical reason the converse does however not hold. This technicality will be a large part of the content of the first chapter \ref{chapter: sites, sheaves, topoi and condensed set}. Let us ignore this complication for now. Since we are considering sheaves only on the site of profinite spaces (vs. say all topological spaces), we cannot expect this functor to be fully faithful (contrary to the standard Yoneda embedding). This is a deliberate choice and reflects the usual practice of restricting attention to a \emphquote{convenient category of spaces} in topology. Indeed, with the above functor we can only hope to distinguish spaces whose topology is generated by profinite spaces. As it turns out, any compact Hausdorff space is a quotient of a profinite space, see corollary \ref{corollary: compact hausdorff spaces are quotients of extremally disconnected spaces}. Thus, these \quote{distinguishable} spaces precisely form the well known class of \emph{compactly generated spaces}, \ie those spaces whose topology can be characterized in terms of continuous images of compact Hausdorff spaces. We will see in theorem \ref{theorem: some correspondences}, that indeed the above functor restricts (up to a minor complication) to a faithful functor
	$$\{X\in \Top \colon X\text{ is compactly generated}\} \to \condensedSets$$
that is fully faithful on spaces that are additionally weak Hausdorff. Thus, condensed set give the desired \emphquote{topos-like} category of \emphquote{space-like} objects.

Then, in the second part of the thesis we will see an application of the theory. We will be able to reprove coherent duality for a large class of morphism of affine schemes -- from the perspective of condensed mathematics, see theorem \ref{theorem: solid coherent duality}. Classically, in the projective case one chooses an embedding into projective space and proves duality there by hand, then pulls back the result to any projective scheme -- this is for example the approach taken in Hartshorne's \emph{Algebraic Geometry} \cite[Chap. III. §7.1 and §7.6]{algebraic-geometry} when proving Serre duality (which is a special case of coherent duality). This, to some degree, obfuscates \eg the naturality of the associated trace map. To remedy this, one (and of course originally Grothendieck) adopts a top-down approach that aims to construct for any morphism of schemes $f\colon X\to Y$ a pair of adjoint functors $f_!$ and $f^!$, called the \emph{exceptional direct image} and \emph{exceptional inverse image} respectively, such that $f_!$ agrees with the \emph{derived direct image} $\RightD f_\ast$ in case that the morphism $f$ is proper. Furthermore, one requires these functors to interact well with the surrounding theory -- this compatibility can be formally captured in what one calls a \emph{$6$-functor formalism}, see for example Scholze's lecture notes \cite{sixfunctors} on the topic. However, in the classical setting it is not clear how one can construct the functors $f_!$ and $f^!$ for a large collection of morphisms $f$. Indeed, in this approach to coherent duality the main difficulty stems from showing the existence of $f^!$ in case that $f$ is proper. 

For relief, we turn to Clausen and Scholze's theory. Using condensed techniques we will define for a large class of morphisms $f$ a functor $f_!$ which does not have an analogue in the classical setting: Its definition makes critical use of condensed mathematics! Furthermore, without many additional assumptions, this functor admits a adjoint $f^!$, see theorem \ref{theorem: exceptional direct and inverse image functors, the affine case}. And indeed, one can construct a whole $6$-functor formalism. This then allows us to prove coherent duality by hand for the affine case -- the global case will be covered in chapter \ref{chapter: globalization} as a final outlook.

%
%

	\newpage
	\tableofcontents
	
	\chapter*{Notations, Conventions \& Prerequisites}
	
All rings are assumed commutative and unitary. The category of rings is denoted by $\commutativeRings$. By $\set$ and $\Set$ we denote the categories of finite respectively all sets. The notation $\eHom$ always indicates an enrichment (\eg of abelian groups or of modules), while $\intHom$ denotes an internal $\Hom$ (usually attached to a closed symmetric monoidal structure). By $L\ladj R$ we denote a pair of adjoint functors where $L$ is left adjoint to $R$. The unit and counit of an adjunction $L\ladj R$ are always denoted by $\eta\colon \id\to R\circ L$ and $\epsilon\colon L\circ R\to \id$ respectively. If the adjunction is enriched over a monoidal category $\category{V}$ we write $F\ladj_\category{V} G$. The reader is assumed familiar with a good wealth of category theory, in particular with the language of objects internal to some category, \eg internal groups, internal abelian groups, internal rings, and so on.

We will usually call a morphism between \emph{space}-like objects a \emph{map}. \Eg a \emph{map of topological spaces} or a \emph{map of sites}.

The reader is assumed familiar with the basics of set-theory in particular with some knowledge of ordinals and cardinals.

Further, we fix the following notation:
\begin{tabbing}
	\hspace*{3.6cm} \= \kill

	$\category C$, $\category D$, ... \> Categories\\
	$\site C$, $\site D$, ... \> Sites\\
	$\preAnaRing{A}$, $\preAnaRing{B}$, ... \> Ppre-)analytic rings\\
	
	$\category C/X$ \> The over category of the category $\category C$ and the object $X\in\category C$\\
	$X/\category C$ \> The under category of the category $\category C$ and the object $X\in\category C$\\
	
	$\limit D = \limit_{i\in \category{I}} D(i)$ \> The limit of the diagram $D\colon \category{I}\to \category C$\\
	$\colimit D = \colimit_{i\in \category{I}} D(i)$ \> The colimit of the diagram $D\colon \category{I} \to \category C$\\
		
	$\ast$ \> A point, suitably interpreted, \eg a one-point topological space,\\
	\> a category with one element, ...\\
	

	
\end{tabbing}

	\chapter{Sites, Sheaves, Topoi \& Condensed Sets}
	\label{chapter: sites, sheaves, topoi and condensed set}
	
\section*{Introduction}
As mentioned in the introduction to this thesis, condensed sets will be modeled by sheaves on the site $\proFiniteSets$ of profinite sets. To understand what this means we will introduce the basics of the very general theory of sites and sheaves in section \ref{section: sites and sheaves}. This will lead us to investigate general categories of sheaves, so-called \emph{Grothendieck topoi}, in section \ref{section: topoi}. A generalization of topoi, called \emph{(infinitary) pretopoi}, are introduced in \ref{section: pretopoi}. Afterwards, in section \ref{section: condensed sets}, we will start to develop the theory of condensed sets and how they relate to traditional topological spaces.

\section{Sites \& Sheaves}
\label{section: sites and sheaves}

A site is a generalization of the notion of open sets and their covers in a topological space. A site provides the minimal amount of data needed to construct a theory of sheaves. Open sets and their inclusions are modeled by a category, while covers of objects are specified as extra datum.

Much of this section follows Johnstone's \emph{Sketches of an Elephant, Volume II} \cite[C2]{elephant2}.

\begin{definition}[Coverages \& Sites, {
	\cite[Definition C2.1.1]{elephant2}}]
	\label{definition: coverages and sites}
	
	Let $\site C$ be a category. A \emph{coverage} on $\site C$ assigns to each $U \in \site C$ a collection $\covers U$ of families $(f_i\colon U_i \to U)_{i\in I}$ over $U$, called \emph{covers}, such that:
	
	If $(f_i\colon U_i \to U)_{i\in I}$ is a cover and $g\colon V\to U$ is any morphism in $\site C$, then there exists a cover $(h_j\colon V_j \to V)_{j\in J}$ of $V$, such that for each $j\in J$ the morphism $g\circ h_j$ factors through some $f_i$ with $i \in I$. That is every solid diagram 
	\begin{center}
		\begin{tikzcd}
			V_j\ar[r, dotted]\ar[d, "h_j"] &U_i\ar[d, "f_i", dotted]\\
			V\ar[r, "g"] &U
		\end{tikzcd}
	\end{center}
	can be completed by a pair of dotted arrows for some $i\in I$.
	
	A category equipped with a coverage is called a \emph{site}. A \emph{(locally) small site} is a site whose underlying category is (locally) small.
\end{definition}

\begin{remark}[Terminology]
	A \emph{cover} (recouvrement) $(f_i\colon U_i \to U)_{i\in I}$ of some object $U$ should be thought of as an open cover of $U$ by the $U_i$. A \emph{coverage} is a collection of covers. These notions should be distinguished from the unrelated \emph{covering} (revêtement) as \eg the universal covering $\R \to S^1$ of the circle $S^1$ in topology. Differentiating between the terms \emph{cover}, \emph{coverage} and \emph{covering} has provided the author with some headache.
\end{remark}

\begin{remark}
	The existence condition in definition \ref{definition: coverages and sites} essentially asserts that a kind of base change is possible -- for every cover of $U$ and any morphism $g\colon V\to U$, there is a cover of $V$ compatible with $g$. Indeed, if $\site C$ admits fiber products such that pullbacks of covers are covers, then this condition is met immediately.
\end{remark}

\begin{example}[Examples of sites]
	We now give two important (classes of) examples of sites.
	\begin{enumerate}
		\item
		Every topological space $X$ with its associated category $\open X$ of open sets naturally admits the structure of a site. Indeed, for every open $U \in \open X$ call a family $(U_i \to U)_{i \in I}$ a cover of $U$ if and only if $U = \bigcup_{i \in I} U_i$. 
		
		\item	
		The \emph{(small) étale site} $\etalesite X$ of a scheme $X$: Denote by $\Schemes$ the category of schemes. As a category, $\etalesite X$ is the full subcategory of $\Overcategory \Schemes X$ of étale morphisms over $X$. A cover of such an étale morphism $U\to X$ is a family $(f_i \colon U_i \to U)_{i\in I}$ over $X$ such that each $f_i$ is étale and the $f_i$ are jointly surjective, \ie the images of the $f_i$ cover $U$.
	\end{enumerate}
\end{example}

It is often important to consider sub-sites, that is subcategories of the underlying category equipped with the induced coverage.
\begin{definition}[Sub-sites]
	Let $\site C$ be a site and $\site C'$ be a subcategory. Declare a family $(U_i\to U)_{i\in I}$ in $\site C'$ to be a cover if it is one in $\site C$. Then $\site C'$ equipped with this coverage is called a \emph{sub-site of $\site C$}. If $\site{C'}$ is a full subcategory, then $\site{C}'$ is called a \emph{full sub-site of $\site{C}$}.
\end{definition}

We can now define presheaves and sheaves on a site. Inspired by the example of a topological space and its category of open sets, a presheaf should map each object of the underlying site (\ie an open set) to an object of a certain category (\eg $\Set$) while morphisms (\ie inclusions of open sets) should be mapped to morphisms in a contravariant manner (\ie to \emph{restrictions} between these associated objects). Hence, presheaves should be contravariant functors on the site.
\begin{definition}[Presheaves]
	Consider a site $\site C$ and a category $\category D$. A \emph{$\category D$-valued presheaf on $\site C$} is a functor $\site C^\op \to \category D$. Morphisms of presheaves are natural transformations. Hence, the category $\Presheaves{\site C}{\category D}$ of $\category D$-valued presheaves on $\site C$ is simply the functor category $\functorCategory{\site C^\op}{\category D}$.
\end{definition}

A presheaf of sets is now called a sheaf if it respects covers. In essence if $(U_i\to U)_i$ is a cover then the values of the (pre)sheaf on $U$ should be precisely determined by \emph{compatible} values on each of the $U_i$.
\begin{definition}[Compatible families \& Sheaves of sets, {
		\cite[Definition C2.1.2]{elephant2}}]
	\label{definition: compatible families and sheaves of sets}
	
	Consider a site $\site C$.
	
	Suppose $(f_i\colon U_i\to U)_{i\in I}$ is a cover in $\site C$ and $\sheaf F$ is a presheaf of sets on $\site C$. A family of sections $s_i\in \sheaf F(U_i)$ is called \emph{compatible} if for every pair of maps $U_i \xleftarrow{g} V \xrightarrow{h} U_j$ with $f_i\circ g = f_j\circ h$ in $\site C$ we have $\sheaf F(g)(s_i) = \sheaf F(h)(s_j)$, \ie the restrictions of $s_i$ and $s_j$ to $\sheaf F(V)$ agree.
	
	The presheaf of sets $\sheaf F$ is called a \emph{sheaf} on $\site C$ if for every cover $(f_i\colon U_i\to U)_{i\in I}$ and every family $(s_i)_{i\in I}$ compatible with the cover, there exists a unique $s\in \sheaf F(U)$ such that $\sheaf F(f_i)(s)=s_i$, \ie the family of sections glues to a unique section over $U$. We obtain the full subcategory
	$$\Sheaves{\site C}{\Set} \subseteq \Presheaves{\site C}{\Set}$$
	of sheaves of sets on $\site C$.
\end{definition}

We will also consider sheaves with values in categories other than sets -- however only in the special case where values are sets with extra structure. What this should mean we shall formalize in the following definition.
\begin{definition}[Very concrete categories]
	\label{definition: very concrete categories}
	
	A category $\category{D}$ is called \emph{very concrete} if it admits all filtered colimits and there is a continuous, faithful and conservative functor $\category{D} \to \Set$ that commutes with filtered colimits.
\end{definition}

\begin{remark}
	Definition \ref{definition: very concrete categories} has two purposes: First, it allows for the general theory of sheaves of sets to translate well to sheaves with values in the very concrete category, see definition \ref{definition: sheaves in general}. Secondly, it provides a rigid setting such that remark \ref{remark: transitioning in general} will apply -- this is fundamental to work with \emph{condensed objects in the very concrete category} (we will see later what this means).
\end{remark}

\begin{example}[Some very concrete categories]
	Very concrete categories are for example the category of groups, its subcategory of abelian groups, the category of rings but not the category of topological spaces (\eg $[0, 2\pi) \to S^1$ is bijective but not a homeomorphism, so the forgetful functor is not conservative).
\end{example}

We obtain the following straightforward generalization of definition \ref{definition: compatible families and sheaves of sets}.
\begin{definition}[Sheaves in general, {
	\cite[Definition C2.1.2]{elephant2}}]
	\label{definition: sheaves in general}
	
	Consider a site $\site C$ and a very concrete category $\category D$. Then a $\category D$-valued presheaf $\sheaf F$ is called a \emph{sheaf}, if the associated underlying presheaf of sets is a sheaf. We obtain a full subcategory
		$$\Sheaves{\site C}{\category D} \subseteq \Presheaves{\site C}{\category D}$$
	of $\category D$-valued sheaves on $\site C$.
\end{definition}

If the site admits fiber products (which one should regard as the analogue of the intersection of open sets), then compatibility of sections is given by \emph{agreement on intersections} and one can write down certain equalizer conditions that are sufficient and necessary for being a sheaf.
\begin{lemma}[Being a sheaf is an equalizer condition]
	\label{lemma: being a sheaf is an equalizer condition}

	Let $\site C$ be a site with fiber products and $\category D$ be a very concrete category. A presheaf $\sheaf F$ is a sheaf if and only if for each $U \in \site C$ and each cover $(f_i \colon U_i \to U)_{i\in I}$ of $U$ the natural morphism
	$$\begin{tikzcd}
		\sheaf F(U) \ar[r] &\equalizer\Bigg(\prod_{i\in I} \sheaf F(U_i) \ar[r, shift left=1.1] \ar[r, shift right=1.1] &\prod_{i, j \in I} \sheaf F(U_i \times_U U_j)\Bigg)
	\end{tikzcd}$$
	is an isomorphism, where the two parallel maps are induced by the assignments $(s_i)_{i\in I}\mapsto \sheaf F(U_i\times_U U_j \to U_i)(s_i)$ and $(s_i)_{i\in I}\mapsto \sheaf F(U_i\times_U U_j \to U_j)(s_j)$.
\end{lemma}
\begin{proof}
	Let $(f_i\colon U_i\to U)_i$ be a cover. Clearly $\prod_{i\in I} \sheaf F(U_i)$ is the set of families of sections $(s_i)_{i\in I}$ where $s_i\in \sheaf F(U_i)$.
	
	Suppose such a family $(s_i)_{i\in I}$ is compatible with the cover. Then by definition of compatibility, for all $i,j\in I$ the restrictions of $s_i$ and $s_j$ under the two maps $\sheaf F(U_i)\rightrightarrows \sheaf F(U_i\times_U U_j)$ agree. Hence, every compatible family lies in the equalizer.
	
	Now suppose we are given a family $(s_i)_{i\in I}$ in this equalizer. We will show that the family is compatible. Suppose that $U_i\xleftarrow g V\xrightarrow h U_j$ is any pair of morphisms with $f_i\circ g=f_j\circ h$. Then by the universal property of the fiber product there is a unique compatible morphism $u\colon V\to U_i\times_U U_j$, such that the diagram
	$$\begin{tikzcd}
		V\ar[dr, "u"]\ar[drr, bend left=15, "g"]\ar[ddr, bend right=15, "h", swap]&&\\
			&U_i\times_U U_j \ar[r, "p_i", swap] \ar[d, "p_j"] &U_i\ar[d, "f_i", swap]\\
			&U_j \ar[r, "f_j"] & U
	\end{tikzcd}$$
	commutes. Now as by definition of the equalizer we have $\sheaf F(p_i)(s_i)=\sheaf F(p_j)(s_j)$ it follows, that
		$$\sheaf F(g)(s_i)=(\sheaf F(u)\circ \sheaf F(p_i))(s_i)=(\sheaf F(u)\circ \sheaf F(p_j))(s_j)=\sheaf F(h)(s_j)$$
	proving compatibility, since $i, j$ were arbitrary.
	
	Hence, the equalizer is precisely the set of compatible families. The claim now follows immediately. Indeed, surjectivity of the induced map to the equalizer ensures that gluing is possible, while injectivity ensures that gluing is unique.
\end{proof}

As we are used to, there exists a canonical way of turning a presheaf of sets into a sheaf of sets -- this is of course sheafification.
\begin{definition}[Sheafification]
	\label{definition: sheafification}
	
	For a small site $\site C$ the inclusion $\Sheaves{\site C}{\Set} \to \Presheaves{\site C}{\Set}$ of sheaves into presheaves admits a left adjoint
		$$\sheafification{-}\colon \Presheaves{\site C}{\Set} \to \Sheaves{\site C}{\Set}$$
	that commutes with finite limits, called \emph{sheafification}.
\end{definition}
\begin{proof}
	This follows from \cite[Prop. C2.2.6]{elephant2} after observing that the colimits in the definition of the \emph{separation} $\sheaf F\mapsto\sheaf F^+$ (and of the sheafification $\sheafification{-}\colon \sheaf{F}\mapsto \sheaf{F}^{++}$) are filtered and hence commute with finite limits.
\end{proof}

\begin{remark}[Complications of sheafification in a general setting]
	We will see more general versions of definition \ref{definition: sheafification} later on. One is often interested in sheaves with values in some category $\category D$ that is not $\Set$. This in general turns out to be somewhat problematic and one needs to impose some conditions on the category $\category D$. A sufficient (but rather technical, hence we will not untangle them) set of conditions on $\category D$ asks, that $\category D$ is bicomplete, that in $\category D$ small filtered colimits are exact and that $\category D$ satisfies the so-called \emph{IPC-property}.
	
	We will not need that kind of generality as we will only consider (pre)sheaves of sets, abelian groups and rings. It is not very hard to see that presheaves of abelian groups/rings are equivalent to abelian/ring objects internal to the category of presheaves of sets. Since sheafification preserves finite limits one naturally obtains the structure of an abelian group/ring on the sheafification of the underlying presheaf. We will repeat this argument once we introduce abelian sheaves in chapter \ref{chapter: rings, modules and condensed algebra}.
\end{remark}

We should now ask ourselves how to relate different sites, in particular how to relate sheaves thereon. Once again, we turn to topological spaces for our defining motivation. If $f\colon X \to Y$ is a map of topological spaces, then by taking preimages under $f$ we obtain a functor $f^{-1}\colon \open Y \to \open X$ mapping unions to unions and covers to covers. As we are dealing with the sites associated to $X$ and $Y$, we only have access to open sets, not elements thereof, so this map $f^{-1}$ alone should determine a map of sites. However, a site should be considered as a generalized space, so the functor $f^{-1}$ \quote{points in the wrong direction}, as we ideally would like any continuous map $f\colon X\to Y$ to have a corresponding map of sites $\open X\to \open Y$ in the same direction. One solves this problem rather pragmatically by declaring a map of sites $f\colon \site{C}\to\site{D}$ to be a functor $f^\top\colon\site{D}\to \site{C}$ in the other direction, satisfying some additional properties. Certainly $f^\top$ should preserve covers. Now in the case of topological spaces, $f^{-1}$ maps intersection to intersections and more generally limits in $\open{Y}$ to limits in $\open{X}$. We will require the same of $f^\top$ if(!) $\site{D}$ is finitely complete. If $\site{D}$ does not admit all finite limits, then one has to be more careful.
\begin{definition}[Maps of sites]
	\label{definition: maps of sites}
	
	Let $\site{C}$ and $\site{D}$ be sites. Suppose that $\site{D}$ is finitely complete. Then a map $f\colon \site C \to \site D$ of sites is a functor $f^\top\colon \site D \to \site C$, commuting with finite limits and mapping covers to covers.
\end{definition}

\begin{remark}[Maps of sites in full generality]
	\label{remark: maps of sites in full generality}
	
	If the site $\site{D}$ in definition \ref{definition: maps of sites} is not finitely complete, one needs an alternative definition of \emph{maps of sites}. A generalization is for example given in \cite[Def. 4.10]{shulman2012exactcompletionssmallsheaves} by Shulman and asks for $f^\top$ to be \emph{covering flat}, essentially requiring that $f^\top$ \emph{preserves all finite limits that should have existed}. We will however only need this generalization in case that $\site{D}$ is a full subcategory of $\site{C}$. In this case it follows by \cite[Def 11.1 \& Thm. 11.2 (b)]{shulman2012exactcompletionssmallsheaves} that the inclusion satisfies this condition, hence giving us a map of sites $\site{C} \to \site{D}$. For a discussion of the notion of covering flatness see  \href{https://golem.ph.utexas.edu/category/2011/06/flat_functors_and_morphisms_of.html}{Shulman's post in the $n$-café}.
\end{remark}

Now having introduced maps of sites we can of course ask to transport sheaves along them. As with topological spaces, if $f\colon X\to Y$ is a continuous map and $\sheaf F$ is a sheaf on $X$, we easily obtain a sheaf on $Y$ by simply setting $f_\ast\sheaf F(V) \coloneqq \sheaf F(f^{-1}(V))$, \ie by precomposing with $f^{-1}$. This works since taking preimages maps open sets to open sets and preserves inclusions. Since in the general case taking preimages of open sets is modeled by the functor $f^\top$, it is clear how to proceed.
\begin{lemma}[Direct images]
	
	Let $f\colon \site C\to \site C'$ be a map of sites and $\category D$ a very concrete category. The functor
	\begin{align*}
		f_\ast \colon \Presheaves{\site C}{\category D} &\to \Presheaves{\site C'}{\category D},\\
		\sheaf F &\mapsto \sheaf F\circ f^\top \ldef \left((\site C')^\op \xrightarrow{f^\top}\site C^\op\xrightarrow{\sheaf F} \category D\right)
	\end{align*}
	is called the \emph{direct image} and restricts to a functor $f_\ast \colon \Sheaves{\site C}{\category D} \to \Sheaves{\site C'}{\category D}$.
\end{lemma}
\begin{shortproof}
	This is very similar to \cite[Lem. C2.3.3]{elephant2} keeping in mind \cite[Rem. C2.3.7]{elephant2}.
\end{shortproof}

Similarly, one can ask to transport sheaves in the other direction. As with topological spaces, this is initially problematic. If one tries to mimic the definition of the direct image by setting $f^{-1}\sheaf F(U) \coloneqq \sheaf F(f(U))$ one immediately runs into the problem that $f(U)$ might not be open and hence $\sheaf F(f(U))$ is ill-defined. Abstractly the solution comes from simply requesting a functor $f^{-1} \ladj f_\ast$ left adjoint to the direct image. More concretely, one can try to approximate $f(U)$ from above by open sets. For this take open sets $V$ such that $U\subseteq f^{-1}(V)$ (which implies $f(U)\subseteq V\cap \image f \subseteq V$) and glue the values $\sheaf F(V)$ suitably. In the general case this is described by a filtered colimit given by the comma-category of maps $U\to f^\top(V)$.
\begin{proposition}[Inverse images]
	\label{Proposition: Existence of inverse image}
	
	Let $f\colon \site C\to \site C'$ be a map of small sites and $\category D$ a very concrete category such that sheafification of $\category D$-valued presheaves exists. Assume that $\site C$ is small and that $\site C'$ is essentially small. There exists a functor $f^{-1}\colon \Sheaves{\site C'}{\category D}\to \Sheaves{\site C}{\category D}$ left adjoint to the direct image $f_\ast$, called the \emph{inverse image}. 
	
	More precisely for each $U \in \site C$ define the category $I_U$ with objects $\{(V, \phi)\setseparator V\in\site C', \phi\colon U \to f^\top(V)\}$ and morphism $\Hom_{I_U}((V, \phi), (V', \phi')) = \{g\colon V\to V'\setseparator g\circ \phi = \phi'\}$. For a sheaf $\sheaf F$ on $\site C'$, we obtain the inverse image as the sheafification
		$$f^{-1}\sheaf F = \sheafification{U\mapsto \colimit_{(V, \phi) \in I_U} \sheaf F(V)}.$$
	Notice how each morphism $g\colon U\to U'$ in $\site C$ functorially induces a morphism $I_U \to I_{U'}, (V, \phi)\mapsto (V, \phi\circ g)$, so the rule above similarly extends to morphisms, hence prior to sheafification we really obtain a presheaf.
\end{proposition}
\begin{shortproof}
	This is very similar to \cite[Cor. C2.3.4]{elephant2} keeping in mind \cite[Rem. C2.3.7]{elephant2}.
\end{shortproof}

In the special case where one site is a full sub-site of the other (see remark \ref{remark: maps of sites in full generality}), the situation is a bit more straightforward. For topological spaces this corresponds to pulling back from a coarser topology to a finer one, along the identity. In this case we cover $f(U)=U$ by open sets $V$ in the coarser topology.
\begin{remark}[Inverse images for inclusions]
	\label{remark: inverse images for inclusions}
	Suppose $f\colon\site C \to \site C'$ is the map of sites corresponding to a full inclusion $\site C'\subseteq \site C$ of underlying categories. Then the categories $I_U$ considered in proposition \ref{Proposition: Existence of inverse image} drastically simplify and we are justified in writing
		$$f^{-1}\sheaf F = \sheafification{U\mapsto \colimit_{\substack{U\to V\\V\in \site C'}} \sheaf F(V)}$$
	for any sheaf $\sheaf F$.
\end{remark}

\section{Topoi}
\label{section: topoi}

We will now study the properties of sheaf categories of small sites.
\begin{definition}[Grothendieck Topoi, {\cite[Definition C2.2.9]{elephant2}}]
	A \emph{(Grothendieck) topos} is a category $\category T$ equivalent to the category of sheaves $\Sheaves{\site C}{\Set}$ on some small site $\site C$.
\end{definition}

\begin{remark}[Being a topos is a property]
	\label{remark: giraud's theorem}
	
	In \ref{theorem: girauds axioms} we will see a theorem due to Giraud that characterizes a topos $\category T$ by certain axioms internal to the category $\category T$. Hence, \emph{being a topos} is a property of a given category, not extra structure. One of these axioms (or rather a relaxation thereof) will play an important role in the theory of condensed mathematics. The axiom states: $\category T$ admits a generating set. We will introduce the notion of generation in \ref{section: generation} and use it to provide the mentioned relaxation in \ref{section: pretopoi}.
\end{remark}

We will now define a suitable notion of morphism between topoi. Recall that for any map of sites $f\colon \site C\to\site D$ there is an adjoint pair $f^{-1}\ladj f_\ast$ between the associated topoi. We mimic such pairs and their properties when defining morphisms of topoi.
\begin{definition}[Morphisms of topoi, {\cite[Definition A4.1.1]{elephant1}}]
	\label{definition: Morphisms of topoi}
	
	A \emph{(geometric) morphism} of topoi $\category T\to \category T'$ is a pair $f^{-1} \ladj f_\ast$ of adjoint functors such that $f^{-1}\colon \category{T'}\to\category{T}$ commutes with finite limits. In particular, isomorphisms of topoi are precisely adjoint equivalences. We will denote such a morphism by
		$$\category T \xleftrightarrows{f^{-1}}{f_\ast} \category T'.$$
\end{definition}

\begin{remark}
	Due to the adjunction $f^{-1} \ladj f_\ast$ in definition \ref{definition: Morphisms of topoi}, we additionally know that $f^{-1}$ and $f_\ast$ commute with all colimits and all limits respectively. Hence, in total $f^{-1}$ is cocontinuous and finitely continuous, while $f_\ast$ is continuous.
\end{remark}

Now one can ask when a given sub-site gives rise to an equivalent category of sheaves. The analogue for topological spaces are of course bases. Recall that that sheaves on a topological space $X$ are the same as sheaves on a base of $X$ -- this is for example used to construct the structure sheaf of an affine scheme $\Spec A$ from its basis of distinguished opens (where one can explicitly give the sections on $D(f)$ with $f\in A$ as the localization $A_f$). We will restrict out attention to \emph{full sub-sites}, \ie sub-sites whose underlying category is full in the underlying category of the larger site.
\begin{definition}[Dense sub-site]
	Let $\site C$ be a site. A sub-site $\site B$ of $\site C$ is said to be \emph{dense}, if it is full and every $U\in \site C$ admits a cover $(f_i\colon U_i\to U)_{i\in I}$ where $U_i\in \site B$ for all $i\in I$.
\end{definition}

We now obtain the desired result which formalizes the idea, that giving a sheaf is the same as giving a sheaf on a base of that topology.
\begin{lemma}[Comparison lemma]
	\label{lemma: comparison lemma}
	Let $\site B$ be a small dense sub-site of a locally small site $\site C$. The map $\site C\to\site B$ of sites (given by the inclusion of the underlying categories) gives rise to an adjoint equivalence
		$$\Sheaves{\site C}{\Set}\xleftrightarrows{i^{-1}}{i_\ast} \Sheaves{\site B}{\Set}$$
	of topoi.
\end{lemma}
\begin{shortproof}
	This is \cite[C2.2.3]{elephant2} or alternatively \cite[Thm. 11.8]{shulman2012exactcompletionssmallsheaves} with suitable choice of $\kappa$.
\end{shortproof}

This lets us introduce some handy terminology. If the underlying category of a site is large, but the site has a small dense sub-site, we obtain from the comparison lemma \ref{lemma: comparison lemma} that their sheaf categories are equivalent. 
\begin{definition}[Essentially small sites]
	\label{definition: essentially small sites}
	
	A locally small site $\site C$ is called \emph{essentially small} if it admits a small dense sub-site.
\end{definition}

\begin{remark}[Terminological clash]
	Notice that there is some clash of terminology when talking about essentially small sites. When referring to an essentially small $\site C$, we might (contrary to the previous definition in definition \ref{definition: essentially small sites}) interpret this as saying that the underlying category of $\site C$ is essentially small (\ie its isomorphism classes are in bijection with a set). Indeed, some authors define an essentially small site to be a site whose underlying category is essentially small.
	
	Clearly, if the underlying category is essentially small, one can obtain a dense sub-site of $\site C$ by considering a skeleton of the underlying category with induced covers. This is then a small dense sub-site, so the site $\site C$ is essentially small. It is however not clear to the author that the converse holds true as well. To him there seems to be no reason that admitting a small dense sub-site tells anything about the isomorphism classes of the underlying category.
\end{remark}

\section{Generation}
\label{section: generation}

As mentioned in remark \ref{remark: giraud's theorem}, we will need to study the notion of generation. When defining condensed sets, we will need to work with sheaves on a large site. On such a site however, it is dangerous to consider the full category of sheaves, since there is no guarantee that it is locally small: For two given sheaves, morphisms between them are natural transformations -- these are collections of morphisms (the components of the transformation) indexed by the objects of the site in question! This is a priori a large amount of data! Indeed, the full category of sheaves on a large site in general is neither a topos nor locally small!

If one imposes certain regularity on the two sheaves, then one can biject the morphisms between any two regular sheaves to a set. Hence, the full subcategory of such regular sheaves turns out to be locally small a posteriori. Condensed sets indeed will be defined as regular sheaves on a large site! This regularity is what we will call \emph{smallness}, see definition \ref{definition: small (pre)sheaves}.

The question now remains: What are the properties of this new category of \emph{regular/small} sheaves? It turns out that for the category of condensed sets the answer is rather satisfactory! It satisfies all but one of Giraud's axioms: It does not admit a set of generators, but instead one obtains only a proper class of them! Such a category will be called a \emph{pretopos} -- we will introduce pretopoi in \ref{section: pretopoi}. 

Let us now work towards the notion of generation. A natural starting point is the theory of stalks. Recall that for any point $x$ of a topological space $X$, one obtains the stalk functor $(-)_x\colon \Sheaves{X}\Set\to \Set$. Let us model this in the setting of topoi. The inclusion $\ast \to X$ of the point induces a morphism of sites $i\colon\open{\ast} \to \open{X}$ which in turn induces a morphism of topoi $$\Set \iso \Sheaves{\ast}{\Set} \xleftrightarrows{i^{-1}}{i_\ast} \Sheaves{X}{\Set}.$$
Here $i^{-1}\sheaf F = \sheaf F_x$ for any sheaf $\sheaf F$ recovers the stalk of $\sheaf F$ at $x$. Such a geometric morphism will be called a \emph{point} of the topos. Recall that a morphism $f\colon \sheaf F\to\sheaf F'$ is an isomorphism if and only if $f_x\colon \sheaf F_x\to\sheaf F'_x$ is an isomorphism for all points $x\in X$. While in a general topos this will not be true, it is still instructive to consider points and stalks for a general topos nonetheless.
\begin{definition}[(Having enough) points in a topos]
	\label{definition: points in a topos}
	
	Let $\category T$ be a topos. A \emph{point} $x=(x^{-1}, x_\ast)$ of $\category T$ is a geometric morphism
		$$\Set\iso \Sheaves \ast \Set \xleftrightarrows{x^{-1}}{x_\ast}\category T.$$
	For $\sheaf F\in \category T$ the object $\sheaf F_x \ldef x^{-1}\sheaf F$ is called the \emph{stalk of $\sheaf F$ at $x$}. The topos $\category T$ is said to \emph{have enough points} if any morphism $f\colon \sheaf F\to \sheaf F'$ in $\category T$ is an isomorphism if and only if the morphism $f_x\ldef x^{-1}f$ of stalks is an isomorphism for each point $x$ of $\category T$.
\end{definition}

Unfortunately not all topoi have enough points. An example from logic is provided by Sina Hazratpour \href{https://sinhp.github.io/scribbling/22-11-2017-a-topos-without-points}{on his website}. However, even if they do, those points might be inexplicit and hard to work with. Furthermore, if one leaves the realm of topoi and works with sheaves on a large site instead, then it seems (see \cite{pointsPretopoiLargeSite}) that defining \emph{points} as in definition \ref{definition: points in a topos} is no longer appropriate.

We will thus abandon the notion of points and search for a suitable replacement: It turns out that the idea of considering functors that \emph{jointly} give information about the category in question is a fruitful approach on even a general category. For this we will ask for a collection of objects $G$ such that the functors $\Hom(G, -)$ are \emph{jointly faithful}. Such objects are called \emph{generators} or \emph{separators}. Clausen and Scholze make extensive use of generators in \cite{condenseddotpdf}: We will for example see in \ref{proposition: abelian sheaves form an abelian category} that the category of condensed abelian groups is abelian and has enough projectives -- the projectives being generators!

Now for a topos (\ie sheaves on a \emph{small} site) the collection of representable sheaves forms (essentially by Yoneda's lemma) a generating set. In particular one usually ask for a \emph{small} collection of generators. This however is problematic for applications in condensed mathematics (or more generally sheaves of a large site). We will thus need a generalization to possibly non-small collections. Unfortunately, Scholze and Clausen mention no such generalization in \eg \cite{condenseddotpdf} and the author is also not aware of any other standard reference. Instead, the author would like to thank \href{https://zll22.user.srcf.net/profile/}{Zhen Lin Low} for providing such a generalization in \cite{generationZhenLin}.
\begin{definition}[Generation]
	\label{definition: generation}
	
	Let $\category C$ be a category. A class $\mathcal G$ of objects of $\category C$ is called \emph{generating}, if for any object $X\in \category C$ there is a small(!) collection (\ie a set) $\mathcal G'\subseteq \mathcal G$, such that for any parallel pair of morphisms $f,g \colon X\to Y$ in $\category C$ with the property that $f\circ h = g\circ h$ for any map $h\colon G\to X$ from some $G\in\mathcal G'$, we already have $f=g$. Such a set $\mathcal G'$ is called \emph{separating for $X$}.
\end{definition}

\begin{remark}
	Suppose that $\category C$ is a category and $\mathcal G$ a generating class of $\category C$. Much like stalks, this generating class can make statements about a morphism $f\colon X\to Y$ being an isomorphism. Suppose we have a candidate morphism $g\colon Y\to X$ that we suspect is an inverse of $f$. This can be reformulated as asking the two questions $\id_X \overset{?}= g\circ f$ and $\id_Y\overset ? = f\circ g$, which can be answered using the generating class $\mathcal G$. Indeed, by checking $(\id_X)_\ast\overset ? = (g\circ f)_\ast \colon \Hom(G, X)\to \Hom(G, X)$ and $(\id_Y)\overset ? = (f\circ g)_\ast\colon \Hom(G, Y)\to\Hom(G, Y)$ for any $G\in\mathcal G'$ in some suitable separating set $\mathcal G'\subseteq \mathcal G$, we can answer these two questions!
\end{remark}

In case the category in question carries more structure, definition \ref{definition: generation} can be reformulated in a manner of ways. Of particular importance to us is the case of an abelian category.
\begin{lemma}[Generating classes in abelian categories]
	\label{lemma: generating classes in abelian categories}
	
	Let $\category A$ be an abelian category with arbitrary direct sums. Then for a class $\mathcal G$ of objects of $\category A$ the following are all equivalent:
	\begin{enumerate}[(i)]
		\item
		\label{generating sets in abelian categories: i}
		$\mathcal G$ is a generating class.
		
		\item
		\label{generating sets in abelian categories: epi}
		Every object $X\in\category A$ admits an epimorphism $\bigoplus_{i\in I} G_i \twoheadrightarrow X$ where $(G_i)_{i\in I}$ is a small(!) family of objects of $\mathcal G$.
		
		\item
		\label{generating sets in abelian categories: presentable}
		Every object can be realized as a coequalizer of a pair of morphisms between small direct sums of objects in $\mathcal G$.
	\end{enumerate}
\end{lemma}
\begin{proof}
	\ref{generating sets in abelian categories: i} $\Rightarrow$ \ref{generating sets in abelian categories: epi}:\newline
	Let $\mathcal G'\subseteq \mathcal G$ be a small separating set for $X$. Consider the morphism $p\colon \bigoplus_{\substack{G\to X\\G\in\mathcal G'}} G \to X$. We claim that this is an epimorphism. Indeed, if $f,g\colon X\to Y$ is a pair of parallel morphism in $\category A$ for which $f\circ p = g\circ p$, then also $f\circ h = g\circ h$ for any $h\colon G\to X$ with $G\in\mathcal G'$ since $h$ admits a factorization $G\to \bigoplus_{\substack{G\to X\\G\in\mathcal G'}} G \to X$ by definition of the direct sum. Thus by choice of $\mathcal G'$, we must have $f=g$. Since $f$ and $g$ were arbitrary, $p$ must be epic.
	
	\ref{generating sets in abelian categories: epi} $\Rightarrow$ \ref{generating sets in abelian categories: i}:\newline
	Let $X$ be any object in $\category A$. Choose an epimorphism $p\colon \bigoplus_{i\in I} G_i \twoheadrightarrow X$ and define $\mathcal G'\ldef \{G_i\setseparator i\in I\}$. We claim that the set $\mathcal G'$ is separating for $X$. For this let $f,g\colon X\to Y$ be a pair of parallel morphisms in $\category A$ such that for any $h\colon G\to X$ with $G\in\mathcal G'$ we have $f\circ h=g\circ h$. Then in particular for $i\in I$ with associated $h_i\colon G_i \to \bigoplus_{i\in I} G_i \twoheadrightarrow X$ we have $f\circ h_i = g\circ h_i$, so also $f\circ p = g\circ p$. But $p$ is epic, hence $f=g$.
	
	\ref{generating sets in abelian categories: epi} $\Rightarrow$ \ref{generating sets in abelian categories: presentable}:\newline
	Consider an epimorphism $p\colon \bigoplus_{i\in I} G_i \twoheadrightarrow X$. Pick a similar epimorphism $p'\colon \bigoplus_{j\in J} G_j' \twoheadrightarrow \ker p$ and define $q=i\circ p'$ where $i$ is the inclusion of the kernel. We obtain the commutative diagram
	$$\begin{tikzcd}
		0 &\bigoplus_{j\in J} G_j' 	&\bigoplus_{i\in I} G_i &\coKernel q 	&0\\
		0 &\kernel p 				&\bigoplus_{i\in I} G_i &X 				&0
		\ar[from=2-1, to=2-2]
		\ar[from=2-2, to=2-3, "i"]
		\ar[from=2-3, to=2-4, "p", two heads]
		\ar[from=2-4, to=2-5]
		\ar[from=1-1, to=1-2]
		\ar[from=1-2, to=1-3, "q"]
		\ar[from=1-3, to=1-4]
		\ar[from=1-4, to=1-5]
		\ar[from=1-2, to=2-2, two heads, "p'"]
		\ar[from=1-3, to=2-3, equal]
		\ar[from=1-4, to=2-4, dashed]
	\end{tikzcd}$$
	with exact rows and induced morphism $\coKernel q\to X$ by the universal property of $\coKernel q$ since $p \circ q = 0$. Hence, by the 5-lemma (with the induced morphism in the center), this induced morphism is an isomorphism, \ie $X$ is a cokernel of $q$.
	
	\ref{generating sets in abelian categories: presentable} $\Rightarrow$ \ref{generating sets in abelian categories: epi}:\newline
	Clear since any coequalizer is just a cokernel.
\end{proof}

\section{Pretopoi}
\label{section: pretopoi}

In this section we will provide the definition of pretopoi, where \emph{being a pretopoi} is a certain relaxation to \emph{being a topos}. As previously mentioned this relaxation hinges on a theorem of Giraud, stated in theorem \ref{theorem: girauds axioms}, which characterizes a topos by properties internal to the category in question. 

We will need a few additional concepts to state Giraud's theorem. These properties are true in any topos and are heavily inspired by the behavior of $\Set$. 

\begin{definition}[Universally disjoint coproducts]
	Let $\category T$ be a category with all finite limits and all small colimits.
	
	A coproduct $X\ldef \coprod_{i\in I} X_i$ in $\category T$ is called \emph{disjoint}, if the associated morphisms $X_i\to X$ are monomorphisms and for each pair $i\ne j$ in $I$ the fiber product $X_i\times_X X_j$ (the intersection of $X_i$ and $X_j$ inside $X$) is an initial object in $\category T$ (\ie the intersection is empty). If additionally for every morphism $X'\to X$ the pullbacks of the morphisms $X_i\to X$ induce an isomorphism
		$$\coprod_{i \in I} (X'\times_X X_i) \to X',$$
	then the coproduct is called \emph{universally disjoint}.
\end{definition}

\begin{definition}[Equivalence relations \& quotients]
	Let $\category T$ be a category with all finite limits and all small colimits.
	
	Let $X$ be an object of $\category T$. A \emph{relation on $X$} is a monomorphism $\iota\colon R\to X\times X$ from some object $R$ of $\category T$. Denote its components by $p_1, p_2\colon R\to X$, so that $\iota=(p_1, p_2)$. Then $\iota$ is called an \emph{equivalence relation}
	\begin{enumerate}[(i)]
		\item
		if the diagonal map $X\to X\times X$ factors over $\iota$ (that is $\iota$ is \emph{reflexive}),
		
		\item
		if the morphism $(p_2, p_1)\colon R\to X\times X$ factors over $\iota$ (that is $\iota$ is \emph{symmetric}),
		
		\item 
		and if the morphism $(p_1\circ q_1, p_2\circ q_2)\colon R\times_X R \to X\times X$ factors over $\iota$ as well, where $R\times_X R$, $q_1$ and $q_2$ are given by the pullback
			$$\begin{tikzcd}
				R\times_X R\ar[r, "q_2"]\ar[d, "q_1"] &R \ar[d, "p_1"]\\
				R \ar[r, "p_2"] &X
			\end{tikzcd}$$
		(that is $\iota$ is \emph{transitive}).
	\end{enumerate}
	If $\iota\colon R\to X\times X$ is an equivalence relation then its \emph{quotient} is the coequalizer
		$$R\rightrightarrows X\twoheadrightarrow E/R,$$
	and it is always an epimorphism.
\end{definition}

\begin{definition}[Kernel pairs]
	Let $\category T$ be a category with all finite limits and all small colimits.
	
	If $u\colon X\to X'$ is any morphism, then a \emph{kernel pair} for $u$ is a parallel pair of morphisms $p_1, p_2\colon R\to X$ such that
		$$\begin{tikzcd}
			R\ar[r, "p_2"]\ar[d, "p_1"] &X\ar[d, "u"]\\
			X\ar[r, "u"] &X'
		\end{tikzcd}$$
	is a pullback square.
\end{definition}

\begin{definition}[Stably exact forks]
	Let $\category T$ be a category with all finite limits and all small colimits.
	
	A diagram in $\category{T}$ of the form
		$$\begin{tikzcd}
			R\ar[r, "p_1", shift left=.25em]\ar[r, "p_2", shift right=.25em, swap] &X \ar[r, "q"] &Q
		\end{tikzcd}$$
	is called a \emph{fork}. A fork is said to be \emph{exact} if $q$ is a coequalizer of $p_1$, $p_2$ and if $p_1$, $p_2$ form a kernel pair for $q$. If additionally any base change
		$$\begin{tikzcd}
			R\times_Q Q' \ar[r, shift left=.25em]\ar[r, shift right=.25em] &X\times_Q Q' \ar[r] &Q \times_Q Q'\cong Q'
		\end{tikzcd}$$
	along some morphism $Q'\to Q$ stays exact, then the fork is called \emph{stably exact}.
\end{definition}

We can now state Giraud's axioms.
\begin{theorem}[Giraud's axioms]
	\label{theorem: girauds axioms}
	
	A locally small category $\category T$ with all finite limits and all small colimits is a topos if and only if
	\begin{enumerate}[(i)]
		\item
		coproducts in $\category T$ are universally disjoint,
		
		\item
		every epimorphism in $\category T$ is a coequalizer,
		
		\item
		every equivalence relation in $\category T$ is a kernel pair,
		
		\item
		every exact fork $R\rightrightarrows X\to Q$ in $\category T$ is stably exact,
		
		\item
		and $\category T$ admits a small generating class.
	\end{enumerate}
\end{theorem}
\begin{shortproof}
	This is proven in \cite[Appendix §3]{sheaves-in-geometry-and-logic}.
\end{shortproof}

This theorem now allows us to define what a \emph{pretopos} should be. As previously mentioned the idea is to relax Giraud's axioms by allowing \emph{large} classes of generators.
\begin{definition}[(Infinitary) pretopoi]
	\label{definition: pretopoi}
	
	A locally small category $\category T$ with all finite limits and all small colimits is called a \emph{(infinitary) pretopos} if
	\begin{enumerate}[(i)]
		\item
		\label{definition: pretopoi, coproducts are universally disjoint}
		coproducts in $\category T$ are universally disjoint,
		
		\item
		\label{definition: pretopoi, every epimorphism is a coequalizer}
		every epimorphism in $\category T$ is a coequalizer,
		
		\item
		\label{definition: pretopoi, every equivalence relation is a kernel pair}
		every equivalence relation in $\category T$ is a kernel pair,
		
		\item
		\label{definition: pretopoi, every exact fork is stably exact}
		every exact fork $R\rightrightarrows X\to Q$ in $\category T$ is stably exact,
		
		\item
		\label{definition: pretopoi, existence of a generating class}
		and $\category T$ admits a generating class.
	\end{enumerate}
	
	A \emph{(geometric) morphism of pretopoi} $\category T\to\category T'$ is a pair $f^{-1}\ladj f_\ast$ of adjoint functors, such that $f^{-1}$ commutes with finite limits. As with morphisms of topoi, we will denote such a morphism of pretopoi by
		$$\category T \xleftrightarrows{f^{-1}}{f_\ast} \category T'.$$
\end{definition}

\section{Regular Projectivity \& Compactness}

The category of condensed sets will be a pretopos, so admits a class of generators. This class of generators can be chosen to consist of very well-behaved objects. In this section we will introduce the terminology that formalizes what this is supposed to mean.
\begin{definition}[Reflexive coequalizers]
	A \emph{reflexive pair} in a category $\category T$ is a diagram of the form
		$$\begin{tikzcd}
			\cdot \ar[r, shift left=0.35em]\ar[r, shift right=0.35em] &\cdot\ar[l]
		\end{tikzcd}$$
	where the arrow to the left is a common section of the two arrows to the right. A colimit of a reflexive pair is called a \emph{reflexive coequalizer}.
\end{definition}

\begin{remark}[Reflexive pairs from reflexive relations]
	\label{remark: reflexive pairs from reflexive relations}
	
	Notice that if $X$ is some object and $\iota\colon R\to X\times X$ is any reflexive relation on $X$ then by definition there is some morphism $X\to R$ providing a factorization of $X\to X\times X$ along $\iota$. This morphism is a common section of the two components $R\to X$ of $\iota$. Hence, to every reflexive relation we can naturally associate a reflexive pair. If the relation is even an equivalence relation then the associated quotient is a reflexive coequalizer for this associated reflexive pair.
\end{remark}

\begin{definition}[Regular projectivity \& compactness]
	\label{definition: regular projectivity and compactness}
	
	Let $\category T$ be a cocomplete category and $P\in \category T$. Then $P$ is called
	\begin{itemize}
		\item
		\emph{(reflexive) projective} if $\Hom_\category T(P, -)$ commutes with reflexive coequalizers,
		
		\item 
		\emph{compact} if $\Hom_\category T(P, -)$ commutes with filtered colimits.
		
		\item 
		\emph{compact projective} if it is both compact and reflexive projective.
	\end{itemize}
\end{definition}

\begin{remark}[Comparison of reflexive projectivity to normal projectivity]
	The definition of reflexive projectivity in definition \ref{definition: regular projectivity and compactness} is closely related to the usual notion of projectivity and provides a relaxation to it. Instead of requiring the lifting against all epimorphisms, being a reflexive projective requires only the lifting against reflexive epimorphisms, that is morphisms that arise as coequalizers of reflexive pairs! By remark \ref{remark: reflexive pairs from reflexive relations} any quotient of an equivalence relation fits that description. Hence, being a reflexive projective allows to lift against \emph{gluing maps}.
	
	If the ambient category is a pretopos, then every epimorphism $q\colon X\to Y$ is already the quotient of the equivalence relation $X\times_Y X\to X\times X$, so is a reflexive epimorphism. Thus, in a pretopos reflexive projectivity and projectivity agree! In particular, we may simply speak of \emph{projectivity} throughout this thesis.
\end{remark}

\section{Condensed Sets}
\label{section: condensed sets}

Having now introduced the necessary sheaf- and (pre)topos-theoretic language, we are now well-equipped to define condensed sets.

Recall that the main application of condensed mathematics is to do algebra when the underlying sets of the algebraic objects in question carry a topology. Usually this leads to categories with poor properties, \eg the category $\TopAb$ of topological abelian groups that is not even abelian. Clausen and Scholze however provide a solution via their theory of condensed mathematics. In a Yoneda-like manner, they associate to every topological space $X$ a restriction of the functor $\Hom_\Top(-, X)$, yielding a certain sheaf of sets on the site $\proetalesite{\ast}$. If now $X$ carries additional algebraic structure compatible with the topology, \eg is a topological abelian group, then for every $S\in \proetalesite{\ast}$ the set of sections $\Hom_\Top(S, X)$ naturally carries the same kind of algebraic structure. This yields for example a functor
\begin{align*}
	\TopAb &\longrightarrow \Sheaves{\proetalesite{\ast}}{\Ab},\\
	X&\longmapsto \eHom_\Top(-, X),
\end{align*}
that is fully faithful on a large class of spaces. As $\Sheaves{\proetalesite{\ast}}{\Ab}$ is abelian this allows us to apply many important tools, like homological algebra, to topological abelian groups -- something that is not as easy for $\TopAb$ on its own. Hence, the solution to the above problem essential comes from enlarging the concept of what one calls a topological space.

This new concept of space is what Clausen and Scholze call a \emph{condensed set}. Ignoring set-theoretic issues for now, they define condensed sets as sheaves of sets on $\proetalesite{\ast}$, the \emph{pro-étale site of a point}. As mentioned in the introduction, the pro-étale site (on a general scheme) was originally introduced by Bhatt and Scholze in \cite{bhatt2014proetaletopologyschemes} to give a new definition of the derived category of constructible $\ell$-adic sheaves. In this text however, only the pro-étale site of a point will be of interest. It turns out that this site (with a rather technical definition) can equivalently be given by the site of \emph{profinite sets} with \emph{finite, jointly surjective} families of continuous maps as covers. Let us define what a profinite set is.
\begin{definition}[Total disconnectedness \& (the site of) profinite sets]
	A topological space $X$ is called \emph{totally disconnected} if every connected component of $X$ is a point.
	
	A totally disconnected, compact Hausdorff space is called a \emph{Stone space}, a \emph{profinite space} or a \emph{profinite set}. In this thesis we will call them profinite sets.
	
	We write $\proFiniteSets$ for the full subcategory of $\Top$ of profinite sets. This category admits the structure of a site by declaring a cover of some profinite set $S$ to be a finite family of jointly surjective maps, that is a family $(f_i\colon S_i\to S)_{i\in I}$ of continuous maps $f_i$, where $I$ is finite and $\bigcup_{i\in I} f_i(S_i) = S$ \ie the images of the $f_i$ jointly cover $S$.
\end{definition}

\begin{remark}[The different names of profinite sets]
	Technically Stone spaces, profinite spaces and profinite sets all form different (but equivalent) categories. As the name hints at, profinite sets are (formal) projective limits of finite sets, so the objects of the pro-category $\ProCategory{\set}$ (see construction \ref{Construction: Pro-categories} in the appendix). Similarly, profinite spaces are given by the pro-category on the full subcategory of $\Top$ of all finite discrete spaces (which is of course equivalent to $\set$).
	
	By \eg \cite[\href{https://stacks.math.columbia.edu/tag/08ZY}{Tag 08ZY}]{stacks-project} any totally disconnected compact Hausdorff space can be written as a projective limit over finite discrete spaces, and conversely any projective limit of finite discrete spaces is totally disconnected and compact Hausdorff. Thus, Stone spaces identify with pro-discrete spaces and thus all three categories agree.
\end{remark}

\begin{example}[(Non-)examples of profinite spaces]
	We have the following (non-)examples of profinite sets.
	\begin{itemize}
		\item 
		A discrete space is totally disconnected if and only if it is finite. Indeed, an infinite discrete space cannot be compact.
		
		\item 
		Arbitrary direct products $\prod_I X$ of some finite discrete space $X$. Notice that this space is not discrete (\eg singletons are not open), but that it is compact by Tychonoff's theorem.
		
		\item
		The Cantor set.
		
		\item
		The \emph{sequence space} $\hat{\N} = \N\cup \{\infty\}$ obtained as the one-point-compactification of the natural numbers $\N$.
		
		\item
		The $p$-adic integers $\Z_p$ or more generally the power series ring $\series{\Z}{T}$ with its $\ideal{T}$-adic topology.
		
		\item
		The field $\adicsQ{p}$ of $p$-adic numbers is not profinite. While it is totally disconnected and Hausdorff, it is not compact.
	\end{itemize}
\end{example}

One problem remains: The site of profinite sets is (essentially) large. As mentioned in the previous sections, we will thus have to restrict to a certain subclass of \emph{regular} sheaves to obtain a locally small category of condensed sets. After all, this is why we introduced the notion of pretopoi.

This \emph{regularity} will be a certain \emph{smallness} condition (see definition \ref{definition: small (pre)sheaves}). We only like to consider sheaves on the large site, that are essentially determined by a small amount of data. To get there, we will impose size restrictions -- not on the sheaves, but on the profinite sets making up the defining site. Using suitably chosen cardinals $\kappa$, we thus obtain (increasingly \emph{finer}) categories of so-called $\kappa$-condensed sets, when increasing $\kappa$. These categories of $\kappa$-condensed sets can then be glued suitably, to obtain the desired category that is independent of any single chosen cardinal.
\begin{convention}
	Throughout this thesis let $\kappa$ denote an uncountable strong limit cardinal. That is $\kappa$ is uncountable and for any cardinal $\lambda < \kappa$ also $2^\lambda < \kappa$. By construction \ref{construction: Beth-numbers and existence of strong limit cardinals} they are ubiquitous. A set is called \emph{$\kappa$-small} if its cardinality is strictly smaller than $\kappa$. For an overview of the set-theoretic conventions used in this thesis, please refer to chapter \ref{chapter: cardinals} of the appendix.
\end{convention}

Using such cardinals as upper bounds, we approximate the large site from below.
\begin{definition}[The site of $\kappa$-small profinite spaces]
	The site $\boundedProFiniteSets{\kappa}$ is the full subcategory of $\Top$ consisting of all $\kappa$-small profinite space with covers given by finite, jointly surjective families of continuous maps.
\end{definition}

Using this restricted site we can now define $\kappa$-condensed sets.
\begin{definition}[The topos of $\kappa$-condensed sets]
	\label{definition: the topos of kappa-condensed sets}
	
	A \emph{$\kappa$-condensed set} is a $\Set$-valued sheaf on the site $\boundedProFiniteSets{\kappa}$ of $\kappa$-small profinite sets. The topos
		$$\boundedCondensedSets{\kappa}\ldef \Sheaves {\boundedProFiniteSets\kappa}{\Set}$$
	is called the \emph{topos of $\kappa$-condensed sets}.
\end{definition}

Due to the following lemma it is in practice not too hard to check if a presheaf on the site $\boundedProFiniteSets{\kappa}$ is actually a sheaf.
\begin{lemma}[Simplified sheaf conditions for $\boundedProFiniteSets{\kappa}$]
	\label{lemma: simplified sheaf conditions for profinite sets}
	
	A presheaf $T\colon \boundedProFiniteSets{\kappa}^\op \to \Set$ is a sheaf if and only if
	\begin{enumerate}[(i)]
		\item
		\label{lemma: simplified sheaf conditions for profinite sets, i}
		For any finite disjoint union $\coprod_{i\in I} S_i$ the natural map $T(\coprod_{i\in I} S_i)\to \prod_{i\in I} T(S_i)$ is a bijection. That is $T$ maps finite coproducts to products
		
		\item
		\label{lemma: simplified sheaf conditions for profinite sets, ii}
		For any surjection $S'\to S$ with projections $p_1, p_2\colon S'\times_S S'\rightrightarrows S'$ the natural map
			$$T(S) \to \{x\in T(S')\setseparator p_1^\ast(x)=p_2^\ast(x) \in T(S'\times_S S')\}=\equalizer(p_1^\ast, p_2^\ast)$$
		is a bijection. That is $T$ maps coequalizers of the form $S'\times_S S'\rightrightarrows S'\to S$ to equalizers.
	\end{enumerate}
	
	These conditions apply verbatim to the full site $\proFiniteSets$.
\end{lemma}
\begin{proof}
	We will prove the equivalence for unrestricted profinite sets. A quick inspection of the following proof shows, that the arguments are still valid for $\boundedProFiniteSets{\kappa}$.
	
	If $T$ is a sheaf, then in both cases $T$ respects the posited colimits. Indeed, the disjoint union has a cover $(S_j\to \coprod_{i\in I} S_i)_{j\in I}$ where all fiber products of $S_i$ and $S_j$ with $i\ne j$ are empty. The map of \ref{lemma: simplified sheaf conditions for profinite sets, i} is thus precisely the natural map from the equalizer condition in lemma \ref{lemma: being a sheaf is an equalizer condition}, hence is a bijection as $T$ is a sheaf. Similarly, the surjection $S'\to S$ itself constitutes a cover of $S$ and the map in \ref{lemma: simplified sheaf conditions for profinite sets, ii} is precisely the one in the equalizer condition for being a sheaf, so a bijection.
		
	It remains to show the converse. So let $T$ be a sheaf satisfying \ref{lemma: simplified sheaf conditions for profinite sets, i} and \ref{lemma: simplified sheaf conditions for profinite sets, ii}. Again we use the equalizer condition. For any cover $(S_i \to S)_{i\in I}$ of $S$ we have a coequalizer $\coprod_{i,j\in I}S_i\times_S S_j \rightrightarrows \coprod_{i\in I} S_i \to S$. Thus, $T$ (since it respects coproducts) induces a diagram $T(S)\to \prod_{i\in I}T(S_i)\rightrightarrows\prod_{i, j \in I} T(S_i\times_S S_j)$. Now consider the natural map $S'\coloneqq \coprod_i S_i \to S$ induced by the cover. Its associated coequalizer $S'\times_S S'\rightrightarrows S'\to S$ is isomorphic to the previous coequalizer (in the sense that there is a compatible isomorphism $\coprod_{i,j\in I} S_i\times_S S_j\to S'\times_S S'$). In particular a presheaf $T$ respecting this latter coequalizer must also respect the former. But by assumption $T$ respects the latter coequalizer, hence the former and thus the induced $T(S)\to \prod_{i\in I}T(S_i)\rightrightarrows\prod_{i, j \in I} T(S_i\times_S S_j)$	is an equalizer. In particular the induced map
		$$T(S)\to \equalizer\left(\prod_{i\in I}T(S_i)\rightrightarrows\prod_{i, j \in I} T(S_i\times_S S_j)\right)$$
	is an equalizer as well and $T$ satisfies the sheaf property for the given cover. Since the cover was arbitrary, $T$ is a sheaf.
\end{proof}

\begin{definition}[Associated $\kappa$-condensed sets]
	\label{definition: associated kappa-condensed sets}
	
	There is a functor
	\begin{align*}
		\Top &\longrightarrow \boundedCondensedSets{\kappa},\\
		X&\longmapsto \associated{X}\ldef \Hom_\Top(-, X),
	\end{align*}
	that assigns to each topological space $X$ its \emph{associated $\kappa$-condensed set} $\associated{X}$. Indeed, $\Hom_\Top(-, X)$ respects all colimits so in particular satisfies the simplified sheaf conditions of \ref{lemma: simplified sheaf conditions for profinite sets}.
\end{definition}

When is this functor $X\mapsto \associated{X}$ fully faithful? We can only ever hope that it faithfully represent a topological space $X$, when the topology of $X$ is generated by continuous images of $\kappa$-small profinite sets. If this is too strict of a condition on a given topology (and at a first glance it really might seem so, as profinite sets really are not very geometric in nature) then we are out of luck. It turns out however, that the class of topological spaces generated by such continuous images is actually the well-known class of \emph{($\kappa$-)compactly generated spaces}!
\begin{definition}[Compactly generation \& k-ification]
	A topological space $(X, \topology T)$ is called \emph{$\kappa$-compactly generated} if $\topology T$ is the final topology (see definition \ref{definition: final topologies}) induced by all continuous maps $K \to X$ from $\kappa$-small compact Hausdorff spaces $K$. That is a $U\subseteq X$ is open if and only if $f^{-1}(U)$ is open in $K$ for any such $f\colon K\to X$. The topological space $(X,\topology T)$ is called \emph{compactly generated} if it is $\kappa$-compactly generated for any $\kappa$.
 
	A compactly generated space is also called a \emph{$k$-space}, hence we will denote the full subcategory of $\Top$ of compactly generated spaces by $\kSpaces$. The full subcategory of $\Top$ of $\kappa$-compactly generated spaces is denoted by $\boundedkSpaces{\kappa}$.
	
	For an exposition on compact generation see 
	\cite[Chap. 3.1]{compact-generation}. The author was however unable to locate an analogous reference that covers $\kappa$-compact generation.
\end{definition}

\begin{remark}[Alternative characterization of ($\kappa$-)compact generation]
	If one were to ignore set theoretic issues, then a topological space $(X, \topology T)$ is $\kappa$-compactly generated if and only if $\topology{T}$ is the quotient topology of the surjective map
		$$\smashoperator{\coprod_{\substack{K\to X\\ S\in\boundedCHaus{\kappa}}}} K \to X.$$
\end{remark}

\begin{example}[First-countable spaces are compactly generated]
	Any first-countable topological space, \ie a space in which every point admits a countable neighborhood basis, is compactly generated. This is shown in \cite[Remark 1.6]{condenseddotpdf}. This class includes all metrizable topological spaces, all locally Euclidean spaces (\eg manifolds) and all analytic spaces over $\R$ and $\C$.
\end{example}

There is a canonical way to turn any space into a $\kappa$-compactly generated one.
\begin{lemma}[k-ification]
	The inclusions $\boundedkSpaces{\kappa}\to \Top$ admits a right-adjoint $(-)^\boundedkificationsymbol{\kappa}\colon \Top\to\boundedkSpaces{\kappa}$ called the \emph{k-ification} (omitting the $\kappa$ from the name for simplicity). For a topological space $X$, the k-ification $X^\boundedkificationsymbol{\kappa}$ has the same underlying set as $X$ and carries the final topology induced by all continuous maps $K\to X$ from a $\kappa$-small compact Hausdorff spaces $K$. 
	
	We obtain an analogous result for the inclusion $\kSpaces\to\Top$. We will call its right-adjoint $(-)^\kificationsymbol\colon \Top\to\kSpaces$ the k-ification as well. It is usually clear from context which k-ification is meant.
	
	
\end{lemma}
\begin{proof}
	In case of $\kSpaces\to \Top$ this is \cite[Prop. 3.2]{compact-generation}. For the case of $\boundedkSpaces{\kappa}\to \Top$ the author could not locate a suitable reference -- the reasoning however should be analogous to the first case as $\kappa$ is strong limit.
\end{proof}

We will now work towards the fact that compactly generated spaces, \ie spaces whose topology is determined by continuous images of compact Hausdorff spaces, are the same as spaces whose topology is determined by continuous images of profinite sets, so \emph{totally disconnected} compact Hausdorff spaces -- this will be shown in proposition \ref{proposition: equivalent characterizations of compact generation}. This class of spaces is then well-represented in condensed sets by proposition \ref{proposition: adjointness and faithfulness}. Both of these facts hinge on the \emph{Stone--\v Cech compactification}.
\begin{proposition}[Stone--\v Cech compactification, {\cite[Thm 38.2 \& Thm. 38.4]{munkres}}]
	\label{proposition: stone-cech compactification}
	
	The inclusion $\CHaus \to \Top$ admits a left adjoint $\stonecechsymbol\colon \Top \to \CHaus$, called the \emph{Stone--\v Cech compactification}. Explicitly, the Stone--\v Cech compactification of a topological space $X$ yields a compact Hausdorff space $\stonecech X$ together with a map $X \to \stonecech X$ such that there is exactly one dotted map completing each solid diagram
	$$\begin{tikzcd}
		X 	&\stonecech X\\
		&K,		
		\arrow[from=1-1, to=1-2]
		\arrow[from=1-2, to=2-2, "\exists!", dotted]
		\arrow[from=1-1, to=2-2, "\forall", swap]
	\end{tikzcd}$$
	for all maps $X\to K$ where $K$ is a compact Hausdorff space.
	
	If $X$ is a discrete topological space, then the Stone--\v Cech compactification of $X$ is a profinite set and is given by the topological space of all ultrafilters on $X$ (see \ref{definition: ultrafilters, principal ultrafilters and the stone-topology}).
\end{proposition}

These Stone--\v Cech compactifications of discrete spaces are a very special type of profinite set and are of vital importance for condensed mathematics. They are \emph{extremally disconnected}.
\begin{definition}[Extremal disconnectedness]
	A topological space $X$ is called \emph{extremally disconnected} if the closure of any open set in $X$ is again open.
\end{definition}

Now indeed, Stone--\v Cech compactifications of discrete spaces are \emph{extremally disconnected}.
\begin{lemma}[{\cite[\href{https://stacks.math.columbia.edu/tag/090C}{Tag 090C}]{stacks-project}}]
	
	If $X$ is a discrete space, then $\stonecech{X}$ is an extremally disconnected profinite set.
\end{lemma}

\begin{convention}
	In this thesis we will never have to work with extremally disconnected spaces that are not also profinite. Hence, from now on we will call such spaces simply \emph{extremally disconnected sets}. Additionally, we will use the adjective \quote{extremally disconnected} to always mean \quote{extremally disconnected and profinite}. Using this convention, we denote by $\extremallyDisconnectedSets$ the full subcategory of $\Top$ of extremally disconnected sets.
\end{convention}

\begin{remark}[About extremal disconnectedness]
	Extremally disconnected sets are rather strange and unintuitive. Clausen and Scholze state in \cite[Warning 2.6]{condenseddotpdf} the following properties:
	\begin{itemize}
		\item
		Any extremally disconnected set is a retract of a Stone--\v Cech compactification of a discrete space.
		
		\item
		The product of two infinite extremally disconnected sets is never again extremally disconnected.
		
		\item
		Any converging sequence of points in an extremally disconnected set is already eventually constant.
	\end{itemize}
	Furthermore, since the points of $\stonecech{X}$ for some discrete space $X$ are by construction ultrafilters on $X$ (with the points of $X$ corresponding to principal ultrafilters), we find that in Zermelo--Fraenkel without some additional axiom (like the axiom of choice) one cannot even be sure that there is any point in $\stonecech{X} \setminus X$, \confer \cite[\href{https://mathoverflow.net/q/15872/546808}{this MathOverflow post}]{non-principal-ultrafilters}.
\end{remark}

We can now state a very important corollary of \ref{proposition: stone-cech compactification}.
\begin{corollary}
	\label{corollary: compact hausdorff spaces are quotients of extremally disconnected spaces}
	
	Any compact Hausdorff space $X$ is a quotient of an extremally disconnected set of cardinality $2^{2^{\vert X\vert}}$.
\end{corollary}
\begin{proof}
	Consider the discrete space $X^\delta$ with the same underlying set as $X$. The universal property of the Stone--\v Cech compactification $\stonecech X^\delta$ applied to the map $\id\colon X^\delta \to X$ yields a surjective map $\stonecech{X^\delta} \to X$ of compact Hausdorff spaces. By proposition \ref{proposition: stone-cech compactification} $\stonecech{X^\delta}$ is extremally disconnected compact Hausdorff and can be constructed as the space of all ultrafilters on $S$, hence $\vert \stonecech{X^\delta}\vert = 2^{2^{\vert X\vert}}$.
	
	It remains to show that $\pi\colon \stonecech{X^\delta} \to X$ is a quotient map. Consider the relation $\sim$ on $\stonecech{X^\delta}$ given by $s\sim s'$ if and only if $\pi(s)=\pi(s')$. The induced map $\stonecech{X^\delta}/{\sim} \to X$ is clearly bijective. As a continuous image of the compact space $\stonecech{X^\delta}$, the quotient $\stonecech{X^\delta}/{\sim}$ is compact. Thus, the induced map is a bijective map from a compact space into a Hausdorff space, so even a homeomorphism. Hence, $\pi$ was a quotient map.
\end{proof}

Now what makes these spaces so important? Due to Gleason in \cite{gleason} we have the following characterization.
\begin{theorem}[Extremally disconnected spaces are projective]
	\label{theorem: extremally disconnected spaces are projective}
	
	The projective objects of $\CHaus$ are precisely the extremally disconnected sets.
\end{theorem}
\begin{shortproof}
	This is proven in \cite[\href{https://stacks.math.columbia.edu/tag/08YN}{Tag 08YN}]{stacks-project} but is based on Gleason's \cite[Thm. 2.5]{gleason}.
\end{shortproof}

\begin{remark}[Enough projectives]
	By \ref{corollary: compact hausdorff spaces are quotients of extremally disconnected spaces} any compact Hausdorff space is a quotient of an extremally disconnected set. Even more, due to Gleason, these spaces are projective. So in essence $\CHaus$ \emph{has enough projectives}, that is any compact Hausdorff space admits an epimorphism from (even more: is a quotient of) a projective object.
	
	When working with condensed abelian groups later on, these extremally disconnected sets will give rise to the projective objects -- then analogously, the category of condensed abelian groups will have enough projectives -- in the sense of homological algebra!
\end{remark}

We are now equipped to show the alternative characterization of compact generation.
\begin{proposition}[Equivalent characterizations of compact generation]
	\label{proposition: equivalent characterizations of compact generation}
	
	For a topological space $(X, \topology T)$ the following are all equivalent:
	\begin{enumerate}[(i)]
		\item
		$(X, \topology T)$ is $\kappa$-compactly generated, \ie $\topology T$ is the final topology induced by all continuous maps $S\to X$ where $S$ is $\kappa$-small compact Hausdorff,
		
		\item
		$\topology T$ is the final topology induced by all continuous maps $S\to X$ where $S$ is $\kappa$-small totally disconnected compact Hausdorff,
		
		\item
		$\topology T$ is the final topology induced by all continuous maps $S\to X$ where $S$ is $\kappa$-small extremally disconnected compact Hausdorff.
	\end{enumerate}
\end{proposition}

\begin{proof}
	Let $\topology T_1$, $\topology T_2$ and $\topology T_3$ be the three final topologies. Clearly $\topology T_1\subseteq \topology T_2 \subseteq \topology T_3$ since the condition for a subset of $X$ to be open gets weaker the less maps one considers. It remains to show that $\topology T_3 \subseteq \topology T_1$. To this end let $U\in \topology T_3$. Consider any $f\colon S\to X$ where $S$ is a $\kappa$-small compact Hausdorff space. By corollary \ref{corollary: compact hausdorff spaces are quotients of extremally disconnected spaces} there is a quotient map $\pi\colon S'\to S$ where $S'$ is an extremally disconnected compact Hausdorff space. By definition of $U$ we know that the preimage $\pi^{-1}(f^{-1}(U))$ under the function $f\circ\pi\colon S'\to X$ is open. Hence, also $f^{-1}(U)$ is open as $\pi$ is a quotient map. Thus, $U\in \topology T_1$.
\end{proof}

These alternative characterizations of compact generation now allow us to study which topological spaces are well-represented by their associated $\kappa$-condensed sets.
\begin{remark}[Recovering a space from its associated condensed set]
	\label{remark: recovering a space from its associated condensed set}
	
	Suppose $X$ is a topological space and $\associated{X}$ is its associated $\kappa$-condensed set. Let us try to reconstruct $X$ from $\associated{X}$. Up to bijection, one can recover the underlying set of $X$ as the set of sections $\associated{X}(\ast)=\Hom_\Top(\ast, X)$ on a point. Similarly, for any $\kappa$-small profinite set $S$ one recovers the continuous functions $S\to X$ as the sections in $\associated{X}(S)=\Hom_\Top(S, X)$, which by Yoneda's lemma are precisely the morphisms of $\kappa$-condensed sets $\associated S\to \associated{X}$. Using this information one can construct a natural topology on the underlying set $\associated{X}(\ast)$, namely the final topology of all maps of sets $S\shortisorightarrow \associated{S}(\ast)\to \associated{X}(\ast)$ ranging over all morphisms $\associated{S}\to \associated{X}$ from some $\kappa$-small profinite sets $S$. Denote the so obtained space by $\topologizedUnderlying{\associated{X}}$.
	
	Now if $X$ was $\kappa$-compactly generated, then by proposition \ref{proposition: equivalent characterizations of compact generation} we have $\topologizedUnderlying{\associated{X}}\cong X$; if not, then $\topologizedUnderlying{\associated{X}}$ is the k-ification $X^\boundedkificationsymbol{\kappa}$ of $X$ instead.
\end{remark}
	
We can generalize this construction to any $\kappa$-condensed set.
\begin{definition}[The underlying space of a $\kappa$-condensed set]
	\label{definition: the underlying space of a kappa-condensed set}
	
	Let $T$ be a $\kappa$-condensed set. Then $T(\ast)$ is called the \emph{underlying set of $T$}. For any profinite set $S$, the set of sections $T(S)$ is called the \emph{set of functions on $S$}. Define a topology on the underlying set $T(\ast)$ as the final topology of all maps of sets $S\shortisorightarrow \associated{S}(\ast)\to T(\ast)$ induced by morphisms $\associated{S}\to T$ ranging over all $\kappa$-small profinite sets $S$. Denote the resulting space by $\topologizedUnderlying{T}$ and call it the \emph{underlying space} of $T$.
\end{definition}

\begin{remark}[Alternative description of the underlying topology]
	If one were to ignore set theoretic problems then the topology of the underlying space in definition \ref{definition: the underlying space of a kappa-condensed set} is equivalently the quotient topology of the surjective map
		$$\smashoperator{\coprod_{\substack{\associated{S}\to T\\S\in\boundedProFiniteSets{\kappa}}}} S\to T(\ast).$$
\end{remark}

\begin{remark}[A shift in perspective]
	\label{remark: a shift in perspective}
	
	Remark \ref{remark: recovering a space from its associated condensed set} and definition \ref{definition: the underlying space of a kappa-condensed set} lead to a shift in perspective. Instead of being a space itself, a ($\kappa$-)condensed set is rather a collection of functions into a \emph{virtual} space. Indeed, there are many non-trivial condensed sets with trivial underlying set. Consider \eg the presheaf (of abelian groups) $Q = \associated{\R}/\associated{\R^\delta}$ that assigns to each profinite $S$ the quotient of continuous function $S\to \R$ by locally constant functions $S\to \R$. One can argue that this is indeed a sheaf, but we will not do so. We will see this example again in remark \ref{remark: the relief} later on. Now clearly $Q$ is not trivial, but nevertheless $Q(\ast)$ is a single point. Much like for schemes, a condensed set carries additional \emph{infinitesimal} information about how functions should behave in neighborhoods of the underlying set -- even if the underlying set is merely a single point. In this sense, a condensed set is very similar to the structure sheaf of a scheme, \eg of the non-reduced $\Spec k[x]/\ideal{x^2}$ for some field $k$.
	
	This example was by no means chosen arbitrarily and is exactly the \emph{missing} (from $\TopAb$) quotient $\associated{\R^\delta}\to\associated{\R}\to Q\to 0$ mentioned in the introduction of this thesis, which is trivial in the category of topological abelian groups.
\end{remark}

\begin{proposition}[Adjointness \& Faithfulness]
	\label{proposition: adjointness and faithfulness}
	
	The functor $\Top\to\boundedCondensedSets{\kappa}, X\mapsto\associated{X}$ is faithful and fully faithful when restricted to $\kappa$-compactly generated spaces.
	
	The functor $\boundedCondensedSets{\kappa} \to \Top, T\mapsto \topologizedUnderlying{T}$ is left adjoint to $X\mapsto \associated{X}$ and the counit $\topologizedUnderlying{\associated{X}}\to X$ of the adjunction agrees with the counit $X^\boundedkificationsymbol{\kappa}\to X$ of the adjunction between $\kappa$-compactly generated spaces and all topological spaces. In particular, $\topologizedUnderlying{\associated{X}}\cong X^\boundedkificationsymbol{\kappa}$ is the k-ification of $X$.
\end{proposition}
\begin{proof}
	The claim formally reduces to proving the adjunction and observing that the counit is componentwise epic for faithfulness and componentwise an isomorphism for full faithfulness, be the general theory of adjunctions.
	
	To this end, consider the natural morphism $\Hom(T, \associated{X}) \to \Hom(\topologizedUnderlying{T}, X)$ given by mapping $g\colon T\to \associated{X}$ to $\topologizedUnderlying{T}\xrightarrow{g_\ast} \topologizedUnderlying{\associated{X}} \xrightarrow{\epsilon_X\ldef\id} X$, where $g_\ast$ is the component of the natural transformation $g$ at the object $\ast\in\boundedProFiniteSets{\kappa}$. We claim that this is an isomorphism.
	
	First we will show injectivity. Suppose we have $f,g\colon T\to\associated{X}$ such that their images agree. As $\id\colon \topologizedUnderlying{\associated{X}}\to X$ is a bijection, we must have $f_\ast=g_\ast$. Now let $S$ profinite and $s\in S$ be arbitrary. Consider the unique isomorphism $\ast\to\{s\}$. Then we obtain a diagram
	$$\begin{tikzcd}
		T(\{s\})\ar[r, "f_{\{s\}}", shift left=.25em]\ar[r, "g_{\{s\}}", shift right=.25em, swap] \ar[d, "\sim" {rotate=90, anchor=north}]&\associated{X}(\{s\})\ar[d, "\sim" {rotate=90, anchor=north}]\\
		T(\ast)\ar[r, "f_\ast=g_\ast"] &\associated{X}(\ast)
	\end{tikzcd}$$
	where each of the two embedded squares commute by naturality of $f$ and $g$ respectively. Then clearly $f_{\{s\}}=g_{\{s\}}$. 
	
	The maps $T(S)\to T(\{s\})$ and $\associated{X}(S)\to \associated{X}(\{s\})$ for the various $s\in S$ induce morphisms $T(S)\to\prod_{s\in S} T(\{s\})$ and $\associated{X}(S)\to\prod_{s\in S} \associated{X}(\{s\})$. The latter morphism corresponds to the inclusion of continuous maps $S\to X$ into all maps $S\to X$ and thus is injective. We obtain the following diagram
	$$\begin{tikzcd}
		T(S) \ar[r]\ar[d, "f_S", shift right=.25em, swap]\ar[d, "g_S", shift left=.25em] &\prod_{s\in S} T(\{s\}) \ar[d]\\
		\associated{X}(S)\ar[r, hook] &\prod_{s\in S}\associated{X}(\{s\})
	\end{tikzcd}$$
	where the right vertical morphism is $\prod_{s\in S} f_{\{s\}} = \prod_{s\in S} g_{\{s\}}$ and each of the two embedded squares commutes. In particular by injectivity of the bottom horizontal map we obtain $f_S=g_S$. As $S$ was arbitrary, this shows $f=g$ and hence injectivity.
	
	We will now show surjectivity. So suppose a map $f\colon \topologizedUnderlying{T}\to X$ is given. Consider the morphism $\eta_T\colon \associated{\topologizedUnderlying{T}} \to T$ given componentwise by
		$$\eta_{T, S}\colon T(S)\cong \Hom(\associated{S}, T)\to \Hom(\topologizedUnderlying{\associated{S}}, \topologizedUnderlying{T})\cong \Hom(S, \topologizedUnderlying{T})=\associated{\topologizedUnderlying{T}}(S)$$
	for $S$ profinite. Define the composite morphism $\hat f\colon T\xrightarrow{\eta_T} \associated{\topologizedUnderlying{T}}\longxrightarrow{\associated{f}} \associated{X}$. Then by chasing along the commutative diagram
	$$\begin{tikzcd}
		\topologizedUnderlying{T} \ar[r, "\eta_{T,\ast}"] \ar[dr, "\id"] &\topologizedUnderlying{\mleft(\associated{\topologizedUnderlying{T}}\mright)} \ar[r, "\associated{f}_\ast"] \ar[d, "\epsilon_{\topologizedUnderlying{T}}"] &\topologizedUnderlying{\associated{X}} \ar[d, "\epsilon_X"]\\
		&\topologizedUnderlying{T} \ar[r, "f"] &X
	\end{tikzcd}$$
	we obtain that $\epsilon_X\circ \hat{f}_\ast = f\circ \id = f$, \ie $\hat{f}$ is a preimage of $f$.
	
	This shows the claimed adjunction. One easily checks that $\eta$ is the unit and $\epsilon$ is the counit. Now indeed for any $X$ we have $\epsilon_X\colon \topologizedUnderlying{\associated{X}}\to X$ which agrees with the counit $X^\boundedkificationsymbol{\kappa}\to X$. In particular, it is always epic (even more it is always a bimorphism) and even an isomorphism if $X$ is $\kappa$-compactly generated, proving the claimed faithfulness and full faithfulness respectively.
\end{proof}

\begin{notation}[Simplifying notation]
	Let $X$ be a $\kappa$-compactly generated space. As a consequence of proposition \ref{proposition: adjointness and faithfulness} we will, if no confusion is possible, simply write $X$ for the associated $\kappa$-condensed set $\associated{X}$.
\end{notation}

Through trying to understand how topological spaces map to condensed sets, we have now seen three important (types of) classes of spaces. The classes of ($\kappa$-small) compact Hausdorff spaces, profinite sets and extremally disconnected sets -- each class being more restrictive than the last. In proposition \ref{proposition: equivalent characterizations of compact generation} we have shown that each of these classes leads to the same notion of topology under generation. Similarly, we will now see that the site of $\kappa$-small profinite sets is only one of several different sites that can be used to define $\kappa$-condensed sets. More precisely, the site $\boundedProFiniteSets{\kappa}$ in definition \ref{definition: the topos of kappa-condensed sets} can be replaced by either the site $\boundedCHaus{\kappa}$ of $\kappa$-small compact Hausdorff spaces or by the site $\boundedExtremallyDisconnectedSets{\kappa}$ of $\kappa$-small extremally disconnected sets, where in both cases covers are given again by finite, jointly surjective families of maps. These equivalences are the content of the following proposition.
\begin{proposition}[Independence to the defining site]
	\label{proposition: independence to the defining site}
	
	Restriction and extension of sheaves define geometric equivalences of topoi
		$$\Sheaves{\boundedCHaus{\kappa}}{\Set} \xleftrightarrows{}{} \Sheaves{\boundedProFiniteSets{\kappa}}{\Set}\xleftrightarrows{}{} \Sheaves{\boundedExtremallyDisconnectedSets{\kappa}}{\Set}.$$
\end{proposition}
\begin{proof}
	This is a straightforward application of the comparison lemma \ref{lemma: comparison lemma}. All three sites involved are (essentially) small and we will show that the two sites $\boundedProFiniteSets{\kappa}$ and $\boundedExtremallyDisconnectedSets{\kappa}$ are dense sub-sites of $\boundedCHaus{\kappa}$ and $\boundedProFiniteSets{\kappa}$ respectively. Now the coverages of all three sites are compatible and all inclusions are full, yielding full sub-sites. It remains to show density. For this, let $S\in\boundedCHaus{\kappa}$ or $S\in \boundedProFiniteSets{\kappa}$ be arbitrary. By corollary \ref{corollary: compact hausdorff spaces are quotients of extremally disconnected spaces} there exists a surjection $S'\to S$ from some $S' \in \boundedExtremallyDisconnectedSets{\kappa}$. Hence, the singleton family given by $S'\to S$ constitutes a cover in $\boundedCHaus{\kappa}$ and $\boundedProFiniteSets{\kappa}$ respectively. This shows density for both inclusions.
\end{proof}
	
While extremally disconnected sets are quite unintuitive in nature, at least they make life a little easier when trying to verify that a presheaf is a sheaf. Akin to lemma \ref{lemma: simplified sheaf conditions for profinite sets} we obtain an equivalent characterization for being a sheaf on the site $\boundedExtremallyDisconnectedSets{\kappa}$.
\begin{lemma}[Simplified sheaf condition for $\boundedExtremallyDisconnectedSets{\kappa}$]
	\label{lemma: simplified sheaf condition for extremally disconnected sets}
	
	A presheaf $T\colon \boundedExtremallyDisconnectedSets{\kappa}^\op \to \Set$ is a sheaf if and only if for any finite disjoint union $\coprod_{i\in I} S_i$, the natural map $T(\coprod_{i\in I} S_i)\to \prod_{i\in I} T(S_i)$ is a bijection. That is $T$ maps finite coproducts to products.
	
	This condition applies verbatim to the full site $\extremallyDisconnectedSets$.
\end{lemma}
\begin{proof}
	We will prove the equivalence for unrestricted extremally disconnected sets. A quick inspection of the following proof shows, that the arguments are still valid for $\boundedExtremallyDisconnectedSets{\kappa}$.
	
	As in the proof of \ref{lemma: simplified sheaf conditions for profinite sets}, if $T$ is a sheaf the simplified sheaf condition holds.
	
	Now suppose that $(f_i\colon S_i\to S)_{i\in I}$ is a finite cover of an extremally disconnected set $S$ by other extremally disconnected sets $S_i$. Unfortunately, limits of extremally disconnected sets need not be extremally disconnected, so we cannot use lemma \ref{lemma: being a sheaf is an equalizer condition} which restates the sheaf condition for the given cover as an equalizer. There is however a modification. Due to corollary \ref{corollary: compact hausdorff spaces are quotients of extremally disconnected spaces} for every pair $i, j\in I$ there is a extremally disconnected set $E_{i,j}$ along with a surjection $E_{i, j} \twoheadrightarrow S_i\times_S S_j$. Now suppose we have a pair of morphisms $S_i\leftarrow V\rightarrow S_j$ between extremally disconnected sets over $S$. By the universal property of $S_i\times_S S_j$ there is a unique compatible morphism $u\colon V\to S_i \times_S S_j$. Now by theorem \ref{theorem: extremally disconnected spaces are projective}, $V$ is a projective object, so the morphism $u$ lifts against the surjection $E_{i, j} \twoheadrightarrow S_i\times_S S_j$ giving a morphism $V\to E_{i, j}$. Very similar to the proof of \ref{lemma: being a sheaf is an equalizer condition} we thus obtain an alternative equalizer condition. Namely, that $T$ satisfies the sheaf condition for the given cover if and only if the induced
	$$\begin{tikzcd}
		T(S) \ar[r] &\equalizer\Bigg(\prod_{i\in I} T(S_i) \ar[r, shift left=1.1] \ar[r, shift right=1.1] &\prod_{i, j \in I} T(E_{i, j})\Bigg)
	\end{tikzcd}$$
	is an isomorphism. As $T$ maps finite coproducts to products this is equivalently that
	$$\begin{tikzcd}
		T(S) \ar[r] &\equalizer\Bigg(T\mleft(\coprod_{i\in I} S_i\mright) \ar[r, shift left=1.1] \ar[r, shift right=1.1] &T\mleft(\coprod_{i, j \in I} E_{i, j}\mright)\Bigg)
	\end{tikzcd}$$
	is an isomorphism. In particular, this is satisfied if $T$ maps coequalizers to equalizers if the coequalizer is of the form $S''\rightrightarrows S' \twoheadrightarrow S$ with $S''\rightrightarrows S'$ being induced by a surjection $S''\twoheadrightarrow S'\times_S S'$. But as surjections of extremally disconnected sets split by theorem \ref{theorem: extremally disconnected spaces are projective}, this coequalizer is split as well and hence preserved by any functor. Thus, this condition is automatic and $T$ satisfies the sheaf condition for the given cover. As the cover was arbitrary it finally follows that $T$ is a sheaf.
\end{proof}

\begin{remark}[The trade-off]
	Consider again the chain of inclusions $\boundedExtremallyDisconnectedSets{\kappa}\subseteq \boundedProFiniteSets{\kappa}\subseteq \boundedCHaus{\kappa}$ and recall the lemmata \ref{lemma: simplified sheaf conditions for profinite sets} and \ref{lemma: simplified sheaf condition for extremally disconnected sets} which provide simplified sheaf conditions for the sites $\boundedProFiniteSets{\kappa}$ and $\boundedExtremallyDisconnectedSets{\kappa}$ respectively. While $\boundedCHaus{\kappa}$ leads to a very good geometric interpretation of $\kappa$-condensed sets, making it much clearer that we obtain a good embedding of $\kappa$-compactly generated spaces into $\kappa$-condensed sets, it is in general rather hard to check the sheaf conditions for this site. At the other extreme, the spaces of $\boundedExtremallyDisconnectedSets{\kappa}$ are utterly ingeometric (all converging sequences are eventually constant; without the axiom of choice it is essentially impossible to pin down points) but they have the added benefit, that the sheaf condition is almost trivial and that by theorem \ref{theorem: extremally disconnected spaces are projective}, extremally disconnected set are projective in $\CHaus$. In between sits the site $\boundedProFiniteSets{\kappa}$ of profinite spaces. While they still enjoy some of the good properties of extremally disconnected spaces (\eg an easier sheaf condition compared to compact Hausdorff spaces) they are also still somewhat familiar in their geometry. Additionally, profinite spaces (by virtue of their definition) allow breaking arguments down to finite sets. In general, we see that there is a trade-off between these sites and depending on the problem we can choose the one that fits best.
\end{remark}

We will now work towards constructing a general category of condensed sets, independent of some cutoff cardinal $\kappa$. Recall that the idea to realize this is to glue the various categories of $\kappa$-condensed sets. For this it is necessary to ask, \emph{how} one should glue. This is at least partially answered in the following proposition.
\begin{proposition}[Transitioning between different cardinals]
	\label{proposition: transitioning between different cardinals}
	
	Suppose $\kappa' > \kappa$ is a second uncountable strong limit cardinal. There is a functor $(-)^\enlargementsymbol{\kappa}{\kappa'}\colon \boundedCondensedSets{\kappa} \to \boundedCondensedSets{\kappa'}$ given by
		$$T \mapsto \sheafification{\tilde S\mapsto \colimit_{\tilde S\to S} T(S)},$$
	where the colimit is $\cofinality\kappa$-filtered and is taken over all $\kappa$-small profinite sets $S$ with a map $\tilde S\to S$. This functor $(-)^\enlargementsymbol{\kappa}{\kappa'}$ is fully faithful, commutes with all colimits and with all $\cofinality{\kappa}$-small limits.
\end{proposition}

\begin{proof}
	We first show the $\cofinality{\kappa}$-filteredness of the colimit. Fix $T$ a $\kappa$-condensed set and $\tilde S$ a $\kappa'$-small profinite set. The index category of the above colimit is $\tilde S/\boundedProFiniteSets{\kappa}$ and want to show that this colimit is $\cofinality{\kappa}$-filtered. Since $T$ is contravariant functor, we thus have to show that the index category is $\cofinality{\kappa}$-cofiltered. So suppose we have a $\cofinality{\kappa}$-small index category $I$ and a diagram $D\colon I\to \tilde S/\boundedProFiniteSets{\kappa}$ composed of compatible maps $\tilde S\to S_i$. Since $\CHaus$ is complete, the limit $L = \limit D$ in $\tilde S/\CHaus$ exists and by corollary \ref{corollary: compact hausdorff spaces are quotients of extremally disconnected spaces} $L$ admits a surjection $S\to L$ from a (extremally disconnected) profinite set $S$ of cardinality $2^{2^{\vert L\vert}}$. Hence, if $L$ is $\kappa$-small, so is $S$ (as $\kappa$ is a strong limit cardinal) and precomposing the limiting cone of $L$ with this surjection yields a cone of $D$ in $\tilde S/\boundedProFiniteSets{\kappa}$. So it remains to show that $L$ indeed is $\kappa$-small. Let $\lambda \coloneqq \vert I\vert$, so that $\lambda < \cofinality{\kappa}$ as $D$ is a $\cofinality{\kappa}$-small diagram by definition. By the definition of cofinality, the category $\boundedSet{\kappa}$ is $\lambda$-cocomplete. In particular $\mu\coloneqq \sup_{i\in I} \vert S_i\vert = \vert \coprod_{i\in I} S_i\vert < \kappa$ and $\mu\cdot \lambda < \kappa$. Finally,
	$$\vert L\vert \le \left\vert \prod_{i\in I} S_i\right\vert \le \prod_{i\in I} \mu = \mu^\lambda \le 2^{\mu \cdot \lambda} < \kappa$$
	where the last inequality follows since $\mu\cdot \lambda <\kappa$ and $\kappa$ is a strong limit cardinal. This shows the desired filteredness. Notice how the tip of the constructed cone is extremally disconnected, showing in particular that categories of the form $\tilde S/\boundedExtremallyDisconnectedSets{\kappa}$ for some extremally disconnected $\kappa'$-small compact Hausdorff space are $\cofinality{\kappa}$-cofiltered as well.
	
	Consider the (up to natural isomorphism) commutative square
	\begin{center}
		\begin{tikzcd}[row sep=35, column sep=40]
			\Sheaves{\boundedExtremallyDisconnectedSets{\kappa}}\Set &\Sheaves{\boundedExtremallyDisconnectedSets{\kappa'}}\Set\\
			\Sheaves{\boundedProFiniteSets{\kappa}}\Set &\Sheaves{\boundedProFiniteSets{\kappa'}}\Set
			\ar[from=1-1, to=1-2, "i^{-1}"]
			\ar[from=1-1, to=2-1, "\sim" {sloped, near start}]
			\ar[from=1-2, to=2-2, "\sim" {sloped, near start}]
			\ar[from=2-1, to=2-2, "(-)^\enlargementsymbol{\kappa}{\kappa'}", swap]
			\ar[from=1-2, to=2-1, "\sim" {sloped, near start}, Rightarrow, shorten =2.5ex]
		\end{tikzcd}
	\end{center}
	induced by the various inclusions of the underlying sites by taking inverse images. Notice that $i^{-1}$ is the pointwise sheafification of the functor
		$$i_p^{-1}\colon T\mapsto \left(\tilde S \mapsto \colimit_{\tilde S\to S} T(S)\right)$$
	where the colimit runs over all $\kappa$-small extremally disconnected compact Hausdorff spaces $S$ with a map $\tilde S\to S$. By lemma \ref{lemma: simplified sheaf condition for extremally disconnected sets} the presheaf $\tilde S \mapsto \colimit_{\tilde S\to S} T(S)$ is already a sheaf on $\boundedExtremallyDisconnectedSets{\kappa'}$ so the natural map $i_p^{-1} \to i^{-1}$ induced by sheafification is an isomorphism of functors, in particular $i_p^{-1} \ladj i_\ast$ as well. Now the vertical morphisms of the square are equivalences by proposition \ref{proposition: independence to the defining site}, hence $(-)^\enlargementsymbol{\kappa}{\kappa'}$ is fully faithful, commutes with all colimits and with all $\cofinality{\kappa}$-small limits if and only if $i^{-1}$ does, so if and only if $i_p^{-1}$ does. We will now show these three properties for $i_p^{-1}$.
	\begin{enumerate}
		\item
		As a left adjoint, $i_p^{-1}$ is fully faithful if and only if its unit $\id \Rightarrow i_\ast\circ i_p^{-1}$ is an isomorphism. Consider $T$ a sheaf and $\tilde S$ a $\kappa$-small extremally disconnected set. Then the corresponding component
			$$T(\tilde S)\cong \colimit_{\id\colon\tilde S\to S} T(\tilde S)\to (i_\ast i_p^{-1} T)(\tilde S) \coloneqq (i_p^{-1}T)(\tilde S) \coloneqq \colimit_{\tilde S\to S} T(S)$$
		of the unit is the map associated to the inclusion $\{\id\}\subseteq \{\tilde S\to S\}$ of index categories, which is an isomorphism since $\id\colon \tilde S\to\tilde S$ is an initial object of the index category, so $\{\id\}$ is a coinitial subset of the indexing poset. Since $\tilde S$ and $T$ were arbitrary, also the unit is an isomorphism and the functor fully faithful.
	
		\item
		Again $i_p^{-1}$ is a left adjoint, so it commutes with all colimits.
		
		\item
		It remains to show the commutativity with $\cofinality{\kappa}$-small limits. As argued above, the colimit in the definition of $i_p^{-1}$ is $\cofinality{\kappa}$-filtered. Now since cofinalities are always regular, this colimit commutes with all $\cofinality\kappa$-small limits of sheaves by remark \ref{remark: interaction of filtered colimits and limits}. In particular, we obtain commutativity pointwise, hence globally, that is the functor $i_p^{-1}$ will commute with all $\cofinality\kappa$-small limits of sheaves.
	\end{enumerate}
\end{proof}

\begin{remark}[Difficulties in gluing]
	Having now introduced $\kappa$-condensed sets to avoid set theoretic issues, and their associated transition maps, we can now try to define what an unrestricted condensed set should be. A very natural choice seems to be to just glue these categories $\boundedCondensedSets{\kappa}$ indexed by uncountable strong limits cardinals $\kappa$, along the (fully faithful) transition maps $(-)^{\enlargementsymbol{\kappa}{\kappa'}}$. Indeed, Clausen and Scholze define in their lecture notes \cite[Def. 2.11]{condenseddotpdf} the category of condensed sets to be
		$$\colimit_{\kappa} \boundedCondensedSets{\kappa}$$
	where the colimit runs over the poset of all uncountable strong limit cardinals $\kappa$. This definition however has a hidden complication unaddressed in the lecture notes. Notice that this colimit is supposed to be taken over a diagram (\ie functor) $D$ following the assignment rule $\kappa \mapsto \boundedCondensedSets{\kappa}$ and $\kappa<\kappa' \mapsto (-)^{\enlargementsymbol{\kappa}{\kappa'}}$. This assignment however fails to be functorial, that is $D$ is not a valid diagram of categories! Indeed, for a composable pair $\kappa < \kappa'$ and $\kappa'<\kappa''$ we in general only have a natural isomorphism
		$$D(\kappa <\kappa'') = (-)^{\enlargementsymbol{\kappa}{\kappa''}} \cong (-)^{\enlargementsymbol{\kappa'}{\kappa''}}\circ (-)^{\enlargementsymbol{\kappa}{\kappa'}}=D(\kappa'<\kappa'')\circ D(\kappa<\kappa'),$$
	not true equality, so $D$ is not compatible with composition. Hence, a priori taking a colimit over the assignment $D$ is ill-defined. As a response to the author's \cite[\href{https://mathoverflow.net/q/484165/546808}{MathOverflow post}]{non-functoriality} about this problem, the user James Hanson submitted an \cite[\href{https://mathoverflow.net/a/492047/546808}{answer}]{condensed-sets-in-zfc} on how to formalize an alternative constructing in $\mathbf{ZFC}$ -- Scholze affirms this construction in the comments. In his answer, Hanson enlarges each $\boundedCondensedSets{\kappa}$ to an equivalent category on which corresponding functors $(-)^\enlargementsymbol{\kappa}{\kappa'}$ commute on the nose (they are simply inclusions of full subcategories), giving a proper diagram of categories. The colimit of this diagram is then shown to be well-defined in $\mathbf{ZFC}$ (it is essentially just the union of the individual categories).
	

	We will however take the approach already given by Clausen and Scholze immediately after \cite[Definition 2.11]{condenseddotpdf}. They state that the category of condensed sets is equivalently the full subcategory of $\Sheaves{\extremallyDisconnectedSets}{\Set}$ of all sheaves $T\colon \extremallyDisconnectedSets^\op \to \Set$ that are already determined by their restriction $T\resTo{\boundedExtremallyDisconnectedSets{\kappa}}$ for some $\kappa$. We will simply make this our definition of condensed sets in \ref{definition: condensed sets}.
\end{remark}
	
To this end, consider the following variant of proposition \ref{proposition: transitioning between different cardinals}.
\begin{proposition}[Unbounded transitioning]
	\label{proposition: unbounded transitioning}
	
	For each uncountable strong limit cardinal $\kappa$, denote by
		$$T\mapsto T\resTo{\kappa}\ldef T\vert_{\boundedExtremallyDisconnectedSets{\kappa}^\op}$$
	the restriction of sheaves on the site $\extremallyDisconnectedSets$ to sheaves on the site $\boundedExtremallyDisconnectedSets{\kappa}$. This restriction admits a left adjoint
		$$(-)^\enlargementsymbol{\kappa}{\infty}\colon\Sheaves{\boundedExtremallyDisconnectedSets{\kappa}}{\Set} \to \Sheaves{\extremallyDisconnectedSets}{\Set},$$
	such that $\kappa$-condensed set $T$ gets mapped to the sheaf
		$$T^\enlargementsymbol{\kappa}{\infty} \ldef \left(\tilde S\mapsto \colimit_{\tilde S\to S} T(S)\right),$$
	where the colimit is $\cofinality{\kappa}$-filtered and runs over all $\kappa$-small extremally disconnected $S$ with a map $\tilde S\to S$. This functor $(-)^\enlargementsymbol{\kappa}{\infty}$ is fully faithful, commutes with all colimits and all $\cofinality{\kappa}$-small limits.
\end{proposition}
\begin{proof}
	Observe that $(-)^\enlargementsymbol{\kappa}{\infty}$ and $(-)\resTo{\kappa}$ are the inverse and direct image of presheaves(!) associated to the map of sites $\extremallyDisconnectedSets\to\boundedExtremallyDisconnectedSets{\kappa}$ induced by the inclusion of the underlying categories. As each $\kappa$-condensed set $T$ preserves finite coproducts by the simplified sheaf conditions \ref{lemma: simplified sheaf condition for extremally disconnected sets} we find that also $T^\enlargementsymbol{\kappa}{\infty}$ does. Hence for any $T$, its image $T^\enlargementsymbol{\kappa}{\infty}$ satisfies the same sheaf conditions so already is a sheaf (before sheafification!). Thus, there is no need to sheafify and the remaining proof is now analogous to the proof of proposition \ref{proposition: transitioning between different cardinals} after observing one does not need to change underlying sites.
\end{proof}

\begin{remark}[Transitioning in general]
	\label{remark: transitioning in general}
	
	Clausen and Scholze state in \cite[Rem. 2.10]{condenseddotpdf} that, except for commutation with $\cofinality{\kappa}$-small limits, the propositions \ref{proposition: transitioning between different cardinals} and \ref{proposition: unbounded transitioning} apply verbatim to sheaves with values in any category $\category{D}$ which admits filtered colimits. They state further that if one wishes for commutation with such limits one has be more careful -- but that it is certainly true if $\category{D}$ admits a conservative forgetful functor to $\Set$ that commutes with all limits and all filtered colimits. Notice that this is readily satisfied if $\category{D}$ is very concrete, see definition \ref{definition: very concrete categories}.
	
	Furthermore, the statements have obvious analogues for presheaves, either for presheaves of sets or for presheaves with values in $\category{D}$.
\end{remark}

We can now give a working definition of unrestricted condensed sets.
\begin{definition}[Condensed sets]
	\label{definition: condensed sets}
	
	A \emph{condensed set} is a sheaf $T\colon \extremallyDisconnectedSets^\op \to \Set$ such that there exists some uncountable strong limit cardinal $\kappa$ for which the induced morphism
		$$(T\resTo{\kappa})^\enlargementsymbol{\kappa}{\infty} \to T$$
	is an isomorphism. The category $\condensedSets$ then is the full subcategory of $\Sheaves{\extremallyDisconnectedSets}{\Set}$ consisting of condensed sets.
\end{definition}

\begin{remark}[Local smallness]
	\label{remark: local smallness}
	
	We have now used precisely what we tried to avoid! The category of all sheaves $\Sheaves{\extremallyDisconnectedSets}{\Set}$ for which we can not guarantee that it is locally small! So is the full subcategory $\condensedSets$ locally small? It is! Consider condensed sets $T$ and $T'$ with respective cardinals $\kappa$ and $\kappa'$ such that $(T\resTo{\kappa})^\enlargementsymbol{\kappa}{\infty}\isorightarrow T$ and $(T'\resTo{\kappa'})^\enlargementsymbol{\kappa'}{\infty}\isorightarrow T'$ are isomorphisms. Consider any third cardinal $\kappa''$ such that $\kappa'' \ge \kappa$ and $\kappa'' \ge \kappa'$. Then also
		$$(T\resTo{\kappa''})^\enlargementsymbol{\kappa''}{\infty}\isorightarrow T\text{ and }(T'\resTo{\kappa''})^\enlargementsymbol{\kappa''}{\infty}\isorightarrow T'$$
	are isomorphisms, and we obtain that
		$$\Hom(T, T')\cong \Hom(T\resTo{\kappa''}, T'\resTo{\kappa''})$$
	by full faithfulness of $(-)^\enlargementsymbol{\kappa''}{\infty}$. Hence, the collection of morphisms between $T$ and $T'$ is small. The effort of approximating from below was not in vain!
\end{remark}

\begin{remark}[Shulman's small sheaves \& independence of the defining site]
	\label{remark: shulman's small sheaves and independence of the defining site}
	
	Consider again proposition \ref{proposition: unbounded transitioning} and definition \ref{definition: condensed sets}. In order to obtain a good transition from the $\kappa$-small world to the unrestricted one, we have to restrict to extremally disconnected sets. Indeed, when trying to mimic the definition of $(-)^\enlargementsymbol{\kappa}{\infty}$ for either the site of profinite sets or compact Hausdorff spaces, one immediately runs into the problem that the inverse image of presheaves does not produce sheaves (as one would expect for a general site). One would have to sheafify the produced presheaves. This however is problematic. On essentially large sites there is no guarantee for the existence of sheafification. Indeed, Waterhouse constructs in \cite{waterhouse} a presheaf on the fpqc-site which does not admit a sheafification. Similarly, it is not clear that analogue of proposition \ref{proposition: independence to the defining site} is viable in the unrestricted case. While one can restrict sheaves from $\CHaus$ to $\proFiniteSets$ to $\extremallyDisconnectedSets$, it is not clear that we obtain good inverse functors without sheafification.
	
	A way out of this conundrum is the theory of \emph{small sheaves} introduced by Shulman in \cite{shulman2012exactcompletionssmallsheaves}. In his work, Shulman constructs for certain $2$-categories a very general \emph{exact completion}. As a special case he defines a category of so-called \emph{small sheaves} on even a large site. In this setting one can then show that condensed sets are precisely small sheaves on the site $\extremallyDisconnectedSets$. Furthermore, Shulman provides an analogue of the comparison lemma \ref{lemma: comparison lemma}, that ensures that categories of small sheaves are invariant under considering dense sub-sites. In particular one obtains the desired analogue of proposition \ref{proposition: independence to the defining site}, that condensed sets are small sheaves on either $\CHaus$, $\proFiniteSets$ or $\extremallyDisconnectedSets$. While much more elegant than considering cutoff cardinals \quote{by hand} and trying to glue the resulting topoi, Shulman's work is quite technical and would need tremendous amount of work to reintroduce in this thesis.
\end{remark}

\begin{definition}[Small (pre)sheaves]
	\label{definition: small (pre)sheaves}
	
	Let $\category{D}$ be a very concrete category and recall that by remark \ref{remark: transitioning in general} transitioning between different cardinals is possible for $\category{D}$-valued (pre)sheaves.
	
	A presheaf $T\colon \extremallyDisconnectedSets^\op\to \category{D}$ is called \emph{small} if there is some $\kappa$ such that the natural $(T\resTo{\kappa})^\enlargementsymbol{\kappa}{\infty}\to T$ is an isomorphism. A sheaf $T\colon \extremallyDisconnectedSets^\op \to \category{D}$ is called \emph{small} if it is small as a presheaf. Denote the resulting full categories of presheaves and sheaves by $\smallPresheaves{\extremallyDisconnectedSets}{\category{D}}$ and $\smallSheaves{\extremallyDisconnectedSets}{\category{D}}$ respectively. Analogous to remark \ref{remark: local smallness} these categories are locally small.
\end{definition}

\begin{remark}[Conflicting terminology]
	\label{remark: conflicting terminology}
	
	The previous definition \ref{definition: small (pre)sheaves} is not entirely faithful to the language introduced in Shulman's \cite{shulman2012exactcompletionssmallsheaves}. He defines small sheaves differently. While the notions agree in this specific case, on a general site they seem to not be identical! On a general site the underlying presheaf of a \emph{Shulman-small} sheaf might not be small (in the sense analogous to the previous definition)! The author wants to thank professor Shulman for pointing this out to him. Certainly, in case of a general large site, \emph{Shulman-smallness} is more appropriate.
\end{remark}

\begin{lemma}[Sheafification of small presheaves]
	\label{lemma: sheafification of small presheaves}
	
	The forgetful functor $\condensedSets \to \smallPresheaves{\extremallyDisconnectedSets}{\Set}$ admits a left adjoint $\sheafification{-}\colon \smallPresheaves{\extremallyDisconnectedSets}{\Set} \to \condensedSets$ that commutes with finite limits, called \emph{sheafification}
\end{lemma}

\begin{proof}
	For a small presheaf $T$ consider the smallest cardinal $\kappa(T)$ such that $(T\resTo{\kappa(T)})^\enlargementsymbol{\kappa(T)}{\infty} \shortisorightarrow T$. Define the condensed set
	$$T^\sheafificationsymbol\ldef \left(\sheafification{T\resTo{\kappa(T)}}\right)^\enlargementsymbol{\kappa(T)}{\infty},$$
	observing that $T^\sheafificationsymbol\resTo{\kappa(T)}^\enlargementsymbol{\kappa(T)}{\infty} \isorightarrow T^\sheafificationsymbol$.
	
	For any two small presheaves $T$ and $T'$ consider $\kappa = \max\{\kappa(T), \kappa(T')\}$. Then, as $\kappa\ge\kappa(T)$ and $\kappa\ge \kappa(T')$ we have that
		$$(-)\resTo{\kappa}\colon \Hom(T,T') \to \Hom(T \resTo{\kappa}, T'\resTo{\kappa})$$
	is a bijection.	Using this, define
	$$\Hom(T, T')\to \Hom(T^\sheafificationsymbol, (T')^\sheafificationsymbol)$$
	by mapping any $f\colon T\to T'$ to the unique map $f^\sheafificationsymbol\colon T^\sheafificationsymbol \to (T')^\sheafificationsymbol$ such that its restriction $f^\sheafificationsymbol\resTo{\kappa}\colon T^\sheafificationsymbol\resTo{\kappa}\to (T')^\sheafificationsymbol\resTo{\kappa}$ corresponds -- under the isomorphisms $T^\sheafificationsymbol \resTo{\kappa} \cong \sheafification{T\resTo{\kappa}}$ and $(T')^\sheafificationsymbol \resTo{\kappa} \cong \sheafification{T'\resTo{\kappa}}$ -- to the map $\sheafification{f\resTo{\kappa}}\colon \sheafification{T\resTo{\kappa}}\to \sheafification{T'\resTo{\kappa}}$. It is clear that $\sheafification{\id_T} = \id_{T^\sheafificationsymbol}$ for any small presheaf $T$. Compatibility with composition of maps $T\to T' \to T''$ can be shown by choosing a cardinal $\kappa$ such that $(T\resTo{\kappa})^\enlargementsymbol{\kappa}{\infty}\shortisorightarrow T$, $(T'\resTo{\kappa})^\enlargementsymbol{\kappa}{\infty}\shortisorightarrow T'$ and $(T''\resTo{\kappa})^\enlargementsymbol{\kappa}{\infty}\shortisorightarrow T''$ and using the full faithfulness of the various functors $(-)^\enlargementsymbol{-}{-}$. We have thus defined a functor $\sheafification{-}\colon \smallPresheaves{\extremallyDisconnectedSets}{\Set}\to\condensedSets$. It is similarly verified that this functor is indeed a left adjoint to the forgetful functor $\condensedSets\to\smallPresheaves{\extremallyDisconnectedSets}{\Set}$. It remains to see that $\sheafification{-}$ commutes with finite limits. Since any fixed finite limit can be computed for some $\kappa$ this follows from sheafification of $\sheafification{-}\colon \Presheaves{\boundedExtremallyDisconnectedSets{\kappa}}{\Set} \to \boundedCondensedSets{\kappa}$ commuting with finite limits.
\end{proof}

Let us now consider some categorical properties of the category $\condensedSets$ of condensed sets. As a first exposition, let us show that it admits all small limits.
\begin{remark}[$\condensedSets$ is complete]
	\label{remark: the category of condensed sets is complete}
	
	Suppose $D\colon I\to \condensedSets$ is a diagram of condensed sets. For any $i\in I$ there is some uncountable strong limit cardinal $\kappa(i)$ such that
	$$(D(i)\vert_{\boundedExtremallyDisconnectedSets{\kappa(i)}^\op})^\enlargementsymbol{\kappa(i)}{\infty} \isorightarrow D(i),$$
	as all $D(i)$ are condensed sets. Choose some uncountable strong limit cardinal $\kappa$ such that $\kappa(i) < \kappa$ for all $i\in I$ and for which the indexing category $I$ is $\cofinality{\kappa}$-small. Then for each $i\in I$ also
	$$(D(i)\resTo{\kappa})^\enlargementsymbol{\kappa}{\infty}\isorightarrow D(i),$$
	as $(-)^\enlargementsymbol{\kappa}{\infty}\circ (-)^\enlargementsymbol{\kappa(i)}{\kappa} \cong (-)^\enlargementsymbol{\kappa(i)}{\infty}$ and the isomorphism is compatible with the counit.
	
	Now $\boundedCondensedSets{\kappa}$ is a topos, so is complete. Thus the restricted diagram $D'\ldef (i\mapsto D(i)\resTo{\kappa})$ admits a limiting cone, with tip say $L$. As $(-)^\enlargementsymbol{\kappa}{\infty}$ commutes with $\cofinality{\kappa}$-small limits and $D'$ is $\cofinality{\kappa}$-small by definition of $\kappa$, the image cone with tip $L^\enlargementsymbol{\kappa}{\infty}$ is limiting as well. Thus, $L^\enlargementsymbol{\kappa}{\infty}$ is a limit for the image diagram $(D')^\enlargementsymbol{\kappa}{\infty} = (i\mapsto D'(i)^\enlargementsymbol{\kappa}{\infty})$. But clearly $(D')^\enlargementsymbol{\kappa}{\infty}\isorightarrow D$, so $L^\enlargementsymbol{\kappa}{\infty}$ is also a limit for $D$. As $D$ was arbitrary, $\condensedSets$ is complete.
	
	%
	%
	%
	%
\end{remark}

We now find that the category $\condensedSets$ behaves extremely well: It is a pretopos.
\begin{theorem}[The pretopos of condensed sets]
	\label{theorem: the pretopos of condensed sets}
	
	The category $\condensedSets$ is a pretopos for which $\mathcal G\ldef \{\associated{S}\setseparator S\in \extremallyDisconnectedSets\}$ forms a class of compact projective generators.
\end{theorem}
\begin{proof}
	By remark \ref{remark: local smallness}, $\condensedSets$ is locally small. By remark \ref{remark: the category of condensed sets is complete}, $\condensedSets$ is complete and by a similar argument it is also cocomplete. To show that $\condensedSets$ is a pretopoi it thus remains to show the five conditions in definition \ref{definition: pretopoi}.
	
	Any coproduct $X=\coprod_{i\in I} X_i$ in $\condensedSets$ and the base changes of the natural morphisms $X_i\to X$ along some fixed morphism $X'\to X$ can be realized in some $\boundedCondensedSets{\kappa}$ for a suitably large $\kappa$. Since $\boundedCondensedSets{\kappa}$ is a topos and hence satisfies Giraud's axioms from theorem \ref{theorem: girauds axioms}, this new coproduct is disjoint and its base change is disjoint as well. By proposition \ref{proposition: unbounded transitioning} base changes and being disjoint are stable under passage $\enlargementsymbol{\kappa}{\infty}$ after possibly enlarging $\kappa$, so the original coproduct and its base change are disjoint too. Thus, the original coproduct is universally disjoint, since the base change was arbitrary. This shows condition \ref{definition: pretopoi, coproducts are universally disjoint} of definition \ref{definition: pretopoi}
	
	Analogously, $\condensedSets$ satisfies conditions \ref{definition: pretopoi, every epimorphism is a coequalizer}, \ref{definition: pretopoi, every equivalence relation is a kernel pair} and \ref{definition: pretopoi, every exact fork is stably exact} of definition \ref{definition: pretopoi}.
	
	Thus, we have shown all but condition \ref{definition: pretopoi, existence of a generating class} -- the existence of a generating class -- for being a pretopoi.  Denote by $\mathcal G$ the class of associated sheaves $\associated{S}$ for $S$ extremally disconnected. We will now show that $\mathcal G$ is a generating class. Observe that any topological space is homeomorphic to a topological space whose underlying set is a cardinal, \ie a well-ordered set. With this in mind, define for each $\kappa$ the class $E_\kappa$ given by
	$$E_\kappa \ldef \{S \in \boundedExtremallyDisconnectedSets{\kappa}\setseparator \text{$S$ has underlying set a cardinal $\lambda < \kappa$}\}$$
	and the subclass $\mathcal G_\kappa\subsetneq \mathcal G$ given by
	$$\mathcal G_\kappa\ldef \{\associated{S} \in \condensedSets \setseparator S\in E_\kappa\}.$$
	In particular for any extremally disconnected $\kappa$-small $S$, there is some $S'\in E_\kappa$ with $S'\cong S$ and $\associated{S'}\cong \associated{S}$. Now the class of cardinals $\lambda < \kappa$ is of size $\vert\powerset{\lambda}\vert$, hence small. Similarly, for any given cardinal $\lambda$ the class of $S$ with underlying set $\lambda$ is small as well, as a topology is just a certain subset of the powerset, so there are at most $\vert\powerset{\powerset{\lambda}}\vert$ many topologies on the set $\lambda$. Thus, $E_\kappa$ and $\mathcal G_\kappa$ are small unions of smalls and hence small. 
	
	We will now show generation. Let $T$ be a condensed set and $\kappa$ such that $(T\resTo{\kappa})^\enlargementsymbol{\kappa}{\infty}\shortisorightarrow T$. Then for any second condensed set $T'$ we have
	$$\Hom(T, T')\cong \Hom((T\resTo{\kappa})^\enlargementsymbol{\kappa}{\infty}, T')\cong \Hom(T\resTo{\kappa}, T'\resTo{\kappa}),$$
	and any $\kappa$-condensed set is of the form $T'\resTo{\kappa}$ for a suitable $T'$. Thus, $\mathcal G_\kappa$ is separating for $T$ if and only if $\mathcal G_\kappa\resTo{\kappa}=\{\associated{S}\in \boundedCondensedSets{\kappa}\setseparator S\in E_\kappa\}$ is separating for $T\resTo{\kappa}$. Clearly the inclusion $E_\kappa \to \boundedExtremallyDisconnectedSets{\kappa}$ is an equivalence and thus gives rise to a dense sub-site. By the comparison lemma \ref{lemma: comparison lemma}, restricting even further from $\boundedExtremallyDisconnectedSets{\kappa}$ to $E_\kappa$ hence is an equivalence of topoi. Thus, the section functors $\Gamma(S, -)$ for $S\in E_\kappa$ are jointly faithful, both for $\Sheaves{E_\kappa}{\Set}$ and $\boundedCondensedSets{\kappa}=\Sheaves{\boundedExtremallyDisconnectedSets{\kappa}}{\Set}$. As $\Hom(\associated{S}, -)\cong \Gamma(S, -)$ by Yoneda's lemma, we thus have that the functors $\Hom(\associated{S}, -)$ are jointly faithful and hence that $\mathcal G_\kappa\resTo{\kappa}$ is a separating set for all $\kappa$-condensed sets, in particular for $T\resTo{\kappa}$. Then $\mathcal G_\kappa$ is separating for $T$ and as $T$ was arbitrary, $\mathcal G$ is a generating class.
	
	%
	%
	
	It remains to show that each $\associated{S}$ for $S$ extremally disconnected is compact projective. Suppose $D\colon I\to \condensedSets$ is either a reflexive pair or a filtered diagram of condensed sets. Any colimit $\colimit_\preSheavesSymbol D$ of $D$ as a diagram of presheaves is calculated componentwise and hence it is straightforward to verify that $\colimit_\preSheavesSymbol D$ respects finite disjoint unions of extremally disconnected sets, so satisfies the simplified sheaf conditions from lemma \ref{lemma: simplified sheaf condition for extremally disconnected sets}. Thus $\colimit_\preSheavesSymbol D$ is already a sheaf and the natural map
	$$\colimit\nolimits_\preSheavesSymbol D \to \sheafification{\colimit\nolimits_\preSheavesSymbol D}\rdef \colimit D$$
	is an isomorphism. Then the natural map $\colimit_{i\in I} \Gamma(S, D_i)\to\Gamma(S, \colimit D)$ is an isomorphism as well, so $\Hom(\associated S, -)\cong \Gamma(S, -)$ commutes with all reflexive coequalizers (so $S$ is reflexive projective) and all filtered colimits (so $S$ is compact), proving the theorem.
\end{proof}

We should now start to reinvestigate all properties and definitions familiar to $\kappa$-condensed sets and try to reimagine them in the glued world. For the remaining part of this chapter we will omit proofs and simply reference them in Clausen's and Scholze's notes \cite{condenseddotpdf}. While the following claims are very important to the foundations of condensed mathematics as a whole, they are not as relevant to the rather \emph{discrete} world of algebraic geometry -- the main application in this thesis. As the proofs would be quite involved, the author chooses to omit them. A direct hindrance to translating all $\kappa$-condensed results verbatim is given in the following remark. 
\begin{remark}[Non-closed points and condensed sets]
	Unfortunately, the natural functor $\Top \to \Sheaves{\extremallyDisconnectedSets}{\Set}, X\mapsto \associated{X}$ does not land in condensed sets. The problem are spaces that are not $\tOne$, \ie spaces in which there exist non-closed points. Any such space contain a copy of the \emph{Sierpinski-space} $X=\{0, 1\}$ with topology $\topology{T} = \{\emptyset, \{1\}, X\}$ in which $\{1\}$ clearly is not a closed set. Already the Sierpinski space is ill-represented: The associated sheaf $\associated{X}$ of $X$ is not small, \confer \cite[Warning 2.14]{condenseddotpdf}. This is essentially due to the fact that continuous maps $S\to X$ are equivalently closed subsets of $S$. Thus, $\associated{X}$ is equivalently the functor that assigns to each $S$ the set of closed subsets -- this is however \quote{too much data to form a small sheaf}. As Clausen and Scholze mention in this warning, if there were to be a suitable $\kappa$ then any closed subset of any extremally disconnected $S$ (both however large) would need to be the intersection of at most $\kappa$ clopen subsets of $S$.
\end{remark}

This is however all that can happen.
\begin{lemma}[$\tOne$-spaces give rise to condensed sets]
	\label{lemma: T1-spaces give rise to condensed sets}
	
	For each $\tOne$-space $X$, the sheaf $\associated{X}$ is a condensed set.
\end{lemma}
\begin{shortproof}
	This is \cite[Prop. 2.15]{condenseddotpdf}.
\end{shortproof}

With this we can now define the analogue of definition \ref{definition: associated kappa-condensed sets}.
\begin{definition}[Associated condensed sets]
	Denote by $\tOneSpaces$ the full subcategory of $\Top$ of $\tOne$-spaces. By lemma \ref{lemma: T1-spaces give rise to condensed sets} there is a functor
	\begin{align*}
		\tOneSpaces &\longrightarrow \boundedCondensedSets{\kappa},\\
		X&\longmapsto \associated{X}\ldef \Hom_\Top(-, X),
	\end{align*}
	that assigns to each $\tOne$-topological space $X$ its \emph{associated condensed set} $\associated{X}$.
\end{definition}

And in the other direction, we also obtain the analogue of \ref{definition: the underlying space of a kappa-condensed set}.
\begin{definition}[The underlying space of a condensed set]
	\label{definition: the underlying space of a condensed set}
	
	Let $T$ be a condensed set. Then $T(\ast)$ is called the \emph{underlying set of $T$}. Define a topology on the underlying set $T(\ast)$ as the final topology of all maps of sets $S\shortisorightarrow \associated{S}(\ast)\to T(\ast)$ induced by morphisms $\associated{S}\to T$ ranging over all profinite sets $S$. Denote the resulting space by $\topologizedUnderlying{T}$ and call it the \emph{underlying space} of $T$.
\end{definition}

\begin{remark}[Alternative description of the underlying topology]
	If one were to ignore set theoretic issues then the topology of the underlying space in definition \ref{definition: the underlying space of a condensed set} is equivalently the quotient topology of the surjective map
		$$\smashoperator{\coprod_{\substack{\associated{S}\to T\\S\in\proFiniteSets}}} S\to T(\ast).$$
\end{remark}

To ensure a good adjunction with right-adjoint $X\mapsto\associated{X}$ and left-adjoint $T\mapsto \topologizedUnderlying{T}$, we have to restrict our attention to condensed sets $T$, such that the underlying space $\topologizedUnderlying{T}$ is $\tOne$. For this (and the following theory) we need the notions of quasi-compactness and quasi-separatedness in a (pre)topos. 
\begin{definition}[Quasi-compactness \& quasi-separatedness]
	\label{definition: quasi-compactness and quasi-separatedness}
	
	Let $\category T$ be a topos or pretopos.
	\begin{itemize}
		\item
		An object $X\in \category T$ is called \emph{quasi-compact (in short qc)} if for any epimorphism $\coprod_{i\in I} X_i \twoheadrightarrow X$ there exists a finite subset $J\subseteq I$ such that $\coprod_{i\in I} X_i \to X$ is still an epimorphism.
		
		\item
		A map $f\colon X\to S$ is called \emph{quasi-compact} if $f\in \category T/S$ is quasi-compact.
		
		If $\category T$ is generated by quasi-compact objects, then quasi-compactness of $f$ is equivalent to the fiber product $X\times_S Y$ being quasi-compact for all morphisms $Y\to S$ and $Y$ quasi-compact -- \emphquote{base changes of qc objects over $S$ along $f$ are qc}.
		
		\item
		An object $X\in \category T$ is called \emph{quasi-separated (in short qs)} if for all quasi-compact $K, L\in \category T$ with maps $K\to X$ and $L\to X$ the fiber product $K\times_X L$ is quasi-compact as well -- \emphquote{intersections of qc objects in $X$ are qc}.
		
		\item
		A morphism $f\colon X\to S$ in $\category T$ is called \emph{quasi-separated} if it is quasi-separated in $\category T/S$.
		
		If $\category T$ is generated by qcqs objects, then quasi-separatedness of $f$ is equivalent to the fiber product $X\times_S Y$ being quasi-separated for all morphisms $Y\to S$ where $Y$ is qcqs -- \emphquote{base changes of qcqs objects over $S$ along $f$ are qs}.
	\end{itemize}
\end{definition}

\begin{remark}[Slice (pre)topoi \& the Fundamental Theorem of Topos Theory]
	In the setting of definition \ref{definition: quasi-compactness and quasi-separatedness}, the slice category $\category T/S$ is again a (pre)topos by what one usually calls the \emph{fundamental theorem of topos theory}. See \eg \cite[A2.3]{elephant1} for a proof in the case of (not necessarily Grothendieck, so-called \emph{elementary}) topoi. In any case, even without the knowledge that $\category T/S$ again is a (pre)topos, one can mimic the definitions of quasi-compactness and quasi-separatedness to allow for more general categories.
\end{remark}

Quasi-compactness of points now gives the analogue of being $\tOne$ in the condensed world.
\begin{lemma}[$\tOne$-spaces and quasi-compact points]
	\label{lemma: T1-spaces and quascompact points}
	
	If $X$ is a $\tOne$-space, then for every point $x\in X$, the corresponding map $\associated{\ast}\to \associated{X}$ is quasi-compact. Conversely, if $T$ is a condensed set such that every point $\associated{\ast}\to T$ is quasi-compact, then $\topologizedUnderlying{T}$ is compactly generated and $\tOne$.
\end{lemma}
\begin{shortproof}
	This is \cite[Prop. 2.15]{condenseddotpdf}.
\end{shortproof}

\begin{remark}[Maps from points and compact fibers]
	If we accept as fact that $\condensedSets$ is generated by qc objects (indeed the compact, extremally disconnected spaces $S$ stay (quasi)compact when passing to $\condensedSets$), we see that quasi-compactness of $\associated{\ast} \to T$ implies that for any quasi-compact $T'$ and any $f\colon T'\to T$ the fiber $T'\times_T \associated{\ast}$ is quasi-compact as well, \ie for $T$, \emphquote{fibers over the point $\associated{\ast}\to T$ inside quasi-compact condensed sets are quasi-compact}. This is precisely what one would expect of a closed point of a topological space, where fibers then are closed (and hence compact if the ambient space is compact as well). So maps from points being quasi-compact is really the right analogue of $\tOne$-ness for condensed sets.
\end{remark}

For such spaces and condensed sets we obtain the following relation.
\begin{corollary}[The adjunction]
	There is an adjunction with right-adjoint $X\mapsto \associated{X}$ and left-adjoint $T\mapsto \topologizedUnderlying{T}$ between $\tOne$-topological spaces and condensed sets for which all maps from points are quasi-compact.
\end{corollary}

This adjunction now gives rise to several correspondences between subclasses of topological spaces and condensed sets. One of these is the class of \emph{compactly generated weak Hausdorff spaces}.
\begin{definition}[Weak Hausdorff spaces]
	A topological space $X$ is said to be \emph{weak Hausdorff} if for any compact Hausdorff space $K$ and any continuous map $K\to X$, its image is closed in $X$. In particular all points of $X$ are closed, so $X$ is a $\tOne$-space.
	
	Recall the well known fact that any continuous map from a compact space to a Hausdorff space is closed. Thus, any Hausdorff space is indeed also weak Hausdorff.
	
	A topological space $X$ is said to be \emph{compactly generated weak Hausdorff} (\CGWH for short) if it is both compactly generated and weak Hausdorff.
\end{definition}

We can now state the claimed correspondences.
\begin{theorem}[Some correspondences]
	\label{theorem: some correspondences}
	
	Under the adjunction from lemma \ref{lemma: T1-spaces and quascompact points} we obtain the following statements:
	\begin{enumerate}[(i)]
		\item
		The functor $X\to\associated{X}$ induces an equivalence between $\CHaus$ and qcqs condensed sets.
		
		\item 
		A space $X$ is \CGWH if and only if $\associated{X}$ is quasi-separated. For every quasi-separated condensed set $T$, the space $\topologizedUnderlying{T}$ is \CGWH
	\end{enumerate}
	
	In particular $X\to \associated{X}$ and its left-adjoint $T\mapsto \topologizedUnderlying{T}$ restrict to an adjunction between \CGWH spaces and quasi-separated condensed sets. Furthermore, $X\to\associated{X}$ is fully faithful on \CGWH spaces.
\end{theorem}
\begin{shortproof}
	This is \cite[Thm. 2.16]{condenseddotpdf}.
\end{shortproof}

\begin{notation}[Simplifying notation]
	Let $X$ be a \CGWH space. Due to theorem \ref{theorem: some correspondences} we are thus justified in simply writing $X$ for the associated condensed set $\associated{X}$, if no confusion is possible.
\end{notation}

	\chapter{Rings, Modules \& Condensed Algebra}
	\label{chapter: rings, modules and condensed algebra}
	
Having now introduced condensed sets -- the foundational object of condensed mathematics -- we can start to build algebra on top of them. For this we will first discuss the general theory of abelian sheaves, sheaves of rings and sheaves of modules over such rings in sections \ref{section: abelian sheaves}, \ref{section: rings, ringed spaces and ringed sites} and \ref{section: modules under morphisms}. Then afterwards, we will start to build the theory of condensed algebra in section \ref{section: condensed algebraic structures}. We will study the surrounding derived theory in chapter \ref{chapter: cohomology and derived condensed algebra}.

To construct a good theory of sheaves of modules one usually has to restrict to essentially small sites, as sheafification is a fundamental part of the theory and is only guaranteed to exist in that case. However almost all results stated solely for essentially small sites hold more generally. As long as one restricts attention to a suitable subcategory of presheaves (and the corresponding category of sheaves with such underlying presheaves) such that sheafification (as the adjoint to the forgetful functor) is defined and commutes with all limits, then virtually all the desired results hold. This is then readily satisfied by small presheaves $\smallPresheaves{\extremallyDisconnectedSets}{\Set}$ and their corresponding sheaves $\condensedSets = \smallSheaves{\extremallyDisconnectedSets}{\Set}$, where the forgetful functor admits such a left adjoint by lemma \ref{lemma: sheafification of small presheaves}. For each claim about essentially small sites, we will give a short remark about the translation to small sheaves on $\extremallyDisconnectedSets$ -- this setting is of course our main application as we wish to \emph{algebra} with condensed sets.

\section{Abelian Sheaves}
\label{section: abelian sheaves}

An important stepping stone are \emph{abelian (pre)sheaves}. Recall that there are two common ways of defining them.
\begin{definition}[Abelian (pre)sheaves]
	Let $\site{C}$ be a site. An \emph{abelian presheaf} (resp. \emph{abelian sheaf}) on $\site{C}$ is an abelian group objects internal to the category $\Presheaves{\site{C}}{\Set}$ (resp. $\Sheaves{\site{C}}{\Set}$). Equivalently, an \emph{abelian (pre)sheaf} is a (pre)sheaf of abelian groups on $\site{C}$.
	
	A \emph{morphism of abelian (pre)sheaves} is a morphism of abelian group objects. Equivalently, this is a morphism of (pre)sheaves of abelian sheaves.
\end{definition}

Of course sheafification of sheaves of sets translates to the abelian setting.
\begin{lemma}[Sheafification of abelian presheaves]
	\label{lemma: sheafification of abelian presheaves}
	
	Let $\site{C}$ be an essentially small site. Sheafification of sheaves of sets gives rise to a left adjoint $\sheafification{-}\colon \Presheaves{\site{C}}{\Ab}\to\Sheaves{\site{C}}{\Ab}$ of the forgetful functor $\Sheaves{\site{C}}{\Ab} \to \Presheaves{\site{C}}{\Ab}$ that commutes with finite limits. This left adjoint is called \emph{sheafification (of abelian presheaves)}.
\end{lemma}
\begin{proof}
	This easily seen if one considers abelian sheaves to be abelian group objects internal to sheaves of sets. The claim then follows from sheafification commuting with finite limits as per definition \ref{definition: sheafification}.
\end{proof}

\begin{remark}[Sheafification of small abelian presheaves]
	\label{remark: sheafification of small abelian presheaves}
	
	Lemma \ref{lemma: sheafification of abelian presheaves} applies verbatim to small (pre)sheaves of abelian groups on the site $\extremallyDisconnectedSets$, as sheafification $\sheafification{-}\colon \smallPresheaves{\extremallyDisconnectedSets}{\Set} \to \smallSheaves{\extremallyDisconnectedSets}{\Set}$ exists and commutes with finite limits by lemma \ref{lemma: sheafification of small presheaves}.
\end{remark}

\begin{proposition}[Abelian sheaves form an abelian category]
	\label{proposition: abelian sheaves form an abelian category}
	
	Let $\site{C}$ be an essentially small site. The category $\Sheaves{\site{C}}{\Ab}$ of abelian sheaves is abelian.
\end{proposition}
\begin{proof}
	This is a direct application \cite[\href{https://stacks.math.columbia.edu/tag/03A3}{Tag 03A3}]{stacks-project} after observing that sheafification of abelian presheaves commutes with finite limits by lemma \ref{lemma: sheafification of abelian presheaves}.
\end{proof}

\begin{remark}[Small abelian sheaves form an abelian category]
	Proposition \ref{proposition: abelian sheaves form an abelian category} applies verbatim to small abelian sheaves on the site $\extremallyDisconnectedSets$. Indeed, \cite[\href{https://stacks.math.columbia.edu/tag/03A3}{Tag 03A3}]{stacks-project} used in the proof of proposition \ref{proposition: abelian sheaves form an abelian category} applies as $\smallPresheaves{\extremallyDisconnectedSets}{\Ab}$ is abelian (all $\smallPresheaves{\boundedExtremallyDisconnectedSets{\kappa}}{\Ab}$ are abelian and the transition functors between different cardinal-cutoffs preserve the relevant structure by remark \ref{remark: transitioning in general}) and sheafification $\sheafification{-}\colon \smallPresheaves{\extremallyDisconnectedSets}{\Ab}\to \smallSheaves{\extremallyDisconnectedSets}{\Set}$ exists and is exact by remark \ref{remark: sheafification of small abelian presheaves}. Hence $\smallSheaves{\extremallyDisconnectedSets}{\Ab}$ is abelian as well.
\end{remark}

\section{Rings, Ringed Spaces \& Ringed Sites}
\label{section: rings, ringed spaces and ringed sites}

Recall the notion of ringed sites, the best-known instances of which probably being schemes.
\begin{definition}[Ringed spaces, {
		\cite[\href{https://stacks.math.columbia.edu/tag/01HA}{Tag 01HA}]{stacks-project}}]
	
	A \emph{ringed space} is a pair $(X, \struct X)$ where $X$ is a topological space and $\struct X$ is a sheaf of commutative unitary rings on $X$, called the \emph{structure sheaf}. If additionally $\stalkStruct X x$ is a local ring for any $x\in X$, then $(X, \struct X)$ is called a \emph{locally ringed space}.
	
	Let $(X,\struct X)$ and $(Y, \struct Y)$ be ringed spaces. A \emph{morphism of ringed spaces} $(X,\struct X) \to (Y, \struct Y)$ is a pair $(f, f^\sharp)$ where $f\colon X\to Y$ is continuous and $f^\sharp \colon \struct Y \to f_\ast \struct X$ is a morphism of sheaves of rings. If $(X, \struct X)$ and $(Y, \struct Y)$ are locally ringed spaces we call the pair $(f, f^\sharp)$ a \emph{morphism of locally ringed spaces} if additionally for any $x\in X$ the induced map
	$$\stalkStruct Y {f(x)} \xrightarrow{f^\sharp_{f(x)}} (f_\ast \struct X)_{f(x)} \to \stalkStruct X x$$
	on stalks is a \emph{local morphism}, \ie maps the maximal ideal to the maximal ideal.
\end{definition}

There is a straightforward generalization (once one is comfortable replacing spaces by sites), that will prove useful when discussing condensed rings and modules over such a ring. Instead of considering sheaves on spaces, we simply generalize to sheaves on sites. First one needs a notion of sheaf of rings.
\begin{definition}[Sheaves of rings]
	Let $\site{C}$ be a site. A \emph{sheaf of rings} or simply a \emph{ring} on $\site{C}$ is a (commutative unitary) ring object internal to the category $\Sheaves{\site{C}}{\Ab}$. Equivalently, a \emph{sheaf of rings} is a sheaf of (commutative unitary) rings on $\site{C}$.
	
	A \emph{morphism of sheaves of rings} is then a morphism of ring objects. Equivalently, it is a morphism in $\Sheaves{\site{C}}{\commutativeRings}$.
\end{definition}

\begin{example}[Many topological rings give rise to sheaves of rings]
	Every topological ring $A$ whose underlying space is $\tOne$ gives rise to a sheaf of rings $\associated{A}$ (where addition and multiplication of each $\associated{A}(S)=\eHom_\Top(S, A)$ is pointwise) on the site $\extremallyDisconnectedSets$ by theorem \ref{theorem: some correspondences}. This of course includes all discrete rings, that is \quote{rings without a topology}, but also topological rings like the $p$-adic integers $\Z_p$ or more generally the ring of formal power series $\series{\Z}{T}$ with its $\ideal{T}$-adic topology.
\end{example}

With this we can provide the needed generalization of ringed spaces.
\begin{definition}[Ringed sites, {
		\cite[\href{https://stacks.math.columbia.edu/tag/03AD}{Tag 03AD}]{stacks-project}}]
	
	A \emph{ringed site} is a pair $(\site C, \structsymbol)$ where $\site C$ is a site and $\structsymbol$ is a sheaf of rings on $\site C$, called the \emph{structure sheaf}.
	
	Let $(\site C,\structsymbol)$ and $(\site C', \structsymbol')$ be ringed sites. A \emph{morphism of ringed sites} $(\site C,\structsymbol) \to (\site C', \structsymbol')$ is a pair $(f, f^\sharp)$ where $f\colon \site C\to \site C'$ is a map of sites and $f^\sharp \colon \structsymbol' \to f_\ast \structsymbol$ is a morphism of sheaves of rings, where $f_\ast \structsymbol$ carries the obvious ring structure.
%
\end{definition}

Of course ringed spaces can naturally be viewed as ringed sites.
\begin{remark}[Ringed spaces as ringed sites]
	Any ringed space $(X,\struct X)$ induces a ringed site $(\open{X}, \struct X)$. Similarly, for any morphism of ringed spaces $(f,f^\sharp)\colon (X,\struct X)\to(Y,\struct Y)$, the continuous map $f\colon X\to Y$ induces a morphism of sites $\open X\to\open Y$ and $f_\ast$ (as known from the theory of sheaves on a topological space) agrees with the direct image associated to this morphism. Hence, this morphism of sites together with the map $f^\sharp$ constitutes a morphism of ringed sites $(\open X, \struct X) \to (\open Y, \struct Y)$.
\end{remark}

\begin{remark}[Open subspaces of ringed spaces]
	Let $(X,\struct X)$ be a ringed space and $U\subseteq X$ open. Then $\struct U\ldef \struct X\restriction_U$ is a sheaf of rings on $U$ and hence $(U, \struct U)$ naturally is a ringed space. Even more, there is a natural map of ringed spaces $(i,i^\sharp)\colon (U, \struct U)\to (X, \struct X)$ where $i\colon U\to X$ is the inclusion and $i^\sharp_V \colon \struct X(V)\to (i_\ast\struct U)(V) = \struct X(U\cap V)$ is the restriction map for all $V$ open in $X$.
\end{remark}

\section{(Pre)sheaves of Modules}
\label{section: (pre)sheaves of modules}

\begin{definition}[(Pre)sheaves of modules, {
	\cite[\href{https://stacks.math.columbia.edu/tag/03CV}{Tag 03CV}]{stacks-project}
	\&
	\cite[\href{https://stacks.math.columbia.edu/tag/03CT}{Tag 03CT}]{stacks-project}}]
	\label{definition: (pre)sheaves of modules}
	
	Let $(\site C, \structsymbol)$ be a ringed site. A \emph{presheaf of $\structsymbol$-modules} $\sheaf F$ is a module object over the sheaf of rings $\structsymbol$, internal to the category $\Presheaves{\site C}{\Ab}$. A \emph{sheaf of $\structsymbol$-modules} is a presheaf of $\structsymbol$-modules whose underlying presheaf is a sheaf.
	
	A \emph{morphism of presheaves of $\structsymbol$-modules} $\sheaf F\to\sheaf G$ is a morphism of the corresponding module objects internal to $\Presheaves{\site C}{\Ab}$. A \emph{morphism of sheaves of $\structsymbol$-modules} is a morphism of the corresponding presheaves of $\structsymbol$-modules.
	
	We denote by $\pMod{\structsymbol}$ and $\mod\structsymbol$ the categories of presheaves and sheaves of $\structsymbol$-modules respectively. Sets of morphisms in $\pMod\structsymbol$ and $\mod\structsymbol$ are denoted by $\Hom_\structsymbol$.
\end{definition}

As definition \ref{definition: (pre)sheaves of modules} is quite abstract, it makes sense to try to gain a more explicit understanding of (pre)sheaves of modules and their morphisms. We will do so by inspecting them on the level of sections.
\begin{remark}[(Pre)sheaves of modules, more explicit]
	The condition in definition \ref{definition: (pre)sheaves of modules} for a sheaf $\sheaf{F}$ to be a (pre)sheaf of $\structsymbol$-modules can be given more explicitly. Indeed, we require, that for any $U \in \site C$ the sections $\sheaf F(U)$ carry the structure of an $\structsymbol(U)$-module and for any morphism $V\subseteq U$ the restriction $\sheaf F(U)\to \sheaf F(V)$ is $\structsymbol(U)$-linear where $\sheaf F(V)$ is considered as an $\structsymbol(U)$-module via restriction along the morphism $\structsymbol(U)\to\structsymbol(V)$. Concretely, for any section $s\in \sheaf F(U)$ and any $r\in \structsymbol(U)$ we require $(r\cdot s)\restriction_V = r\restriction_V \cdot s\restriction_V$.
	
	Similarly, the condition on being a morphism of (pre)sheaves of $\structsymbol$-modules can be rephrased as follows: An additive morphism $f\colon \sheaf F\to \sheaf G$ between (pre)sheaves of $\structsymbol$-modules is a morphism of (pre)sheaves of $\structsymbol$-modules if
	\begin{center}
		\begin{tikzcd}
			\structsymbol \times \sheaf F &\sheaf F\\
			\structsymbol \times \sheaf G &\sheaf G
			\ar[from=1-1, to=1-2]
			\ar[from=2-1, to=2-2]
			\ar[from=1-1, to=2-1, "\id \times f", swap]
			\ar[from=1-2, to=2-2, "f"]
		\end{tikzcd}
	\end{center}
	commutes, \ie for any $U\subseteq X$, any $r\in \structsymbol(U)$ and any $s\in \sheaf F(U)$ we require $f_U(r\cdot s) = r\cdot f_U(s)$.
\end{remark}

As one expects, there is a canonical way of turning a presheaf of modules into a sheaf of modules. This is formalized in the following lemma.
\begin{lemma}[Sheafification of presheaves of modules]
	\label{lemma: sheafification of presheaves of modules}
	
	Let $(\site C, \structsymbol)$ be an essentially small ringed site. For each presheaf of $\structsymbol$-modules $\sheaf F$ there is a unique $\structsymbol$-module structure $\structsymbol\times \sheaf F^\sheafificationsymbol \to \sheaf F^\sheafificationsymbol$ on $\sheaf F^\sheafificationsymbol$ such that
	\begin{center}
		\begin{tikzcd}
			\structsymbol \times \sheaf F &\sheaf F\\
			\structsymbol \times \sheaf F^\sheafificationsymbol & \sheaf F^\sheafificationsymbol
			\ar[from=1-1, to=1-2]
			\ar[from=1-1, to=2-1]
			\ar[from=1-2, to=2-2]
			\ar[from=2-1, to=2-2]
		\end{tikzcd}
	\end{center}
	commutes. Additionally, if $\sheaf G$ is a sheaf of $\structsymbol$-modules, then any morphism $\sheaf F\to\sheaf G$ of presheaves of $\structsymbol$-modules factors uniquely as $\sheaf F\to \sheaf F^\sheafificationsymbol\to\sheaf G$ where $\sheaf F^\sheafificationsymbol \to \sheaf G$ is a morphism of sheaves of $\structsymbol$-modules. Hence sheafification of abelian sheaves gives rise to an adjunction
		$$\mod\structsymbol \xleftrightarrows{\sheafification{-}}{} \pMod\structsymbol.$$
\end{lemma}
\begin{proof}
	This follows from \cite[\href{https://stacks.math.columbia.edu/tag/03CY}{Tag 03CY}]{stacks-project} using that the unit $\structsymbol\to\structsymbol ^\sheafificationsymbol$ is an isomorphism.
\end{proof}

\begin{remark}[Sheafification of small presheaves of modules]
	The previous lemma \ref{lemma: sheafification of presheaves of modules} applies the small presheaves of modules on the site $\extremallyDisconnectedSets$ for any small sheaf of rings. Indeed, this is a direct consequence of the existence of sheafification and its commutation with finite limits.
\end{remark}

We also find that the category of sheaves of modules are abelian.
\begin{proposition}[The abelian category of sheaves of modules]
	\label{proposition: the abelian category of sheaves of modules}
	
	Let $(\site C, \structsymbol)$ be an essentially small ringed site. The category $\mod\structsymbol$ of sheaves of $\structsymbol$-modules is an abelian category. The forgetful functor $\mod\structsymbol \to \Sheaves{\site C}{\Ab}$ is exact, hence kernels, cokernels and exactness of sequences of sheaves of $\structsymbol$-modules agree with their corresponding notions of abelian sheaves.
\end{proposition}
\begin{shortproof}
	This is \cite[\href{https://stacks.math.columbia.edu/tag/03DA}{Tag 03DA}]{stacks-project}.
\end{shortproof}

\begin{remark}
	While proposition \ref{proposition: the abelian category of sheaves of modules} translates to small sheaves as well, it is more instructive to state the corresponding result later on in theorem \ref{theorem: condensed modules form an abelian category} as there is much more to say about such categories. 
\end{remark}

\begin{lemma}[Tensor products \& internal $\Hom$s]
	\label{lemma: tensor products and internal Homs}
	
	Let $(\site{C}, \structsymbol)$ be a ringed site. Then $\pMod{\structsymbol}$ admits a symmetric monoidal tensor product given on presheaves $\sheaf{F}$ and $\sheaf{G}$ by
		$$\sheaf F \otimes_{p, \structsymbol} \sheaf{G} \ldef \mleft(U\mapsto \sheaf{F}(U)\otimes_{\structsymbol(U)} \sheaf{G}(U)\mright).$$
	This tensor product admits a partial adjoint given on presheaves $\sheaf{F}$ and $\sheaf{G}$ by
		$$\intHom_\structsymbol(\sheaf{F}, \sheaf{G})\ldef \mleft(U\mapsto \Hom_\structsymbol(\sheaf{F}\resTo{U}, \sheaf{G}\resTo{U})\mright),$$
	where $\sheaf{F}\resTo{U}$ and $\sheaf{G}\resTo{U}$ are the induced presheaves on the slice-category $\site{C}/U$. This makes $\pMod{\structsymbol}$ a closed symmetric monoidal category.
	
	If $\site{C}$ is essentially small, the category $\mod{\structsymbol}$ of sheaves of $\structsymbol$-modules is closed symmetric monoidal as well, with tensor product $-\otimes_\structsymbol- \ldef \sheafification{-\otimes_{p, \structsymbol}-}$ and internal $\Hom$ given by $\intHom_\structsymbol$.
\end{lemma}
\begin{proof}
	The existence of the tensor product and the internal Hom (as well as their properties)
	can be found at \cite[\href{https://stacks.math.columbia.edu/tag/03EK}{Tag 03EK}]{stacks-project} and \cite[\href{https://stacks.math.columbia.edu/tag/04TT}{Tag 04TT}]{stacks-project} respectively. For the claimed adjunctions consider \cite[\href{https://stacks.math.columbia.edu/tag/03EO}{Tag 03EO}]{stacks-project}.
\end{proof}

The partial adjunctions even enrich to the internal $\Hom$.
\begin{lemma}[(Enriched) tensor-Hom adjunction]
	\label{lemma: enriched tensor-Hom adjunction}
	
	Let $(\site C, \structsymbol)$ be an essentially small ringed site. There are isomorphisms
	$$\intHom_\structsymbol(\sheaf F\otimes_\structsymbol\sheaf G, \sheaf H) \cong \intHom_\structsymbol(\sheaf F, \intHom_\structsymbol(\sheaf G, \sheaf H))$$
	natural in the sheaves $\sheaf F$, $\sheaf G$ and $\sheaf H$ of $\structsymbol$-modules.
\end{lemma}
\begin{shortproof}
	This is \cite[\href{https://stacks.math.columbia.edu/tag/03EO}{Tag 03EO}]{stacks-project}.
\end{shortproof}

\begin{remark}[Tensor products and internal $\Hom$s of small sheaves]
	Again, by the usual argument, lemma \ref{lemma: tensor products and internal Homs} and lemma \ref{lemma: enriched tensor-Hom adjunction} translate to small (pre)sheaves of modules on $\extremallyDisconnectedSets$ for any small sheaf of rings. 
\end{remark}

\section{Modules under Morphisms}
\label{section: modules under morphisms}

\begin{lemma}[Direct and inverse images of (pre)sheaves of modules]
	\label{lemma: direct and inverse images of (pre)sheaves of modules}
	
	Let $f\colon \site C\to\site C'$ be a morphism of essentially small sites and let $\structsymbol$ and $\structsymbol'$ be sheaves of rings on $\site C$ and $\site C'$ respectively. Let $\sheaf F$ be a (pre)sheaf of $\structsymbol$-modules. There is a natural morphism $f_\ast\structsymbol \times f_\ast \sheaf F\to f_\ast\sheaf F$ that turns $f_\ast\sheaf F$ into a (pre)sheaf of $f_\ast\structsymbol$-modules. Similarly, for any (pre)sheaf $\sheaf G$ of $\structsymbol'$-modules there is a natural morphism $f^{-1}\structsymbol' \times f^{-1}\sheaf G\to f^{-1}\sheaf G$ which turns $f^{-1}\sheaf G$ into a (pre)sheaf of $f^{-1}\structsymbol'$-modules. Both constructions are functorial, that is we obtain functors
		$$\maybePMod\structsymbol \to \maybePMod{f_\ast\structsymbol}$$
	and
		$$\maybePMod{\structsymbol'}\to\maybePMod{f^{-1}\structsymbol'}.$$
\end{lemma}
\begin{proof}
	This is shown in \cite[\href{https://stacks.math.columbia.edu/tag/03D1}{Tag 03D1}]{stacks-project} and \cite[\href{https://stacks.math.columbia.edu/tag/03D2}{Tag 03D2}]{stacks-project} and boils down to the fact that both $f_\ast$ and $f^{-1}$ commute with finite products of sheaves of sets.
\end{proof}

\begin{remark}[Direct and inverse images of (pre)sheaves of modules]
	\label{remark: direct and inverse images of (pre)sheaves of modules}
	
	One has to be careful when translating lemma \ref{lemma: direct and inverse images of (pre)sheaves of modules} to small sheaves. First, it is necessary to distinguish the cases of morphisms of sites $\extremallyDisconnectedSets\to\site{C'}$ and $\site{C} \to \extremallyDisconnectedSets$ and even endomorphisms $\extremallyDisconnectedSets \to \extremallyDisconnectedSets$ (essentially by virtue that we defined smallness only for the site $\extremallyDisconnectedSets$, not in general). Secondly, it is not clear the usual definitions of $f_\ast$ and $f^{-1}$ simply translate. In an exchange with \href{https://home.sandiego.edu/~shulman/}{Professor Mike Schulman}, he mentioned that in general (contrary to what one might initially expect) a morphism of sites always induces an inverse image $f^{-1}$ of small sheaves (small in the sense of his paper \cite{shulman2012exactcompletionssmallsheaves}, \confer remarks \ref{remark: shulman's small sheaves and independence of the defining site} and \ref{remark: conflicting terminology}), but what fails is the existence of the direct image! The author is however not aware of any reference of this fact, as the literature on small sheaves (again, in the sense of Shulman) is rather thin. He could not locate the result in \eg \cite{shulman2012exactcompletionssmallsheaves}.
	
	Either way, in this thesis we will only need the very special case with $f^\top = \id_\extremallyDisconnectedSets$ being the identity. In this case of course $f_\ast$ and $f^{-1}$ are simply the identity and the result follows. This will prove useful in defining scalar restriction and extension of modules on $\extremallyDisconnectedSets$ later on.
\end{remark}

\begin{lemma}[Pullback and pushforward of sheaves of modules, {
	\cite[\href{https://stacks.math.columbia.edu/tag/03D6}{Tag 03D6}]{stacks-project}}]
	\label{lemma: pullback and pushforward of sheaves of modules}
	
	Let $(f,f^\sharp)\colon (\site C, \structsymbol) \to (\site C',\structsymbol')$ be a morphism of essentially small ringed sites. We define two functors
		$$f_\ast\colon \mod\structsymbol\rightarrow\mod{\structsymbol'}$$
	and
		$$f^\ast\colon \mod\structsymbol\leftarrow \mod{\structsymbol'}.$$
	\begin{itemize}
		\item
		Let $\sheaf F$ be a sheaf of $\structsymbol$-modules. The \emph{pushforward of $\sheaf F$ along $f$} is the sheaf of $\structsymbol'$-modules with underlying abelian sheaf $f_\ast \sheaf F$ and module structure given by the map $f^\sharp\colon \structsymbol'\to f_\ast \structsymbol$ and the $f_\ast\structsymbol$-module structure $f_\ast\structsymbol\times f_\ast\sheaf F\to f_\ast\sheaf F$ of lemma \ref{lemma: direct and inverse images of (pre)sheaves of modules}.
		
		\item
		Let $\sheaf G$ be a sheaf of $\structsymbol'$-modules. The \emph{pullback of $\sheaf G$ along $f$} is the sheaf of $\structsymbol$-modules with underlying abelian sheaf
			$$f^\ast\sheaf G\ldef \structsymbol\otimes_{f^{-1}\structsymbol'} f^{-1}\sheaf G$$
		where $\structsymbol$ is a sheaf of $f^{-1}\structsymbol'$-modules by virtue of the map $f^{-1}\structsymbol'\to\structsymbol$ corresponding to $f^\sharp\colon \structsymbol'\to f_\ast\structsymbol$ under the adjunction $f^{-1}\ladj f_\ast$. The module structure again is induced by lemma \ref{lemma: direct and inverse images of (pre)sheaves of modules} after applying $\structsymbol \otimes_{f^{-1}\structsymbol'} -$ and using the isomorphism $\structsymbol \to \structsymbol\otimes_{f^{-1}\structsymbol'} f^{-1}\structsymbol$.
	\end{itemize}
	
	By \cite[\href{https://stacks.math.columbia.edu/tag/03D7}{Tag 03D7}]{stacks-project}, the pushforward and pullback form and adjunction
		$$\mod{\structsymbol}\xleftrightarrows{f^\ast}{f_\ast}\mod{\structsymbol'},$$
	that is, the pullback $f^\ast$ along $f$ is left adjoint to the pushforward $f_\ast$ along $f$.
\end{lemma}

\begin{corollary}[Change of rings]
	\label{corollary: change of rings}
	
	Let $\site C$ be an essentially small site and $\phi\colon \structsymbol \to \structsymbol'$ be a morphism of rings on $\site C$. There is an adjoint pair of functors
		$$\mod {\structsymbol'} \xleftrightarrows{\scalarExtension{\phi}}{\scalarRestriction{\phi}} \mod\structsymbol$$
	where $\scalarExtension{\phi}\sheaf G \ldef \structsymbol' \otimes_\structsymbol \sheaf G$. The functors $\scalarExtension{\phi}$ and $\scalarRestriction{\phi}$ are called \emph{scalar extension} and \emph{scalar restriction} respectively.
\end{corollary}
\begin{proof}
	This trivially follows from lemma \ref{lemma: pullback and pushforward of sheaves of modules} observing that $(\id_\site C, \phi)\colon (\site C,\structsymbol')\to(\site C, \structsymbol)$ is a morphism of ringed sites.
\end{proof}

\begin{remark}[Pullback and pushforward of small sheaves]
	As in remark \ref{remark: direct and inverse images of (pre)sheaves of modules}, for translating lemma \ref{lemma: pullback and pushforward of sheaves of modules} to small sheaves, we restrict to the special case of $f^\top = \id_\extremallyDisconnectedSets$ being the identity. In this case the proof immediately follows from \cite[\href{https://stacks.math.columbia.edu/tag/03CZ}{Tag 03CZ}]{stacks-project}. This is then also precisely corollary \ref{corollary: change of rings} for the site $\extremallyDisconnectedSets$.
\end{remark}

It will prove useful later on to know that scalar restriction is always exact.
\begin{lemma}[Scalar restriction is exact]
	Let $\site{C}$ be an essentially small site and $\phi\colon \structsymbol \to \structsymbol'$ be a morphism of rings on $\site{C}$. Then scalar restriction $\scalarRestriction{\phi}\colon \mod{\structsymbol'}\to\mod{\structsymbol}$ is exact.
\end{lemma}
\begin{proof}
	As a right adjoint, $\scalarRestriction{\phi}$ certainly is left exact. It remains to show that $\scalarRestriction{\phi}$ is right exact. Recall from proposition \ref{proposition: the abelian category of sheaves of modules} that both categories of modules are abelian and that exactness can be verified on the level of abelian sheaves. The claim then follows from \cite[\href{https://stacks.math.columbia.edu/tag/04DB}{Tag 04DB}]{stacks-project} after observing that the geometric morphism of topoi $\Sheaves{\site{C}}{\Set} \xleftrightarrows{}{} \Sheaves{\site{C}}{\Set}$ induced by $\id_\site{C}\colon \site{C}\to \site{C}$ trivially satisfies the four equivalent conditions. 
\end{proof}

\begin{lemma}[Free sheaves]
	\label{lemma: free sheaves}
	
	Let $\site C$ be an essentially small site. 
	\begin{enumerate}
		\item
		The forgetful functor $\Sheaves{\site C}{\Ab} \to \Sheaves{\site C}{\Set}$ admits a left adjoint $\free {\constantSheaf{\Z}} -$, such that
			$$\free {\constantSheaf{\Z}} T = \sheafification{U \mapsto \free \Z {T(U)}},$$
		for any sheaf of sets $T$. We call $\free {\constantSheaf{\Z}} T$ the \emph{free abelian group on $T$}.
		
		\item
		Let $\structsymbol$ be a sheaf of rings on $\site C$. The forgetful functor $\mod\structsymbol \to \Sheaves{\site C}{\Set}$ admits a left adjoint $\free \structsymbol -$, such that
			$$\free \structsymbol T = \structsymbol \otimes_{\constantSheaf\Z} \free {\constantSheaf \Z} T$$
		is the extension of scalars along the unique map $\constantSheaf{\Z} \to \structsymbol$, for any sheaf of sets $T$ . We call $\free \structsymbol T$ the \emph{free $\structsymbol$-module on $T$}.
	\end{enumerate}
\end{lemma}
\begin{proof}
	Consider the functor $\free{\associated{\Z}^p}{-}\colon\Presheaves{\site{C}}{\Set} \to \Presheaves{\site{C}}{\Ab}$ given by $\sheaf{F}\mapsto \mleft(U\mapsto \free{\Z}{\sheaf{F}(U)}\mright)$. Let $T$ be a presheaf of sets and $A$ be a presheaf of abelian groups. Then there is a natural morphism
		$$\Hom_\Presheaves{\site{C}}{\Ab}(\free{\associated{\Z}^p}{T}, A) \to \Hom_\Presheaves{\site{C}}{\Set}(T, A)$$
	given by mapping a natural transformation $\eta$ to the natural transformation with component $T(U)\hookrightarrow\free{\Z}{T(U)} \xrightarrow{\eta_U} A(U)$ at $U\in\site{C}$. By the universal property of $\free{\Z}{-}$ this natural morphism admits and inverse, hence $\free{\associated{\Z}^p}{-}$ is a left adjoint to the forgetful functor $\Presheaves{\site{C}}{\Ab} \to \Presheaves{\site{C}}{\Set}$. Define $\free{\associated{\Z}}{-}\ldef \free{\associated{\Z}^p}{-}^\sheafificationsymbol\colon \Sheaves{\site{C}}{\Set}\to\Sheaves{\site{C}}{\Ab}$. Now if $T$ is a sheaf of sets and $A$ is an abelian sheaf, then
		$$\Hom_\sheavesSymbol(\free{\associated{\Z}}{T}, A) \cong \Hom_\preSheavesSymbol(\free{\associated{\Z}^p}{T}, A) \cong \Hom_\preSheavesSymbol(T, A)=\Hom_\sheavesSymbol(T, A)$$
	showing that $\free{\associated{\Z}}{-}$ is the desired adjoint. Denote by $\phi$ the unique morphism of rings $\associated{\Z}\to \structsymbol$. Since the forgetful functor from sheaves of $\structsymbol$-modules to sheaves of sets factors over the forgetful functor from sheaves of abelian groups to sheaves of sets, it is also clear that $\free{\structsymbol}{-} \ldef \scalarExtension{\phi}\free{\associated{\Z}}{-} = \structsymbol\otimes_{\associated{\Z}}\free{\associated{\Z}}{-}$ is left adjoint to $\mod{\structsymbol}\to\Sheaves{\site{C}}{\Set}$.
\end{proof}

\begin{remark}[Free small sheaves]
	\label{remark: free small sheaves}
	
	Lemma \ref{lemma: free sheaves} translates well to small sheaves, as the proof only relies on the existence of sheafification as an adjoint to the forgetful functor. It is also clear that for a small (pre)sheaf $\sheaf F$, the presheaf $\free{\associated{\Z}^p}{\sheaf{F}}$ is again small, as the analogous functors for $\boundedExtremallyDisconnectedSets{\kappa}$ commute with the transition functors $(-)^\enlargementsymbol{\kappa}{\kappa'}$.
\end{remark}

\section{Condensed Algebraic Structures}
\label{section: condensed algebraic structures}

%

We can now start to investigate how condensed sets synergize with algebraic structures. As already seen for abelian (pre)sheaves and (pre)sheaves of rings, we find that internal groups/rings/\dots in (pre)sheaves of sets precisely correspond to sheaves of groups/rings/\dots on the same site (at least for reasonably well-behaved classes of internal objects). This motivates the following definition.
\begin{definition}[Condensed groups/rings/\dots]
	\label{definition: condensed groups/rings/...}
	
	For $\category C$ the category of groups/rings/\dots, the category of \emph{condensed groups/rings/\dots} $\condensed{\category C}$ is the category of internal group/ring/\dots objects of $\condensedSets$ where the morphisms are those of the underlying condensed sets, that are compatible with the internal structures.
	
	Equivalently this is the full subcategory of the sheaf category $\Sheaves{\extremallyDisconnectedSets}{\category C}$ consisting of those sheaves $T\colon \extremallyDisconnectedSets^\op\to \category C$ for which there is some uncountable strong limit cardinal $\kappa$, such that at $T$ the component of the counit $( T\resTo{\kappa})^\enlargementsymbol{\kappa}{\infty}\to T$ is an isomorphism.
\end{definition}

\begin{definition}[Associated condensed groups/rings/\dots]	
	
	Let $X$ be a topological group/ring/\dots, \ie an internal group/ring/\dots in the category $\Top$, whose underlying topological space is $\tOne$. Then $\associated X$ is a condensed set by theorem \ref{theorem: some correspondences} and inherits the structure of an internal group/ring/\dots in $\condensedSets$, as the functor $X\mapsto \associated{X}$ has a left adjoint, hence commutes with all limits (and internal objects are most commonly given by certain maps from limits). Indeed, for any extremally disconnected $S$ the set of sections $\associated{X}(S) = \Hom_\Top(S, X)$ is the set of continuous functions $S\to X$, which inherits the structure of a group/ring/\dots from $X$ by pointwise addition/multiplication/\dots of functions. Akin to the topological case, if no confusion is possible, we will simply write $X$ for its condensed analogue $\associated{X}$.
\end{definition}

We obtain the following special case.
\begin{definition}[Discrete objects]
	If $X$ is a group/ring/\dots, then $X$ viewed as a discrete space is naturally a $\tOne$-topological group/ring/\dots. Thus, the associated condensed set $\associated{X}$ to this discrete space is naturally a condensed group/ring/\dots. Any condensed group/ring/\dots isomorphic to such a condensed object will be called \emph{discrete}.
\end{definition}

We are now well-equipped to work with condensed rings and modules thereof.
\begin{remark}[Condensed rings and modules]
	We obtain from definition \ref{definition: condensed groups/rings/...} the notion of \emph{condensed rings} -- these are precisely sheaves of rings on the site $\extremallyDisconnectedSets$. Given such a condensed ring $A$ we then obtain the ringed site $(\extremallyDisconnectedSets, A)$ and the corresponding notion of sheaves of $A$-modules. From now on out we shall call them \emph{(condensed) $A$-modules}. With our efforts from sections \ref{section: rings, ringed spaces and ringed sites}, \ref{section: (pre)sheaves of modules} and \ref{section: modules under morphisms} we thus obtain a well-rounded theory of condensed $A$-modules. Recall 
	\begin{itemize}
		\item
		that there is a sheafification from presheaves of $A$-modules to condensed $A$-modules left adjoint to the forgetful functor,
		
		\item
		that the category $\mod{A}$ is closed symmetric monoidal, giving us a tensor product $-\otimes_A -$ and an internal $\intHom_A$ of condensed $A$-modules (with the unit being $A$),
		
		\item 
		that the tensor-$\Hom$-adjunction enriches, giving isomorphisms
			$$\intHom_A(L \otimes_A M, N) \cong \intHom_A(L, \intHom_A(M, N))$$
		natural in the condensed $A$-modules $L$, $M$ and $N$,
		
		\item
		that for any morphism $\phi\colon A\to A'$ of condensed rings there is a notion of scalar restriction $\scalarRestriction{\phi}$ and scalar extension $\scalarExtension{\phi}$ that are part of an adjunction
			$$\mod{A'}\xleftrightarrows{\scalarExtension{\phi}}{\scalarRestriction{\phi}}\mod{A}$$
		where scalar restriction $\scalarRestriction{\phi}$ is exact,
		
		\item
		and that there are free condensed $A$-modules, giving a functor $\free{A}{-}\colon \condensedSets \to \mod{A}$ left adjoint to the forgetful functor $\mod{A}\to\condensedSets$.
	\end{itemize}
\end{remark}

An important special case is the following.
\begin{definition}[Condensed abelian groups]
	Observe that condensed abelian groups as per definition \ref{definition: condensed groups/rings/...} are abelian group objects internal to condensed sets, \ie small abelian sheaves. As these are the fundamental building blocks of algebra their category deserves its own name. From now on denote by $\condensedAb$ the category of condensed abelian groups. Observe further that $\condensedAb$ and $\mod{\associated{\Z}}$ are equivalent (in fact isomorphic) and hence that the category $\condensedAb$ enjoys all the good properties of a category of condensed modules. In particular, we will see in theorem \ref{theorem: condensed modules form an abelian category} that it is an abelian category generated by the class of free condensed abelian groups $\free{\associated{\Z}}{S}$ (for $S$ extremally disconnected) which are compact projectives.
\end{definition}

We can further investigate the internal $\intHom$ and notice that it indeed provides an enrichment of the usual $\Hom$.
\begin{remark}[Internal $\Hom$s are enriched]
	\label{remark: internal Homs are enriched}
	
	Let $A$ be a condensed ring and $S$ extremally disconnected. For any two condensed $A$-modules $M$ and $N$ the $S$-sections of $\intHom_A(M, N)$ can be described using Yoneda's lemma and some adjunctions. Indeed, we obtain
	\begin{align*}
		\intHom_A(M, N)(S) &\cong \Hom(\associated{S}, \intHom_A(M, N))\\
		&\cong \Hom_A(\free A {\associated{S}}, \intHom_A(M, N))\\
		&\cong \Hom_A(\free A  {\associated S}\otimes_A M, N)
	\end{align*}
	and in particular the underlying set of $\intHom_A(M, N)$ is
		$$\intHom_A(M, N)(\ast)\cong \Hom_A(M, N),$$
	as $\free A  \ast \cong A$. Thus $\intHom_A$ provides a natural enrichment of $\Hom_A$.
\end{remark}

A useful result that allows us to understand many $\intHom$s is given by the following lemma.
\begin{lemma}[Internal $\Hom$ -- but explicit]
	\label{lemma: internal hom -- but explicit}
	
	Let $M$ and $N$ be topological modules over the discrete ring $A$ and endow the $A$-module $\intHom_A(M, N)$ with the compact open topology.
	
	If $M$ and $N$ are Hausdorff and $M$ is compactly generated, then there is a natural isomorphism
		$$\intHom_{\associated{A}}(\associated{M}, \associated{N}) \cong \associated{\intHom_A(M, N)}$$
	of condensed $\associated{A}$-modules.
\end{lemma}
\begin{proof}
	For condensed abelian groups (\ie $R=\Z$) this is \cite[Prop. 4.2]{condenseddotpdf}. The general case is however analogous after observing that $\free{\associated{R}}{-}$ is (strong) monoidal. 
\end{proof}

\begin{remark}[Internal $\Hom$s of locally compact Hausdorff modules]
	Recall that a topological space $X$ is called \emph{locally compact} if any point $x\in X$ has a basis of compact Hausdorff neighborhoods. If $M$ and $N$ are locally compact and Hausdorff, the requirements of lemma \ref{lemma: internal hom -- but explicit} are satisfied, as any locally compact space is compactly generated. This is in particular the case if $M$ and $N$ are discrete, as any discrete space is locally compact and Hausdorff.
\end{remark}

\begin{notation}[Simplifying notation]
	Akin to writing $T$ for the condensed set $\associated{T}$ associated to some topological space or topological abelian group, due to lemma \ref{lemma: internal hom -- but explicit} we are in many cases justified in simply writing $\intHom_A(M, N)$ instead of $\intHom_{\associated{A}}(\associated{M}, \associated{N})$ or $\associated{\intHom_A(M, N)}$ when no confusion is possible.
\end{notation}



We now have the following major theorem, which characterizes the properties of categories of condensed modules.
\begin{theorem}[Condensed modules form an abelian category]
	\label{theorem: condensed modules form an abelian category}
	
	Suppose that $A$ is a condensed ring. The category $\mod{A}$ of condensed $A$-modules is abelian and satisfies Grothendieck's axioms \grothendieckaxiom{$3$}, \grothendieckaxiom{$3^\ast$}, \grothendieckaxiom{$4$} \grothendieckaxiom{$4^\ast$}, \grothendieckaxiom{$5$} and \grothendieckaxiom{$6$}. That is, in order:
	\begin{itemize}
		\item all colimits exist,
		\item all limits exist,
		\item arbitrary direct sums exist and are exact,
		\item arbitrary direct products exist and are exact,
		\item filtered colimits are exact,
		\item and for any set $J$, filtered categories $I_j$ for $j\in J$ and diagrams $i\mapsto M_i$ in $\mod{A}$ the natural map
			$$\colimit_{i\in\prod_{j\in J} I_j} \prod_{j\in J} M_{i_j} \to \prod_{j\in J} \colimit_{i_j \in I_j} M_{i_j}$$
		is an isomorphism.
	\end{itemize}
	Furthermore, the collection of free modules $\free{A}{S}$ with $S$ extremally disconnected forms a class of compact projective generators of $\mod{A}$.
\end{theorem}
\begin{proof}
	Analogous to \cite[\href{https://stacks.math.columbia.edu/tag/03DA}{Tag 03DA}]{stacks-project} and \cite[\href{https://stacks.math.columbia.edu/tag/03DB}{Tag 03DB}]{stacks-project}, the category $\mod{A}$ is abelian, admits all limits and colimits and the forgetful functor $\mod{A}\to\smallSheaves{\extremallyDisconnectedSets}{\Ab}=\condensedAb$ is exact and commutes with both limits and colimits. Thus, to show that $\mod{A}$ satisfies the posited axioms, it is enough to show them in case of $\condensedAb \cong \mod{\associated{\Z}}$.
	
	Observe that in $\Ab$ any finite product is already a direct sum, so a coproduct. Let $i\mapsto T_i$ be a diagram of condensed set of shape $I$ and denote by $T = \colimit_{\preSheavesSymbol, i} T_i$ the colimit of $i\mapsto T_i$ considered as a diagram of presheaves. Then for any finite disjoint union $S \ldef \coprod_{j\in J} S_j$, the natural map $T(\coprod_{j\in J} S_j) \to \prod_{j\in J} T(S_j)$ is precisely the map
		$$(\colimit\nolimits_{\preSheavesSymbol, i} T_i)\left(\coprod\nolimits_{j\in J} S_j\right) \isorightarrow \colimit_i \bigoplus_{j\in J} T_i(S_j) \isorightarrow \bigoplus_{j\in J} \colimit_i T_i(S_j)=\prod_{j\in J} (\colimit\nolimits_{\preSheavesSymbol, i} T_i)(S_j)$$
	which is an isomorphism since each $T_i$ satisfies the simplified sheaf condition from lemma \ref{lemma: simplified sheaf condition for extremally disconnected sets} and colimits commute with coproducts. In particular $T$ already satisfies the same simplified sheaf condition so is already a colimit of $i\mapsto T_i$ of sheaves. That is to say both limits and colimits in $\condensedAb$ are calculated pointwise. Hence, all of the axioms follow at once from the corresponding statements in $\Ab$.
	
	Thus, it remains to show that the free modules $\free{A}{S}$ for extremally disconnected sets $S$ form a collection of compact projective generators of $\mod{A}$. By adjointness and the Yoneda lemma we obtain the isomorphism
		$$\Hom_A(\free{A}{S}, M)\cong \Hom(S, M) \cong M(S)$$
	natural in the condensed $A$-module $M$. Thus, the functor $\Hom_A(\free{A}{S}, -)$ commutes with a given (co)limit if and only if $\Gamma(S, -)\colon M\mapsto M(S)$ does. But as previously shown, limits and colimit of condensed abelian groups and hence of condensed $A$-modules can be computed pointwise, implying that $\Gamma(S, -)$ and hence $\Hom_A(\free{A}{S}, -)$ commutes with all limits and colimits! Thus, the objects $\free{A}{S}$ are compact projective. It remains to show that they generate.
	
	For this, let $M$ be any condensed $A$-module. By theorem \ref{theorem: the pretopos of condensed sets}, associated condensed sets $\associated{S}$ of extremally disconnected sets $S$ generate the pretopos of condensed sets. Choose a separating set $\mathcal G'$ of such associates for the condensed set $M$. Consider the morphism of condensed $A$-modules
		$$p\colon \hat M\ldef \bigoplus_{\substack{\associated{S}\to M\\\associated{S}\in\mathcal G'}} \free{A}{\associated{S}}\to M$$
	where any $\associated{S}\to M$ induces a morphism $\free{A}{\associated{S}} \to M$ by adjointness. Since $\mathcal G'$ is small and the morphisms $\associated{S}\to M$ are (by Yoneda) in bijection with the set $M(S)$, this direct sum is small. If we can show that $p$ is an epimorphism, then by lemma \ref{lemma: generating classes in abelian categories} generation follows.
	
	Suppose that $f, g\colon M\to M'$ are morphisms such that $f\circ p = g\circ p$. Observe that for any $\associated{S}\in\mathcal G'$ and $h\colon \associated{S}\to M$ we have the commuting triangle
	$$\begin{tikzcd}
		\associated{S}\ar[rr, "\eta_S"]\ar[dr, "h", swap] &[-1em]&[-1em] \free{A}{\associated{S}}\ar[dl, "\hat h"]\\
		&M
	\end{tikzcd}$$	
	where $\eta$ is the unit of the adjunction between $\free{A}{-}$ and the forgetful functor $\mod{A}\to\condensedSets$ and $\hat h$ is the unique extension of $h$ to the free $A$-module. Denote by $i\colon \free{A}{\associated{S}}\to \hat M$ the inclusion of $\free{A}{\associated{S}}$ into the direct sum given by $h$. Then by definition of $p$ the triangle
	$$\begin{tikzcd}
		\free{A}{\associated{S}}\ar[rr, "i"]\ar[dr, "\hat h", swap] &[-1em]&[-1em] \hat M\ar[dl, "p"]\\
		&M
	\end{tikzcd}$$
	commutes as well. Overall we thus have
		$$f\circ h = f\circ \hat h\circ \eta_S=f\circ p \circ i \circ \eta_S = g\circ p \circ i \circ \eta_S = g\circ \hat h\circ \eta_S = g\circ h.$$
	But as $h\colon \associated{S}\to M$ and $\associated{S}\in\mathcal G'$ were arbitrary, we must have $f=g$ since $\mathcal G'$ is separating for $M$. Hence, for any parallel pair $f, g$ with $f\circ p = g\circ p$ we already have $f=g$, so $p$ is an epimorphism.
\end{proof}

\begin{remark}[About the Grothendieck-axioms]
	It is quite unusual for direct products in sheaf categories to be exact. For example \cite[Theorem 1.1]{exactness-direct-products} asserts that for a (divisorial Noetherian) scheme $(X, \struct{X})$ the category of (quasi-coherent) sheaves has exact direct products if and only if $(X, \struct{X})$ is affine. Clausen and Scholze state further, in \cite[Remark 1.11]{condenseddotpdf}, that abelian sheaves on any (essentially small) site form an abelian category satisfying \grothendieckaxiom{$3$}, \grothendieckaxiom{$4$}, \grothendieckaxiom{$5$} and \grothendieckaxiom{$3^\ast$}, but that the category only rarely satisfies \grothendieckaxiom{$4^\ast$} and \grothendieckaxiom{$6$} -- however these two axioms are a direct consequence of the existence of a class of compact projective generators.
\end{remark}

\begin{remark}[Not enough injectives]
	\label{remark: not enough injectives}
	
	We now know, that for any condensed ring $A$, the category $\mod{A}$ admits a class of compact projective generators, in particular that it has enough projectives. What about injectives? It turns out that $\condensedAb$ (and hence by extension $\mod{A}$) admits no non-zero injective objects, as shown in \cite[\href{https://mathoverflow.net/a/356261/546808}{this}]{injectives} MathOverflow answer of Scholze! It is interesting to note that for each cardinal $\kappa$ the abelian category of sheaves $\Sheaves{\boundedExtremallyDisconnectedSets{\kappa}}{\Ab}$ does admit enough injectives (this is a standard fact for abelian objects in a topos), however that these injectives are not stable under the transitions $(-)^\enlargementsymbol{\kappa}{\kappa'}$ -- the relevant extension property true in the $\kappa$-condensed world fails in the $\kappa'$-condensed world for large enough $\kappa'$, as one needs to lift against many more monomorphism not visible on the $\kappa$-condensed level.
	
	These missing injectives indeed provide an initial hindering in defining right derived functors in a condensed algebraic world. The author is not sure if and when right derived functors exist in general, but for this thesis only the various cases of right derived $\Hom$ functors are relevant. In this special case one can solve this rather pragmatically. One simply resolves the first argument by projectives while keeping the second argument fixed. One obtains a bifunctor this way. It is however unclear to the author if this yields the correct derived functor also for the second argument. We will come back to this issue in remark \ref{remark: about the notation and choices}.
\end{remark}

We will now revisit remark \ref{remark: a shift in perspective}.
\begin{remark}[The relief]
	\label{remark: the relief}
	
	Recall that one motivation for introducing condensed algebra was that \quote{algebra does not play well with topology}. A concrete instance of the phenomenon was that the category $\TopAb$ of topological abelian groups (equivalently the category of abelian group objects of $\Top$) is not abelian. Indeed, in $\TopAb$ the continuous homomorphism
		$$\id_\R\colon \R^\delta \to \R$$
	from the discrete to the metric topology is not an isomorphism, however both its kernel and cokernel are trivial: the respective universal properties are directly verified. For any map $f$ with $f\circ \id = 0$ also $f=0$ and similarly if $g\circ \id = 0$ then $g=0$ -- both maps must factor over the zero-module. Now a morphism not being an isomorphism and yet still having trivial kernel and cokernel is something that is impossible in an abelian category. Conclusion: $\TopAb$ is not abelian.
	
	If we now pass to condensed abelian groups (which do form an abelian category!) and consider the corresponding morphism $\associated{\R^\delta} \to \associated{\R}$ (both spaces are \CGWH and hence well-represented) then it is still not an isomorphism (and surely it should not be, the spaces are not homeomorphic) but(!) we obtain a short exact sequence
		$$0\to \associated{\R^\delta}\to \associated{\R} \to Q \to 0$$
	where $Q$ is the sheaf (recall from the proof of theorem \ref{theorem: condensed modules form an abelian category} that limits and colimit can be calculated as presheaves and hence pointwise) given by mapping $S\mapsto \associated{\R}(S)/\associated{\R^\delta}(S)$. Here $\associated{\R}(S)$ is the abelian group of continuous maps $S\to \R$ and $\associated{\R^\delta}(S)$ is the abelian group of locally constant maps $S\to \R$. As $\R^\delta\to \R$ is not a homeomorphism the induced morphism $\associated{\R^\delta} \to \associated{\R}$ is not an isomorphism and hence it is clear that $Q$ is non-trivial. However, the underlying space $\topologizedUnderlying{Q}$ of $Q$ is! Thus, $Q$ \quote{remembers some infinitesimal information about the map $\id_\R$ albeit having a trivial underlying space}.
	
%
\end{remark}

	\chapter{Cohomology \& Derived Condensed Algebra}
	\label{chapter: cohomology and derived condensed algebra}
	
\section*{Introduction}

In this chapter we will try to develop some of the derived theory needed to handle condensed phenomena later on, \eg in the definition of analytic rings in definition \ref{definition: analytic rings}. For the derived theory of sheaves of modules on a scheme (a setting with plenty of injectives), we generically refer to the \cite[The Stacks Project]{stacks-project} which covers all needed formalism extensively.

Why does the author write \emphquote{try to}? In the main reference for the material covered in this thesis -- the lecture notes \cite{condenseddotpdf} -- Clausen and Scholze are quite lenient in their use of derived algebra. They seldomly show that objects (like derived functors) exist and why they satisfy certain properties. While they occasionally give short hints, \eg that one can construct (left) derived tensor products and (right) derived internal $\Hom$s by \quote{[u]sing projective resolutions} \cite[p. 13 (iii)]{condenseddotpdf}, these hints are certainly not comprehensive. It is to note that Clausen and Scholze (in their lectures) seem to mainly work in the language of $\infty$-categories, implicit in the first chapters but explicit in \cite[Lectures IX, X \& XI]{condenseddotpdf} where the derived categories of interest are explicitly considered as \emph{stable $\infty$-categories} (which then allows one to glue them). The author is sure that all claims made by Clausen and Scholze (either implicitly or explicitly) are \quote{obvious enough} when one works (and is used to) $\infty$-categorical language -- this however does not apply to the author. 

As there are (too) many things to say about $\infty$-land for one thesis, the author decided on a more straightforward approach. We will try to introduce all necessary objects and prove their relevant properties on the fly, using more \quote{traditional} theory. This however did not work in all cases! This will be indicated (obviously) at each such place. We will however still use these results later on, even if we do not have a proof. Then with the peace in mind that Clausen and Scholze can provide a proof in $\infty$-world.


\section{Prerequisites}

We assume the reader familiar with the basic notions and definitions of (bounded -- either to the left, to the right, or in both directions) derived algebra. This includes
\begin{itemize}
	\item 
	the definition of \emph{co-chain complexes} in some abelian category $\category{A}$ (which we will simply call \emph{chain complexes} or \emph{complexes} as we will solely work with such) and the definition of \emph{chain maps} as their morphisms, constituting a category $\complexes{\category{A}}$ (see \cite[\href{https://stacks.math.columbia.edu/tag/010V}{Tag 010V}]{stacks-project} for both),
	
	\item 
	the notion of \emph{homotopy of chain maps} (see \cite[\href{https://stacks.math.columbia.edu/tag/0112}{Tag 0112}]{stacks-project}) which leads to the formation of the \emph{homotopy category of complexes} $\homotopyCategoryOfComplexes{\category{A}}$ (see \cite[\href{https://stacks.math.columbia.edu/tag/05RN}{Tag 05RN}]{stacks-project}),
	
	\item
	the definition of \emph{triangulated categories} (see \cite[\href{https://stacks.math.columbia.edu/tag/05QK}{Tag 05QK}]{stacks-project}) with example $\homotopyCategoryOfComplexes{\category{A}}$ (see \cite[\href{https://stacks.math.columbia.edu/tag/014S}{Tag 014S}]{stacks-project}), as well as the notion of \emph{exact/triangulated functor} (see \cite[\href{https://stacks.math.columbia.edu/tag/014V}{Tag 014V}]{stacks-project}) between such categories,
	
	\item
	the notion of \emph{acyclicity}, the definition of \emph{quasi-isomorphisms} (see \cite[\href{https://stacks.math.columbia.edu/tag/010Z}{Tag 010Z}]{stacks-project} for both) and the general theory of \emph{localization of triangulated categories} at a multiplicative system compatible with the triangulated structure (see \cite[\href{https://stacks.math.columbia.edu/tag/05R1}{Tag 05R1}]{stacks-project}), leading to the definition of the \emph{derived category} $\derived{\category{A}}$ as the localization of the triangulated category $\homotopyCategoryOfComplexes{\category{A}}$ at the multiplicative system $\Sigma$ of all quasi-isomorphisms of complexes in $\category{A}$ (see \cite[\href{https://stacks.math.columbia.edu/tag/05RU}{Tag 05RU}]{stacks-project}),
	
	\item
	the definitions of \emph{left \& right resolutions} and the notions of \emph{having enough projectives} and \emph{having enough injectives} (see \cite[\href{https://stacks.math.columbia.edu/tag/0643}{Tag 0643}]{stacks-project} for the projective case, the injective case is formally dual),
	
	\item
	the definition of \emph{left and right derived functors} as the natural extension of an exact functor between triangulated categories to a localization of its domain (see \cite[\href{https://stacks.math.columbia.edu/tag/05S7}{Tag 05S7}]{stacks-project}) and how well chosen left and right resolutions allow to \emph{compute} (see \cite[\href{https://stacks.math.columbia.edu/tag/05SX}{Tag 05SX}]{stacks-project} and \cite[\href{https://stacks.math.columbia.edu/tag/06XN}{Tag 06XN}]{stacks-project}) the values of derived functors, \eg by left or right resolutions of projectives and injectives respectively.
\end{itemize}

\section{\texorpdfstring{$K$}{K}-projectivity}
\label{section: K-projectivity}

In this thesis our main application of the derived theory are of course categories of condensed objects, \eg condensed abelian group or condensed modules over a condensed ring. As already mentioned these usually do not admit enough (or rather any non-zero) injectives (see remark \ref{remark: not enough injectives}). Consequentially, we will restrict our attention to projectives (of which there are plenty since all these categories are generated by compact projectives). It is to note that the main reference for this section, \cite[The Stacks Project]{stacks-project}, works with a very restrictive class of abelian categories, see \cite[\href{https://stacks.math.columbia.edu/tag/09PA}{Tag 09PA}]{stacks-project} -- into which our categories of desire do not fit. As mentioned in the tag, if one wishes to work with a \emph{large} abelian category, one must ensure that the derived category, as a localization of the homotopy category of complexes $\homotopyCategoryOfComplexes{\category{A}}$, is locally small. This will be the content of corollary \ref{corollary: local smallness of the derived category}. Furthermore, in this thesis we will predominantly work with \emph{unbounded} derived categories. While in right bounded case its is enough to work with right bounded projective resolutions, in the unbounded case one has to be more careful. This leads to the notion of \emph{$K$-projectivity}, originally introduced by Spaltenstein in \cite{spaltenstein}. The author wants to thank \href{https://zmao-math.github.io/}{Zhouhang Mao} for introducing him to Spaltenstein's work. We will mainly follow The Stacks Project \cite{stacks-project} and only reference Spaltenstein's \cite{spaltenstein} if necessary.
\begin{definition}[K-projectivity, {\cite[\href{https://stacks.math.columbia.edu/tag/070H}{Tag 070H}]{stacks-project}}]
	
	Let $\category{A}$ be an abelian category. A complex $P^\bullet$ in $\category{A}$ is called \emph{$K$-projective} if for any acyclic complex $N^\bullet$ in $\category{A}$ we have that $\Hom_\homotopyCategoryOfComplexes{\category{A}}(P^\bullet, N^\bullet)\cong 0$.
\end{definition}

We obtain the following very useful characterization of K-projectivity.
\begin{lemma}[Morphisms out of a $K$-projective complex, {\cite[\href{https://stacks.math.columbia.edu/tag/070I}{Tag 070I}]{stacks-project}}]
	\label{lemma: morphism out of a K-projective complex}
	
	Let $\category{A}$ be an abelian category. A complex $P^\bullet$ in $\category{A}$ is $K$-projective if and only if for all complexes $N^\bullet$ in $\category{A}$ the natural map
	$$\Hom_\homotopyCategoryOfComplexes{\category{A}}(P^\bullet, N^\bullet) \to \Hom_\derived{\category{A}}(P^\bullet, N^\bullet)$$
	is a bijection.
\end{lemma}

An immediate consequence is that if the abelian category $\category{A}$ \quote{has enough $K$-projective resolutions} then the derived category $\derived{\category{A}}$ is locally small.
\begin{corollary}[Local smallness of the derived category]
	\label{corollary: local smallness of the derived category}
	
	Let $\category{A}$ be an abelian complex such that every complex admits a left resolution by a $K$-projective complex. Then $\derived{\category{A}}$ is locally small.
\end{corollary}
\begin{proof}
	Let $M^\bullet$ and $N^\bullet$ be arbitrary complexes. By assumption there is a left resolution $P^\bullet \to M^\bullet$ of $M^\bullet$ by a $K$-projective complex $P^\bullet$. Then by lemma \ref{lemma: morphism out of a K-projective complex} we obtain that
	$$\Hom_\derived{\category{A}}(M^\bullet, N^\bullet) \isorightarrow \Hom_\derived{\category{A}}(P^\bullet, N^\bullet) \isoleftarrow \Hom_\homotopyCategoryOfComplexes{\category{A}}(P^\bullet, N^\bullet)$$
	and thus that $\Hom_\derived{\category{A}}(M^\bullet, N^\bullet)$ is small.
\end{proof}

Furthermore, $K$-projective complexes \emph{compute} left derived functors.
\begin{lemma}[Existence of left derived functors, {\cite[\href{https://stacks.math.columbia.edu/tag/070K}{Tag 070K}]{stacks-project}}]
	\label{lemma: existence of left derived functors}
	
	Let $\category{A}$ be an abelian category such that every complex admits a left resolution by a $K$-projective complex and let $F\colon \homotopyCategoryOfComplexes{\category{A}} \to \category{D}'$ be an exact functor of triangulated categories. Then $\LeftD F\colon \derived{\category{A}} \to \category{D}'$ exists and for any $K$-projective complex $P^\bullet$ we obtain that the natural morphism $\LeftD F(P^\bullet) \to F(P^\bullet)$ is a quasi-isomorphism.
\end{lemma}

Spaltenstein now gives the following proposition that ensures the existence of \emph{enough $K$-projective complexes}. Recall that a \emph{sequential} (co)limit is a (co)limit over a diagram indexed by the partially ordered set $(\N, \le)$.
\begin{proposition}[Enough projectives on steroids]
	\label{proposition: enough projectives on steroids}
	
	If $\category{A}$ is an abelian category with enough projectives such that filtered colimits are exact, then every complex in $\category{A}$ admits a left resolution by a $K$-projective complex that is a sequential colimit of bounded to the right, term-wise projective complexes.
\end{proposition}
\begin{proof}
	By \cite[Example. 3.2 (a)]{spaltenstein} every bounded to the right complex in $\category{A}$ admits a left resolution by projectives and any such complex of projectives is $K$-projective. Denote by $\mathfrak P$ the class of bounded to the right complexes of projectives. Denote further the closure of $\mathfrak{P}$ under \emph{special direct limits} (see \cite[Def. 2.6]{spaltenstein}) by $\mathfrak P^\rightarrow$.
	
	Now since filtered colimits are exact, theorem \cite[Thm. 3.4]{spaltenstein} and its corollary \cite[Cor. 3.5]{spaltenstein} assert that every complex of $\category{A}$ admits a left $\mathfrak P^\rightarrow$-resolution of the desired form and that any such complex is $K$-projective.
\end{proof}

A useful result is the following.
\begin{lemma}[Preservation of $K$-projectivity, {\cite[\href{https://stacks.math.columbia.edu/tag/08BJ}{Tag 08BJ}]{stacks-project}}]
	\label{lemma: preservation of K-projectivity}
	
	Let $F\colon \category{A} \to \category{B}$ and $G\colon \category{B} \to \category{A}$ be a pair of additive functors such that $F$ is exact and left adjoint to $G$. Then $G$ preserves $K$-projectivity, \ie for any $K$-projective complex $P^\bullet$ the complex $G(P^\bullet)$ is $K$-projective as well.
\end{lemma}

\section[Derived Tensor Products and Derived Homs]{Derived Tensor Products and Derived $\Hom$s}

We will now work towards the definition and existence of left derived tensor products and right derived $\Hom$s. We first start off with introducing the tensor products of complexes.
\begin{definition}[Tensor products of complexes]
	\label{definition: tensor products of complexes}
	
	Let $\category{A}$ be a symmetric monoidal abelian category with all direct sums. Define a functor
	$$-\otimes_\category{A}^\bullet - \colon \homotopyCategoryOfComplexes{\category{A}} \times \homotopyCategoryOfComplexes{
		\category{A}}\to \homotopyCategoryOfComplexes{\category{A}}$$
	given on complexes $X^\bullet$ and $Y^\bullet$ in degree $n\in \Z$ by
	$$X^\bullet\otimes_\category{A}^n Y^\bullet \ldef \bigoplus_{p+q=n} X^p \otimes_\category{A} Y^q$$
	with differential
	$$d^n \ldef \bigoplus_{p+q=n} d_X^p\otimes \id_{Y^q} + (-1)^p \cdot \id_{X^p}\otimes d_Y^q.$$
	Denoting by $\mathbf{1}$ the unit object of $\category{A}$ it is clear that $-\otimes_\category{A}^\bullet -$ makes $\homotopyCategoryOfComplexes{\category{A}}$ into a symmetric monoidal category with unit given by $\mathbf{1}[0]$.
\end{definition}

We need an analogous construction for $\Hom$-sets. As there are several enrichments floating around at any given time, we give a general definition. The main specializations for us being the usual enrichment over abelian groups (as is standard for any abelian category), but also over condensed abelian groups and over condensed $A$-modules for some condensed ring $A$, recalling that both $\condensedAb$ and $\mod{A}$ are closed monoidal, \ie have an internal $\intHom$ enriching the usual $\Hom$.
\begin{definition}[Complexes of $\Hom$s]
	\label{definition: complexes of Homs}
	
	Let $\category{A}$ be an abelian category enriched over the monoidal abelian category $\category{V}$. Assume that $\category{V}$ admits all direct products. Define a functor
		$$\eHom_\category{A}^\bullet\colon \homotopyCategoryOfComplexes{\category{A}}^\op \times \homotopyCategoryOfComplexes{\category{A}}\to \homotopyCategoryOfComplexes{\category{B}}$$
	given on complexes $X^\bullet$ and $Y^\bullet$ in degree $n\in \Z$ by
		$$\eHom_\category{A}^\bullet(X^\bullet, Y^\bullet)\ldef \prod_{p+q=n} \eHom_\category{A}(X^p, Y^q)$$
	with differential
		$$d^n \ldef \prod_{p+q=n} (d_Y^q)_\ast + (-1)^{n+1} \cdot (d_X^{p-1})^\ast$$
	where as usual, $(-)_\ast$ and $(-)^\ast$ denote the images of morphisms under $\eHom_\category{A}$ when fixing the first respectively second argument.
\end{definition}

Tensor products of complexes and complexes of internal $\intHom$s inside a closed symmetric monoidal $\category{A}$ interact well and ensure that the homotopy category of complex $\homotopyCategoryOfComplexes{\category{A}}$ is closed monoidal as well.
\begin{remark}[Chain complexes are closed monoidal]
	\label{definition: chain complexes are closed monoidal}
	
	Let $\category A$ be closed symmetric monoidal abelian category with arbitrary direct sums and arbitrary direct products. Denote by $\mathbf{1}$ the unit object of $\category{A}$. Since all partial adjunctions $-\otimes_\category{A}^\bullet Y \ladj \intHom_\category{A}^\bullet(Y, -)$ for any complex $Y$ are compatible, we obtain that $\otimes_\category{A}^\bullet$ and $\intHom_\category{A}^\bullet$ make $\homotopyCategoryOfComplexes{\category{A}}$ into a closed symmetric monoidal category with unit given by $\mathbf{1}[0]$. Furthermore, we obtain an even stronger statement, namely that there is an isomorphism
	$$\intHom_\category{A}^\bullet(X \otimes_\category{A}^\bullet, Y, Z) \cong \intHom_\category{A}^\bullet(X, \intHom_\category{A}^\bullet(Y, Z))$$
	natural in all three complexes $X$, $Y$ and $Z$.
\end{remark}

We can now define what the left derived tensor product should be.
\begin{definition}[Left derived tensor products]
	\label{definition: left derived tensor products}
	
	Recall the setting of definition \ref{definition: tensor products of complexes}. Observe that for any morphism of complexes $Y\to Y'$ there is a natural transformation of functors $-\otimes_\category{A}^\bullet Y \Rightarrow -\otimes_\category{A}^\bullet Y'$ inducing a natural transformation of left derived functors $\LeftD(-\otimes_\category{A}^\bullet Y) \Rightarrow \LeftD(-\otimes_\category{A}^\bullet Y')$ in case that both exist. The \emph{left derived tensor product} is the functor
	$$-\Dotimes_\category{A} - \colon \derived{\category{A}} \times \derived{\category{A}} \to \derived{\category{A}}$$
	such that for any complex $Y$ the partial functor $-\Dotimes_\category{A} Y$ is the left derived functor $\LeftD(-\otimes_\category{A}^\bullet Y)$.
\end{definition}

In a similar fashion, we obtain the definition of right derived $\Hom$s.
\begin{definition}[Right derived $\Hom$s]
	\label{definition: right derived Homs}
	
	Recall the setting of definition \ref{definition: complexes of Homs}. Observe that for any morphism $Y\to Y'$ of complexes there is a natural transformation $\eHom_\category{A}^\bullet(-, Y) \Rightarrow \eHom_\category{A}^\bullet(-, Y')$ of functors $\homotopyCategoryOfComplexes{\category{A}}^\op \to \homotopyCategoryOfComplexes{\category{A}}$, inducing a natural transformation of right derived functors $\RightD\eHom_\category{A}^\bullet(-, Y)\Rightarrow \RightD\eHom_\category{A}^\bullet(-, Y')$ in case that both exist. The \emph{right derived $\Hom$} is the functor
		$$\RightD\eHom_\category{A}\colon \derived{\category{A}}^\op\times\derived{\category{A}}\to \derived{\category{B}}$$
	such that for any complex $Y$ the partial functor $\RightD\eHom_\category{A}(-, Y)$ is the right derived functor $\RightD\eHom_\category{A}^\bullet(-, Y)$.
\end{definition}

\begin{remark}[About the notation and choices]
	\label{remark: about the notation and choices}
	
	Recall the settings of definitions \ref{definition: left derived tensor products} and \ref{definition: right derived Homs}. One should now ask the following two questions: First, why not name the derived functors $\otimes^{\LeftD\bullet}_\category{A}$ and $\RightD\eHom_\category{A}^\bullet$ instead of dropping the $\bullet$? And secondly, why fix the second argument and not the first? After all that is an arbitrary choice.
	
	To answer the first question one has to notice that if either the first argument or second argument of $\otimes_\category{A}^\bullet$ or $\eHom_\category{A}^\bullet$ is an object of $\category{A}$ concentrated in degree $0$, let us assume the first argument $X\in\homotopyCategoryOfComplexes{\category{A}}$ is general, and the second argument is some $Y\in\category{A}$, then
	$$X\otimes_\category{A}^\bullet Y[0] \cong X\otimes_\category{A} Y[0]\text{ and }\eHom_\category{A}^\bullet(X, Y[0])\cong \eHom_\category{A}(X, Y[0])$$
	agree with the natural extensions of the functors $-\otimes_\category{A} Y$ and $\eHom_\category{A}(-, Y)$ to functors $\homotopyCategoryOfComplexes{\category{A}}\to\homotopyCategoryOfComplexes{\category{A}}$ and $\homotopyCategoryOfComplexes{\category{A}}^\op \to \homotopyCategoryOfComplexes{\category{B}}$. So in this case $- \Dotimes_\category{A} Y$ and $\RightD\eHom_\category{A}(-, Y)$ really are the derived functors expected from the notation. Because of this and since we are not really interested in $\otimes_{\category{A}}^\bullet$ and $\eHom_\category{A}^\bullet$ themselves, we simply drop the $\bullet$.
	
	To understand why we fix the second argument instead of the first, we should recall the case of modules over a discrete ring $R$. The resulting category $\mod{R}$ has both enough projectives and enough injectives. Then one might as well fix the first argument and resolve the second one by a \emph{(K-)injective} complex: Using a simple dimension shifting argument one shows by hand that the two notions agree and hence that either choice of argument yields isomorphic functors. As a consequence for either choice of argument the constructed derived functor must be also the derived functor of the other argument. An argument like this is not possible in a condensed world as there are simply not enough injectives. The author was unable to provide an alternative proof so we must make this an \assumptionPreWord.
\end{remark}


Using out previous results on $K$-projectivity we can now ensure the existence of the left derived tensor product.
\begin{corollary}[Existence of the left derived tensor product]
	Assume further in the setting of definition \ref{definition: left derived tensor products} that $\category{A}$ has enough projectives and that filtered colimits in $\category{A}$ are exact. Then the left derived tensor product $\Dotimes_\category{A}$ exists.
\end{corollary}
\begin{proof}
	This is direct consequence of lemma \ref{lemma: existence of left derived functors} once one knows that there are enough $K$-projective resolutions due to proposition \ref{proposition: enough projectives on steroids}.
\end{proof}

Similarly, we obtain the existence of right derived $\Hom$s.
\begin{corollary}[Existence of right derived $\Hom$s]
	\label{corollary: existence of right derived Homs}
	
	Assume further in the setting of definition \ref{definition: right derived Homs} that $\category{A}$ has enough projectives and that filtered colimits in $\category{A}$ are exact. Then the right derived functor $\RightD\eHom_\category{A}$ exists.
\end{corollary}
\begin{proof}
	Observe that for any complex $Y$ the right derived functor of $\eHom_\category{A}^\bullet(-, Y)\colon \homotopyCategoryOfComplexes{\category{A}}^\op \to \homotopyCategoryOfComplexes{\category{B}}$ is equivalently the left derived functor of $\eHom_\category{A}^\bullet(-, Y)\colon \homotopyCategoryOfComplexes{\category{A}} \to \homotopyCategoryOfComplexes{\category{B}}^\op\equiv \homotopyCategoryOfComplexes{\category{B}^\op}$. Then once again one knows that there are enough $K$-projective resolutions due to proposition \ref{proposition: enough projectives on steroids} giving the existence of the left derived functor by \ref{lemma: existence of left derived functors}.
\end{proof}

Denoting by $\eHom_\category{A}$ then enrichment of $\Hom_\category{A}$ to abelian groups, we see that $\eRHom_\category{A}$ actually enriches $\eHom_\derived{\category{A}}$. In particular in the derived category, to each $\Hom$-set there is a corresponding complex of abelian groups.
\begin{lemma}[$\eRHom$ enriches $\Hom$-sets in the derived category]
	\label{lemma: eRHom enriches Hom-sets in the derived category}
	
	Let $\category{A}$ be an abelian category that admits all direct products, has enough projectives and such that filtered colimits are exact. Then $\eRHom_\category{A}\colon \derived{\category{A}}^\op \times\derived{\category{A}}\to\derived{\Ab}$ exists by corollary \ref{corollary: existence of right derived Homs}, and there is a natural isomorphism
	$$H^0\eRHom_\category{A} \cong \eHom_\derived{\category{A}},$$
	providing a natural enrichment of $\eHom_\derived{\category{A}}$ by $\eRHom_\category{A}$.
\end{lemma}
\begin{proof}
	Let $X^\bullet$ and $Y^\bullet$ be complexes. Choosing a left resolution of $X^\bullet$ by a $K$-projective complex $P^\bullet$ we see that
	$$H^0\eRHom_\category{A}(X^\bullet, Y^\bullet) \cong H^0 \eHom_\category{A}^\bullet(P^\bullet, Y^\bullet) \cong \eHom_\homotopyCategoryOfComplexes{\category{A}}(P^\bullet, Y^\bullet).$$
	As $P^\bullet$ is $K$-projective we obtain by lemma \ref{lemma: morphism out of a K-projective complex} that
	$$\eHom_\homotopyCategoryOfComplexes{\category{A}}(P^\bullet, Y^\bullet)\cong \eHom_\derived{\category{A}}(P^\bullet, Y^\bullet) \cong \eHom_\derived{\category{A}}(X^\bullet, Y^\bullet),$$
	proving the claim.
\end{proof}

\section{Deriving Adjunctions}

In this section we will study when adjunctions are stable under derivation, in particular enriched adjunctions. If there is no enrichment we have the following general result.
\begin{lemma}[Adjunctions are stable under derivation, {\cite[\href{https://stacks.math.columbia.edu/tag/09T5}{Tag 09T5}]{stacks-project}}]
	\label{lemma: adjunctions are stable under derivation}
	
	Let $\category{A}$ and $\category{B}$ be abelian categories and $F\colon \category{A} \to \category{B}$ be left adjoint to $G\colon \category{B}\to\category{A}$. If $\LeftD F$ and $\RightD G$ exist, then there are isomorphisms
	$$\Hom_\derived{\category{B}}(\LeftD F(X^\bullet), Y^\bullet) \cong \Hom_\derived{\category{A}}(X^\bullet, \RightD G(Y^\bullet))$$
	natural in the complexes $X^\bullet \in \derived{\category{A}}$ and $Y^\bullet \in \derived{\category{B}}$.
\end{lemma}

An immediate consequence is that a closed monoidal structure on a (well-behaved) abelian category translates to the derived setting.
\begin{lemma}[Deriving a closed monoidal structure]
	\label{lemma: deriving a closed monoidal structure}
	
	Let $\category{A}$ be a closed symmetric monoidal abelian category. Assume that $\category{A}$ is \grothendieckaxiom{$5$} and has enough projectives. Then there is an isomorphism
		$$\Hom_\derived{\category{A}}(X\Dotimes_\category{A}Y, Z)\cong \Hom_\derived{\category{A}}(X, \intRHom_\category{A}(Y, Z))$$
	natural in the complexes $X$, $Y$ and $Z$ making $\derived{\category{A}}$ into a closed symmetric monoidal category. Furthermore, there even is a natural isomorphism
		$$\intRHom_\category{A}(X\Dotimes_\category{A}Y, Z)\cong \intRHom_\category{A}(X, \intRHom_\category{A}(Y, Z))$$
	of enrichments of $\Hom_\derived{\category{A}}$.
\end{lemma}
\begin{proof}
	By lemma \ref{lemma: adjunctions are stable under derivation} we obtain that for any complex $X$ the two functors $- \Dotimes_\category{A} X$ and $\intRHom_\category{A}(X, -)$ are adjoint. As these adjunctions for various $X$ are compatible it is easily verified that $\derived{\category{A}}$ is closed symmetric monoidal with unit $\mathbf{1}[0]$, the unit of $\category{A}$ concentrated in degree $0$.
	
	It is now purely formal that the adjunctions $-\Dotimes_\category{A}X\ladj \intRHom_\category{A}(X, -)$ enrich. Indeed, in any closed monoidal category $\category{V}$ and $X, Y, Z\in \category{V}$ we find that
	\begin{equation}
		\scalebox{0.95}{\parbox{\textwidth}{
			\begin{align*}
				\Hom(-, \intHom(X\otimes Y, Z))&\cong \Hom(-\otimes X\otimes Y, Z)\\
				&\cong \Hom(-\otimes X, \intHom(Y, Z))\cong \Hom(-, \intHom(X, \intHom(Y, Z)))
			\end{align*}
		}}
	\end{equation}
	implying by Yoneda's lemma that there is a natural isomorphism $\intHom(X\otimes Y, Z)\cong \intHom(X,\intHom(Y, Z))$ implying that $-\otimes X\ladj_\category{V} \intHom(X, -)$ for all $X \in \category{V}$.
\end{proof}

More generally we are also interested in translating enriched adjunctions $F\ladj_\category{V} G$ to the derived category. In this thesis the main application of this is the case of scalar extension and restriction where the adjunction is enriched over condensed abelian groups (or over the respective module categories). This however is somewhat precarious and the author only managed to prove that an enriched adjunction translates under quite strong assumptions (tailored to this specific use case). In the condensed setting we find that scalar extension maps compact projective generators to compact projective generators all while scalar restriction is exact. This helps a lot in applying Spaltenstein's techniques.
\begin{lemma}[Enriched adjunctions are stable under derivation]
	\label{lemma: enriched adjunctions are stable under derivation}
	
	Let $\category{A}$ and $\category{B}$ be abelian categories, both enriched over the monoidal abelian category $\category{V}$. Assume that all three categories are \grothendieckaxiom{$5$} with enough projectives. Let $F\colon \category{A} \to \category{B}$ and $G\colon \category{B}\to\category{A}$ be a pair of additive functors such that $F\ladj_\category{V} G$ is a $\category{V}$-enriched adjunction, such that $F$ maps projectives to projectives and such that $G$ is exact. Then in $\derived{\category{V}}$ there are isomorphisms
	$$\eRHom_\category{B}(\LeftD F(X^\bullet), Y^\bullet) \cong \eRHom_\category{A}(X^\bullet, \RightD G(Y^\bullet))$$
	natural in the complexes $X^\bullet \in \derived{\category{A}}$ and $Y^\bullet \in \derived{\category{B}}$.
\end{lemma}
\begin{proof}
	From $\eHom_\category{B}(F(-), -)\cong \eHom_\category{A}(-, G(-))$ it follows formally that
	$$\eHom_\category{B}^\bullet(\tilde F(-), -)\cong \eHom^\bullet_\category{A}(-, \tilde G(-))$$
	where $\tilde F$ and $\tilde G$ are the natural extensions of $F$ and $G$ to accommodate complexes. Denote by $\mathfrak{P}_\category{A}$ the class of bounded to the right complexes in $\category{A}$ which are term-wise projective and denote similarly $\mathfrak{P}_\category{B}$. Recall from \cite[Def. 2.6]{spaltenstein} the notion of a \emph{$\mathfrak{P}_\category{A}$-special direct system}. Suppose that $(P_n^\bullet)_{n\in E}$ is such a $\mathfrak{P}_\category{A}$-special direct system with colimit $P^\bullet$. Then $\tilde F(P^\bullet_n) \in \mathfrak{P}_\category{B}$ for any $n\in E$ as $F$ maps projectives to projectives by assumption. Now as $F$ is additive, it preserves split exact sequences in $\category{A}$ and hence $\tilde F$ preserves \emph{semi-split sequences} (see \cite[Def. 0.5]{spaltenstein}) of complexes in particular $\tilde F(P^\bullet_{n-1})\to\tilde F(P^\bullet_n)$ is injective as well, for any $n\in E$ which has a predecessor $n-1$. We thus obtain that $(\tilde F(P^\bullet_n))_{n\in E}$ is a $\mathfrak{P}_\category{B}$-special direct system. As $\tilde F$ is left adjoint to $\tilde G$ and thus preserves colimits, we obtain that $\tilde F(P^\bullet) = \tilde F(\colimit_n P^\bullet_n) \cong \colimit_n \tilde F(P^\bullet_n)$ is a colimit of a $\mathfrak{P}_\category{B}$-special direct system and hence must be $K$-projective.
	
	Let $X^\bullet \in \derived{\category{A}}$ and $Y^\bullet \in \derived{\category{B}}$. Choose a left resolution $P^\bullet \to X^\bullet$ by a $K$-projective complex $P^\bullet$ that is the colimit of a $\mathfrak{P}_\category{A}$-special direct system. Then by the above argument $\tilde F(P^\bullet)$ is $K$-projective and we obtain that
	\begin{align*}
		\eRHom_\category{B}(\LeftD F(X^\bullet), Y^\bullet)
		&\cong\eRHom_\category{B}(\tilde F(P^\bullet), Y^\bullet)\\
		&\cong\eHom_\category{B}^\bullet(\tilde F(P^\bullet), Y^\bullet)\\
		&\cong \eHom_\category{A}^\bullet(P^\bullet, \tilde G(Y^\bullet)) = \eRHom_\category{A}(X^\bullet, \RightD G(Y^\bullet))
	\end{align*}
	proving the claim.
\end{proof}

While the previous results cover many cases of interest, occasionally we need the following more general assumption.
\begin{assumption}[Derived adjunctions are enriched]
	\label{assumption: derived adjunctions are enriched}
	
	Let $\category{A}$ and $\category{B}$ be abelian categories enriched over the monoidal abelian category $\category{V}$. Assume that for each of the three categories the associated derived category is locally small. If $F\colon \derived{\category{A}} \to \derived{\category{B}}$ if left adjoint to $G\colon \derived{\category{B}}\to\derived{\category{A}}$ then in $\derived{\category{V}}$ there are isomorphisms
		$$\eRHom_\category{B}(F(X^\bullet), Y^\bullet) \cong \eRHom_\category{A}(X^\bullet, G(Y^\bullet))$$
	natural in the complexes $X^\bullet \in \derived{\category{A}}$ and $Y^\bullet \in \derived{\category{B}}$ providing a natural enrichment of the adjunction.
\end{assumption}

\newpage
\section{(Commuting with) Homotopy (Co-)Limits}
 
It is often very convenient to break down an argument about an unbounded complex to a left or right bounded one. Recall the several notions of truncation of complexes from \cite[\href{https://stacks.math.columbia.edu/tag/0118}{Tag 0118}]{stacks-project}.
\begin{lemma}[A complex as a (co)limit of its truncations]
	\label{lemma: a complex as a (co)limit of its truncations}
	
	Let $\category{A}$ be an abelian category and $C^\bullet \in \complexes{\category{A}}$ any complex.
	\begin{itemize}
		\item
		The stupid truncations $\sTruncation{\le n}C^\bullet$ with the natural projections $\sTruncation{\le m} \to \sTruncation{\le n}$ for $n\le m$ form a sequential inverse system. If $\category{A}$ has sequential limits then $C^\bullet \to \limit_{n\in \N} \sTruncation{\le n} C^\bullet$ is an isomorphism.
		
		\item
		The stupid truncations $\sTruncation{\le n}C^\bullet$ with the natural inclusions $\sTruncation{\le n} \to \sTruncation{\le m}$ for $n\le m$ form a sequential direct system. If $\category{A}$ has sequential colimits then $\colimit_{n\in \N} \sTruncation{\le n}C^\bullet \to C^\bullet$ is an isomorphism.
		
		\item
		The canonical truncations $\cTruncation{\le n}C^\bullet$ with the natural inclusions $\cTruncation{\le n} \to \cTruncation{\le m}$ for $n\le m$ form a sequential direct system. If $\category{A}$ has sequential colimits then $\colimit_{n\in \N} \cTruncation{\le n}C^\bullet \to C^\bullet$ is an isomorphism.
		
		\item
		The stupid truncations $\sTruncation{\ge -n}C^\bullet$ with the natural inclusions $\sTruncation{\ge -n}C^\bullet \to \sTruncation{\ge -m}C^\bullet$ for $n\le m$ for a sequential direct system. If $\category{A}$ has sequential colimits them $\colimit_{n\in \N} \sTruncation{\ge -n} C^\bullet \to C^\bullet$ is an isomorphism.
	\end{itemize}
\end{lemma}
\begin{shortproof}
	Clear since in each fixed degree $n\in \Z$ each diagram is eventually constant.
\end{shortproof}

Recall the definition of (sequential) homotopy limits from \cite[\href{https://stacks.math.columbia.edu/tag/08TC}{Tag 08TC}]{stacks-project} and of (sequential) homotopy colimits from \cite[\href{https://stacks.math.columbia.edu/tag/090Z}{Tag 090Z}]{stacks-project}.
\begin{lemma}[Ordinary sequential colimits are homotopy colimits, {\cite[\href{https://stacks.math.columbia.edu/tag/0949}{Tag 0949}]{stacks-project}}]
	\label{lemma: oridnary sequential colimits are homotopy limits}
	
	Let $\category{A}$ be an abelian category such that sequential colimits exist and are exact (for example if $\category{A}$ is  \grothendieckaxiom{$5$}). For any sequential diagram $n\mapsto X_n^\bullet$ of complexes in $\category{A}$ the colimit of complexes $\colimit_{n\in \N} X_n^\bullet$ is a homotopy colimit of $n\mapsto X_n^\bullet$ in $\derived{\category{A}}$.
\end{lemma}
This implies that any of the colimits in lemma \ref{lemma: a complex as a (co)limit of its truncations} is actually a homotopy colimit. For non-sequential diagrams the definition of homotopy (co)limits is more inexplicit.
\begin{definition}[Homotopy limits and colimits]
	Let $\category{A}$ be an abelian category and $\category{I}$ be a small category. Consider the abelian category $\functorCategory{\category{I}}{\category{A}}$ of $\category{I}$-shaped diagrams in $\category{A}$. Observe that we have additive functors
		$$\limit\colon \functorCategory{\category{I}}{\category{A}} \to \category{A}\text{ and }\colimit\colon \functorCategory{\category{I}}{\category{A}}\to\category{A}$$
	if $\category{A}$ is complete respectively cocomplete. Let $D\colon \category{I} \to \complexes{\category{A}}$ be a diagram. In case that the right derived functor $\RightD\limit$ is defined at $D$ we call its image
		$$\homotopyLimit D \ldef (\RightD\limit)(D) \in \derived{\category{A}}$$
	the \emph{homotopy limit of $D$}. Similarly, if the left derived functor $\LeftD\colimit$ is defined at $D$ we call its image
		$$\homotopyColimit D\ldef (\LeftD \colimit)(D) \in \derived{\category{A}}$$
	the \emph{homotopy colimit of $D$}.
\end{definition}


\begin{remark}
	If the abelian category $\category{A}$ is \grothendieckaxiom{$5$} and $\category{I}$ is filtered then by definition $\colimit\colon \functorCategory{\category{I}}{\category{A}} \to \category{A}$ is exact and hence $\homotopyColimit D \cong \colimit D$ for any diagram $D\colon \category{I}\to\complexes{\category{A}}$.
\end{remark}

We now find that derived $\Hom$s maps (sequential) homotopy colimits to homotopy limits in the first argument, as expected.
\begin{lemma}[Right derived Hom \& homotopy colimits in the first argument]
	\label{lemma: right derived Hom and homotopy colimits in the first argument}
	
	Let $\category{A}$ be an abelian category enriched over the monoidal abelian category $\category{V}$. Assume that both categories are \grothendieckaxiom{$5$} and have enough projectives. Let $Y$ be a complex and assume that $\eRHom_\category{A}(-, Y)$ maps direct sums to direct products. For any sequential diagram $n\mapsto X_n$ of complexes the complex $\eRHom_\category{A}(\colimit_{n\in \N} X_n, Y)$ is a homotopy limit of $n\mapsto \eRHom_\category{A}(X_n, Y)$ in $\derived{\category{V}}$. We thus may write
		$$\eRHom_\category{A}(\colimit_{n\in \N} X_n, Y) \cong \homotopyLimit_{n\in \N} \eRHom_\category{A}(X_n ,Y).$$
\end{lemma}
\begin{proof}
	Since $\colimit_{n\in \N} X_n$ is a homotopy colimit of $n\mapsto X_n$ there is a distinguished triangle
		$$\bigoplus_{n\in \N} X_n\to \bigoplus_{n\in \N} X_n\to \colimit_{n\in \N} X_n \to \bigoplus_{n\in \N} X_n[1].$$
	Application of $\eRHom_\category{A}(-, Y)$ yields that
		$$\eRHom_\category{A}(\colimit_{n\in \N} X_n, Y) \to \eRHom_\category{A}(\bigoplus_{n\in \N} X_n, Y)\to \eRHom_\category{A}(\bigoplus_{n\in \N} X_n, Y)\to (...)[1]$$
	is a distinguished triangle as well. As $\eRHom_\category{A}(-, Y)$ maps direct sums to products we obtain a morphism of distinguished triangles
	$$\begin{tikzcd}[column sep=1em, scale cd=0.9]
		\eRHom_\category{A}(\colimit_{n\in \N} X_n, Y) \ar[r]\ar[d, equal] &\eRHom_\category{A}(\bigoplus_{n\in \N} X_n, Y) \ar[r]\ar[d] &\eRHom_\category{A}(\bigoplus_{n\in \N} X_n, Y) \ar[r]\ar[d] &(...)[1]\\
		\eRHom_\category{A}(\colimit_{n\in \N} X_n, Y) \ar[r] & \prod_{n\in \N} \eRHom_\category{A}(X_n, Y) \ar[r] &\prod_{n\in \N} \eRHom_\category{A}(X_n, Y) \ar[r] &(...)[1].
	\end{tikzcd}$$
	Hence $\eRHom_\category{A}(\colimit_{n\in \N} X_n ,Y)$ is a homotopy limit of $n\mapsto \eRHom_\category{A}(X_n, Y)$.
\end{proof}

Completely analogous is of course the commutation with homotopy limits in the second argument.
\begin{lemma}[Right derived Hom \& homotopy limits in the second argument]
	\label{lemma: right derived Hom and homotopy limits in the second argument}
	
	Let $\category{A}$ be an abelian category enriched over the monoidal abelian category $\category{V}$. Assume that both categories are \grothendieckaxiom{$5$} and have enough projectives. Let $X$ be a complex and assume that $\eRHom_\category{A}(X, -)$ commutes with direct products. For any sequential diagram $n\mapsto Y_n$ of complexes we find that
		$$\eRHom_\category{A}(X, \homotopyLimit_{n\in \N} Y_n) \cong \homotopyLimit_{n\in \N} \eRHom_\category{A}(X, Y_n)$$
	in $\derived{\category{V}}$.
\end{lemma}
\begin{proof}
	Since $\homotopyLimit_{n\in \N} Y_n$ is a homotopy limit of $n\mapsto Y_n$ there is a distinguished triangle
		$$\homotopyLimit_{n\in \N} Y_n\to \prod_{n\in \N} Y_n\to  \prod_{n\in \N} Y_n \to (\cdots)[1].$$
	Application of $\eRHom_\category{A}(X, -)$ yields that
		$$\eRHom_\category{A}(X, \homotopyLimit_{n\in \N} Y_n) \to \eRHom_\category{A}(X, \prod_{n\in\N} Y_n)\to \eRHom_\category{A}(X, \prod_{n\in \N} Y_n)\to (...)[1]$$
	is a distinguished triangle as well. As $\eRHom_\category{A}(X, -)$ maps products to products we obtain a morphism of distinguished triangles
	$$\begin{tikzcd}[column sep=1em, scale cd=0.9]
		\eRHom_\category{A}(X, \homotopyLimit_{n\in \N} Y_n) \ar[r]\ar[d, equal] &\eRHom_\category{A}(X, \prod_{n\in \N} Y_n) \ar[r]\ar[d] &\eRHom_\category{A}(X, \prod_{n\in \N} Y_n) \ar[r]\ar[d] &(...)[1]\\
		\eRHom_\category{A}(X, \homotopyLimit_{n\in \N} Y_n) \ar[r] & \prod_{n\in \N} \eRHom_\category{A}(X, Y_n) \ar[r] &\prod_{n\in \N} \eRHom_\category{A}(X, Y_n) \ar[r] &(...)[1]
	\end{tikzcd}$$
	implying that $\eRHom_\category{A}(X, \homotopyLimit_{n\in \N} Y_n)$ is a homotopy limit of $n\mapsto \eRHom_\category{A}(X, Y_n)$.
\end{proof}

\begin{remark}[Compatibility of enriched $\Hom$s with direct sums and products]
	\label{remark: compatibility of enriched Homs with direct sums and products}
	
	That $\eRHom_\category{A}$ maps direct sums to direct products is often true. Assume for example that the category $\category{D}$ is closed symmetric monoidal (as is the case for a closed symmetric monoidal abelian category $\category{A}$ where $\Dotimes_\category{A}$ and $\intRHom_\category{A}$ exist and are (partially) adjoint, see \ref{lemma: deriving a closed monoidal structure}). Then for $X_i$ and $Y$ in $\category{D}$ we obtain that
	\begin{align*}
		\Hom(-,\intHom(\bigoplus_i X_i, Y)) &\cong \Hom(-\otimes \bigoplus_i X_i, Y)\\
		&\cong \Hom(\bigoplus_i X_i, \intHom(-, Y))\\
		&\cong \prod_i \Hom(X_i, \intHom(-, Y))\\
		&\cong \prod_i \Hom(-, \intHom(X_i, Y)) \cong \Hom(-, \prod_i \intHom(X_i, Y))
	\end{align*}
	as $\Hom$ certainly maps colimits to limits in the first argument. By Yoneda we thus find that $\intHom(\bigoplus_i X_i, Y) \cong \prod_i \intHom(X_i, Y)$ as desired -- although technically one should check that this is the correct comparison map. Similarly, one shows that $\intHom(X, \prod_i Y_i) \cong \prod_i \intHom(X, Y_i)$ as well.
\end{remark}

The commutation with general homotopy colimits is not so clear and for now we assume it given.
\begin{assumption}[Right derived $\Hom$ maps colimits to limits]
	\label{assumption: right derived Hom maps colimits to limits}
	
	Let $\category{A}$ be an abelian category enriched over the monoidal abelian category $\category{V}$ such that both categories are \grothendieckaxiom{$5$} with enough projectives. Suppose that $i\mapsto X_i^\bullet$ is a filtered system of complexes, then the natural morphism
	$$\eRHom_\category{A}(\colimit_i X_i^\bullet, Y^\bullet) \to \homotopyLimit_i \eRHom_\category{A}(X_i^\bullet, Y^\bullet)$$
	in $\derived{\category{V}}$ is an isomorphism.
\end{assumption}

\section{Derived Compactness}

We will need slight reformulations of compactness and generation when dealing with triangulated categories. Recall that a triangulated category that has all direct sums already has all (homotopy) colimits, this can be found in \cite[\href{https://stacks.math.columbia.edu/tag/0A5K}{Tag 0A5K}]{stacks-project}. The reformulation of compactness for triangulated categories will only talk about direct sums. 
\begin{definition}[(Derived) compactness, {\cite[\href{https://stacks.math.columbia.edu/tag/07LS}{Tag 07LS}]{stacks-project}}]
	
	Let $\category{D}$ be a triangulated category with arbitrary direct sums. An object $K\in \category{D}$ is called \emph{(derived) compact} if for all direct sums $\bigoplus_{i \in I} X_i$ in $\category{D}$ the natural morphism
	$$\bigoplus\nolimits_{i\in I} \Hom_\category{D}(K, X_i) \to \Hom_\category{D}(K, \bigoplus\nolimits_{i\in I} X_i)$$
	is an isomorphism.
\end{definition}

\begin{example}[Bounded complexes of compact projectives are derived compact]
	\label{example: bounded complexes of compact projectives are derived compact}
	
	Let $\category{A}$ be an \grothendieckaxiom{$5$} abelian category generated by compact projectives. Observe that for any compact projective $P$ of $\category{A}$ the complex $P[0]$ is $K$-projective so that for every complex $X$ we find by lemma \ref{lemma: morphism out of a K-projective complex} that
		$$\Hom_\homotopyCategoryOfComplexes{\category{A}}(P[0], X)\cong \Hom_\derived{\category{A}}(P, X).$$
	But as $P$ is compact this directly implies that $\Hom_\derived{\category{A}}(P[0], -)$ commutes with arbitrary direct sums as we can verify this in the homotopy category of complexes. Hence, for each compact projective $P$ the complex $P[0]$ is derived compact. Furthermore, using an induction on the length (the difference between the largest degree after which there is no entry and the smallest degree before which there is no entry), one easily checks that this implies that any \emph{bounded} complex $P^\bullet$ of which each term is compact projective is derived projective as well.
\end{example}

While in general, contrary to notion of compactness introduced in definition \ref{definition: regular projectivity and compactness}, we might not commute with all filtered (homotopy) colimits (at least it is not known to the author if it holds true), we still can be sure that we commute with sequential colimits.
\begin{lemma}[Commutation with sequential colimits]
	\label{lemma: commutation with sequential colimits}
	
	Let $\category{A}$ be an abelian category in which sequential colimits exist and are exact (for example if $\category{A}$ is \grothendieckaxiom{5}). Suppose that $\N \to \complexes{\category{A}}, n\mapsto X_n^\bullet$ is a sequential diagram in $\complexes{\category{A}}$ and that $K^\bullet \in \derived{\category{A}}$ is derived compact. Then the natural morphism
	$$\colimit_{n\in \N} \Hom_\derived{A}(K^\bullet, X_n^\bullet) \to \Hom_\derived{A}(K^\bullet, \colimit_{n\in \N} X_n^\bullet)$$
	is an isomorphism. Even more, the natural morphism
	$$\colimit_{n\in \N} \eRHom_\category{A}(K^\bullet, X_n^\bullet) \to \eRHom_\category{A}(K^\bullet, \colimit_{n\in \N} X_n^\bullet)$$
	in $\derived{\Ab}$ is an isomorphism, providing a natural enrichment.
\end{lemma}
\begin{proof}
	The first isomorphism is precisely \cite[\href{https://stacks.math.columbia.edu/tag/094A}{Tag 094A}]{stacks-project} after noticing that (sequential) homotopy colimits and ordinary colimits agree by lemma \ref{lemma: oridnary sequential colimits are homotopy limits}. For the second isomorphism, we can check this in each degree $m\in \Z$. Indeed, as 
	$$H^m \eRHom_\category{A}(-, -)\cong H^0\eRHom_\category{A}(-, -[m]) \cong \eHom_\derived{\category{A}}(-, -[m])$$
	this follows from the first case, as then in degree $m$ the induced morphism in cohomology
	$$\colimit_{n\in \N} \eHom_\derived{A}(K^\bullet, X_n^\bullet[m])\to \eHom_\derived{A}(K^\bullet, \colimit_{n\in \N} X_n^\bullet[m])$$
	is an isomorphism since $\colimit_{n\in \N}$ is exact and hence commutes with $H^n$.
\end{proof}

%

In case that the triangulated category in question is the derived category of a monoidally enriched abelian category we can generalize the notion of derived compactness.
\begin{definition}[Enriched compactness]
	\label{definition: enriched compactness}
	
	Let $\category{A}$ be an abelian category enriched over the monoidal abelian category $\category{V}$. Assume that both categories are \grothendieckaxiom{$5$} and have enough projectives. An object $K \in \derived{\category{A}}$ is called \emph{($\category{V}$-)enriched (derived) compact} if for all direct sums $\bigoplus_{i\in I} X_i$ in $\derived{\category{A}}$ the natural morphism
		$$\bigoplus_{i\in I} \eRHom_\category{A}(K, X_i) \to \eRHom_\category{A}(K, \bigoplus_{i\in I} X_i)$$
	in $\derived{\category{V}}$ is an isomorphism.
\end{definition}

For enriched compacts we can now show commutation with sequential homotopy colimits.
\begin{lemma}[Enriched compact objects preserve sequential homotopy colimits]
	\label{lemma: enriched compact objects preserve sequential homotopy colimits}
	
	Let $\category{A}$ be an abelian category enriched over the monoidal abelian category $\category{V}$. Assume that both categories are \grothendieckaxiom{$5$} and have enough projectives. If $K$ is enriched compact and $n\mapsto X_n^\bullet$ is a sequential diagram of complexes then the natural morphism
		$$\colimit_{n\in \N} \eRHom_\category{A}(K^\bullet, X_n^\bullet) \to \eRHom_\category{A}(K^\bullet, \colimit_{n\in \N} X_n^\bullet)$$
	in $\derived{\category{V}}$ is an isomorphism.
\end{lemma}
\begin{proof}
	Consider the defining distinguished triangle
		$$\bigoplus_{n\in \N} X_n^\bullet \to \bigoplus_{n\in \N} X_n^\bullet \to \colimit_{n\in \N} X_n^\bullet \to \bigoplus_{n\in \N} X_n^\bullet[1]$$
	of the homotopy colimit $\colimit_{n\in \N} X_n^\bullet$. Applying the triangulated functor $\eRHom_A(K^\bullet, -)$ we obtain that
		$$\eRHom_A(K^\bullet, \bigoplus_{n\in \N} X_n^\bullet) \to \eRHom_A(K^\bullet, \bigoplus_{n\in \N} X_n^\bullet) \to \eRHom_A(K^\bullet,\colimit_{n\in \N} X_n^\bullet) \to (\cdots)[1]$$
	is distinguished as well. By enriched compactness of $K$ we have a morphism of distinguished triangles
	$$\begin{tikzcd}[column sep=1em, scale cd=.9]
		\eRHom_A(K^\bullet, \bigoplus_{n\in \N} X_n^\bullet)\ar[r]\ar[d]& \eRHom_A(K^\bullet, \bigoplus_{n\in \N} X_n^\bullet)\ar[r]\ar[d]& \eRHom_A(K^\bullet,\colimit_{n\in \N} X_n^\bullet) \ar[r]\ar[d]& (\cdots)[1]\\
		\bigoplus_{n\in \N} \eRHom_A(K^\bullet, X_n^\bullet)\ar[r]& \bigoplus_{n\in \N} \eRHom_A(K^\bullet,  X_n^\bullet)\ar[r]& \colimit_{n\in \N}\eRHom_A(K^\bullet, X_n^\bullet) \ar[r]& (\cdots)[1]
	\end{tikzcd}$$
	where two out of three vertical morphisms are isomorphism. By \cite[\href{https://stacks.math.columbia.edu/tag/014A}{Tag 014A}]{stacks-project} so is the third and we obtain that
		$$\colimit_{n\in \N}\eRHom_A(K^\bullet, X_n^\bullet) \to \eRHom_A(K^\bullet, \colimit_{n\in \N} X_n^\bullet)$$
	is an isomorphism.
\end{proof}

\section{Compact Generation \& Adjoints}

If a triangulated category is \emph{generated} by derived compact objects then many useful consequences follow. One of these consequences is a very powerful adjoint functor theorem based on Brown representability, originally due to Neeman in \cite{neeman}. Neeman uses this theorem to study Grothendieck duality, essentially providing a proof of the theorem using methods from homotopy theory. We state a slight variation of the adjoint functor theorem in \ref{theorem: an adjoint functor theorem via brown representability} that can be found in \cite[\href{https://stacks.math.columbia.edu/tag/0A8E}{Tag 0A8E}]{stacks-project}. Furthermore, we have to assume that the theorem translates the case where the triangulated category in question is generated by a large(!) class of objects. This is \assumptionPreWord \ref{assumption: another adjoint functor theorem}. 

As in the condensed setting we are often dealing with a large class of compacts we have to be careful in the definition of a compactly generated triangulated category. The standard definition assumes a small collection of generators $(G_i)_{i\in I}$ and requires that the single object $\bigoplus_{i\in I} G_i$ generates the category. The right generalization to a proper class of generators seems to be the analogue of definition \ref{definition: generation}. In the language of \cite[\href{https://stacks.math.columbia.edu/tag/09SI}{Tag 09SI}]{stacks-project} we make the following definition.
\begin{definition}[Compactly generated triangulated categories]
	\label{definition: compactly generated triangulated categories}
	
	Let $\category{D}$ be a triangulated category. A class of $\mathcal G$ of objects of $\category{D}$ is called \emph{generating} if for every object $X\in\category{D}$ there is a small $\mathcal{G'}\subseteq \mathcal G$ such that $X$ is contained in $\langle \mathcal{G'}\rangle$ -- the \emph{smallest strictly full, saturated, triangulated subcategory} of $\category{D}$ containing $\mathcal{G'}$. The triangulated category $\category{D}$ is called \emph{compactly generated} if it admits a class of derived compact generators. We call $\category{D}$ \emph{classically compactly generated} if it can be generated by a small class of derived compact objects.
\end{definition}

Now having defined compactly generated triangulated categories, we can state Brown representability. As already hinted at by the name, this result ensures that certain functors out of the triangulated category in question are representable. This is used in the proof of the adjoint functor theorem to ensure that for each given object one can construct the value of the adjoint functor at that object -- the object being provided by Brown representability. 
\begin{proposition}[Brown representability {\cite[\href{https://stacks.math.columbia.edu/tag/0A8F}{Tag 0A8F}]{stacks-project}}]
	\label{proposition: brown representability}
	
	Let $\category{D}$ be classically compactly generated triangulated category admitting arbitrary direct sums. Let $H\colon \category{D}\to \Ab$ be a contravariant functor which transforms direct sums into products. Then $H$ is representable.
\end{proposition}

This now gives the desired adjoint functor theorem.
\begin{theorem}[An Adjoint Functor Theorem via Brown representability]
	\label{theorem: an adjoint functor theorem via brown representability}
	
	Let $\category{D}$ be a classically compactly generated triangulated category admitting arbitrary direct sums. If $F\colon \category D \to \category D'$ is an exact functor of triangulated categories which preserves direct sums, then $F$ has a right adjoint.
\end{theorem}
\begin{proofsketch}[This is {\cite[\href{https://stacks.math.columbia.edu/tag/0A8G}{Tag 0A8G}]{stacks-project}}]
	
	Let $Y$ be any object of $\category{D}'$ and consider the contravariant functor $\category{D}\to \Ab, X\mapsto \Hom_{\category{D'}}(F(X), Y)$. This functor is cohomological as $F$ is exact and transforms directs sums into products (as $F$ preserves direct sums). Hence, Brown representability \ref{proposition: brown representability} guarantees that this functor is representable, \ie that there exists some object $R(Y)\in\category{D}$ such that $\Hom_\category{D}(-, R(Y)) \cong \Hom_{\category{D'}}(F(-), Y)$. These objects $R(Y)$ then assemble into a functor $R\colon \category{D'}\to \category{D}$.
\end{proofsketch}

The question now is: When does this theorem extend to non-classically compactly generated triangulated categories. Neither is the author aware of such a generalization in the literature nor was he able to prove one. In their notes, Clausen and Scholze always use an(!) adjoint functor theorem \quote{implicitly} when they for example state that \quote{[b]y definition $f_!$ commutes with all direct sums (...) [t]his implies formally that $f_!$ admits a right adjoint $f^!$} at \cite[p. 57]{condenseddotpdf} -- this certainly implies that they have this (or a very similar) adjoint functor theorem in mind.
\begin{assumption}[Another adjoint functor theorem]
	\label{assumption: another adjoint functor theorem}

	Let $\category{D}$ be a compactly generated triangulated category admitting arbitrary direct sums. If $F\colon \category D \to \category D'$ is an exact functor of triangulated categories which preserves direct sums, then $F$ has a right adjoint functor.
\end{assumption}

If one has a pair of functors between two compactly generated triangulated categories, one can say even more. 
\begin{proposition}[{\cite[Thm. 5.1]{neeman}}]
	\label{proposition: preserves compact iff right adjoint preserves coproducts}
	
	Let $F\colon \category{D} \to \category{D'}$ be an exact functor of compactly generated triangulated categories with right adjoint $G$. Then $F$ preserves compact objects if and only if $G$ commutes with direct sums.
\end{proposition}

\begin{remark}
	It seems to the author that the proof of proposition \ref{proposition: preserves compact iff right adjoint preserves coproducts} given by Neeman applies verbatim to the case where $\category{D}$ and $\category{D'}$ are generated by a large class of compacts.
\end{remark}

This implies in particular the following criterion for enriched compactness.
\begin{corollary}[A sufficient condition for enriched compactness]
	\label{corollary: a sufficient condition for enriched compactness}
	
	Let $\category{A}$ be a closed monoidal abelian category that is \grothendieckaxiom{$5$} and is generated by compact projectives. Let $K \in \derived{\category{A}}$ be derived compact. If $-\Dotimes_\category{A} K$ preserves derived compactness then $K$ is enriched compact. This is in particular the case if inside $\category{A}$ the tensor product of any two compact objects is again compact.
\end{corollary}
\begin{proof}
	Observe that $\derived{\category{A}}$ is compactly generated. Since the right adjoint of $-\Dotimes_\category{A} K$ is $\intRHom_\category{A}(K, -)$ by lemma \ref{lemma: deriving a closed monoidal structure} it follows that the right adjoint preserves direct sums (so $K$ is enriched compact) if and only if the left adjoint preserves derived compactness. For the final claim observe that as $\derived{\category{A}}$ is generated by the compact projectives of $\category{A}$ it is enough to check the preservation of derived compactness by $-\Dotimes_\category{A} K$ on these generators since the tensor product commutes with arbitrary direct sums (and hence homotopy colimits). But we can replace $K$ by a complex of such compact projectives and hence a tensor product of the replaced $K$ with a compact generators will be derived compact too.
\end{proof}

\section{A Detour on Cohomology of Topological Spaces}

A sample consequence of the power of derived condensed algebra is that one can recover the sheaf cohomology of any compact Hausdorff space (and hence by extension of any space that is reasonably glued from such compacta). First we need a general definition of cohomology of an object in a (pre)topos.
\begin{definition}[Cohomology internal to a pretopos \& condensed cohomology]
	Let $\category T$ be a pretopos, such that $\Ab(\category T)$ is an abelian category where $\Ext$s exist and such that $\Ab(\category T)\to \category T$ admits a left-adjoint $\Z[-]$. For an object $X\in\category T$ and an abelian group objects $M\in\Ab(\category T)$, define
	$$H_\category T^\bullet(X, M) \ldef \Ext^\bullet_{\Ab(\category T)}(\Z[X], M) \in \derived{\Ab},$$
	the \emph{cohomology of $X$ with coefficients in $M$}. This applies verbatim to condensed sets. Thus, for a condensed set $X$ and a condensed abelian group $M$, define
	$$H_\text{cond}^\bullet(X, M)\ldef H_\condensedSets^\bullet(X, M),$$
	the \emph{condensed cohomology of $X$ with coefficients in $M$}.
\end{definition}

We then have the claimed isomorphism of cohomological theories.
\begin{theorem}[Comparing condensed and sheaf cohomology, {\cite[Thm. 3.2]{condenseddotpdf}}]
	Let $M$ be a discrete topological abelian group. There is an isomorphism
	$$H_\text{sheaf}^\bullet(S, \associated{M}) \ldef \RightD\Gamma(S, \associated{M})\cong \eRHom_{\associated{\Z}}(\freeAbelian{\associated{S}}, \associated{M}) \rdef H_\text{cond}^\bullet(\associated S, \associated M)$$
	natural in $S\in\CHaus$. 
\end{theorem}

	\chapter{Scheme Theory \& \texorpdfstring{$6$}{6}-Functor Formalisms}
	\label{chapter: scheme theory and $6$-functor formalisms}
	
\section*{Introduction}
Before delving into the theory of analytic rings and their associated notion of \emph{completeness} in chapter \ref{chapter: analytic rings and completeness}, we will in this chapter first recall the needed geometric foundations. In section \ref{section: schemes} we will recall the definition of spectra and schemes as ringed spaces that are locally such spectra. In section \ref{section: (quasi-)coherentness} we recall the notions of quasi-coherentness and coherentness of sheaves of modules. Section \ref{section: properness} recalls the notion of proper morphism between schemes as these are fundamental to duality. Then \ref{section: yoga} is concerned with introducing the basics of so-called \emph{$6$-functor formalisms} and their implications for duality in section \ref{section: coherent duality}.

\section{Schemes}
\label{section: schemes}

Recall the definition of affine spectra -- the basic building blocks of schemes.
\begin{definition}[Affine spectra, {
		\cite[\href{https://stacks.math.columbia.edu/tag/00DY}{Tag 00DY}]{stacks-project}
		\&
		\cite[\href{https://stacks.math.columbia.edu/tag/01HU}{Tag 01HU}]{stacks-project}}]
	\label{definition: affine spectra}
	
	Let $R$ be a commutative unitary ring. The \emph{(affine) spectrum of $R$} is the (locally) ringed space with underlying set $\Spec(R)$, the set of prime ideals of $R$. We endow this set with the Zariski topology, \ie closed sets are precisely those of the form
	$$V(I) \ldef \{P \in \Spec R \setseparator \forall f\in I\colon f\in P\}$$
	for some $I\subseteq R$. For each $f\in R$ call $D(f)\ldef \Spec R\setminus V(f)$ \emph{distinguished open}. These sets $D(f)$ form a basis of the Zariski topology, and the assignment $D(f)\mapsto R_f$ where $R_f$ is the localization of $R$ at $f$, forms a sheaf of rings on this basis. Hence, there is a sheaf of rings on $\Spec R$ denoted $\struct {\Spec R}$ extending this sheaf on the basis, \ie such that there is a natural isomorphism
	$$R_f \to \struct {\Spec R}(D(f))$$
	for any $f\in R$.
\end{definition}

For convenience, we define an affine scheme to be any locally ringed space that is isomorphic (as a locally ringed space) to an affine spectrum.
\begin{definition}[Affine schemes, {
		\cite[\href{https://stacks.math.columbia.edu/tag/01HW}{Tag 01HW}]{stacks-project}}]
	
	An \emph{affine scheme} is a locally ringed space $(X, \struct X)$ isomorphic to an affine spectrum as a locally ringed space.
\end{definition}

\begin{definition}[Schemes and their morphisms, {
		\cite[\href{https://stacks.math.columbia.edu/tag/01IJ}{Tag 01IJ}]{stacks-project}}]
	
	A \emph{scheme} is a locally ringed space $(X, \struct X)$ in which every point admits an open neighborhood $U$ such that $(U, \struct U)$ is an affine scheme. A morphism of schemes is a morphism of locally ringed spaces.
\end{definition}

One is usually only interested in morphisms of schemes whose fibers are not infinite-dimensional.
\begin{definition}[Being locally of finite type]
	\label{definition: being locally of finite type}
	
	A morphism $f\colon X\to Y$ of schemes is said to be \emph{locally of finite type} if for any point $x\in X$ there is an affine open neighborhood $\Spec A\cong U\subseteq Y$ of $x$ and an affine open neighborhood $\Spec R\cong V\subseteq Y$ of $f(x)$ such that $f(U)\subseteq V$ and such that the associated morphism of rings $R\to A$ induced by $f$ is \emph{of finite type}, \ie that $A$ is isomorphic as an $R$-algebra to a quotient of some polynomial algebra $R[X_1,\dots, X_n]$ over $R$.
\end{definition}

\begin{remark}[Interpreting being locally of finite type]
	In the setting of definition \ref{definition: being locally of finite type} if we write $A \cong R[X_1,\dots, X_n]/I$ for some ideal $I$ then we obtain that $U\cong \Spec A \cong \Spec R[X_1,\dots, X_n]/I \rdef V(I)$ is a closed subscheme of the affine space $\A_R^n$ over the space $\Spec R$. \Ie if $f\colon X\to Y$ is locally of finite type then it locally can be realized as a finite dimensional subscheme of affine space over (a part of) the base.
\end{remark}

There is an even stronger condition that essentially requires the fibers to be finite.
\begin{definition}[Being locally finite]
	\label{definition: being locally finite}
	
	A morphism $f\colon X\to Y$ of schemes is said to be \emph{locally finite} if for every point $x\in X$ there is an affine open neighborhood $\Spec A\cong U\subseteq Y$ of $x$ and an affine open neighborhood $\Spec R\cong V\subseteq Y$ of $f(x)$ such that $f(U)\subseteq V$ and such that the associated morphism of rings $R\to A$ induced by $f$ is \emph{finite}, \ie that $A$ is a finitely generated $R$-module.
\end{definition}

\begin{remark}[Interpreting finiteness]
	In the setting of definition \ref{definition: being locally finite} we find that since $R\to A$ is finite that by \cite[\href{https://stacks.math.columbia.edu/tag/05DR}{Tag 05DR}]{stacks-project} the fibers of $f\colon U\cong \Spec A\to \Spec R=V$ are finite, as desired. Thus, $f\colon X\to Y$ is best imagined as a branched cover of $Y$ (or rather of its image inside $Y$ -- the inclusion of a closed subset is for example finite as well, but many fibers are empty).
\end{remark}

A very well-behaved (for these we will obtain boundedness of the dualizing complex, see theorem \ref{theorem: exceptional direct and inverse image functors, the affine case}) class of (affine) schemes is given by \emph{complete intersections}.
\begin{definition}[Complete intersections]
	A morphism $f\colon \Spec A\to \Spec R$ of affine schemes is said the be a \emph{complete intersection} if there is a regular sequence $(f_1,\dots,f_r)$ in $R$ such that the natural morphism $R/\ideal{f_1,\dots, f_r}\to A$ of $R$-algebras induced by $R\to A$ is an isomorphism. Indeed, then $f$ is isomorphic (as schemes over $\Spec R$) to $V(f_1,\dots, f_r) = \Spec R/\ideal{f_1,\dots, f_r}$ into $\Spec R$.
\end{definition}

\section{(Quasi-)coherentness}
\label{section: (quasi-)coherentness}

Next we will introduce the notions of quasi-coherent and coherent sheaves on some scheme $X$. These essentially should correspond to sheaves of modules which locally look like modules respectively finitely generated modules. While this is not true in general (at least in the case of coherence) it is still a good way to think about them.

\begin{definition}[Quasi-coherent sheaves, {
		\cite[\href{https://stacks.math.columbia.edu/tag/01BE}{Tag 01BE}]{stacks-project}}]
	
	Let $(X, \struct X)$ be a ringed space. A sheaf $\sheaf F\in\structMod X$ is called \emph{quasi-coherent} if for every $x\in X$ there exists an open neighborhood $x\in U\subseteq X$ for which there is a presentation
	$$\struct U^{\oplus \vert J\vert} \to \struct U^{\oplus \vert I\vert} \to \sheaf F\restriction_U \to 0$$
	of $\sheaf F\restriction_U$, where $I$ and $J$ are possibly infinite sets, \ie if \emph{$\sheaf F$ is locally presentable}. We denote by $\qCoh X$ the full subcategory of quasi-coherent sheaves of $\structMod X$. If additionally for any $x$ the set $I$ and $J$ can be chosen finite, then $\sheaf F$ is called \emph{locally of finite presentation}.
\end{definition}

\begin{remark}[Terminology]
	Some authors (including the authors of The Stacks Project, \confer \cite[\href{https://stacks.math.columbia.edu/tag/01BM}{Tag 01BM}]{stacks-project}) call \emphquote{locally of finite presentation} just \emphquote{of finite presentation}. To the author this seems like a misnomer.
\end{remark}

\begin{definition}[Associated modules, {
		\cite[\href{https://stacks.math.columbia.edu/tag/01BH}{Tag 01BH}]{stacks-project}}]
	
	Let $(X, \struct X)$ be a ringed space and $\alpha\colon R\to \Gamma(X, \struct X)$ be any ring morphism. There is a functor $\associatedSheaf{(-)} \colon \mod R \to \structMod{X}$ left adjoint to $\Gamma(X, -)\colon \structMod{X} \to \mod R$ where for some sheaf $\sheaf F$ the $R$-module structure on $\Gamma(X,\sheaf F)$ is induced by $\alpha$. For an $R$-module $M$ the sheaf $\associatedSheaf{M}$ is called the \emph{sheaf associated to $M$ and $\alpha$}. If $R=\Gamma(X, \struct X)$ and $\alpha = \id_R$, then $\associatedSheaf{M}$ is simply called the \emph{sheaf associated to $M$}.
\end{definition}

\begin{definition}[Fundamental systems of neighborhoods]
	Let $X$ be a topological space and $x\in X$. A collection $\mathcal N\subseteq \mathcal N_x$ of neighborhoods of $x$ is called a \emph{fundamental system of neighborhoods of $x$} or \emph{neighborhood basis of $x$}, if for every neighborhood $N$ of $x$ there exists some neighborhood $N' \in \mathcal N$ with $N'\subseteq N$. Hence, if we define the partial order $\le$ by $N \le N' \Leftrightarrow N'\subseteq N$, then $\mathcal N$ is a fundamental system of neighborhoods for $x$ if and only if $\mathcal N$ is cofinal in $\mathcal N_x$.
\end{definition}

\begin{lemma}[Quasi-coherent sheaves are locally associated modules, {\cite[\href{https://stacks.math.columbia.edu/tag/01BK}{Tag 01BK}]{stacks-project}}]
	
	Let $(X,\struct X)$ be a ringed space and $x\in X$ a point with a fundamental system of quasi-compact neighborhoods. If $\sheaf F$ is a quasi-coherent $\struct X$-module, then there exists an open neighborhood $U$ of $x$ such that $\sheaf F\restriction_U$ is isomorphic to the sheaf $\associatedSheaf{M}$ associated to some $\Gamma(U, \struct X)$-module $M$.
\end{lemma}

For schemes this is always true.
\begin{corollary}[Quasi-coherent sheaves of schemes]
	Since in a scheme $X$ every point has a fundamental system of quasi-compact neighborhoods, namely the system of affine neighborhoods, we obtain that any quasi-coherent sheaf on $X$ is locally an associated module.
\end{corollary}

If the scheme in question is even affine, we obtain an even stronger statement. Quasi-coherent sheaves are nothing more than modules of the corresponding ring.
\begin{remark}[Quasi-coherent sheaves of affine schemes, {\cite[\href{https://stacks.math.columbia.edu/tag/01IB}{Tag 01IB}]{stacks-project}}]
	\label{remark: quasi-coherent sheaves of affine schemes}
	
	If $X=\Spec R$ is affine, then 
	$$\mod{R} \xleftrightarrows{\Gamma(X, -)}{\associatedSheaf{(-)}} \qCoh{\struct{X}}$$
	is an equivalence of categories.
\end{remark}

Thus, just like schemes are glued from affine schemes, quasi-coherent sheaves (on a scheme) are glued from modules associated to the local patches. A question arises: What about (quasi-coherent) sheaves that a locally \emph{finitely presented} (or equivalently \emph{finitely generated} if the scheme in question is locally Noetherian)? This will lead us to the notion of \emph{coherentness}. First we introduce the analogue of \emph{being finitely generated}.
\begin{definition}[Sheaves of finite type, {
		\cite[\href{https://stacks.math.columbia.edu/tag/01B5}{Tag 01B5}]{stacks-project}}]
	
	Let $(X, \struct X)$ be a ringed space. A sheaf $\sheaf F$ of $\struct X$-modules  is called \emph{of finite type} if for every $x\in X$ there is an open neighborhood $x\in U\subseteq X$ and a surjection
	$$\struct U^{\oplus n} \to \sheaf F\restriction_U$$
	for some $n\in \N$, \ie if $\sheaf F$ is locally generated by finitely many sections.
\end{definition}

Being \emph{coherent} is now given by requiring certain local maps to have finite type kernel. This is a priori stronger than requiring the sheaf to be locally of finite presentation as one needs many kernels to be of finite type. We will see in lemma \ref{lemma: coherentness vs being locally of finite presentation} that for a wide class of examples \emph{being coherent} and \emph{being locally of finite presentation} agree.
\begin{definition}[Coherent sheaves of modules, {
		\cite[\href{https://stacks.math.columbia.edu/tag/01BV}{Tag 01BV}]{stacks-project}}]
	
	Let $(X, \struct X)$ be a ringed space. A sheaf of $\struct X$-modules $\sheaf F$ is called \emph{coherent} if it is of finite type and for any open $U\subseteq X$ and any finite collection $s_1,\dots,s_n\in \sheaf F(U)$ the kernel of the corresponding morphism $\struct U^{\oplus n} \to \sheaf F\restriction_U$ is of finite type. The full subcategory of coherent sheaves of $\struct X$-modules is denoted $\Coh X$.
\end{definition}

A direct consequence is that any coherent sheaf is already quasi-coherent.
\begin{remark}[Coherent sheaves are quasi-coherent]
	Any coherent sheaf $\sheaf F$ on a locally ringed space $(X, \struct X)$ is locally of finite presentation and hence quasi-coherent. Indeed, since $\sheaf F$ is of finite type, choose for each $x\in X$ an open neighborhood $U$ and a surjection $\pi\colon\struct U^{\oplus n} \to \sheaf F\restriction_U$. By coherentness its kernel must be of finite type and thus itself admit a surjection $\struct V^{\oplus m} \to \kernel(\pi)\restriction_V$ for some $x\in V\subseteq U$. Hence, we obtain a local presentation
	$$\struct V^{\oplus m} \to \struct V^{\oplus n} \to \sheaf F\restriction_V\to 0$$
	showing that $\sheaf F$ is quasi-coherent.
\end{remark}

As mentioned before, if one imposes some regularity, then coherentness is nothing more than being locally of finite presentation.
\begin{lemma}[Coherentness vs being locally of finite presentation, {\cite[\href{https://stacks.math.columbia.edu/tag/01BZ}{Tag 01BZ}]{stacks-project}}]
	\label{lemma: coherentness vs being locally of finite presentation}
	
	Let $(X, \struct X)$ be a ringed space for which the structure sheaf $\struct X$ is coherent and let $\sheaf F$ be a sheaf of $\struct X$-modules. Then $\sheaf F$ is coherent if and only if it is locally of finite presentation.
\end{lemma}

\begin{remark}[For which schemes $X$ is $\struct X$ coherent?]
	The question now is: To which schemes $X$ does lemma \ref{lemma: coherentness vs being locally of finite presentation} apply? \Ie for which schemes $X$ is the structure sheaf $\struct{X}$ coherent. This is certainly true for any scheme that is locally Noetherian, \confer \cite[\href{https://stacks.math.columbia.edu/tag/01XZ}{Tag 01XZ}]{stacks-project}. This covers many examples of interest and certainly all varieties over fields. In general there is the notion of a coherent ring in the affine case (a ring $R$ whose associated module $\struct{\Spec R}$ is coherent, equivalently that all its finitely generated ideals are finitely presented).
\end{remark}

\section{Properness}
\label{section: properness}

Very important to the theory of duality is the notion of a proper morphism.

\begin{definition}[Closed immersions, {
		\cite[\href{https://stacks.math.columbia.edu/tag/01HK}{Tag 01HK}]{stacks-project}}]
	
	Let $(i, i^\sharp)\colon (Z,\struct Z)\to (X,\struct X)$ be a morphism of locally ringed spaces. We call $(i, i^\sharp)$ a \emph{closed immersion} if
	\begin{itemize}
		\item
		the map $i$ is a homeomorphism onto a closed subset of $X$,
		
		\item
		the morphism $i^\sharp\colon\struct X\to i_\ast \struct Z$ is an epimorphism, and
		
		\item 
		the $\struct X$-module $\sheaf I\ldef \ker i^\sharp$ is locally generated by sections.  
	\end{itemize}
\end{definition}

\begin{example}[Interpreting closed immersions]
	The most basic closed immersion is the case of the morphism of schemes $V(I)\ldef \Spec A/I \to \Spec A$ where $A$ is some ring and $I$ is some ideal of $A$. Every other closed immersion is built from these examples as by \cite[\href{https://stacks.math.columbia.edu/tag/01QO}{Tag 01QO}]{stacks-project} a morphism of schemes $(i, i^\sharp)\colon Z\to X$ is a closed immersion if and only if for each affine open $\Spec A\cong U \subseteq X$ the scheme $i^{-1}(U)$ over $\Spec A$ is isomorphic to $\Spec A/I \to \Spec A$ for some ideal $I$ of $A$, \ie $(i, i^\sharp)\colon Z\to X$ is a closed immersion if $Z$ inside $X$ is \emphquote{locally a vanishing set of some ideal}.
\end{example}

\begin{definition}[Quasi-compactness, {
		\cite[\href{https://stacks.math.columbia.edu/tag/01K3}{Tag 01K3}]{stacks-project} \& \cite[\href{https://stacks.math.columbia.edu/tag/005A}{Tag 005A}]{stacks-project}}]
	
	A morphism of schemes $f\colon X\to Y$ is said to be \emph{quasi-compact}, if the (set-theoretic) preimage $f^{-1}(V)$ of any compact open $V\subseteq Y$ is itself compact.
\end{definition}

\begin{definition}[(Quasi-)separatedness, {
		\cite[\href{https://stacks.math.columbia.edu/tag/01KK}{Tag 01KK}]{stacks-project}}]
	Let $f\colon X\to Y$ be a morphism of schemes and consider the diagonal morphism $\Delta\colon X\to X\times_Y X$. If $\Delta$ is a closed immersion, then $f$ is called \emph{separated}. If $\Delta$ is quasi-compact then $f$ is called \emph{quasi-separated}. By \cite[\href{https://stacks.math.columbia.edu/tag/01K7}{Tag 01K7}]{stacks-project} any separated morphism is also quasi-separated.
\end{definition}

\begin{remark}
	Imposing quasi-separatedness ensures that geometrically pathological cases like the affine line with doubled origin $(\A_K^1\sqcup \A_K^1)/(\A_K^1\setminus\{0\}) \to \Spec K$ is excluded from consideration. Including such examples would be an obstruction to developing a good local theory as in non-(quasi-)separated spaces weird local artifacts occur: The two origins are distinct, yet there are no two \quote{essentially} disjoint closed subsets separating these two points. In an ordinary variety (over an algebraically closed field) one would expect for any two given distinct points to find two (up to a codimension $1$ subvariety) disjoint curves, one passing through the first point one passing through the second one. This is not the case in the affine line with doubled origin. Even worse, there are two ways to complete the inclusion of the smooth locally closed subvariety $\A_K^1\setminus\{0\}$ to a closed smooth subvariety. We essentially loose \quote{uniqueness of limits}. 
\end{remark}


\begin{definition}[Universally closed morphisms]
	Let $f\colon X\to Y$ be a morphism of schemes. Then $f$ is called \emph{universally closed} if the map of underlying topological spaces is closed and any base change $X\times_Y Z\to Z$ of $f$ along a morphism $h\colon Z\to Y$ of schemes is closed as well.
\end{definition}

\begin{definition}[Properness, {
		\cite[\href{https://stacks.math.columbia.edu/tag/01W1}{Tag 01W1}]{stacks-project}}]
	
	A morphism $f\colon X\to Y$ of schemes is called \emph{proper} if it is separated, universally closed and of finite type.
\end{definition}

\begin{example}[Interpretation of properness]
	There is an alternative characterization of properness that is quite geometrically pleasing. Suppose that $f\colon X\to Y$ is a quasi-separated morphism of finite type. Then by \cite[\href{https://stacks.math.columbia.edu/tag/0BX4}{Tag 0BX4}]{stacks-project} the morphism $f$ is proper if and only if for every valuation ring $A$ with field of fractions $K$ fitting into a commutative solid diagram
	$$\begin{tikzcd}
		\Spec K \ar[r]\ar[d] &X\ar[d]\\
		\Spec A\ar[r]\ar[ur, dotted] &Y
	\end{tikzcd}$$
	there is a unique dotted arrow making the diagram commute. Observe that $\Spec K$ is the generic point of $\Spec A$ and that $\Spec A\to Y$ can be imagined as a parameterized (by the one-dimensional space $\Spec A$) smooth curve in $Y$. Each such solid diagram then essentially provides a lift of the generic point (so of \quote{almost all} of the curve) to $X$. The existence and uniqueness of the dotted arrow then assert that the missing points in the lifted curve can be uniquely imputed to obtain a lift of the curve in $Y$ to a smooth curve in $X$. So properness is a kind of \emph{(relative) completeness} for schemes.
	
	This for example is clear in case of a closed immersion $X\hookrightarrow Y$: If \quote{almost all} of a curve $S\subseteq Y$ already lies in $X$ then so does the whole curve: $X$ is closed in $Y$ so already contains all its limiting points. Properness only gets properly interesting in cases where the underlying map is surjective. \Eg the projection $\A_\C^2\to \A_\C^1$ to the first component fails to be proper: Denote by $x$ and $y$ the coordinate functions of the coordinate ring of $\A^2_\C$. The \quote{almost} lift of the curve $\A^1_\C\to \A^1_\C, t\mapsto t$ to the curve $\A^1_\C \to \A^2_\C, t\mapsto (t, 1/t)$ with image in $V(xy - 1)$ already fails to lift. This of course is the starting point for the cherished concept of a projective space.
\end{example}

Proper morphisms of schemes are both quasi-separated (by definition) and quasi-compact. Indeed, a proper morphism is universally closed by definition and hence quasi-compact by \cite[\href{https://stacks.math.columbia.edu/tag/04XU}{Tag 04XU}]{stacks-project}. Along qcqs morphisms it is possible to transport quasi-coherent sheaves.

\begin{lemma}[Direct and inverse images of quasi-coherent sheaves]
	\label{lemma: direct and inverse images of quasi-coherent sheaves}
	
	Let $f\colon X\to Y$ be a morphism of ringed spaces. If $\sheaf G \in \qCoh{Y}$, then also $f^\ast \sheaf G \in \qCoh X$. If $f$ is a qcqs morphism of schemes and $\sheaf F\in\qCoh X$, then also $f_\ast\sheaf F \in \qCoh Y$. Hence, if $f$ is qcqs morphism of schemes, direct and inverse image restrict to an adjoint pair
	$$\qCoh X \xleftrightarrows{f^\ast}{f_\ast} \qCoh Y.$$
\end{lemma}
\begin{proof}
	For inverse images this is \cite[\href{https://stacks.math.columbia.edu/tag/01BG}{Tag 01BG}]{stacks-project}, for direct images this is \cite[\href{https://stacks.math.columbia.edu/tag/01LC}{Tag 01LC}]{stacks-project}.
\end{proof}

\begin{remark}
	As any morphism$ f\colon \Spec B \to \Spec A$ of affine schemes is qcqs by \cite[\href{https://stacks.math.columbia.edu/tag/01S7}{Tag 01S7}]{stacks-project}, the previous lemma \ref{lemma: direct and inverse images of quasi-coherent sheaves} directly applies. This is not very surprising 
	because if $f$ corresponds to $\phi \colon A\to B$ then restricted to quasi-coherent modules, $f_\ast$ and $f^\ast$ correspond to scalar restriction and scalar extension of modules along $\phi$.
\end{remark}

\section{Yoga}
\label{section: yoga}

The natural starting point for coherent duality (and many related duality theories) is the idea of a \emph{$6$-functor formalism}. In general, a $6$-functor formalism ascribes to each object of a category of geometric objects $C$ (be it schemes or some category of locally compact topological spaces) a category (best thought of as a derived category of sheaves) and then relates them by $6$ classes of functors -- two of which (the \emphquote{exceptional} ones) are only defined on a pre-specified subclass $E$ of morphisms of $C$. Manipulating these $6$-functors is also called \emph{the yoga of $6$-functors}.

There is a rich theory surrounding this very general formalism -- this is for example covered in Scholze's notes \cite{sixfunctors} -- but we are only interested in one specific application (namely coherent duality). Thus, let us not introduce the general theory and instead consider by hand what this formalism should be in our case.


For us, the category of geometric objects is of course a subcategory $X$ of $\Schemes$, the category of schemes. Classically one attaches to each scheme $X$ a derived category $\derived{X}$ of sheaves of $\struct X$-modules, \eg coherent or quasi-coherent sheaves. Let us postpone the discussion of the appropriate $E$ for now and let us specify instead the $6$ classes of functors. The first $2$ classes are given by
\begin{enumerate}
	\item
	the (derived) tensor products $-\otimes^\LeftD_{\struct X} -$
	
	\item
	and their (partial) right adjoints $\intRHom_{\struct X}$ the (derived) sheaf of homomorphisms,
\end{enumerate}
making each $\derived{X}$ closed symmetric monoidal. The next $2$ classes of functors associate to any morphism $f\colon X\to Y$ of schemes the
\begin{enumerate}
	\item[3.]
	(derived) direct image $f_\ast\colon \derived{X}\to \derived{Y}$ of $\struct X$-modules and
	
	\item[4.]
	its left adjoint $f^\ast\colon \derived{Y}\to \derived{X}$, the (derived) inverse image of $\struct Y$-modules.
\end{enumerate}
These are the derived functors of the usual direct and inverse image $f_\ast$ and $f^\ast$ and hence should maybe be more appropriately notated by $\RightD f_\ast$ and $\LeftD f^\ast$. It is often however more convenient to conflate the names and the symbols used to describe these functors. In particular, for a general $6$-functor formalism these do not even have to be derived functors nor do the categories $\derived{X}$ need to be derived (in Scholze's \cite{sixfunctors} they are even \quote{generic} $\infty$-categories). The final two classes associate to any morphism $f\colon X\to Y$ of $E$ the
\begin{enumerate}
	\item[5.]
	\emph{exceptional direct image} $f_!\colon \derived{X} \to \derived{Y}$ and
	
	\item[6.]
	its right-adjoint $f^!\colon \derived{Y} \to \derived{X}$, the \emph{exceptional inverse image}.
\end{enumerate}
These $6$ classes of functors now are subjected to several constraints and compatibilities they ought to satisfy. One such constraint is that one requires that $f_!$ agrees with $f_\ast$ in case $f$ is proper. Another example is the \emph{projection formula}: If $f\colon X\to Y$ is proper then for all $\sheaf F\in\derived{X}$ and $\sheaf G\in\derived{Y}$ there is a natural isomorphism
$$f_\ast\sheaf F\otimes_\struct{Y} \sheaf G \to f_\ast(A\otimes_{\struct{X}} f^\ast \sheaf G)$$
in $\derived{Y}$.

Classically there area few choices. If one wants to restrict to (quasi-)coherent sheaves then one needs to be restrictive in the definition of $C$ to ensure that for any morphism $f$ in $C$ the two functors $f_\ast$ and $f^\ast$ are well-defined, see lemma \ref{lemma: direct and inverse images of quasi-coherent sheaves}. Further, to find a large class $E$ usually requires some thought, as it is hard to construct an appropriate functor $f_!$. In case of sheaves on the étale site of a scheme, one obtains such a functor $f_!$ as the \emph{direct image with compact support}, see \cite[\href{https://stacks.math.columbia.edu/tag/0F4W}{Tag 0F4W}]{stacks-project}. In general, it is not so clear how to make a good choice of $f_!$ (and hence of $f^!$).

For us this issue will be resolved using condensed mathematics. Clausen and Scholze introduce in \cite{condenseddotpdf} for every schemes $X$ a category $\derived{\struct{X,\blacksquare}}$ and for a large class of morphisms an exceptional direct image $f_!$. In the affine case of $X=\Spec A$ this category $\derived{\struct{X, \blacksquare}}$ will be the derived category $\derived{A_\blacksquare}$ of so-called \emph{$A_\blacksquare$-complete} condensed $A$-modules (see chapter \ref{chapter: analytic rings and completeness}) while the $f_!$ will be defined for any $f\colon \Spec A\to \Spec R$ where $R$ is Noetherian and $f$ is of finite Tor-dimension (see chapter \ref{chapter: affine duality}). These definitions make critical use of condensed mathematics and provide a new approach to constructing a good derived theory of schemes.

%
%

\section{Coherent Duality}
\label{section: coherent duality}

One of the ways to formulate coherent duality as follows.
\begin{theorem}[Coherent Duality, {\cite[Chapter VII, Corollary 3.4 (a) \& (c)]{residues}}]
	\label{theorem: coherent duality}
	
	Let $f\colon X\to Y$ be a proper morphism of finite type between \emph{Noetherian} (locally Noetherian and quasi-compact) schemes of finite Krull dimension admitting a dualizing complex. Then the right adjoint $f^!\colon \derivedPlus{\Coh{Y}} \to \derivedPlus{\Coh{X}}$ of $\RightD f_\ast\colon  \derivedPlus{\Coh{X}} \to \derivedPlus{\Coh{Y}}$ exists and the natural morphism
	$$\RightD f_\ast \RightD\intHom_{\struct X}^\bullet(F^\bullet, f^! G^\bullet) \to \RightD\intHom_{\struct Y}^\bullet(\RightD f_\ast F^\bullet, G^\bullet)$$
	is a (quasi)isomorphism for all complexes $F^\bullet \in \derivedMinus{\qCoh{X}}$ and $G^\bullet \in \derivedPlus{\Coh{Y}}$.
\end{theorem}

\begin{remark}[Existence of a dualizing complexes, {\cite[V \S10, sufficient conditions]{residues}}]
	When does a dualizing complex exist?
	\begin{itemize}
		\item
		If the scheme $X$ is Gorenstein (see \cite[\href{https://stacks.math.columbia.edu/tag/0AWV}{Tag 0AWV}]{stacks-project}) and of finite Krull dimension then $\struct{X}[0]$ is a dualizing complex for $X$. This covers the case of $X$ being Cohen--Macaulay (see \cite[\href{https://stacks.math.columbia.edu/tag/02IN}{Tag 02IN}]{stacks-project}) and of course the case where $X$ is regular and locally Noetherian (the definition of a Gorenstein scheme and of a Cohen--Macaulay scheme require local Noetherianness).
		
		\item
		If $f\colon X\to Y$ is a morphism of finite type with $Y$ Noetherian then if $Y$ admits a dualizing complex $\canonical_Y^\bullet$ so does $X$. If $f$ is proper with right adjoint $f^!$ then this is $\canonical_X^\bullet \ldef f^!\canonical_Y^\bullet$. This of course covers the case where $Y=\Spec K$ for some field $K$ as $Y$ is regular and (locally) Noetherian.
	\end{itemize}
\end{remark}

In the case where $f\colon X\to \Spec K$ is proper over the field $K$ this can be simplified a bit. Notice that if $f$ is of finite type with $X$ quasi-compact the Noetherianness is automatic.
\begin{corollary}[Coherent duality, a special case]
	\label{corollary: coherent duality, a special case}
	
	Let $f\colon X\to \Spec K$ be a morphism of finite type such that $X$ is quasi-compact and of finite Krull dimension. Denote by $(-)^\vee$ the $K$-dual. Then in $\derived{K}$ there is an isomorphism
	$$\eRHom_\struct{X}(F^\bullet, f^!K)\cong \RightD\Gamma(F^\bullet, X)^\vee$$
	natural in the complexes $F^\bullet \in \derivedMinus{\qCoh{X}}$.
\end{corollary}
\begin{proof}
	First, observe that $f_\ast$ is nothing more than the global sections functor $\Gamma(X, -)\cong \eHom_\struct{X}(\struct{X}, -)$. Furthermore, $\Gamma(X, -) \circ \intHom_\struct{X} = \eHom_\struct{X}$ gives that $\RightD f_\ast\intRHom_\struct{X} \cong \eRHom_\struct{X}$. As $(-)^\vee$ is short for $\eRHom_K(-, K)$ we find that theorem \ref{theorem: coherent duality} implies
	$$\eRHom_\struct{X}(F^\bullet, f^!K)\cong \RightD\Gamma(F^\bullet, X)^\vee$$
	for $G^\bullet = K$.
\end{proof}
This of course includes all quasi-compact \emph{varieties} (reduced, separated schemes of finite type over a(n algebraically closed) field).

\begin{remark}
	\label{remark: upper shriek of structure sheaf is concentrated}
	
	Suppose we are again in the setting of theorem \ref{theorem: coherent duality} and recall the definition of a Cohen--Macaulay morphism of schemes from \cite[\href{https://stacks.math.columbia.edu/tag/045Q}{Tag 045Q}]{stacks-project}.	In case that $f\colon X\to \Spec Y$ is Cohen--Macaulay (this includes the case where $f$ is smooth) such that all fibers have the same dimension $d\in\N$ we find that $f^!\struct{Y}$ is a coherent sheaf concentrated in degree $d$ by \cite[Chapter VII, \S4]{residues}.
\end{remark}

One immediately obtains the following corollary of coherent duality. This is the version of Serre duality presented in \cite[Chap. III, Thm. 7.6]{algebraic-geometry} for projective schemes.
\begin{corollary}[Serre duality for Cohen--Macaulay schemes]
	\label{corollary: serre duality for cohen--macaulay schemes}
	
	If $X\to \Spec K$ is a finite type proper Cohen--Macaulay scheme over some field $K$ of pure dimension $n\in \N$ then there is a coherent sheaf $\canonical_{X/K}$ such that for all $k\in \Z$ there is an isomorphism
		$$\Ext_{\struct X}^k(\sheaf F, \canonical_{X/K}) \iso H^{n-k}(X, \sheaf F)^\dual$$
natural in the coherent sheaf $\sheaf F$.
\end{corollary}
\begin{proof}
Observe that by remark \ref{remark: upper shriek of structure sheaf is concentrated} the complex $f^!K$ is a coherent sheaf concentrated in degree $n$. Denote by $\canonical_{X/K}$ this sheaf such that $f^!K = \canonical_{X/K}[n]$. By corollary \ref{corollary: coherent duality, a special case} we find that
$$\eRHom_\struct{X}(\sheaf F, \canonical_{X/K}[n]) \cong \RightD\Gamma(\sheaf F, X)^\vee$$
and application of $H^{k-n}=H^k\circ (-)[-n]$ yields that
$$H^k\eRHom_\struct{X}(\sheaf{F}, \canonical_{X/K}) \cong (H^{n-k}\RightD\Gamma(\sheaf F, X))^\vee$$
since $\eRHom_\struct{X}$ is compatible with shifts and since $(-)^\vee$ is exact, however reverses the complexes. Thus,
$$\Ext^k_\struct{X}(\sheaf F, \canonical_{X/K})\cong H^{n-k}(X, \sheaf{F})^\vee.$$
\end{proof}

In case that the coherent sheaf $\sheaf F$ in corollary \ref{corollary: serre duality for cohen--macaulay schemes} is locally free, one rephrase the obtained isomorphism in more familiar terms. For this we need the following characterization of internal $\Hom$s out of a locally free sheaf.
\begin{lemma}
\label{lemma: sheaf hom for vector bundles}

Let $X$ be a scheme. If $\sheaf L$ is locally free of finite rank and $\sheaf G$ is any sheaf of $\struct X$-modules, then
$$\intHom_{\struct X}(\sheaf L, \sheaf G) \iso \sheaf L^\dual \structTensor{X} \sheaf G.$$
If $X$ is a scheme over $\Spec K$ for some field $K$, then this isomorphism is $K$-linear.
\end{lemma}
\begin{shortproof}
This can be calculated locally where it is clear due to freeness.
\end{shortproof}

\begin{corollary}
\label{corollary: ext of vector bundle}

In the setting of the previous lemma \ref{lemma: sheaf hom for vector bundles} we obtain that
$$\Ext^\bullet_{\struct X} (\sheaf L, -) \ldef \eStructRHom{X}(\sheaf L, -) \iso \RightD \Gamma(X, \sheaf L^\dual \structTensor{X} -) \rdef H^\bullet(X, \sheaf L^\dual \structTensor{X} -)$$
as functors $\derived{\struct{X}}\to\derived{\Ab}$. If $X$ is a scheme over $\Spec K$ for some field $K$ then this is an isomorphism in $\derived{K}$ as well.
\end{corollary}
\begin{proof}
We have
\begin{align*}
	\eStructHom{X}(\sheaf L, -)&\iso \eStructHom{X}(\sheaf L\structTensor{X} \struct X, -)\\
	&\iso \eStructHom{X}(\struct X, \intHom_{\struct X}(\sheaf L, -))\\
	& \iso \eStructHom{X}(\struct X, \sheaf L^\dual\structTensor{X} -)\\
	&\iso \Gamma(X, \sheaf L^\dual \structTensor{X} -).
\end{align*}
If $X$ is a scheme over $\Spec K$ then $\mod{\struct{X}}$ is $K$-enriched and these isomorphisms are all of $K$-modules. It follows that the derived functors
$$\eStructRHom{X}(\sheaf L, -) \iso \RightD\Gamma(X, \sheaf L^\dual\structTensor{X} -)$$
agree as well as $\sheaf{L}^\vee\otimes_\struct{X}-$ is exact.
\end{proof}

We now can reformulate Serre duality in its classical form.
\begin{corollary}[Serre Duality]
\label{theorem: serre duality}

	If $X\to \Spec K$ is a finite type proper Cohen--Macaulay scheme over some field $K$ of pure dimension $n\in \N$ then there is a coherent sheaf $\canonical_{X/K}$ such that for all $k\in \Z$ there is an isomorphism
	$$H^k(X, \sheaf F) \cong H^{n-k}(X, \sheaf F^\vee \otimes_{\struct{X}} \canonical_{X/K})^\dual$$
	natural in the locally free sheaf $\sheaf F$ of finite rank.
\end{corollary}
\begin{proof}
	This follows by noticing that $\Ext^k_{\struct{X}}(\sheaf F, \canonical_{X/K})\cong H^k(X, \sheaf F^\vee \otimes_{\struct{X}} \canonical_{X/K})$, by applying the (self-inverse) $(-)^\vee$ and replacing $k$ by $n-k$.
\end{proof}

	\chapter{Analytic Rings \& Completeness}
	\label{chapter: analytic rings and completeness}

In the spirit of section \ref{section: yoga} we will work towards constructing a $6$-functor formalism suitable to imply coherent duality. Contrary to the classical approach we will not associate to a scheme $X$ the derived category of (quasi-)coherent $\struct{X}$-modules but instead enlarge this category considerably. We will -- in case of an affine $X=\Spec A$ -- associate to $X$ the derived category $\derived{A_\blacksquare}$ of modules \emph{complete} over the \emph{analytic ring} $A_\blacksquare$. Defining the notion of (pre-)analytic rings the associated notion of completeness is the content of this chapter. The specific (pre-)analytic rings $A_\blacksquare$ will be introduced in section \ref{section: the pre-analytic rings Ablack and (A,R)black}, their analyticity will be shown in sections \ref{section: analyticity of Ablack and (A, R)black} and \ref{section: an analogue of hilbert's basis theorem}.

\section{Pre-analytic \& analytic rings}
\label{section: pre-analytic and analytic rings}

What is a (pre)-analytic ring supposed to be? During their investigation of condensed mathematics, Clausen and Scholze noticed that there is a reoccurring theme: Inside a category of condensed modules one finds a certain subclass of \quote{nice} objects one deems \quote{complete}. This covers for example the case of so-called \emph{solid} modules (the ones we will study in this thesis) but also various categories of analytically interesting modules like $p$-liquid vector spaces inside condensed $\R$-modules (which we will neither define nor investigate). A common theme is that one has the following data: Certain \emph{free complete} modules $\pfree{A}{S}$ associated to any extremally (or even profinite) set $S$ and for any such $S$ a natural morphism $S\to \pfree{A}{S}$ of (condensed) sets. This is what one call a \emph{pre-analytic ring}. Any module $M$ in question is then said to be \emph{complete} if for any $S$ and any specification of how elements of $S$ map into $M$ (\ie a map $S\to M$ of (condensed) sets) there is a unique extension of the map to $\pfree{A}{S}$. Depending on the free complete modules chosen this provides varying flavors of theory. In the solid case the theory is very non-archimedean, however in the $p$-liquid case the theory is well suited to \quote{traditional} functional analysis over the real numbers. As the notion of pre-analytic rings is very general, one has to impose some regularity to obtain a good theory of completeness (for example that the full subcategory of complete modules being abelian). This then leads to the notion of \emph{analyticity}, which for example implies that the free complete modules $\pfree{A}{S}$ are themselves complete (something that is not a priori clear). For an \emph{analytic ring} (a pre-analytic ring that is analytic) a myriad of good results follow at once, \confer theorem \ref{theorem: properties of the category of complete modules}. This is the content of this section.

\begin{definition}[Pre-analytic rings]
	\label{definition: pre-analytic rings}
	
	A \emph{pre-analytic ring} $\preAnaRing A$ is a triple $(\underlying A, \pfree A -, \eta)$, where $\underlying A$ is a condensed ring, called the \emph{underlying ring}, $\pfree A - \colon \extremallyDisconnectedSets \to \mod{\underlying A}$ is a functor preserving finite disjoint unions, called the \emph{free completion} and $\eta$ is a natural transformation, containing the datum of a morphism $\eta_S\colon S \to \pfree A S$ of condensed sets for every $S\in \extremallyDisconnectedSets$.
	
	A \emph{map of pre-analytic rings} $\preAnaRing{A} \to \preAnaRing{B}$ is a pair $(\phi, \epsilon)$ where $\phi$ is a map of underlying rings $\underlying A \to \underlying B$ and $\epsilon$ is a natural transformation $\pfree A - \Rightarrow \scalarRestriction{\phi}\pfree B -$ consisting of maps of $\underlying A$-modules, such that for all extremally disconnected $S$ the induced diagram
	\begin{center}
		\begin{tikzcd}[column sep=0.75em]
								&S\ar[dl]\ar[dr]	&\\
			\pfree A S\ar[rr]	&					&\scalarRestriction{\phi}\pfree B S
		\end{tikzcd}
	\end{center}
	of condensed sets commutes.
	
	Denote the category of pre-analytic rings by $\preAnalyticRings$.
\end{definition}

\begin{example}
	\label{example: examples for pre-analytic rings}
	We have the following list of examples of pre-analytic rings:
	\begin{itemize}
		\item 
		Any condensed ring $A$ naturally admits the structure of a pre-analytic with underlying ring $A$, free completion $\free{A}{-}$ and natural map $S\to \free{A}{S}$ given by adjointness. This pre-analytic ring will be denoted by $A$ as well. We obtain a faithful embedding of condensed rings into $\preAnalyticRings$.
		
		\item
		The pre-analytic ring $\Z_\blacksquare$ with underlying ring $\Z$ and free completion $S=\limit_i S_i\mapsto \free{\Z_\blacksquare}{S}\ldef \limit_i \free{\Z}{S_i}$ along with the natural map $S\to \free{\Z_\blacksquare}{S}$ induced by the cone $S_i\to \free{\Z}{S_i}$.
		
		\item
		More general we will define for any discrete ring $A$ the pre-analytic ring $A_\blacksquare$ with underlying ring $A$ in section \ref{section: the pre-analytic rings Ablack and (A,R)black}. For any discrete $R$-algebra $A$ we will also define a relative version $(A, R)_\blacksquare$.
	\end{itemize}
\end{example}

\begin{definition}[Complete modules]
	\label{definition: completeness}
	
	Let $\preAnaRing{A}$ be a pre-analytic ring. An $\underlying{A}$-module $M$ is called \emph{$\preAnaRing{A}$-complete} if for any extremally disconnected set $S$ the natural map
		$$\Hom_{\underlying{A}}(\pfree A S, M) \to \Hom(\associated{S}, M),$$
	induced by $\associated{S}\to \pfree A S$, is an isomorphism. That is, for each solid diagram
	\begin{center}
		\begin{tikzcd}
			\associated{S} &[-1.5em]&[-1.5em] M\\
			&\pfree A S
			\ar[from=1-1, to=1-3, "\forall"]
			\ar[from=1-1, to=2-2]
			\ar[from=2-2, to=1-3, "\exists!\text{ $\underlying{A}$-linear}", swap, dotted]
		\end{tikzcd}
	\end{center}
	there is a dotted arrow making the triangle commute. Write $\pCompleteMod A$ for the (full) subcategory of $\mod{\underlying{A}}$ of $\preAnaRing{A}$-complete modules.
\end{definition}

\begin{remark}[Why \quote{completeness}?]
	This definition is of course rather abstract. We will see an interpretation of $\preAnaRing{A}$-completeness in the case of $\preAnaRing{A} = \Z_\blacksquare$ in remark \ref{remark: interpreting Zblack-completeness}. In this case the completeness of some (discrete) $\underlying{A}$-module $M$ requires \emph{certain sequences of elements in $M$ to \quote{converge}}. This then of course explains why this notion is termed \emph{completeness}.
\end{remark}

\begin{remark}
	Let $\preAnaRing{A}=(\underlying{A}, \pfree A -, \eta)$ be a pre-analytic ring. Since the functor $\ufree A -$ giving free $\underlying A$-modules is left adjoint to the forgetful functor $\mod{\underlying{A}} \to \condensedSets$, the natural map $\eta$ gives rise to morphisms $\ufree A S\to \pfree A S$ of $\underlying{A}$-modules, natural in the extremally disconnected set $S$.
	
	An $\underlying{A}$-module $M$ is then $\preAnaRing{A}$-complete if and only if the induced morphism
		$$\Hom_\underlying{A}(\pfree{A}{S}, M) \to \Hom_\underlying{A}(\ufree{A}{S}, M)$$
	is an isomorphism for all extremally disconnected $S$.
\end{remark}

We will now start to investigate when the category $\completeMod{\preAnaRing{A}}$ of $\preAnaRing{A}$-complete modules is well-behaved. Informally, if one is able to show that the free completions $\pfree{A}{S}$ are $\preAnaRing{A}$-complete (which is not clear!), then one can argue that all modules presentable by such free completes are $\preAnaRing{A}$-complete as well. 
\begin{definition}[Free and finitely free complete modules]
	Let $\preAnaRing{A}$ be a pre-analytic ring. An $\underlying{A}$-module $M$ is called \emph{finitely-$\preAnaRing{A}$-free} if there is an isomorphism $M\cong \pfree A S$ for some extremally disconnected set $S$. If more generally $M$ is a direct sum of finitely-$\preAnaRing{A}$-free modules then $M$ is called \emph{$\preAnaRing{A}$-free}.
\end{definition}

\begin{definition}[Presentable modules]
	Let $\preAnaRing{A}$ be a pre-analytic ring. An $\underlying{A}$-module $M$ is called \emph{$\preAnaRing{A}$-presentable} if there is an exact sequence
	$$\bigoplus_{i\in I} \pfree A {S'_i} \to \bigoplus_{j\in J} \pfree A {S_j}\to M\to 0$$
	presenting $M$ in $\mod{\underlying A}$ via $\preAnaRing{A}$-free modules. Write $\pPresentableMod A$ for the (full) subcategory of $\mod{\underlying{A}}$ of $\preAnaRing{A}$-presentable modules.
\end{definition}

As mentioned before, $\preAnaRing{A}$-free modules are not automatically $\preAnaRing{A}$-complete. Solely imposing this condition is however not enough to build a good (derived) theory of $\preAnaRing{A}$-complete modules. Instead, one essentially requires that bound to the right complexes of $\preAnaRing{A}$-free modules are derived $\preAnaRing{A}$-complete.
\begin{definition}[(Derived) completeness]
	Let $\preAnaRing{A}$ be a pre-analytic ring. A complex $C^\bullet \in \derived{\underlying{A}}$ is called \emph{(derived) $\preAnaRing{A}$-complete} if for any extremally disconnected $S$ the morphism
		$$\eRHom_\underlying{A}(\pfree{A}{S}, C^\bullet) \to \eRHom_\underlying{A}(\ufree{A}{S}, C^\bullet)$$
	in $\derived{\Ab}$ is an isomorphism.
\end{definition}

We can observe that derived completeness of an object concentrated in degree $0$ implies ordinary completeness.
\begin{remark}[Derived completeness implies ordinary completeness]
	\label{remark: derived completeness implies ordinary completeness}
	
	Suppose that $X$ is an $\underlying{A}$-module such that $X[0]$ is derived $\preAnaRing{A}$-complete. As $H^0\circ \eRHom_\underlying{A} \cong \eHom_\underlying{A}$ on objects of $\mod{\underlying{A}}$ and since $H^0$ preserves isomorphisms, we obtain that for any $S$ extremally disconnected the induced
		$$\eHom_\underlying{A}(\pfree{A}{S}, X) \to \eHom_\underlying{A}(\ufree{A}{S}, X)$$
	is an isomorphism, hence that $X$ is $\preAnaRing{A}$-complete. It is a priori not clear that the converse holds true! If $\preAnaRing{A}$ is however analytic (see definition \ref{definition: analytic rings}, then it is true a fortiori by theorem \ref{theorem: properties of the category of complete modules}.
\end{remark}

We can now state the definition for an analytic ring. As just mentioned, one needs to impose the condition that bound to the right complexes of $\preAnaRing{A}$-free modules are $\preAnaRing{A}$-complete. Even more, we will require the analogous morphism of condensed abelian groups $\eRHom_\underlying{A}$ to be an isomorphism. This will allow us to control the internal $\intHom_\underlying{A}$'s of $\preAnaRing{A}$-complete modules later on in theorem \ref{theorem: properties of the category of complete modules}.
\begin{definition}[Analytic rings]
	\label{definition: analytic rings}
	
	A pre-analytic ring $\preAnaRing{A}$ is called \emph{analytic} if for any complex
	$$C^\bullet = \left(\dots \to C^{-2} \to C^{-1} \to C^0 \to 0\right)$$
	of $\preAnaRing{A}$-free modules and any extremally disconnected $S$ the induced morphism
	$$\eRHom_\underlying{A}(\preAnaRing A[S], C^\bullet) \to \eRHom_\underlying A(\underlying A[S], C^\bullet)$$
	in $\derived{\condensedAb}$ is an isomorphism.
\end{definition}

\begin{remark}[Free completions are complete]
	If $\preAnaRing{A}$ is analytic and $S$ is extremally disconnected then $\pfree{A}{S}$ is $\preAnaRing{A}$-complete. Indeed, the complex $C^\bullet=\pfree{A}{S}[0]$ is derived $\preAnaRing{A}$-complete and the claim follows by remark \ref{remark: derived completeness implies ordinary completeness}.
\end{remark}

As per usual, we need a good notion of map between analytic rings. Such a map should not(!) simply be a map of the underlying pre-analytic rings.
\begin{definition}[Maps of analytic rings]
	\label{definition: maps of analytic rings}
	
	Let $\preAnaRing A=(\underlying{A}, \pfree{A}{-}, \eta)$ and $\preAnaRing{B}=(\underlying{B}, \pfree{B}{-}, \eta')$ be analytic rings. A \emph{map of analytic rings} $\preAnaRing{A}\to\preAnaRing B$ is a map of the underlying condensed rings $\phi\colon \underlying{A}\to \underlying{B}$ such that for each extremally disconnected $S$ the free module $\pfree B S$ considered as an $\underlying{A}$-module is $\preAnaRing{A}$-complete, \ie the scalar restriction $\scalarRestriction{\phi} \pfree B S$ is $\preAnaRing{A}$-complete.
\end{definition}

The following remark shows how any map of analytic rings induces a map of pre-analytic rings.
\begin{remark}[Maps of analytic rings induce maps of pre-analytic rings]
	\label{remark: maps of analytic rings induce maps of pre-analytic rings}
	
	In the setting of definition \ref{definition: maps of analytic rings}, since $\scalarRestriction{\phi}\pfree B S$ is $\preAnaRing{A}$-complete, the map of condensed sets $\associated S\to \scalarRestriction{\phi}\pfree B S$ factors through $\associated S\to \pfree A S$ via a unique map $\pfree A S \to \scalarRestriction{\phi} \pfree B S$. Hence, a map of analytic rings induces a corresponding map of pre-analytic rings.
\end{remark}

Now to define a map of analytic rings $\preAnaRing{A} \to \preAnaRing{B}$ one has to give a map of underlying rings $\phi\colon \underlying{A}\to\underlying{B}$ and ensure that for any extremally disconnected $S$ the scalar restriction $\scalarRestriction{\phi}\pfree{B}{S}$ is $\preAnaRing{A}$-complete -- this is potentially quite inconvenient. In practice however it is often good enough to specify a map of pre-analytic rings instead.
\begin{proposition}[Maps of analytic rings from maps of pre-analytic rings]
	\label{proposition: maps of analytic rings from maps of pre-analytic rings}
	
	Suppose that $\preAnaRing{A}$ and $\preAnaRing{B}$ are analytic rings and that $(\phi, \epsilon)\colon \preAnaRing{A}\to \preAnaRing{B}$ is a map of pre-analytic rings. If either $\underlying{A}$ is discrete or $\scalarRestriction{\phi}\pfree{B}{\ast}$ is $\preAnaRing{A}$-complete, then $\phi\colon \underlying{A}\to\underlying{B}$ is a map of analytic rings $\preAnaRing{A}\to\preAnaRing{B}$.
\end{proposition}
\begin{shortproof}
	This is weakened version of \cite[Prop. 7.14]{condenseddotpdf}, keeping in mind \cite[Rem. 7.15]{condenseddotpdf}.
\end{shortproof}

We can now state the main result of this section: For any analytic ring $\preAnaRing{A}$ the category of $\preAnaRing{A}$-complete modules is an extremely well-behaved subcategory of $\mod{\underlying{A}}$.
\begin{theorem}[The category of complete $\preAnaRing{A}$-modules]
	\label{theorem: properties of the category of complete modules} 
	
	Let $\preAnaRing{A}$ be an analytic ring.
	\begin{enumerate}
		\item
		The category $\pCompleteMod{A}$ of complete $\preAnaRing{A}$-modules is an abelian category stable under all limits, colimits and extensions in $\mod{\underlying{A}}$. The objects $\pfree A S$ for $S$ extremally disconnected form a (proper) class of compact projective generators of $\pCompleteMod{A}$. The inclusion $\pCompleteMod{A} \to \mod{\underlying{A}}$ admits a left adjoint
			$$\completion{-}{\preAnaRing{A}}\ldef \completionAsTensor{\preAnaRing{A}}{\underlying{A}}{-}\colon\mod{\underlying{A}} \to \pCompleteMod{A}$$
		which is the unique colimit preserving extension of the functor $\ufree A S \mapsto \pfree A S$ defined on the image of $\ufree{A}{-}\colon \extremallyDisconnectedSets\to \mod{\underlying{A}}$ and is called the \emph{$\preAnaRing{A}$-completion} or simply \emph{completion}. If $\underlying{A}$ is commutative, there is a unique symmetric monoidal tensor product $\otimes_\preAnaRing{A}$ on $\pCompleteMod{A}$ for which the completion $\completion{-}{\preAnaRing{A}}$ is (strong) symmetric monoidal. Furthermore, $\pPresentableMod{A}=\pCompleteMod{A}$, that is $\preAnaRing{A}$-presentability and  $\preAnaRing{A}$-completeness agree.
		
		\item
		The functor
			$$\derived{\pCompleteMod{A}} \to \derived{\underlying{A}}$$
		induced by the inclusion $\pCompleteMod{A} \to \mod{\underlying A}$ is fully faithful. Its essential image is stable under all limits and colimits in $\derived{\underlying{A}}$ and is given by those complexes $C^\bullet \in \derived{\underlying{A}}$ for which one of the two equivalent conditions is met:
		\begin{itemize}
			\item
			$C^\bullet$ is derived $\preAnaRing{A}$-complete, that is for all extremally disconnected $S$ the induced map
				$$\eRHom_{\underlying A}(\pfree A S, C^\bullet) \to \eRHom_{\underlying A}(\ufree A S, C^\bullet)$$
			in $\derived{\Ab}$ is an isomorphism. In this case, also the $\eRHom_{\underlying{A}}$s in $\derived{\condensedAb}$ agree.
			
			\item
			For all $n\in\Z$ the $\underlying A$-module $H^n(C^\bullet)$ is $\preAnaRing{A}$-complete.
		\end{itemize}
		Furthermore, the inclusion $\derived{\pCompleteMod{A}} \to \derived{\mod{\underlying{A}}}$ admits a left adjoint
			$$\derivedCompletion{-}{\preAnaRing{A}}\ldef \derivedCompletionAsTensor{\preAnaRing{A}}{\underlying{A}}{-}\colon \derived{\mod{\underlying{A}}}\to\derived{\pCompleteMod{A}}$$
		which is the left-derived functor of $\completion{-}{\preAnaRing{A}} = \completionAsTensor{\preAnaRing{A}}{\underlying{A}}{-}$. If $\underlying{A}$ is commutative, then $\derived{\pCompleteMod{A}}$ admits a unique symmetric monoidal tensor product $\Dotimes_\preAnaRing{A}$ for which $\derivedCompletion{-}{\preAnaRing{A}}$ is (strong) symmetric monoidal.
	\end{enumerate}
\end{theorem}
\begin{proof}
	Most of this follows from lemma \ref{lemma: THE LEMMA} -- once we have proven it in the next section. What remains is that the $\eRHom_\underlying{A}$s in $\derived{\condensedAb}$ agree and that $\derived{\completeMod{\preAnaRing{A}}}$ is symmetric monoidal in case of a commutative $\underlying{A}$.
	
	So suppose that $X^\bullet \in \derived{\completeMod{\preAnaRing{A}}}$.
	Clausen and Scholze argue in the proof of \cite[Prop. 7.5, p.46]{condenseddotpdf} that the claim reduces to the case where $X^\bullet$ is bounded to the right. The author is not sure how they argue this, as he sees no reason why in general $\pfree{A}{S}$ should be $\condensedAb$-enriched compact (which would ensure the commutation with the homotopy colimit $\colimit_{n\in \N} \cTruncation{\le n}X^\bullet$) or why $\eRHom_\underlying{A}$ should commute with limits (which would ensure the commutation with the limit $\limit_{n\in \N} \sTruncation{\le n}X^\bullet$) -- only commutation with homotopy limits seems ensured. We thus have to make it an \assumptionPreWord that we can reduce to the case that $X^\bullet$ is bounded to the right. As $\eRHom_A$ is compatible with shifts we can further assume that $X^\bullet$ is concentrated in non-positive degrees. Then $X^\bullet$ admits a left resolution by a complex $C^\bullet$ of $\preAnaRing{A}$-free modules that is concentrated in non-positive degrees too. By definition of analyticity of $\preAnaRing{A}$ we then have
		$$\eRHom_{\underlying{A}}(\pfree{A}{S}, X^\bullet) \cong \eRHom_\underlying{A}(\pfree{A}{S}, C^\bullet) \isorightarrow \eRHom_\underlying{A}(\ufree{A}{S}, C^\bullet) \cong \eRHom_{\underlying{A}}(\ufree{A}{S}, X^\bullet)$$
	in $\derived{\condensedAb}$ which shows the claim.
	
	It remains to argue the existence of the monoidal structure on $\derived{\completeMod{\preAnaRing{A}}}$ for which the derived completion is (strong) symmetric monoidal in case that $\underlying{A}$ is commutative. For this we will simply reference the proof given by Clausen and Scholze for \cite[Prop. 7.5]{condenseddotpdf}.
\end{proof}

\begin{notation}
	If $\preAnaRing{A}$ is an analytic ring we will abbreviate
		$$\mod{\preAnaRing{A}}\ldef \pCompleteMod{\preAnaRing{A}}\text{ and }\derived{\preAnaRing{A}}\ldef \derived{\pCompleteMod{A}}$$
	and refer to these categories as the \emph{category of $\preAnaRing{A}$-modules} and the \emph{derived category of $\preAnaRing{A}$-modules} respectively.
\end{notation}

\begin{remark}[$\Dotimes_\preAnaRing{A}$ and $\otimes_\preAnaRing{A}$]
	\label{remark: when is the tensor product the derived tensor product?}
	
	There is a notational subtlety: For a general analytic ring $\preAnaRing{A}$ the functor $\Dotimes_\preAnaRing{A}$ might not be the left derived functor of $\otimes_\preAnaRing{A}$!
	
	Clausen and Scholze state in \cite[Warning 7.6]{condenseddotpdf} that this however holds true in any instance know to them. Furthermore, they state that one might equivalently check that for any two extremally disconnected sets $S$ and $T$ the complex of $\underlying{A}$-modules $\pfree{A}{S\times T}\ldef \nerve{\pfree{A}{S_\bullet}}$ (obtained by simplicially resolving the non-extremally disconnected set $S\times T$ by extremally disconnecteds, \ie choosing a simplicial hypercover $S_\bullet$ of $S\times T$ of extremally disconnected sets $S_n$, giving a simplicial $\underlying{A}$-module $\pfree{A}{S_\bullet}$ and attaching to it by the Dold--Kan correspondence a connective complex $\nerve{\pfree{A}{S_\bullet}}$) is concentrated in degree $0$.
	
	We will see in remark \ref{remark: free completions for profinite sets} that we can -- for all analytic rings introduced in this thesis -- safely ignore the issue and always assume $\Dotimes_\preAnaRing{A}$ to be the left derived functor of $\otimes_\preAnaRing{A}$.
\end{remark}

Now knowing, that complete modules give rise to extremally well-behaved abelian categories, we should also investigate how one can transport complete modules along a map of analytic rings.
\begin{proposition}[Scalar restriction and extension of complete modules]
	\label{proposition: scalar restriction and extension of complete modules}
	
	Let $\phi\colon \preAnaRing{A}\to \preAnaRing{B}$ be a map of analytic rings.
	\begin{enumerate}
		\item
		The scalar restriction of any $\preAnaRing{B}$-complete module is $\preAnaRing{A}$-complete and the induced functor $\scalarRestriction{\phi}\colon \mod{\preAnaRing{B}}\to\mod{\preAnaRing{A}}$ admits a left adjoint
			$$\completedScalarExtension{\phi}{\preAnaRing{B}}\ldef \preAnaRing{B}\otimes_{\preAnaRing{A}}-\colon \mod{\preAnaRing{A}}\to\mod{\preAnaRing{B}}$$
		that is the unique colimit preserving extension of $\pfree{A}{S}\mapsto \pfree{B}{S}$ and is called the \emph{completed scalar extension}.
		
		
		\item
		The scalar restriction of any complex in $\derived{\preAnaRing{B}}$ is $\preAnaRing{A}$-complete and the induced functor $\scalarRestriction{\phi}\colon \derived{\preAnaRing{B}}\to\derived{\preAnaRing{A}}$ admits a left adjoint
			$$\derivedCompletedScalarExtension{\phi}{\preAnaRing{B}}\ldef \preAnaRing{B}\Dotimes_\preAnaRing{A} -\colon \derived{\preAnaRing{A}}\to\derived{\preAnaRing{B}}$$
		that is the left derived functor of the completed scalar extension $\completedScalarExtension{\phi}{\preAnaRing{B}}$.
		
	\end{enumerate}
\end{proposition}
\begin{proof}
	Denote by $i_\preAnaRing{A}\colon \mod{\preAnaRing{A}}\to\mod{\underlying{A}}$ and $i_\preAnaRing{B}\colon \mod{\preAnaRing{B}}\to\mod{\underlying{B}}$ the exact inclusions.
	\begin{enumerate}
		\item
		By definition of maps of analytic rings, for any extremally disconnected $S$, the scalar restriction $\scalarRestriction{\phi}\pfree{B}{S}$ is $\preAnaRing{A}$-complete. As scalar restriction commutes with arbitrary colimits and $\mod{\preAnaRing{B}}$ is generated, under colimits, by such $\pfree{B}{S}$, we obtain that scalar restriction restricts to a functor $\mod{\preAnaRing{B}}\to\mod{\preAnaRing{A}}$. Hence, we obtain the commutative square
		$$\begin{tikzcd}
			\mod{\preAnaRing{B}} \ar[r, "\scalarRestriction{\phi}", swap]\ar[d, "i_\preAnaRing{B}"] &\mod{\preAnaRing{A}}\ar[d, "i_\preAnaRing{A}", swap]\\
			\mod{\underlying{B}}\ar[r, "\scalarRestriction{\phi}"]&\mod{\underlying{A}}.
		\end{tikzcd}$$
		Observe that the inclusions $i_\preAnaRing{A}$ and $i_\preAnaRing{B}$ as well as the scalar restriction $\scalarRestriction{\phi}\colon \mod{\underlying{B}}\to \mod{\underlying{A}}$ admit left adjoints given by the completions $\completion{-}{\preAnaRing{A}}$ and $\completion{-}{\preAnaRing{B}}$ as well as scalar extension $\scalarExtension{\phi}\colon \mod{\underlying{A}}\to\mod{\underlying{B}}$. Define the functor
			$$\completedScalarExtension{\phi}{\preAnaRing{B}}\colon \mod{\preAnaRing{A}}\to\mod{\preAnaRing{B}}$$
		as the composition $M\mapsto \completion{\scalarExtension{\phi} i_\preAnaRing{A}M}{\preAnaRing{B}}$. Observe, that for any $\preAnaRing{A}$-complete module $M$ and any $\preAnaRing{B}$-complete module $N$ we have a natural isomorphism
			$$\Hom(\completedScalarExtension{\phi}{\preAnaRing{B}}M, N) \cong \Hom(i_\preAnaRing{A}M, \scalarRestriction{\phi}i_\preAnaRing{B} N) = \Hom(i_\preAnaRing{A}M, i_\preAnaRing{A}\scalarRestriction{\phi} N) \cong \Hom(M, \scalarRestriction{\phi} N)$$
		by adjointness, by commutativity of the above square and by full faithfulness of $i_\preAnaRing{A}$. In particular, $\completedScalarExtension{\phi}{\preAnaRing{B}}$ is the unique (up to unique isomorphism) left adjoint of $\scalarRestriction{\phi}\colon \mod{\preAnaRing{B}}\to\mod{\preAnaRing{A}}$. Thus, all functors in the above commutative square admit a left adjoint and the induced diagram
		$$\begin{tikzcd}
			\mod{\preAnaRing{B}} &\mod{\preAnaRing{A}}\ar[l, "\completedScalarExtension{\phi}{\preAnaRing{B}}", swap]\\
			\mod{\underlying{B}}\ar[ur, Rightarrow, "\sim" {sloped, near start}]\ar[u, "{\completion{-}{\preAnaRing{B}} = \preAnaRing{B}\otimes_\underlying{B} -}"] &\mod{\underlying{A}}\ar[l, "\scalarExtension{\phi}"]\ar[u, "{\completion{-}{\preAnaRing{A}} = \preAnaRing{A}\otimes_\underlying{A} -}", swap].
		\end{tikzcd}$$
		of left adjoints must commute up to unique isomorphism as well, as the composition of two left adjoints is a left adjoint of the composition of the original functors. In particular,
			$$\completedScalarExtension{\phi}{\preAnaRing{B}}\pfree{A}{S} \cong  \completedScalarExtension{\phi}{\preAnaRing{B}}\completion{\ufree{A}{S}}{\preAnaRing{A}} \cong \completion{\scalarExtension{\phi}\ufree{A}{S}}{\preAnaRing{B}} = \completionAsTensor{\preAnaRing{B}}{\underlying{B}}{\underlying{B}\otimes_\underlying{A}\ufree{A}{S}} \cong \completionAsTensor{\preAnaRing{B}}{\underlying{B}}{\ufree{B}{S}} \cong \pfree{B}{S}$$
		so $\completedScalarExtension{\phi}{\preAnaRing{B}}$ is the unique colimit preserving extension of $\pfree{A}{S} \mapsto \pfree{B}{S}$ as claimed.
		
		\item 
		Scalar restriction $\scalarRestriction{\phi}$ is exact. Recall that for an analytic ring, a complex in the derived category is complete if and only if all of its cohomology groups are complete. Thus, if $C^\bullet\in\derived{\preAnaRing{B}}$ then for each $n\in \Z$ the $\underlying{B}$-module $H^n(C^\bullet)$ is $\preAnaRing{B}$-complete. But as $\scalarRestriction{\phi}$ is exact we obtain that $H^n(\scalarRestriction{\phi}C^\bullet) \cong \scalarRestriction{\phi}H^n(C^\bullet)$ is $\preAnaRing{A}$-complete by the first part and hence $\scalarRestriction{\phi}C^\bullet \in \derived{\preAnaRing{A}}$, that is $\scalarRestriction{\phi}$ restricts to a functor $\derived{\preAnaRing{B}}\to\derived{\preAnaRing{A}}$. The remaining proof is then analogous to the first part.
	\end{enumerate}
\end{proof}

\begin{remark}[Notational compatibility]
	Recall that any condensed ring $A$ can be considered as a pre-analytic ring, denoted $A$ as well. This pre-analytic ring $A$ is clearly analytic and any module is automatically complete. Then the notation $\mod{A}$ for complete modules, as well as the notation $\preAnaRing{A}\otimes_{A} -$ for any analytic ring $A\to \preAnaRing{A}$ over $A$ and $A\otimes_\preAnaRing{A}-$ for any analytic ring $\preAnaRing{A}\to A$ under $A$ are compatible with this identification.
\end{remark}

We can further observe the following.
\begin{lemma}[Completed scalar extension is monoidal]
	\label{lemma: completed scalar extension is monoidal}
	
	Let $\phi\colon \preAnaRing{A}\to \preAnaRing{B}$ be a map of analytic rings. Then completed scalar extension
		$$\completedScalarExtension{\phi}{\preAnaRing{B}}\colon \mod{\preAnaRing{A}}\to\mod{\preAnaRing{B}}$$
	and derived completed scalar extension
		$$\derivedCompletedScalarExtension{\phi}{\preAnaRing{B}}\colon \derived{\preAnaRing{A}}\to\derived{\preAnaRing{B}}$$
	are (strong) symmetric monoidal.
\end{lemma}
\begin{proof}
	We will prove the claim only for $\completedScalarExtension{\phi}{\preAnaRing{B}}$ as the other case is analogous. Observe first that for any two extremally disconnected $S$ and $T$ we have
		$$\pfree{A}{S}\otimes_\preAnaRing{A} \pfree{A}{T} \cong \ufree{A}{S}^\completionSymbol{\preAnaRing{A}} \otimes_\preAnaRing{A} \ufree{A}{T}^\completionSymbol{\preAnaRing{A}}\cong\completion{\ufree{A}{S}\otimes_\underlying{A}\ufree{A}{T}}{\preAnaRing{A}} \cong \ufree{A}{S\times T}^\completionSymbol{\preAnaRing{A}}$$
	as $\completion{-}{\preAnaRing{A}}$ is (strong) monoidal. Similarly, $\pfree{B}{S}\otimes_\preAnaRing{B} \pfree{B}{T} \cong \ufree{B}{S\times T}^\completionSymbol{\preAnaRing{B}}$ for $\preAnaRing{B}$. Recall further, from the proof of proposition \ref{proposition: scalar restriction and extension of complete modules}, that $\completedScalarExtension{\phi}{\preAnaRing{B}}\circ \completion{-}{\preAnaRing{A}} \cong \completion{-}{\preAnaRing{B}}\circ \scalarExtension{\phi}$. We now obtain
	\begin{align*}
		\completedScalarExtension{\phi}{\preAnaRing{B}}\mleft(\pfree{A}{S}\otimes_\preAnaRing{A}\pfree{A}{T}\mright)
		&\cong \completedScalarExtension{\phi}{\preAnaRing{B}} \mleft(\ufree{A}{S\times T}^\completionSymbol{\preAnaRing{A}}\mright)\\
		&\cong \mleft(\scalarExtension{\phi}\ufree{A}{S\times T}\mright)^\completionSymbol{\preAnaRing{B}}\\
		&\cong \ufree{B}{S\times T}^\completionSymbol{\preAnaRing{B}}\\
		&\cong\pfree{B}{S}\otimes_\preAnaRing{B}\pfree{B}{T} \cong \completedScalarExtension{\phi}{\preAnaRing{B}}\pfree{A}{S} \otimes_\preAnaRing{B} \completedScalarExtension{\phi}{\preAnaRing{B}}\pfree{A}{T}
	\end{align*}
	which proves that $\completedScalarExtension{\phi}{\preAnaRing{B}}$ is (strong) monoidal on generators. As $\completedScalarExtension{\phi}{\preAnaRing{B}}$ as well as $\otimes_\preAnaRing{A}$ and $\otimes_\preAnaRing{B}$ commute with arbitrary colimits, this already shows the claim.
\end{proof}

\section{The Underlying Machinery}
\label{section: the underlying machinery}

In this section we will prove the machinery underlying the main results of chapter \ref{chapter: analytic rings and completeness}. 
\begin{notation}
	Throughout this section, let $\category A$ be a cocomplete abelian category that admits a subcategory $i\colon \category A_0 \subseteq \category A$ of compact projective generators. Furthermore, let $F\colon \category A_0\to \category A$ be a functor equipped with a natural transformation $\eta\colon i\Rightarrow F$.
\end{notation}
Notice that this data mimics indeed mimics the setting of section \ref{section: pre-analytic and analytic rings}. Given some pre-analytic ring $\preAnaRing{A}$, we have a cocomplete abelian category $\category A\ldef\mod{\underlying{A}}$ of $\underlying{A}$-modules generated by the class $\category A_0\ldef \{\free {\underlying{A}} S\setseparator S\text{ extremally disconnected}\}$ of $\underlying{A}$-free modules (which are compact projective). The free completion $\pfree A -$ and the natural morphisms $S\to \pfree A S$ induce a functor $\category A_0\to \category A, \free{\underlying{A}}S\mapsto \pfree A S$ and a natural transformation $i\Rightarrow F$.

Let us translate some terminology know from the previous setting.
\begin{definition}[$F$-Freeness]
	An object $X\in\category A$ is called \emph{finitely $F$-free} if there is an isomorphism $X\cong F(G)$ for some $G\in\category A_0$. Furthermore, $X$ is called \emph{$F$-free} if $X$ is isomorphic to a (not necessarily finite) direct sum of finitely $F$-free modules.
\end{definition}

\begin{definition}[$F$-Presentability]
	An object $X\in\category A$ is called \emph{$F$-presentable} if there is an exact sequence
	$$\bigoplus_{i\in I} F(G'_i) \to \bigoplus_{j\in J} F(G_j)\to X\to 0,$$
	presenting $X$ in $\category A$ via $F$-free modules. Write $\presentableMod F$ for the (full) subcategory of $\category A$ of $F$-presentable modules.
\end{definition}

\begin{definition}[$F$-Completeness]
	An object $X\in\category A$ is called \emph{$F$-complete} if for any generator $G\in\category A_0$ the natural map of abelian groups
		$$\eHom(F(G), X) \to \eHom(G, X)$$
	induced by $\eta_G\colon G\to F(G)$ is an isomorphism. Write $\completeMod F$ for the (full) subcategory of $\category A$ of $F$-complete modules. Similarly, a complex $C^\bullet \in \derived{\category A}$ is called \emph{(derived) $F$-complete} if for any generator $G\in\category A_0$ the induced map
		$$\eRHom(F(G), C^\bullet) \to \eRHom(G, C^\bullet)$$
	in $\derived{\Ab}$ is an isomorphism. Write $\categoryname{D}_F(\category A)$ for the full subcategory of $\derived{\category{A}}$ of all derived $F$-complete complexes.
\end{definition}

\begin{remark}[$F$-complete objects vs. $F$-complete complexes]
	Suppose that $X\in\category{A}$ is an object such that the complex $X[0]$ concentrated in degree $0$ is $F$-complete. Then $X$ itself is $F$-complete, since $H^0$ clearly preserves isomorphisms and there is a natural isomorphism $\RightD^0\eHom(-, X[0]) \cong \eHom(-, X)$ of functors $\category{A}^\op \to \Ab$.
\end{remark}

With these definitions we are well-equipped to state the first main result of this section. Under mild assumptions, it ensures that $F$-completeness and $F$-presentability agree and that the class of such modules behaves very well. Furthermore, this well-behavedness passes to the derived setting.
\begin{lemma}[The Underlying Lemma]
	\label{lemma: THE LEMMA}
	
	Assume that any object $X\in\presentableMod{F}$ viewed as a complex concentrated in degree $0$ is derived $F$-complete, that is for any generator $G\in\category A_0$ the morphism
	$$\eRHom(F(G), X) \to \eRHom(G, X)$$
	in $\derived{\Ab}$ is an isomorphism. Then the following properties hold true:
	\begin{itemize}
		\item
		The category $\completeMod{F}$ is an abelian subcategory stable under all limits, colimits and extensions in $\category A$ and the objects $F(G)$ for $G\in\category A_0$ form a family of compact projective generators. The inclusion $\completeMod F\subseteq \category A$ admits a left-adjoint $\tilde F\colon \category A\to\completeMod{F}$ that is the unique colimit preserving extension of $F\colon \category A_0\to \category A$. Additionally, $\presentableMod{F} = \completeMod{F}$, that is $F$-presentability and $F$-completeness agree.
		
		\item
		Furthermore, the functor $\derived{\completeMod F}\to \derived{\category A}$ is fully faithful and identifies the derived category $\derived{\completeMod{F}}$ with $\categoryname{D}_F(\category A)$. A complex $C^\bullet\in\derived{\category A}$ lies in $\categoryname{D}_F(\category A)$ (\ie is derived $F$-complete) if and only if all $n\in \Z$ the objects $H^n(C^\bullet)\in\category A$ are $F$-complete. The inclusion
		$$\derived{\completeMod{F}}\equiv\categoryname{D}_F(\category A) \to\derived{\category A}$$
		admits a left-adjoint, which is given by left-derived functor of $\tilde F$.
	\end{itemize}
\end{lemma}
\begin{proof}
	Suppose $X=\limit_{i\in I} X_i \in\category A$ is a limit (in $\category A$) of objects $X_i\in\completeMod F$. Since for any generator $G\in\category A_0$, the functors $\eHom(F(G), -)$ and $\eHom(G, -)$ map limits to limits, we obtain the commutative square
	\begin{center}
		\begin{tikzcd}[column sep=1em]
			\eHom(F(G), X) 						&\eHom(G, X)\\
			\limit_{i\in I} \eHom(F(G), X_i) 	&\limit_{i\in I}\eHom(G, X_i)
			\ar[from=1-1, to=1-2]
			\ar[from=1-1, to=2-1, "\sim", sloped]
			\ar[from=1-2, to=2-2, "\sim", sloped]
			\ar[from=2-1, to=2-2]
		\end{tikzcd}
	\end{center}
	where the bottom morphism is an isomorphism by $F$-completeness of the $X_i$. Hence, the top morphism is an isomorphism and $X$ is $F$-complete, so $\completeMod{F}$ is stable under limits in $\category A$.
	
	We will now show that $\completeMod{F}$ is complete under cokernels from $\category A$. For this let $f\colon X\to Y$ be any morphism in $\completeMod{F}$, we will show that the cokernel $\coKernel f$ in $\category A$ is $F$-complete. By generation, choose an epimorphism $p\colon \bigoplus_{i\in I} G_i \to Y$ with $G_i\in\category A_0$ for $i\in I$. Since $Y$ is $F$-complete we have that $\eta_{G_i}^\ast\colon\eHom(F(G_i), Y)\to \eHom(G_i, Y)$ is an isomorphism for each $i\in I$. Hence, we obtain a morphism $p'\colon\bigoplus_{i\in I} F(G_i) \to Y$ for which $p = p'\circ (\eta_{G_i}^\ast)_{i\in I}$. Since $p$ was an epimorphism, so must be $p'$. Write $\widetilde Y\ldef \bigoplus_{i\in I}F(G_i)$. Then since $p'$ is an epimorphism by lemma \cite[\href{https://stacks.math.columbia.edu/tag/08N4}{Tag 08N4}]{stacks-project}, the pullback square
	\begin{center}
		\begin{tikzcd}
			\tilde X 	&\widetilde Y\\
			X 			&Y
			\ar[from=1-1, to=1-2, "\tilde f"]
			\ar[from=1-1, to=2-1]
			\ar[from=1-2, to=2-2, "p'"]
			\ar[from=2-1, to=2-2, "f"]
		\end{tikzcd}
	\end{center}
	is also a pushout. Hence, by lemma \cite[\href{https://stacks.math.columbia.edu/tag/08N3}{Tag 08N3}]{stacks-project} the induced map $\coKernel{\tilde f}\to \coKernel f$ is an isomorphism. By replacing $f$ by $\tilde f$, we now assume without loss of generality that $Y$ is $F$-free. One can choose a similar epimorphism $q\colon\bigoplus_{j\in J}F(G'_j)\to X$ and argue that the cokernels of $f$ and $f\circ q$ agree (since $q$ is epic). Thus, $X$ can be chosen $F$-free as well. Hence, the cokernel is $F$-presentable so $F$-complete both as a complex concentrated in degree $0$ (by assumption) and then as an object in $\category A$ (as $\RightD^0\eHom \cong \eHom$ on objects of $\category{A}$).
	
	In particular, we now obtain $\presentableMod{F}=\completeMod{F}$, since on the one hand $\presentableMod{F}\subseteq\completeMod{F}$ by assumption and on the other hand, every $X\in \completeMod{F}$ can be realized as the cokernel of $0\to X$, which by the above argument is $F$-presentable. It is then clear that $\completeMod{F}=\presentableMod{F}$ is closed under arbitrary direct sums from $\category A$ and hence closed under arbitrary colimits from $\category A$. 
	
	Thus, $\completeMod{F}$ is a full subcategory of the abelian category $\category{A}$ which is stable under all limits and colimits from $\category{A}$ and hence must be an abelian subcategory. It is further clear from $\completeMod{F}=\presentableMod{F}$ that $\completeMod{F}$ is closed under extensions from $\category{A}$
	
	Now since colimits in $\completeMod{F}$ agree with those in $\category{A}$ and since for any compact projective generator $G\in\category{A}_0$ of $\category{A}$ we have $\eHom_{\completeMod{F}}(F(G), -) \cong \eHom_{\category{A}}(G, -)$ by $F$-completeness, the objects $F(G)$ are compact projective as well. As $\completeMod{F} = \presentableMod{F}$ these objects $F(G)$ surely generate as well.
	
	Now let $\tilde F\colon \category{A}\to\completeMod{F}$ be the (necessarily unique) colimit preserving extension of the functor $F\colon \category{A}_0\to \completeMod{F}$ defined on generators of $\category{A}$. Let $X\in \category A$ and $Y\in \completeMod{F}$. As $\category{A}_0$ generates $\category A$ we can write $X\cong \colimit_{i\in I} G_i$ as a colimit with $G_i \in \category{A}_0$. Then
		$$\eHom(\tilde F(X), Y) \cong \limit_i \eHom(F(G_i), Y) \cong \limit_i \eHom(G_i, Y) \cong \eHom(X, Y)$$
	as $Y$ is $F$-complete,	showing that $\tilde F$ is left adjoint to the inclusion $\completeMod{F}\to \category{A}$.
	
	By lemma \ref{lemma: right derived Hom and homotopy colimits in the first argument}, it is enough to check faithfulness on generators. For this let $F(G)\in\completeMod{F}$ with $G\in\category{A}_0$ be such a generator and $C^\bullet \in \derived{\completeMod{F}}$ be general. We will show that the natural morphism
		$$\eRHom_\completeMod{F}(F(G), C^\bullet) \to \eRHom_\category{A}(F(G), C^\bullet) \isorightarrow \eRHom_\category{A}(G, C^\bullet)$$
	is an isomorphism as then also
		$$\Hom_\derived{\completeMod{F}}(F(G), C^\bullet) \to \Hom_\derived{\category{A}}(F(G), C^\bullet)$$
	is an isomorphism by \ref{lemma: eRHom enriches Hom-sets in the derived category}. To do so we will show that this morphism induces an isomorphism in each degree $n\in \Z$. Writing the complex $C^\bullet$ as the iterated sequential colimit of its truncations $C^\bullet_{[n, m]}\ldef \cTruncation{\ge n}\cTruncation{\le m} C^\bullet$, we can reduce the claim to the bounded complexes $C^\bullet_{[n, m]}$ by lemma \ref{lemma: commutation with sequential colimits} as $F(G)$ and $G$ are (derived) compact in their respective (derived) categories. Using the distinguished triangles
		$$C^\bullet_{[n+1, m]}\to C^\bullet_{[n, m]}\to C^n[-n]\to C^\bullet_{[n+1, m]}[1]$$
	for $n<m$ we can further reduce in finitely many steps to the case that $C^\bullet$ is concentrated in a single degree (if $n\ge m$ then either $C^\bullet_{[n,m]}$ is concentrated in degree $0$ or is $0$). As the $\eRHom$s are compatible with shifts, assume $C^\bullet = C[0]$ with $C \in \completeMod{F}$ is concentrated in degree $0$. It then remains to see that for any $i\in \Z$ the induced morphism
		$$\Ext_\completeMod{F}^i(F(G), C) \to \Ext_\category{A}^i(G, C)$$
	is an isomorphism. For $i<0$ both sides are trivial by definition. For $i>0$ both sides vanish as $F(G)$ and $G$ are projective in their respective abelian category. For $i=0$ this morphism is simply
		$$\Hom_\completeMod{F}(F(G), C) \to \Hom_\category{A}(G, C),$$
	which is an isomorphism since $C$ is $F$-complete. This shows full faithfulness.
	
	Denote by $\categoryname{D}'_F(\category A)$ the full subcategory of $\derived{\category{A}}$ of all complexes $C^\bullet$ such that for all $n\in \Z$ the cohomology group $H^nC^\bullet$ is $F$-complete.	By the results already established, $\categoryname{D}'_F(\category A)$ is a triangulated subcategory of $\derived{\category A}$ stable under all direct sums and products. Since $\derived{\completeMod{F}}$ is generated by the objects $F(G)[0]$ for $G\in\category{A}_0$ and each $F(G)[0] \in \categoryname{D}'_F(\category A)$ (as $F(G)$ is $F$-complete), we obtain that $\derived{\completeMod{F}}\subseteq \categoryname{D}'_F(\category A)$. We claim that the converse holds as well. Suppose that $C^\bullet \in \categoryname{D}'_F(\category A)$, \ie that for all $n\in \Z$ the cohomology group $H^nC^\bullet$ is $F$-complete. Similar to the previous argument, the derived $F$-completeness of $C^\bullet$ reduces to the case of a bounded $C^\bullet$ -- which in turn reduces by induction on the length of the bounded complex to the case where $C^\bullet = C[0]$ with $C\in\category{A}$ is concentrated in degree $0$. But in this case $H^0C^\bullet = H^0 C[0] \cong C$ has to be $F$-complete and the result follows. Having shown both inclusions we thus obtain that $\categoryname{D}'_F(\category A) = \derived{\completeMod{D}}$. Now the fully faithful functor $\derived{\completeMod{F}} \to \derived{\category{A}}$ clearly factors over $\derived{\completeMod{F}}\to \categoryname{D}_F(\category A)$ and this functor is an equivalence on complexes concentrated in degree $0$. In fact, the subcategory of complexes concentrated in degree $0$ is precisely $\completeMod{F}$. Now this factored functor is still fully faithful and commutes with arbitrary direct sums and products. But both categories are stable under arbitrary direct sums, direct products and shifts all while $\categoryname{D}_F(\category{A})$ is generated by (the shifts of) its complexes in degree $0$. It follows that this fully faithful functor must be an equivalence.
	
	It remains to see that the inclusion admits a left adjoint. The existence follows formally by the adjoint functor theorem of \assumptionPreWord \ref{assumption: another adjoint functor theorem}. This adjoint takes $G\in\mathcal{A}_0$ to $F(G)$ and hence must be the left derived functor of the completion $\category{A} \to \completeMod{F}$.
\end{proof}

\begin{remark}
	Compare this to lemma \ref{lemma: coherentness vs being locally of finite presentation}, where on a ringed space $(X,\struct X)$ for which the structure sheaf $\struct X$ is coherent, we obtain that being coherent is the same as being locally of finite presentation. Something similar happens here too. Once we know that $F$-free modules are (derived) complete, then completeness is the same as being $F$-presentable.
	
	In particular the condition that $F$-frees are $F$-complete has some similarity with Oka's coherence theorem. It states that the structure sheaf (\ie the most basic free sheaf) on a complex analytic space is coherent. Contrary to local Noetherian schemes and complex analytic spaces, where coherentness of the structure sheaf is automatic, $F$-completeness might not be and hence has to be enforced to obtain a good theory.
\end{remark}

\begin{lemma}[Kernels of $F$-free objects are complete]
	\label{lemma: kernels of F-free objects are complete}
	
	Suppose any complex
		$$C^\bullet = \left(\dots \to C^{-2} \to C^{-1} \to C^0 \to 0\right) \in \derived{\category A}$$
	of $F$-free modules is derived $F$-complete. Then any kernel of a morphism between $F$-presentable objects viewed as a complex concentrated in degree $0$ is already derived $F$-complete.
\end{lemma}
\begin{proof}
	Suppose that $i \colon K \hookrightarrow X$ is a kernel of the morphism $f\colon X\to Y$ between $F$-presentable objects $X$ and $Y$. We need to show that $K[0]$ is derived $F$-complete.
	
	By generation we can choose a resolution of $K$, that is choose a quasi-isomorphism $\epsilon\colon P^\bullet \to K$ where $P^\bullet = (\dots\to G^{-2} \to G^{-1} \to G^0 \to 0)$ is a complex such that each $G^i$ is of the form $\bigoplus_{j\in J_i} G_{i, j}$ where $J_i$ is some set and each $G_{i, j}\in\category A_0$ is a compact projective generator. Denote by $\tilde P^\bullet$ the induced complex with $\tilde P^i = \bigoplus_{j\in J_i} F(G_{i, j})$. Then $\tilde P^\bullet$ is a complex consisting of $F$-free objects and there is a natural map $\tilde \eta\colon P^\bullet\to \tilde P^\bullet$ induced by the natural morphisms $\eta_{G_{i, j}}\colon G_{i,j}\to F(G_{i, j})$.
	
	By assumption the complexes $X[0]$ and $Y[0]$ are $F$-complete, so for any $G_{i,j}$ we obtain isomorphisms $\eRHom(F(G_{i,j}), X)\isorightarrow \eRHom(G_{i,j}, X)$ and $\eRHom(F(G_{i,j}), Y)\isorightarrow \eRHom(G_{i,j}, Y)$. In particular, we obtain natural isomorphisms
		$$\eRHom(\tilde P^\bullet, X)\isorightarrow \eRHom(P^\bullet, X)\text{ and }\eRHom(\tilde P^\bullet, Y)\isorightarrow \eRHom(P^\bullet, Y).$$
	Indeed, by writing $P^\bullet$ and $\tilde P^\bullet$ as a sequential colimits of their truncations $\sTruncation{\ge -n}P^\bullet$ and $\sTruncation{\ge -n}\tilde P^\bullet$ we can reduce this claim to the case that both complexes are bounded by using lemma \ref{lemma: right derived Hom and homotopy colimits in the first argument}. Using the distinguished triangles $\sTruncation{\ge -(n-1)}P^\bullet \to \sTruncation{\ge -n}P^\bullet \to P^{-n}[-n] \to (\cdots)[1]$ for $P^\bullet$ and the analogous one for $\tilde P^\bullet$ the bounded case further reduces to the case that both complexes are concentrated in degree $0$. But then direct sum pulls out as products, and we finally reduce to $P^\bullet = G[0]$ and $\tilde P^\bullet = F(G)[0]$ for some $G$ where the isomorphy is clear.
	
	Consider the map $i\circ\epsilon \colon P^\bullet\to K \to X$. By the isomorphism of $\eRHom$s for $X$ this map extends to a map $\widehat{i\circ\epsilon}\colon \tilde P^\bullet \to X$ such that the diagram
	$$\begin{tikzcd}
		&P^\bullet\ar[d, "\epsilon"]\ar[r, "\tilde\eta", swap] &\tilde P^\bullet\ar[d, "\widehat{i\circ \epsilon}", dotted, swap]\\
		0\ar[r]&K\ar[r, "i"] &X \ar[r, "f"] &Y
	\end{tikzcd}$$
	commutes. Now $g\ldef f\circ \widehat{i\circ \epsilon}$ is a map $\tilde P^\bullet\to Y$ such that $g\circ \tilde\eta = f\circ i\circ \epsilon$. By the isomorphism of $\eRHom$s for $Y$ this map must then be the unique extension $\widehat{f\circ i \circ \epsilon}$. But $f\circ i\circ\epsilon = 0$ so its extension $g$ is $0$ as well. Hence, $f\circ \widehat{i\circ \epsilon} = 0$, so by the universal property of the kernel $\widehat{i\circ \epsilon}$ must factor along the kernel $i$ along some map $h\colon \tilde P^\bullet\to K$. Thus, we have a commuting triangle
	$$\begin{tikzcd}[column sep=1em]
		P^\bullet\ar[rr, "\tilde\eta"]\ar[dr, "\epsilon", swap] &&[-.3em] \tilde P^\bullet\ar[dl, "h"]\\
		&K
	\end{tikzcd}$$
	in both $\complexes{\category A}$ and $\derived{\category{A}}$. Recall that $\epsilon$ is a quasi-isomorphism, so is invertible in $\derived{\category{A}}$. Define $s\ldef \tilde\eta\circ\epsilon^{-1}$, then
	$$h\circ s = h\circ \tilde\eta\circ \epsilon^{-1} = \epsilon\circ\epsilon^{-1}=\id_K$$
	so $s$ is a section of $h$ and $K$ is a retract of $\tilde P^\bullet$ in the derived category. Thus, the two embedded squares in the diagram
	$$\begin{tikzcd}
		\eRHom(F(G), K)\ar[r]\ar[d, "s_\ast", shift left=3pt] &\eRHom(G, K)\ar[d, shift left=3pt, "s_\ast"]\\
		\eRHom(F(G), \tilde P^\bullet)\ar[r]\ar[u, "h_\ast", shift left=3pt] &\eRHom(G, \tilde P^\bullet)\ar[u, "h_\ast", shift left=3pt]
	\end{tikzcd}$$
	obtained by choosing either both $s_\ast$ or both $h_\ast$ are commutative. By functoriality we still have $h_\ast\circ s_\ast = (h\circ s)_\ast = \id$. Now since $\tilde P^\bullet$ is a complex of $F$-free objects, concentrated in non-positive degrees the bottom map is an isomorphism by assumption. But then also the top map is an isomorphism as well with the obvious inverse.
\end{proof}

We can now ask, when the assumption of \ref{lemma: THE LEMMA} is met. Indeed, this is the case if the analogous property of \ref{definition: analytic rings} (for being an analytic ring) is satisfied. This is
\begin{lemma}
	Suppose any complex
		$$C^\bullet = \left(\dots \to C^{-2} \to C^{-1} \to C^0 \to 0\right) \in \derived{\category A}$$
	of $F$-free modules is already derived $F$-complete. Then any $F$-presentable object viewed as a complex concentrated in degree $0$ is already derived $F$-complete and hence lemma \ref{lemma: THE LEMMA} applies.
\end{lemma}
\begin{proof}
	Suppose that $Q$ is $F$-representable with associated right exact sequence 
		$$X\xrightarrow{f} Y\to Q\to 0$$
	where $f$ is a morphism between $F$-free modules $X$ and $Y$. Denote by $i\colon K\to X$ the kernel of $f$, for which $K[0]$ is $F$-complete by \ref{lemma: kernels of F-free objects are complete}. We obtain the exact sequence
		$$0\to K\to X\to Y\to Q\to 0$$
	for which the two induced sequences
		$$0\to K\to X\to \image f \to 0\text{ and }0\to \image f \to Y \to Q \to 0$$
	are exact as well. Let $G\in\category{A}_0$ be arbitrary. The morphism $\eta_G\colon G\to F(G)$ induces a morphism of long exact sequences in cohomology for $\Ext^\bullet_{\category{A}}$, such that for any $n\in\N$ (even $n\in \Z$) we have the commutative diagram
	$$\begin{tikzcd}[column sep=1em, scale cd=.8]
		\Ext_{\category{A}}^n(F(G), K)\ar[r]\ar[d, "\sim", sloped, pos=0.15, swap] &\Ext_{\category{A}}^n(F(G), X)\ar[r]\ar[d, "\sim", sloped, pos=0.15, swap] &\Ext_{\category{A}}^n(F(G), \image f)\ar[r, "\delta^n"]\ar[d] &\Ext_{\category{A}}^{n+1}(F(G), K)\ar[r]\ar[d, "\sim", sloped, pos=0.15, swap] &\Ext_{\category{A}}^{n+1}(F(G), X)\ar[d, "\sim", sloped, pos=0.15, swap]\\		
		\Ext_{\category{A}}^n(G, K)\ar[r] &\Ext_{\category{A}}^n(G, X)\ar[r] &\Ext_{\category{A}}^n(G, \image f)\ar[r, "\delta^n"] &\Ext_{\category{A}}^{n+1}(G, K)\ar[r] &\Ext_{\category{A}}^{n+1}(G, X)
	\end{tikzcd}$$
	with exact rows, where all but the middle vertical morphisms are isomorphisms, since $K[0]$ and $X[0]$ are derived $F$-complete. Thus, by the $5$-lemma also the fifth vertical morphism is an isomorphism. Hence, for any degree $n\in \Z$ the morphism in cohomology induced by
		$$\eta_G^\ast\colon\eRHom(F(G), \image f) \to \eRHom(G, \image f)$$
	is an isomorphism, making $\eta_G^\ast$ a quasi-isomorphism, that is $\image(f)[0]$ is derived $F$-complete.
	
	Since we now know that $\image(f)[0]$ is derived $F$-complete, we obtain with a similar argument applied to the second short exact sequence $0\to \image f\to Y\to Q \to 0$ that $Q[0]$ is derived $F$-complete as well.
\end{proof}

\newpage
\section{\texorpdfstring{The pre-analytic rings $A_\blacksquare$ and $(A, R)_\blacksquare$}{Solid Rings}}
\label{section: the pre-analytic rings Ablack and (A,R)black}

In this section we will introduce the pre-analytic rings $A_\blacksquare$ and $(A, R)_\blacksquare$ and work towards proving their analyticity. On the way we will try to develop an intuition for $\Z_\blacksquare$-completeness. Before defining these pre-analytic rings we however need to make a technical remark that Clausen and Scholze do not address in \cite{condenseddotpdf}.
\begin{remark}[Natural decomposition of a profinite set]
	To define the various free completions we will break down an extremally disconnected set (which we recall is profinite) into finite constituents $S = \limit_i S_i$. To ensure that the resulting completions are functorial, we need to choose such decompositions functorially as well. This can be achieved for any profinite set by inspecting the proof of implication $(2)\Rightarrow (1)$ in \cite[\href{https://stacks.math.columbia.edu/tag/08ZY}{Tag 08ZY}]{stacks-project} and observing that the posited diagram is functorial in $S$. We will however simply write $S=\limit_i S_i$ and assume the mentioned construction tacitly.
\end{remark}

If $A$ is finitely generated, the definition of $A_\blacksquare$ is analogous to $\Z_\blacksquare$ from example \ref{example: examples for pre-analytic rings}, otherwise one has to approximate $A_\blacksquare$ from below, using finitely generated subrings. 
\begin{definition}[The pre-analytic ring $A_\blacksquare$ for finitely generated $A$]
	\label{definition: the pre-analytic ring Ablack for finitely generated A}
	
	Let $A$ be a discrete finitely generated $\Z$-algebra. The pre-analytic ring $A_\blacksquare$ has underlying ring $A$, free completion $S=\limit_i S_i \mapsto \limit_i \free{A}{S_i}$ and corresponding natural transformation induced by the natural maps $S_i\to \free{A}{S_i}$ for finite $S_i$.
\end{definition}

To define $A_\blacksquare$ in full generality we will need a relative version of the previous definition.
\begin{definition}[The pre-analytic ring $(A, R)_\blacksquare$ for finitely generated $R$]
	Let $\phi\colon R\to A$ be a map of discrete $\Z$-algebras, such that $R$ is finitely generated. The pre-analytic ring $(A, R)_\blacksquare$ has underlying ring $A$, free completion $\scalarExtension{\phi}\free{R_\blacksquare}{-} = A\otimes_R \free{R_\blacksquare}{-}$ and corresponding natural transformation induced by the morphisms $S\to \free{R_\blacksquare}{S}$.
\end{definition}

We can now define $A_\blacksquare$ in full generality.
\begin{definition}[The pre-analytic ring $A_\blacksquare$]
	\label{definition: the pre-analytic ring Ablack}
	
	Let $A$ be a discrete ring. For any extremally disconnected $S$ define
		$$\free{A_\blacksquare}{S}\ldef \colimit_{A'\subseteq A} \free{(A, A')_\blacksquare}{S},$$
	where the colimit is filtered and runs over all finitely generated $\Z$-subalgebras $A'$ of $A$. Now $A_\blacksquare$ is the pre-analytic ring with underlying ring $A$, free completion given by $\free{A_\blacksquare}{-}$ and corresponding natural transformation induced by any of the compatible morphisms $S\to \free{(A, A')_\blacksquare}{S}$.
\end{definition}

\begin{remark}[$A_\blacksquare$ as a colimit of pre-analytic rings.]
	In definition \ref{definition: the pre-analytic ring Ablack}, $A_\blacksquare$ can be viewed as the colimit of the pre-analytic rings $(A, A')_\blacksquare$ along the various $A'$, all with the same underlying ring $A$. We will however not really need this characterization and thus opted for the more explicit definition.
\end{remark}

\begin{remark}[Sanity check]
	Suppose that in definition \ref{definition: the pre-analytic ring Ablack} the $\Z$-algebra $A$ is a finitely generated. Then clearly $A'\ldef A$ is cofinal in the partially ordered indexing set of the colimit, ensuring that the natural map
		$$\free{A_\blacksquare^{\text{old}}} S \cong \colimit_{A' = A} \free{A'_\blacksquare} S \to \colimit_{A'\subseteq A} \free{A'_\blacksquare} S \rdef \free{A_\blacksquare} S$$
	is an isomorphism, recovering the definition \ref{definition: the pre-analytic ring Ablack for finitely generated A} of the finitely generated case.
\end{remark}

We can now reintroduce the relative ring $(A, R)_\blacksquare$ for a general base $R$.
\begin{definition}[The pre-analytic ring $(A, R)_\blacksquare$]
	Let $\phi\colon R\to A$ be a morphism of discrete rings. Then $(A, R)_\blacksquare$ is the pre-analytic ring with underlying ring $A$, with free completion $\scalarExtension{\phi}\free{R_\blacksquare}{-} = A\otimes_R \free{R_\blacksquare}{-}$ and corresponding natural transformation induced by the maps $S\to \free{R_\blacksquare}{S}$.
\end{definition}

\begin{remark}[Free completions for profinite sets]
	\label{remark: free completions for profinite sets}
	
	Observe, that for all the pre-analytic rings introduced in this section, the free completions can readily be generalized to an arbitrary profinite set $S$. Furthermore, as a consequence of \cite[Prop. 5.6]{condenseddotpdf} we observe that this generalized definition agrees with the construction mentioned in remark \ref{remark: when is the tensor product the derived tensor product?}. In particular, the complex constructed there is always concentrated in degree $0$ and hence $\Dotimes_\preAnaRing{A}$ really is the left derived functor of $\otimes_\preAnaRing{A}$ for call cases of $\preAnaRing{A}$ relevant to us -- of course only once we know the rings $A_\blacksquare$ and $(A, R)_\blacksquare$ to be analytic.
\end{remark}

\begin{remark}[Solidness]
	These (pre-)analytic rings $A_\blacksquare$ play an essential role in much of the algebraic theory of condensed mathematics, so much so, that $A_\blacksquare$-complete modules deserve a special name. A condensed $A$-module is called \emph{$A$-solid} if it is $A_\blacksquare$-complete. A $\Z$-solid condensed abelian group is simply called \emph{solid}. We will use this terminology occasionally as well.
\end{remark}

Let us now try to work towards an understanding of what \emph{being solid} means for a module. Recall that \emph{solidness} is characterized by existence of certain linear maps from the condensed modules $\free{A_\blacksquare}{S}$ with $S$ extremally disconnected. It is thus a good idea to develop an understanding of these $A_\blacksquare$-free modules themselves. We will try to do so for even for a general profinite set $S$.
\begin{definition}[The abelian group of continuous maps]
	For any profinite set $S$ and any discrete abelian group $M$ denote by $\continuous{S}{M}$ the discrete topological abelian group given by pointwise addition.
\end{definition}

We can now prove a first result regarding the structure of $A_\blacksquare$-free abelian groups.
\begin{lemma}[Finitely-free $A_\blacksquare$-complete modules as $A$-duals]
	\label{lemma: finitely-free A-solid modules as A-duals}
	
	Let $A$ be a discrete finitely generated $\Z$-algebra. For any profinite $S$ there is a natural isomorphism
		$$\free{A_\blacksquare}{S} \cong \associated{\intHom_\Ab(\continuous S \Z, A)}\cong \intHom_\condensedAb(\associated{\continuous{S}{\Z}}, \associated{A})$$
	of $A$-modules.
\end{lemma}
\begin{proof}
	Write $S\cong \limit_i S_i$ where each $S_i$ is finite and discrete. Then $\continuous{S_i}{\Z}\cong \bigoplus_{S_i} \Z$ is free of finite rank and hence there is a natural isomorphism
		$$\free {A} {S_i}\cong \intHom_\Ab(\continuous{S_i}{\Z}, A).$$
	As $\free{A_\blacksquare}{S} \ldef \limit_i \free{A}{S_i}$, we in particular obtain
		$$\free{A_\blacksquare}{S}\cong \limit_i \associated{\intHom_\Ab(\continuous{S_i}{\Z}, A)}\cong \associated{\intHom_\Ab(\continuous{S}{\Z}, A)}.$$
	This shows the first isomorphism. As both $\continuous S \Z$ and $A$ are discrete, additive maps $\continuous S \Z\to A$ as abelian groups and as topological abelian groups agree, so the second isomorphism follows as well.
\end{proof}

Thus, these $A_\blacksquare$-free modules are $A$-duals of abelian groups of continuous maps, at least for finitely generated $A$. One can interpret elements of this dual space as a kind of measure.
\begin{definition}[Pseudo-measures]
	Let $S$ be a set and $\mathcal{G}$ be a $\cup$-stable, $\cap$-stable and complement-stable system of subsets of $S$ such that $\emptyset\in \mathcal{G}$. Let $A$ be a discrete abelian group. An \emph{$A$-valued pseudo-measure} on $S$ is a function
		$$\mu\colon \mathcal{G} \to A$$
	such that $\mu(\emptyset)=0$ and $\mu(U\cup V) = \mu(U) + \mu(V)$ for any two disjoint $U, V\in\mathcal{G}$. Tacitly assuming the system of sets to be fixed, we denote by $\measures{S}{A}$ the set of pseudo-measures on $S$.
\end{definition}

If $S$ is profinite, one such system of subsets is given by the clopen subsets of $S$. Indeed, $\emptyset$ is clopen, the union and intersection of clopen subsets is clopen and certainly the complement of a clopen set is clopen as well. Furthermore, by identifying a clopen subset with its indicator function we obtain the following characterization.
\begin{lemma}
	\label{lemma: indicator functions generate}
	
	Let $S$ be a profinite set. Then
		$$\continuous{S}{\Z} = \langle 1_U\colon S\to \Z\setseparator U\subseteq S\text{ clopen}\rangle_\Z$$
	that is any continuous function $S\to \Z$ is a finite $\Z$-linear combination of indicator functions of clopen subsets of $S$.
\end{lemma}
\begin{proof}
	If $U\subseteq S$ is clopen then $1_U\colon S\to \Z$ is continuous. Indeed, if $V$ is a(n open) subset of $\Z$ then its preimage is either $\emptyset$, the clopen $U$, its complement $S\setminus U$ or the whole space $S$ -- all of which are open in $S$. Certainly, also any finite $\Z$-linear combination of such indicator functions is continuous as well.
	
	For the converse assume that $f\colon S\to \Z$ is continuous. As $S$ is compact so is its image $f(S)$ and since $\Z$ is discrete so is its subspace $f(S)$. But a discrete space is compact if and only if it is finite hence $f$ takes only finitely many values. Let $n\in f(S)$. Then $\{n\}$ is clopen in $\Z$ (by discreteness of $\Z$) and hence so is the fiber $f^{-1}(\{n\})$. It follows that
		$$f = \sum_{n\in f(S)} n \cdot 1_{f^{-1}(\{n\})}$$
	is such a $\Z$-linear combination as claimed.
\end{proof}

We can thus characterize duals of such groups $\continuous{S}{\Z}$.
\begin{lemma}[Finitely-free $A_\blacksquare$-modules as spaces of measures]
	Let $A$ be a discrete finitely generated $\Z$-algebra. For any profinite set $S$ there is a natural isomorphism of $A$-modules
		$$\intHom_\Ab(\continuous S \Z, A) \cong \measures S A$$
	between the $A$-dual of $\continuous S \Z$ and the $A$-valued pseudo-measures on the system of clopen subsets of $S$. In particular, lemma \ref{lemma: finitely-free A-solid modules as A-duals} yields that
		$$\free{A_\blacksquare}{S} \cong \associated{\measures S A}$$
	is the condensed set associated to the space of $A$-valued pseudo-measures on $S$.
\end{lemma}
\begin{proof}
	Recall that by the previous lemma \ref{lemma: indicator functions generate} the abelian group $\continuous{S}{\Z}$ is generated by indicator functions $1_U$ of clopen subsets $U$ of $S$. As stated above, the result then follows from identifying these clopen subsets $U$ of $S$ with their indicator functions $1_U$. More precisely, map an $f\colon \continuous S \Z \to A$ to the pseudo-measure $\mu_f\colon U\mapsto f(1_U)$ and map every pseudo-measure $\mu$ to the additive function $f_\mu$ such that $1_U\mapsto \mu(U)$ for every clopen subset $U$ of $S$. Clearly $\mu_{-}$ and $f_{-}$ are additive. Furthermore, they satisfy
		$$\mu_{f_\mu} = \mu\text{ and }f_{\mu_f}=f$$
	so are inverse of each other.
\end{proof}

With this interpretation of $\Z_\blacksquare$-free abelian groups at hand, we can try to understand solidness.
\begin{remark}[Interpreting $\Z_\blacksquare$-completeness]
	\label{remark: interpreting Zblack-completeness}
	
	Suppose that the $\Z$-module $\associated{M}$ is $\Z_\blacksquare$-complete. If $f\colon S\to M$ is a map from a profinite set, then we obtain by completeness a commuting triangle
	$$\begin{tikzcd}[column sep=1pt]
		S\ar[rr, "f"]\ar[dr] &&M\\
		&\measures{S}{\Z}\ar[ur, "\hat f", dotted, swap]
	\end{tikzcd}$$
	Now if $\mu\in \measures S \Z$ is any pseudo-measure, then we can define $\int f\d \mu\ldef \hat f(\mu)$. Since the map $S\to\measures{S}{\Z}$ corresponds to mapping any $s \in S$ to its corresponding Dirac-pseudo-measure $\delta_s \in \measures{S}{\Z}$, this gives the familiar identity $\int f \d\delta_s = \hat f(\delta_s) = f(s)$. Furthermore, since $f\mapsto \hat f$ is $\Z$ linear we obtain further that $\int (\alpha f + \beta g) \d\mu = \alpha \int f\d\mu + \beta \int g\d\mu$ for all $f, g\colon S\to M$ and integers $\alpha, \beta$. We find that this \emph{integral} really has the expected properties. Thus, $\Z_\blacksquare$-completeness of $M$ can be characterized as saying that \emph{any continuous function $S\to M$ is integrable against any $\Z$-valued pseudo-measure on $M$}.
	
	
	Let us investigate this integrability in the special case of the profinite set $S = \hat\N=\N\cup\{\infty\}$, that is the one-point-compactification of $\N$. Then a continuous map $f\colon \hat \N\to M$ corresponds precisely to a convergent sequence $(a_n)_{n\in \N}$ in $M$ with limit $a_\infty = f(\infty)$. Denoting by $\continuousZero{\N}{M}$ the null-sequences in $M$, we obtain a decomposition $\continuous {\hat \N} M \cong \continuousZero {\N} M \oplus M$ of abelian groups by mapping $f\mapsto (f\vert_\N-f(\infty), f(\infty))$. This induces a decomposition of pseudo-measures $\measures{\hat\N}{\Z} \cong \measures{\N}{\Z}\oplus \Z$, by means of $\mu\mapsto (\mu\vert_\N, \mu(\infty))$. As $\N$ is a discrete space, any subset is clopen hence the system of clopen subsets of $\N$ agrees with the usual $\sigma$-algebra (both are simply the powerset of $\N$). It is clear that any measure on $\N$ is a pseudo-measure on $\N$ (but not conversely as $\sigma$-additivity is not clear). Now, using the linearity of the integral we see that for any $f\in \continuous{S}{\Z}$ and any $\mu \in \measures{S}{\Z}$ we have
		$$\int_{\hat \N} f \d\mu = \int_{\hat \N} f-f(\infty) \d\mu + f(\infty)\mu(\infty).$$
	In particular if $f(\infty) = 0$ then there is an element
		$$\sum_{n\in \N} f(n) \mu(n) \ldef \int_{\hat\N} f\d\mu \in M,$$
	that is \emph{any null-sequence in $M$ is summable against any $\Z$-valued pseudo-measure $\mu$ on $\hat \N$}. By choosing a pseudo-measure $\mu$ with $\mu(n)=1$ for all $n\in \N$ we see that for any null-sequence $f$ in $M$ the sum $\sum_{n\in \N} f(n)$ ought to exist in $M$ -- it follows that solidness is a \emph{non-archimedean} notion. Indeed, in an archimedean setting, this is surely false. \Eg in $\R$ one has the standard example of the (partial sums of the) harmonic series which diverge.
\end{remark}

A rather important result on the structure of $\continuous S \Z$, due to Nöbling in \cite{noebling}, based on results of Specker in \cite{specker}, is given in the following theorem. We will only sketch the proof (using an argument attributed to Bergman in \cite[Thm. 97.2]{fuchs-infinite-abelian-groups}) -- a formally verified proof can be found in \cite{asgeirsson}.
\begin{theorem}[Specker's theorem]
	\label{theorem: specker's theorem}
	
	For any profinite set $S$, the discrete abelian group $\continuous{S}{\Z}$ of continuous integer valued functions on $S$ is free with rank at most $2^{\vert S\vert}$.
	
	In particular, there is some set $I$ with $\vert I\vert\le 2^{\vert S\vert}$ such that for any discrete finitely generated $\Z$-algebra we have an isomorphism of $A$-modules
		$$\free{A_\blacksquare}{S} \cong \prod_I \associated{A}.$$
\end{theorem}
\begin{proofsketch}
	Pick a continuous injection $S\to \prod_\lambda \{0, 1\}$ onto a closed subset, where $\lambda$ is suitably large and assumed an ordinal (here one can take $\lambda$ to be equicardinal to the set of clopen subsets of $S$ such that for each clopen subset the corresponding projection to $\{0, 1\}$ is the indicator function on the subset). In particular, the elements $\mu \in \lambda$ correspond to ordinals $\mu < \lambda$. For such a $\mu < \lambda$ consider the idempotent $e_\mu \in \continuous S \Z$ given by $S\to\prod_\lambda \{0, 1\}\to \{0, 1\}\subseteq \Z$ induced by the projection onto the $\mu$-th factor. For $r\in \N_0$, order the products $e_{\mu_1}\cdot (\cdots) \cdot e_{\mu_r}$ with $\mu_1 > \dots > \mu_r$ lexicographically. Denote by $E$ the set of such products that cannot be written as a linear combination of smaller such products. We claim that $E$ is a basis of $\continuous{S}{\Z}$. We will proceed by induction on the ordinal $\lambda$. Denote by $\tilde S$ the image of $S$ and for any $\mu <\lambda$ denote by $\tilde S_\mu$ the image of $S$ in $\prod_{\mu'<\mu} \{0, 1\}$ after projecting from $\prod_\lambda \{0, 1\}$. 
	
	If $\lambda = 0$ then $S=\{\ast\}$ and the claim is trivial.
	
	If $\lambda$ is a limit ordinal, \ie $\lambda = \bigcup_{\mu < \lambda} \mu$, then $\tilde S \cong \colimit_{\mu < \lambda} \tilde S_\mu$ is the union of the $\tilde S_\mu$. Denoting by $E_\mu$ the analogous subsets of $\continuous{\tilde S_\mu}{\Z}$ -- which are bases by induction, we observe that $E \cong \colimit_{\mu < \lambda} E_\mu$ is a basis for $\continuous{\tilde S}{\Z}$.
	
	Finally assume that $\lambda = \rho + 1$ is a successor ordinal. Then the above map identifies with the map $S \hookrightarrow \tilde S_\rho \times \{0, 1\}$. Denote by $\tilde S_0$ and $\tilde S_1$ the intersections $\tilde S\cap (\tilde S_\rho \times \{0\})$ and $\tilde S\cap (\tilde S_\rho \times \{1\})$ respectively. As intersections of closed subsets, these are closed. Furthermore, we have $\tilde S = \tilde S_0 \cup \tilde S_1$. Denote by $\tilde S'$ the intersection of the images of $\tilde S_0$ and $\tilde S_1$ under the projection $\tilde S_\rho\times\{0, 1\}\to \tilde S_\rho$. We thus have the short exact sequence
		$$0\to \continuous{\tilde S_\rho}\Z\to \continuous {\tilde S} \Z \to \continuous{\tilde S'}{\Z}\to 0$$
	where the left morphism maps $f\colon \tilde S_\rho\to\Z$ to the continuous function $(s, b)\mapsto f(s)$, and the right morphism maps $g\colon \tilde S\to \Z$ to the difference of restrictions $s\mapsto f(s,0) - f(s, 1)$. Now applying the induction hypothesis for both $\tilde S_\rho$ and $\tilde S'$ (mapping to $\prod_{\mu<\rho} \{0, 1\}$) we see that the basis vectors of $E$ that do not start with $e_\rho$ form a basis of $\continuous{\tilde S_\rho}{\Z}$, while the basis vectors that do start with $e_\rho$ project to a basis of $\continuous{\tilde S'}{\Z}$. So $E$ really is a basis by exactness of the above sequence.
	
	For the second claim follows by choosing an isomorphism $\continuous{S}{\Z}\cong \bigoplus_I \Z$ of abelian groups. Then
		$$\free{A_\blacksquare}S \cong \associated{\intHom_\Ab(\continuous{S}{\Z}, A)} \cong \associated{\intHom_\Ab(\bigoplus\nolimits_I \Z, A)}=\prod_I \associated{A},$$
	as $\intHom_\Ab(-, A)$ maps direct sums to products and as $X\mapsto\associated{X}$ commutes with them.
\end{proofsketch}

To get an intuition on how $A_\blacksquare$-completed tensor products behave, we state the following proposition but simply reference its proof.
\begin{proposition}
	\label{proposition: tensor products of compacts, f.g. case}
	
	Let $A$ be a discrete finitely generated $\Z$-algebra. If $M=\prod_I \associated{A}$ and $N=\prod_J \associated{A}$, then
	$$M\Dotimes_{A_\blacksquare} N \cong \prod_{I\times J} \associated{A}.$$
\end{proposition}
\begin{shortproof}
	In case of $A=\Z$ this is \cite[Prop. 6.3]{condenseddotpdf}. The other case is analogous.
\end{shortproof}

A sample consequence of this proposition is given in
\begin{example}[Some completed tensor products]
	One can calculate the following tensor products directly:
	\begin{itemize}
		\item $\Z_p \Dotimes_{\Z_\blacksquare} \R \cong 0$
		\item $\Z_p \Dotimes_{\Z_\blacksquare} \Z_\ell \cong 0$
		\item $\Z_p \Dotimes_{\Z_\blacksquare} \Z_p \cong \Z_p$
		\item $\Z_p \Dotimes_{\Z_\blacksquare} \series{\Z}{T} \cong \series{\Z_p}{T}$
		\item $\series{\Z}{U} \Dotimes_{\Z_\blacksquare} \series{\Z}{T} \cong \series{\Z}{U, T}$
	\end{itemize}
	for primes $p\ne \ell$. The precise reasoning is given in \cite[Example 6.4]{condenseddotpdf}. What one should observe is that the (derived) completed tensor product $\Dotimes_{\Z_\blacksquare}$ \quote{seems to ask both factors for what $I$-adic topology they have, and then combines} them. In particular solidness and more generally $A$-solidness for some ring $A$ are non-archimedean in nature (see also remark \ref{remark: interpreting Zblack-completeness}). Clausen and Scholze however state, that solidness is very well suited for applications in $p$-adic functional analysis. If one wishes to obtain a good condensed foundation for real or complex functional analysis, one has to put in significantly more work. This is the content of the (now completed, but very subtle) \emph{Liquid Tensor Experiment}, see \href{https://xenaproject.wordpress.com/2020/12/05/liquid-tensor-experiment/}{for a short introduction}. 
\end{example}

This allows us to calculate tensor products even in the general case. Recall that in case of an $A$ that is not finitely generated we have $\free{A_\blacksquare}{S} = \colimit_{A'\subseteq A} A\otimes_{A'}\free{A'_\blacksquare}{S}$ for any extremally disconnected $S$. With $I$ as in Specker's theorem \ref{theorem: specker's theorem}, we thus obtain that $\free{A_\blacksquare}{S} \cong \colimit_{A'\subseteq A} A\otimes_{A'} \prod_I A'$. First, observe that all objects of this form are already compact projective.
\begin{lemma}[More compact projectives]
	\label{lemma: more compact projectives}
	
	Let $A$ be a discrete ring such that $A_\blacksquare$ is analytic. Then for any set $I$ the $A$-module
		$$\colimit_{A'\subseteq A} A\otimes_{A'} \prod_I A'$$
	is $A_\blacksquare$-complete and a compact projective object of $\mod{A_\blacksquare}$.
\end{lemma}
\begin{proof}
	Choose a large enough extremally disconnected $S$ such that the set $J$ associated to it by Specker's theorem \ref{theorem: specker's theorem} satisfies $\vert I\vert \le \vert J \vert$. Choosing an injection $I \hookrightarrow J$ we find that for all $A'$ the product $\prod_{I} A'$ is a retract of the product $\prod_J A'$ where the retraction is natural. Hence, the $A$-module in question is a retract of
		$$\colimit_{A'\subseteq A} A\otimes_{A'} \prod_J A'\cong \free{A_\blacksquare}{S}$$
	which is compact projective. As being $A_\blacksquare$-complete and being compact projective are stable under retracts it follows that $\colimit_{A'\subseteq A} A\otimes_{A'} \prod_I A'$ is $A_\blacksquare$-complete and compact projective as well.
\end{proof}

\begin{corollary}
	\label{corollary: tensor products of compacts, general case}
		
	Let $A$ be a discrete ring such that $A_\blacksquare$ is analytic. If $M = \colimit_{A'\subseteq A} A\otimes_{A'} \prod_I A'$ and $M = \colimit_{A'\subseteq A} A\otimes_{A'} \prod_j A'$ for two sets $I$ and $J$, then
		$$M\Dotimes_{A_\blacksquare} N \cong M\otimes_{A_\blacksquare} N\cong \colimit_{A'\subseteq A} A\otimes_{A'} \prod_{I\times J} A'.$$
	In particular, by the previous lemma \ref{lemma: more compact projectives} the tensor product $M\otimes_{A_\blacksquare} N$ is compact projective in $\mod{A_\blacksquare}$. Furthermore, $M\Dotimes_{A_\blacksquare} N$ is derived compact in $\derived{A_\blacksquare}$.
\end{corollary}
\begin{proof}
	By this previous lemma, we also know $M$ and $N$ to be projective. We immediately obtain that
		$$M\Dotimes_{A_\blacksquare} N \cong M\otimes_{A_\blacksquare} N,$$
	and we do not need to work in the derived setting. Since $A_\blacksquare$ is analytic, we know by theorem \ref{theorem: properties of the category of complete modules} that $A_\blacksquare$-completion exists and is strong monoidal, yielding that
		$$M\otimes_{A_\blacksquare} N \cong \completion{\mleft(\colimit_{A'} A\otimes_{A'} \prod_I A'\mright)\otimes_A \mleft(\colimit_{A''} A\otimes_{A''} \prod_J A''\mright)}{A_\blacksquare}.$$
	Since tensor products commute with arbitrary colimits we find that
	\begin{align*}
		M\otimes_{A_\blacksquare} N &\cong \completion{\colimit_{A'}\colimit_{A''} \mleft(A\otimes_{A'} \prod_I A'\mright) \otimes_A \mleft(A\otimes_{A''} \prod_J A''\mright)}{A_\blacksquare}\\
		&\cong \completion{\colimit_{A'}\colimit_{A''} \mleft(\prod_I A'\mright)\otimes_{A'} A \otimes_A A\otimes_{A''} \mleft(\prod_J A''\mright)}{A_\blacksquare}\\
		&\cong \completion{\colimit_{A'}\colimit_{A''} \mleft(\prod_I A'\mright)\otimes_{A'} A \otimes_{A''} \mleft(\prod_J A''\mright)}{A_\blacksquare}.
	\end{align*}
	By injecting the pair $(A', A'')$ into $(A''', A''')$ where $A''' \ldef \Z[A'\cup A'']$ is $\Z$-finitely generated, we find that the diagonal $A'=A''$ is cofinal and hence that
	\begin{align*}
		M\otimes_{A_\blacksquare} N &\cong \completion{\colimit_{A'} \mleft(\prod_I A'\mright)\otimes_{A'} A \otimes_{A'} \mleft(\prod_J A'\mright)}{A_\blacksquare}\\
		&\cong \completion{\colimit_{A'} A\otimes_{A'}\mleft(\prod_I A'\mright)\otimes_{A'}\mleft(\prod_J A'\mright)}{A_\blacksquare}.
	\end{align*}
	Now $\completion{-}{A_\blacksquare}$ is a left adjoint so commutes with arbitrary colimits. Furthermore, if we denote by $\phi\colon A'\to A$ the inclusion then $\completion{-}{A_\blacksquare} \circ \scalarExtension{\phi} \cong \completedScalarExtension{\phi}{A_\blacksquare} \circ \completion{-}{A'_\blacksquare}$ as follows from the proof of proposition \ref{proposition: scalar restriction and extension of complete modules}, \ie $A_\blacksquare \otimes_{A} A\otimes_{A'} - \cong A_\blacksquare\otimes_{A'} - \cong A_\blacksquare\otimes_{A'_\blacksquare} A'_\blacksquare\otimes_{A'} -$. Thus,
	\begin{align*}
		M\otimes_{A_\blacksquare} N &\cong \colimit_{A'} A_\blacksquare\otimes_{A'_\blacksquare}\completion{\mleft(\prod_I A'\mright)\otimes_{A'}\mleft(\prod_J A'\mright)}{A'_\blacksquare}\\
		 &\cong \colimit_{A'} A_\blacksquare\otimes_{A'_\blacksquare}\mleft(\prod_I A'\mright)\otimes_{A'_\blacksquare}\mleft(\prod_J A'\mright).
	\end{align*}
	By proposition \ref{proposition: tensor products of compacts, f.g. case} we thus have
		$$M\otimes_{A_\blacksquare} N \cong \colimit_{A'} A_\blacksquare\otimes_{A'_\blacksquare} \prod_{I\times J} A' \cong \completion{\colimit_{A'} A\otimes_{A'}\prod_{I\times J} A'}{A_\blacksquare}$$
	which is simply
		$$\colimit_{A'} A\otimes_{A'} \prod_{I\times J} A'$$
	as this $A$-module is already $A_\blacksquare$-complete by the previous lemma.
\end{proof}

These previous statements assert that the tensor product of compact objects is still compact. We obtain the following result as a consequence that occasionally proves useful.
\begin{lemma}[Solid derived compacts are enriched compact]
	\label{lemma: solid derived compacts are enriched compact}
	
	Let $A$ be a discrete ring such that $A_\blacksquare$ is analytic. Let $K$ be a derived compact object of $\derived{A_\blacksquare}$. Then $K$ is $\condensedAb$-enriched compact.
	
	Similarly, if $R\to A$ is a morphism of discrete rings such that $(A, R)_\blacksquare$ is analytic, then any derived compact of $\derived{(A, R)_\blacksquare}$ is $\condensedAb$-enriched compact.
\end{lemma}
\begin{proofsketch}
	By theorem \ref{theorem: properties of the category of complete modules} we know that if $Y$ is $A_\blacksquare$-complete and $S$ is extremally disconnected then $\eRHom_A(\free{A_\blacksquare}{S}, Y) \cong \eRHom_A(\free{A}{S}, Y)$ in $\derived{\condensedAb}$. This implies that also $\intRHom_A(\free{A_\blacksquare}{S}, Y) \cong \intRHom_A(\free{A}{S}, Y)$ as $\eRHom_A$ is the underlying condensed abelian group of $\intRHom_A$. By applying $H^0$ we thus obtain that $\intHom_A(\free{A_\blacksquare}{S}, Y) \cong \intHom_A(\free{A}{S})$ in $\mod{A}$. Using this, we can show that $\mod{A_\blacksquare}$ is actually closed monoidal. Indeed, we find for all $X, Y\in\mod{A_\blacksquare}$ the $A$-module $\intHom_A(X, Y)$ to be $A_\blacksquare$-complete. For this let $S$ be extremally disconnected. Then
	\begin{align*}
		\Hom(\free{A_\blacksquare}{S}, \intHom_A(X, Y)) &\cong \Hom(X, \intHom_R(\free{A_\blacksquare}{S}, Y))\\
		& \cong \Hom(X, \intHom_A(\free{A}{S}, Y)) \cong \Hom(\free{A}{S}, \intHom_A(X, Y))
	\end{align*}
	showing the desired isomorphy. One can then check that this really is a partial right adjoint to $\otimes_{A_\blacksquare}$.
	
	Now knowing that $\mod{A_\blacksquare}$ is closed monoidal we find by corollary \ref{corollary: a sufficient condition for enriched compactness} that if $-\Dotimes_{A_\blacksquare} K$ preserves compactness then $K$ is $\mod{A_\blacksquare}$-enriched compact. But the tensor products of the derived compact generators of $\derived{A_\blacksquare}$ are again derived compact by corollary \ref{corollary: tensor products of compacts, general case} so this follows. This shows $\mod{A_\blacksquare}$-enriched compactness of $K$. Now as the $\condensedAb$-enrichment is simply given by the underlying condensed abelian groups of the $\mod{A_\blacksquare}$-enrichment it follows that $K$ is also already $\condensedAb$-enriched compact (since direct sums in all categories in question agree).
	
	The case of $(A, R)_\blacksquare$ is analogous, after one observes that the relevant claim about tensor products of generators being compact again translates.
\end{proofsketch}

Having now investigated the rings $(A, R)_\blacksquare$ and $A_\blacksquare$ on their own, we should then set them in relation to each other. The following lemma provides various natural maps between them.
\begin{lemma}[Some maps of pre-analytic rings]
	\label{lemma: some maps of pre-analytic rings}
	
	We have the following maps of pre-analytic rings
	\begin{enumerate}[(i)]
		\item
		For any discrete ring $A$ a map $A\to A_\blacksquare$ with underlying morphism of rings $\id_A$.
		
		\item
		\label{lemma: some maps of pre-analytic rings, ii}
		For any morphism $\phi\colon R\to A$ of discrete rings a map $R_\blacksquare\to A_\blacksquare$ with morphism of underlying rings $\phi$.
		
		\item
		For any morphism $R\to A$ of discrete rings a map $(A, R)_\blacksquare\to A_\blacksquare$ with morphism of underlying rings $\id_A$.
		
		\item
		For any commutative square
		$$\begin{tikzcd}
			R\ar[r]\ar[d] &A\ar[d, "\phi"]\\
			R'\ar[r] & A'
		\end{tikzcd}$$
		of discrete rings a map $(A, R)_\blacksquare \to (A', R')_\blacksquare$ with morphism of underlying rings $\phi$.
	\end{enumerate}
\end{lemma}
\begin{proof}
	We will only show the existence of the desired natural morphisms between (finitely) free modules and omit the compatibility with the natural transformations of underlying condensed sets.
	\begin{enumerate}[(i)]
		\item
		Suppose that the $\Z$-algebra $A$ is finitely generated and that $S \cong \limit_i S_i$. Then the projections $S\to S_i$ induce natural morphisms $\free{A}{S}\to \free{A}{S_i}$ which in turn induce a natural morphism $\free{A}{S} \to \limit_i \free{A}{S_i} \rdef \free{A_\blacksquare}{S}$, giving the desired $A\to A_\blacksquare$.
		
		Now suppose that $A$ is general. If $A'$ is any finitely generated $\Z$-subalgebra of $A$ then by the previous argument we have a morphism $A'\to A'_\blacksquare$. These induce $A$-linear maps
		$$\free{A}{S} \cong \colimit_{A'\subseteq A} \free{A'}{S}\cong \colimit_{A'\subseteq A} A\otimes_{A'} \free{A'}{S} \to \colimit_{A'\subseteq A} A\otimes_{A'}\free{A'_\blacksquare}{S}$$
		natural in the extremally disconnected $S$, giving the desired $A\to A_\blacksquare$.
		
		\item
		Suppose that the $\Z$-algebras $R$ and $A$ are finitely generated. If $S\cong \limit_i S_i$ then the natural maps $\free{R}{S_i} \to \free{A}{S_i}$ induce a natural $R$-linear map $\free{R_\blacksquare}{S} \ldef \limit_i \free{R}{S_i} \to \limit_i \free{A}{S_i} \rdef \free{A_\blacksquare}{S}$, giving the desired $R_\blacksquare\to A_\blacksquare$.
		
		Now suppose that $A$ is general, but $R$ is still finitely generated. Then the poset of $\Z$-finitely generated $R$-subalgebras $A''$ of $A$ is cofinal in the poset of all finitely generated $\Z$-subalgebras of $A$, as any $\Z$-subalgebra $A'$ can be embedded in the still $\Z$-finitely generated $R$-algebra $A''=\Z[A', R]$. Thus, the natural morphism of colimits
		$$\colimit_{R\to A''\subseteq A} A\otimes_{A''}\free{A''_\blacksquare}{S} \to \colimit_{A'\subseteq A} A\otimes_{A'}\free{A'_\blacksquare}{S}\rdef \free{A_\blacksquare}{S}$$
		is an isomorphism. Then by the previous paragraph we have compatible $R$-linear maps $\free{R_\blacksquare}{S} \to \free{A''_\blacksquare}{S} \to A\otimes_{A''} \free{A''_\blacksquare}{S}$ that then induce a natural $R$-linear map into the first, hence into the second colimit which is $\free{A_\blacksquare}{S}$, giving the desired $R_\blacksquare \to A_\blacksquare$.
		
		Now if both $R$ and $A$ are general, then any finitely generated $\Z$-subalgebra $R'$ gives rise to a map $R'_\blacksquare \to A_\blacksquare$ by the previous paragraph. Hence we obtain an $R$-linear map
		$$\free{R_\blacksquare}{S} \ldef \colimit_{R'\subseteq R} R\otimes_{R'}\free{R'_\blacksquare}{S} \to \free{A_\blacksquare}{S}$$
		natural in the extremally disconnected $S$ from the compatible morphisms $R\otimes_{R'} \free{R'_\blacksquare}{S} \to R\otimes_{R'} \free{A_\blacksquare}{S} \to \free{A_\blacksquare}{S}$ where the second map is $R$-multiplication on the $R$-module $\free{A_\blacksquare}{S}$. This gives the desired $R_\blacksquare \to A_\blacksquare$.
		
		\item
		By \ref{lemma: some maps of pre-analytic rings, ii} we have a map $R_\blacksquare\to A_\blacksquare$. For any extremally disconnected $S$ we thus obtain the $A$-linear map
		$$\free{(A, R)_\blacksquare}{S}\ldef A\otimes_R \free{R_\blacksquare}{S} \to A\otimes_R \free{A_\blacksquare}{S} \to \free{A_\blacksquare}{S}$$
		natural in $S$.
		
		\item 
		By \ref{lemma: some maps of pre-analytic rings, ii} there is a map $R_\blacksquare \to R'_\blacksquare$. Denote its corresponding natural transformation $\free{R_\blacksquare}{-}\Rightarrow \free{R'_\blacksquare}{-}$ by $\eta$. Observe that the morphism $R\to R'$ induces a natural transformation $\otimes_R \Rightarrow \otimes_{R'}$. Thus we obtain $A$-linear maps
		$$\free{(A, R)_\blacksquare}{S} \ldef A\otimes_R \free{R_\blacksquare}{S} \xrightarrow{\phi\otimes \eta_S} A'\otimes_R \free{R'_\blacksquare}{S} \to A'\otimes_{R'} \free{R'_\blacksquare}{S} \rdef \free{(A', R')}{S}$$
		natural in the extremally disconnected $S$, giving the desired $(A, R)_\blacksquare \to (A', R')_\blacksquare$.
	\end{enumerate}
\end{proof}

\newpage
\section{Analyticity of \texorpdfstring{$A_\blacksquare$ and $(A, R)_\blacksquare$}{Solid Rings}}
\label{section: analyticity of Ablack and (A, R)black}

The goal of this chapter is to prove that the pre-analytic rings $A_\blacksquare$ and $(A, R)_\blacksquare$ are analytic. The case of $(A, R)_\blacksquare$ follows essentially formally from knowledge of analyticity of $R_\blacksquare$ by lemma \ref{lemma: partial analyticity of (A, R)black}. Thus the main difficulty lies in $A_\blacksquare$ and we will obtain the result by tackling the analyticity of $A_\blacksquare$ for increasingly more general $A$, the most basic instance being $\Z_\blacksquare$.
\begin{theorem}[Analyticity of $\Z_\blacksquare$]
	\label{theorem: analyticity of Zblack}
	
	The pre-analytic ring $\Z_\blacksquare$ is analytic.
\end{theorem}
While the claim of theorem \ref{theorem: analyticity of Zblack} is very interesting and fundamental to the rest of this thesis (and hence deserves a full proof), the argument given in \cite{condenseddotpdf} by Clausen and Scholze is a rather long and technical calculation. Hence, it will simply be referenced.
\begin{shortproof}[of theorem \ref{theorem: analyticity of Zblack}]
	This is shown in \cite[Lecture VI, pp. 38--41]{condenseddotpdf}
\end{shortproof}

The next level of generality will be analogous to Hilbert's basis theorem.
\begin{theorem}[A kind of Hilbert's basis theorem]
	\label{theorem: a kind of hilbert's basis theorem}
	
	If $R$ is a discrete finitely generated $\Z$-algebra such that $R_\blacksquare$ is analytic, then the pre-analytic ring $R[T]_\blacksquare$ is analytic as well.
\end{theorem}
The proof of this theorem uses many intermediate results which will be needed again at a later time, hence we will devote to this proof its own section \ref{section: an analogue of hilbert's basis theorem}. Assuming the statement proven for now, we immediately obtain the following corollary.
\begin{corollary}[Analyticity of $A_\blacksquare$ for finitely generated free $\Z$-algebras]
	\label{corollary: analyticity of Ablack for finitely generated free Z-algebras}
	
	The pre-analytic ring $\Z[T_1,\dots, T_n]_\blacksquare$ associated to the polynomial ring in $n$ variables over $\Z$ is analytic.
\end{corollary}
\begin{shortproof}
	This follows inductively from theorem \ref{theorem: a kind of hilbert's basis theorem} and theorem \ref{theorem: analyticity of Zblack}.
\end{shortproof}

The next stepping stone will be analyticity of any finitely generated $\Z$-algebra, that is precisely the quotients of the previous rings $\Z[T_1,\dots, T_n]$. The idea is essentially to choose an epimorphism $R\ldef\Z[T_1, \dots, T_n] \twoheadrightarrow A$ and observe that $A_\blacksquare\cong (A, R)_\blacksquare$. Then one obtains analyticity of $A_\blacksquare$ after noticing that analyticity of $(A, R)_\blacksquare$ can always be reduced to $R_\blacksquare$ as shown in lemma \ref{lemma: partial analyticity of (A, R)black}. To prove this reduction, one needs to control the derived scalar extension of (finitely) $R_\blacksquare$-free modules. Clausen and Scholze show this in the case that $A$ and $R$ are both finitely generated $\Z$-algebras in \cite[Intro. Appendix to Lec. VIII]{condenseddotpdf}. In this thesis, we will need the generalization where we do not require $R$ or $A$ to be finitely generated.
\begin{lemma}[Derived scalar extension is concentrated in degree $0$]
	\label{lemma: derived scalar extension is concentrated in degree 0}
	
	Let $\phi\colon R\to A$ be a morphism of discrete rings. Then for any extremally disconnected set $S$ the natural morphism
		$$\LeftD\scalarExtension{\phi}\free{R_\blacksquare}{S}\ldef A\Dotimes_R \free{R_\blacksquare} S\to A\otimes_R \free{R_\blacksquare}S \rdef \scalarExtension{\phi}\free{R_\blacksquare}{S}$$
	in $\derived{\associated{A}}$ is an isomorphism.
\end{lemma}
\begin{proof}
	We show the analogous claim for the (derived) tensor product with an arbitrary $R$-module $M$ instead of $A$. Observe that an $A$-linear map is an isomorphism if and only if it is an isomorphism as an $R$-linear map. We will thus show that the given morphism is an isomorphism of $R$-modules.
	
	Assume that $R$ is a finitely generated $\Z$-algebra. Then by Specker's theorem \ref{theorem: specker's theorem} there is a set $I$ for which $\free{R_\blacksquare}{S}\cong \prod_I R$. Now $R$ is Noetherian and thus coherent, so by \cite[\href{https://stacks.math.columbia.edu/tag/05CZ}{Tag 05CZ}]{stacks-project} the product $\prod_I R$ and hence $\free{R_\blacksquare}{S}$ are $R$-flat. In particular the natural morphism
		$$M\Dotimes_R \free{R_\blacksquare}{S} \to M\otimes_R \free{R_\blacksquare}{S}$$
	is a (quasi-)isomorphism.
	
	Now assume that $R$ is general. Recall that tensor products (as left adjoints) commute with colimits. As colimits of complexes are calculated degree-wise we observe further that also $\otimes_\category{A}^\bullet$ commutes with colimits of complexes. We thus obtain (by replacing $M\Dotimes_\category{A} -$ by $P^\bullet\otimes_\category{A}^\bullet -$ where $P^\bullet$ is a suitable left resolution of $M$) that
		$$M\Dotimes_R \free{R_\blacksquare}{S} \cong \colimit_{R'\subseteq R} M \Dotimes_R (R \otimes_{R'} \free{R'_\blacksquare}{S})$$
	where the colimit is calculated in chain complexes (but considered in $\derived{\associated{R}}$), is filtered and runs over all finitely generated $\Z$-subalgebras $R'$ of $R$. But by the previous case we already know that $R\Dotimes_{R'}\free{R'_\blacksquare}{S} \to R\otimes_{R'}\free{R'_\blacksquare}{S}$ is a (quasi-)isomorphism and hence we obtain by exactness of $\colimit_{R'\subseteq R}$ that
		$$\colimit_{R'\subseteq R} M \Dotimes_R (R \otimes_{R'} \free{R'_\blacksquare}{S})\cong \colimit_{R'\subseteq R} M \Dotimes_R R \Dotimes_{R'} \free{R'_\blacksquare}{S}\cong \colimit_{R'\subseteq R} M \Dotimes_{R'} \free{R'_\blacksquare}{S}.$$
	But again by the previous case (and exactness of $\colimit_{R'\subseteq R}$) we have
		$$\colimit_{R'\subseteq R} M\Dotimes_{R'} \free{R'_\blacksquare}{S} \isorightarrow \colimit_{R'\subseteq R} M\otimes_{R'} \free{R'_\blacksquare}{S}\cong \colimit_{R'\subseteq R} M\otimes_R R\otimes_{R'} \free{R'_\blacksquare}{S} \cong M\otimes_R \free{R_\blacksquare}{S},$$
	concluding the proof.
	
\end{proof}

We need one further technical result for lemma \ref{lemma: partial analyticity of (A, R)black}.
\begin{lemma}[Relative completeness]
	\label{lemma: relative completeness}

	Let $R\to A$ be a morphism of discrete rings such that $R_\blacksquare$ is analytic. Then the scalar restriction of any $(A, R)_\blacksquare$-free module is $R_\blacksquare$-complete.
\end{lemma}
\begin{proof}
	Since scalar restriction commutes with direct sums, we can assume with out of loss of generality that $M$ is of the form $\free{(A, R)_\blacksquare}{S}$ for some extremally disconnected set $S$. Observe that $\scalarRestriction{\phi}\free{(A, R)_\blacksquare}{S} = A\otimes_R \free{R_\blacksquare}{S}$. We will show more generally that for any discrete $R$-module $N$ the tensor product $N\otimes_R \free{R_\blacksquare}{S}$ is complete. For this choose an $R$-presentation $\bigoplus_I R \to \bigoplus_J R\to N \to 0$ of $N$. As $-\otimes_R \free{R_\blacksquare}{S}$ commutes with arbitrary colimits and is right exact, the induced sequence
		$$\bigoplus\nolimits_I \free{R_\blacksquare}{S} \to \bigoplus\nolimits_J \free{R_\blacksquare}{S}\to N \otimes_R \free{R_\blacksquare}{S} \to 0$$
	is exact. Thus $N\otimes_R \free{R_\blacksquare}{S}$ is $R_\blacksquare$-presentable hence $R_\blacksquare$-complete.
\end{proof}

We can now show the claimed partial result on analyticity of $(A, R)_\blacksquare$.
\begin{lemma}[Partial analyticity of $(A, R)_\blacksquare$]
	\label{lemma: partial analyticity of (A, R)black}
	
	Let $\phi\colon R\to A$ be a morphism of discrete rings such that $R_\blacksquare$ is analytic. Then $(A, R)_\blacksquare$ is analytic as well.
\end{lemma}
\begin{proof}
	Let $S$ be any extremally disconnected set and $C^\bullet$ a complex of $(A, R)_\blacksquare$-free modules, concentrated in non-positive degrees. By lemma \ref{lemma: derived scalar extension is concentrated in degree 0} we obtain $\scalarExtension{\phi}\free{R_\blacksquare} S \cong \LeftD\scalarExtension{\phi}\free{R_\blacksquare}S$. As $\scalarExtension{\phi}\ladj \scalarRestriction{\phi}$ and as $\scalarRestriction{\phi}$ is exact we obtain by lemma \ref{lemma: enriched adjunctions are stable under derivation} that in $\derived{\condensedAb}$ we have
	\begin{align*}
		\eRHom_A(\free{(A,R)_\blacksquare}{S}, C^\bullet)
		&= \eRHom_A(\scalarExtension{\phi}\free{R_\blacksquare}{S}, C^\bullet)\\
		&\cong \eRHom_A(\LeftD\scalarExtension{\phi}\free{R_\blacksquare}{S}, C^\bullet)\\
		&\cong\eRHom_R(\free{R_\blacksquare}{S}, \RightD\scalarRestriction{\phi}C^\bullet) \cong \eRHom_R(\free{R_\blacksquare}{S},\scalarRestriction{\phi}C^\bullet)
	\end{align*}
	where the last isomorphism follows since $\RightD\scalarRestriction{\phi} = \scalarRestriction{\phi}$ by exactness of $\scalarRestriction{\phi}$. But the complex $\scalarRestriction{\phi}C^\bullet$ is $R_\blacksquare$-complete by lemma \ref{lemma: relative completeness} and we thus further obtain that
		$$\eRHom_R(\free{R_\blacksquare}{S},\scalarRestriction{\phi}C^\bullet) \cong \eRHom_R(\free R S, \scalarRestriction{\phi}C^\bullet)\cong \eRHom_A(\free A S, C^\bullet)$$
	since clearly $\LeftD\scalarExtension{\phi}\free R S \cong \scalarExtension{\phi}\free R S \cong \free A S$ as $\free R S$ is projective.
\end{proof}

%
%
%

The final ingredient for proving analyticity for any finitely generated $\Z$-algebra is now to show that $(A, R)_\blacksquare$ and $A_\blacksquare$ are isomorphic as pre-analytic rings in case that $R\to A$ is a finite morphism. Clausen and Scholze argue this in \eg \cite[Exposition to Lecture IX]{condenseddotpdf} in case of finitely generated $\Z$-algebras. Again, we will need the generalization to arbitrary rings later on. We already state and prove the generalization here.
\begin{lemma}[Analytic rings associated to finite morphisms]
	\label{lemma: analytic rings associated to finite morphisms}
	
	If $\phi\colon R\to A$ is a finite morphism of discrete rings, then the morphism of pre-analytic rings $(A, R)_\blacksquare \to A_\blacksquare$ with underlying map of rings $\id_A$ from \ref{lemma: some maps of pre-analytic rings} is an isomorphism.
\end{lemma}
\begin{proof}
	Consider first the case that both $R$ and $A$ are finitely-generated $\Z$-algebras. Let $S$ be extremally disconnected and choose $I$ as in Specker's theorem \ref{theorem: specker's theorem}. Since $R\to A$ is finite, the $R$-module $A$ is a finitely generated (and hence finitely presented) and \cite[\href{https://stacks.math.columbia.edu/tag/059K}{Tag 059K}]{stacks-project} asserts that the natural morphism
		$$A\otimes_R \prod_I R \to \prod_I A\otimes_R R \cong \prod_I A$$
	is an isomorphism. But again by Specker, this is precisely the morphism $\free{(A, R)_\blacksquare}{S} \to \free{A_\blacksquare}{S}$, proving that $(A, R)_\blacksquare \to A_\blacksquare$ is an isomorphism (with inverse given by the inverse natural transformation $\free{A_\blacksquare}{-}\Rightarrow \free{(A, R)_\blacksquare}{-}$ and with the same underlying morphism $\id_A$).
	
	Now assume that $R$ and $A$ are arbitrary. Consider again some extremally disconnected $S$. Consider the partially ordered sets $P_R$ and $P_A$ of $\Z$-finitely generated subalgebras of $R$ and $A$ respectively. Denote by $P$ the category of finite morphisms $R'\to A'$ with $R' \in P_R$ and $A' \in P_A$ such that
	$$\begin{tikzcd}
		R' \ar[r]\ar[d] &A'\ar[d]\\
		R \ar[r] & A
	\end{tikzcd}$$
	commutes. By \cite[\href{https://stacks.math.columbia.edu/tag/0BTG}{Tag 0BTG}]{stacks-project} the two projections $P\to P_R$ and $P\to P_A$ are final. Then
	\begin{align*}
		\free{(A, R)_\blacksquare}{S}
		&\ldef A\otimes_R (\colimit_{R'} R\otimes_{R'} \free{R'_\blacksquare}{S})\\
		&\cong \colimit_{R'} A\otimes_{R'} \free{R'_\blacksquare}{S}\\
		&\cong \colimit_{R'\to A'}A\otimes_{R'} \free{R'_\blacksquare}{S}\\
		&\cong \colimit_{R'\to A'}A\otimes_{A'} A' \otimes_{R'} \free{R'_\blacksquare}{S}\\
		&\cong \colimit_{R'\to A'} A\otimes_{A'} \free{A'_\blacksquare}{S}\\
		&\cong \colimit_{A'} A\otimes_{A'} \free{A'_\blacksquare}{S} \rdef \free{A_\blacksquare}{S}
	\end{align*}
	where the first isomorphism is due to the tensor product $\otimes_R$ commuting with filtered colimits, the second and fifth isomorphism being due to finalness of $P\to P_R$ and $P\to P_A$ respectively, and the fourth isomorphism following from the previous case applied to each $R'\to A'$. By naturality of all isomorphisms it follows that $\free{(A, R)_\blacksquare}{S}\to\free{A_\blacksquare}{S}$ is an isomorphism as desired.
\end{proof}

\begin{theorem}[Analyticity of finite type $\Z$-algebras]
	\label{theorem: analyticity of finite type Z-algebras}
	
	For every finitely generated discrete $\Z$-algebra $A$, the ring $A_\blacksquare$ is analytic.
\end{theorem}
\begin{proof}
	Choose a surjection $R\ldef \Z[T_1, \dots, T_n]\twoheadrightarrow A$. Then $R_\blacksquare$ is analytic by corollary \ref{corollary: analyticity of Ablack for finitely generated free Z-algebras} and hence lemma \ref{lemma: partial analyticity of (A, R)black} gives the analyticity of $(A, R)_\blacksquare$. But lemma \ref{lemma: analytic rings associated to finite morphisms} yields that $A_\blacksquare \cong (A, R)_\blacksquare$ as any surjection of rings is finite (indeed a morphism of rings is surjective if and only if it is finite and epic by \cite[\href{https://stacks.math.columbia.edu/tag/04VT}{Tag 04VT}]{stacks-project}). Hence, $A_\blacksquare$ is analytic.
\end{proof}

We immediately obtain the following corollary.
\begin{corollary}[Analyticity of $(A, R)_\blacksquare$ for finitely generated $R$]
	\label{corollary: analyticity of (A, R)black for finitely generated R}
	
	For every map $R\to A$ of discrete $\Z$-algebras for which $R$ is finitely generated, the ring $(A, R)_\blacksquare$ is analytic.
\end{corollary}
\begin{shortproof}
	This immediately follows from proposition \ref{theorem: analyticity of finite type Z-algebras} together with lemma \ref{lemma: partial analyticity of (A, R)black}.
\end{shortproof}

This allows us to prove the final stage, namely that $A_\blacksquare$ is analytic for any discrete ring $A$. While Clausen and Scholze hint at the fact that this should be true (\eg immediately after \cite[Definition 9.1]{condenseddotpdf} -- albeit in slightly different form), they do not provide an argument in their notes.
\begin{theorem}[Analyticity of $A_\blacksquare$ for non-finite $\Z$-algebras $A$]
	\label{theorem: analyticity of Ablack for non-finite Z-algebras A}
	
	For any ring discrete $A$ the pre-analytic ring $A_\blacksquare$ is analytic.
\end{theorem}
\begin{proof}	
	Let $A'\subseteq A$ be a finitely generated subring. Then the pre-analytic ring $(A, A')_\blacksquare$ is analytic by corollary \ref{corollary: analyticity of (A, R)black for finitely generated R}.

	Consider the poset $P\ldef \{A''\subseteq A\setseparator A'' \text{ finitely generated}\}$ of such finitely generated subrings of $A$ and its sub-poset $P'\ldef \{A'' \in P\setseparator A'\subseteq A''\}$ of those subrings that contain $A'$. Denote be $i\colon P'\to P$ the inclusion. Observe that $i$ is a final functor. Indeed, for any given $A''\in P$ the comma category $(A''/i)=\{ A''' \in P'\setseparator A''\subseteq A'''\} = A''/P'$ is connected: This category is non-empty since the subring $\Z[A' \cup A''] \in P'$ contains $A''$. Furthermore, any two $A'''$ and $A''''$ are connected by a zig-zag $A'''\subseteq \Z[A'''\cup A''''] \supseteq A''''$. 
	
	In particular for any extremally disconnected set $S$ we obtain that the natural morphism of $A$-modules
	$$\colimit_{A'\subseteq A''\subseteq A} \free{(A, A'')_\blacksquare}{S} \to \colimit_{A''\subseteq A} \free{(A, A'')_\blacksquare}{S} \rdef \free{A_\blacksquare}{S}$$
	is an isomorphism. 
	
	Now for any such $A'\subseteq A''$ both $(A, A')_\blacksquare$ and $(A, A'')_\blacksquare$ are analytic by \ref{corollary: analyticity of (A, R)black for finitely generated R}. Thus, by proposition \ref{proposition: maps of analytic rings from maps of pre-analytic rings} the existence of the map of pre-analytic rings $(A, A')_\blacksquare \to (A, A'')_\blacksquare$ from lemma \ref{lemma: some maps of pre-analytic rings} implies that $\id_A$ is a map of analytic rings $(A, A')_\blacksquare \to (A, A'')_\blacksquare$. Thus, by proposition \ref{proposition: scalar restriction and extension of complete modules} each $\free{(A, A'')_\blacksquare}{S}$ is already $(A, A')_\blacksquare$-complete, implying that $\free{A_\blacksquare}{S} = \colimit_{A'\subseteq A'' \subseteq A} \free{(A, A'')_\blacksquare}{S}$ as a colimit of $(A, A')_\blacksquare$-complete modules is $(A, A')_\blacksquare$-complete as well. By a similar argumentation we find that the inclusion $A'\subseteq A$ is a map of analytic rings $A'_\blacksquare \to (A, A')_\blacksquare$, further implying that the scalar restriction of $\free{A_\blacksquare}{S}$ to $A'$ is also $A'_\blacksquare$-complete, by proposition \ref{proposition: scalar restriction and extension of complete modules}.
	
	Now suppose that $C^\bullet$ is a complex of $A_\blacksquare$-free modules concentrated in non-positive degrees. Then we have just shown that $C^\bullet$ is $(A, A')_\blacksquare$-complete and its scalar restriction to $A'$ is $A'_\blacksquare$-complete for any $A'\in P$. Now for some extremally disconnected $S$ we have
	\begin{align*}
		&\eRHom_A(\free{A_\blacksquare}{S}, C^\bullet) \\
		=\,&\eRHom_A(\colimit\nolimits_{\phi\colon A'\subseteq A} \free{(A, A')_\blacksquare}{S}, C^\bullet)\\
		\cong\,&\homotopyLimit\nolimits_{\phi\colon A'\subseteq A} \eRHom_A(\free{(A, A')_\blacksquare}{S}, C^\bullet)
	\end{align*}
	by \assumptionPreWord \ref{assumption: right derived Hom maps colimits to limits}. By lemma \ref{lemma: derived scalar extension is concentrated in degree 0} the natural map $\LeftD\scalarExtension{\phi}\free{A'_\blacksquare}{S}\to \scalarExtension{\phi}\free{A'_\blacksquare}{S}$ is a quasi-isomorphism. Recall that $\scalarExtension{\phi}\ladj \scalarRestriction{\phi}$ and that $\scalarRestriction{\phi}$ is exact. Hence, we obtain from lemma \ref{lemma: enriched adjunctions are stable under derivation} that in $\derived{\condensedAb}$ we have
	\begin{align*}
		\eRHom_A(\free{(A, A')_\blacksquare}{S}, C^\bullet)&= \eRHom_A(\scalarExtension{\phi}\free{A'_\blacksquare}{S}, C^\bullet)\\
		& \cong \eRHom_A(\LeftD\scalarExtension{\phi}\free{A'_\blacksquare}{S}, C^\bullet)\\
		&\cong \eRHom_{A'}(\free{A'_\blacksquare}{S}, \RightD\scalarRestriction{\phi}C^\bullet)\\
		&\cong \eRHom_{A'}(\free{A'_\blacksquare}{S}, \scalarRestriction{\phi}C^\bullet)\\
		&\cong \eRHom_{A'}(\free{A'}{S}, \scalarRestriction{\phi}C^\bullet) \cong \eRHom_A(\free{A}{S}, C^\bullet)
	\end{align*}
	since $\RightD\scalarRestriction{\phi} \cong \scalarRestriction{\phi}$, since $\scalarRestriction{\phi}C^\bullet$ is $A'_\blacksquare$-complete and since $\LeftD\scalarExtension{\phi}\free{A'}{S} = \free{A}{S}$. Applying this to the previous chain of isomorphisms we obtain that
	$$\eRHom_A(\free{A_\blacksquare}{S}, C^\bullet) \cong \homotopyLimit\nolimits_{\phi\colon A'\subseteq A} \eRHom_A(\free{A}{S}, C^\bullet) \cong \eRHom_A(\free{A}{S}, C^\bullet)$$
	by another application of \assumptionPreWord \ref{assumption: right derived Hom maps colimits to limits}. Thus, $A_\blacksquare$ is analytic as well.
\end{proof}

We finally obtain the relative case in full generality.
\begin{corollary}[Analyticity of $(A, R)_\blacksquare$]
	\label{corollary: analyticity of (A, R)black}
	
	For any morphism of discrete rings $R\to A$ the pre-analytic ring $(A, R)_\blacksquare$ is analytic.
\end{corollary}
\begin{shortproof}
	By theorem \ref{theorem: analyticity of Ablack for non-finite Z-algebras A} the ring $R_\blacksquare$ is analytic, hence by lemma \ref{lemma: partial analyticity of (A, R)black} so is $(A, R)_\blacksquare$.
\end{shortproof}

We can now translate the morphism of pre-analytic rings from lemma \ref{lemma: some maps of pre-analytic rings} into the analytic world.
\begin{corollary}[Some maps of analytic rings]
	\label{corollary: some maps of analytic rings}
	
	We have the following maps of analytic rings
	\begin{enumerate}[(i)]
		\item
		For any discrete ring $A$ the map $\id_A\colon A\to A_\blacksquare$.
		
		\item
		Any morphism $\phi\colon R\to A$ of discrete rings is the map $\phi\colon R_\blacksquare\to A_\blacksquare$.
		
		\item
		For any morphism $R\to A$ of discrete rings the map $\id_A\colon (A, R)_\blacksquare\to A_\blacksquare$.
		
		\item
		For any commutative square
			$$\begin{tikzcd}
				R\ar[r]\ar[d] &A\ar[d, "\phi"]\\
				R'\ar[r] & A'
			\end{tikzcd}$$
		of discrete rings the map $\phi\colon (A, R)_\blacksquare \to (A', R')_\blacksquare$.
	\end{enumerate}
\end{corollary}
\begin{proof}
	By proposition \ref{proposition: maps of analytic rings from maps of pre-analytic rings} it is enough to exhibit in each of the cases a map of pre-analytic rings with matching morphism of underlying rings, as this map of underlying rings then is a map of analytic rings. But lemma \ref{lemma: some maps of pre-analytic rings} provides exactly such maps.
\end{proof}

\newpage
\section{An Analogue of Hilbert's Basis Theorem}
\label{section: an analogue of hilbert's basis theorem}

The goal of this section is provide a proof of theorem \ref{theorem: a kind of hilbert's basis theorem} -- the analogue of Hilbert's basis theorem.

\begin{convention}
	Throughout this section let $R$ be a finitely generated discrete $\Z$-algebra. Denote by $A\ldef R[T]$ the ring of polynomials and by $A_\infty\ldef \laurent{R}{T^{-1}}$ the ring of \emph{formal Laurent series in $T^{-1}$} (\ie series $\sum_{n\gg \infty} r_n T^{-n}$ with finite principal part), with its natural $\ideal{T^{-1}}$-adic topology.
	
	Assume that $R_\blacksquare$ is analytic. Then by lemma \ref{lemma: partial analyticity of (A, R)black} also $(A, R)_\blacksquare$ is analytic and the identity $\id_A$ is a morphism of analytic rings $R_\blacksquare \to (A, R)_\blacksquare$ by corollary \ref{corollary: some maps of analytic rings}, implying that any $(A, R)_\blacksquare$-complete module is also $R_\blacksquare$-complete by proposition \ref{proposition: scalar restriction and extension of complete modules}.
\end{convention}

\begin{lemma}
	\label{lemma: Ainfty is A-flat, f.g. case}
	
	The $A$-module $A_\infty$ is $A$-flat.
\end{lemma}
\begin{proof}
	Observe that $A_\infty = \laurent{R}{T^{-1}} = \series{R}{T^{-1}}[T]$ is the $\ideal{T^{-1}}$-adic completion of the ring $A'\ldef R[T, T^{-1}]$. Since $R$ is $\Z$-finitely generated and hence Noetherian, Hilbert's basis theorem asserts that $A'\cong R[T, U]/\ideal{UT - 1}$ is Noetherian as well. Thus by \cite[\href{https://stacks.math.columbia.edu/tag/00HT}{Tag 00HT}]{stacks-project}, the completion is $A'$-flat. But $A' = R[T, T^{-1}]\cong R[T]_T = A_T$ is a localization of $A$, hence $A$-flat by \cite[\href{https://stacks.math.columbia.edu/tag/00HT}{Tag 00HT}]{stacks-project}. As flatness is transitive by \cite[\href{https://stacks.math.columbia.edu/tag/00HT}{Tag 00HT}]{stacks-project}, we obtain that $A_\infty$ is $A$-flat.
\end{proof}

\begin{lemma}
	\label{lemma: resolution for Ainfty, f.g. case}
	There is an exact sequence of $A$-modules
		$$0\to A \otimes_R \series R U \xrightarrow{UT - 1} A \otimes_R \series R U \to A_\infty\to 0$$
	where the first two terms are compact projectives of $\mod{(A, R)_\blacksquare}$. In particular $A_\infty$ is $(A, R)_\blacksquare$-complete, $R_\blacksquare$-complete and quasi-isomorphic to a bounded complex of compact projectives of $\mod{(A, R)_\blacksquare}$ hence compact in $\derived{(A, R)_\blacksquare}$.
\end{lemma}
\begin{proof}
	Observe that $A \otimes_R \series R U \cong \mleft(\bigoplus_\N R\mright)\otimes_R \series{R}{U}\cong \bigoplus_\N \series{R}{U}\cong \series R U [T]$ is the ring of polynomials over the ring of power series, hence exactness of the sequence immediately follows as $A_\infty = \laurent{R}{T^{-1}} = \series{R}{T^{-1}}[T]$. It remains to show that $A\otimes_R \series{R}{U}$ is a compact projective object of $\mod{(A, R)_\blacksquare}$. But by (the proof of) lemma \ref{lemma: more compact projectives} the object $A\otimes_R \series{R}{U} \cong A\otimes_R \prod_\N R$ is a compact projective of $\mod{(A, R)_\blacksquare}$ (as the retraction is preserved under $A\otimes_R -$).

\end{proof}

We will need the following statement even for general rings $R$ that are not $\Z$-finitely generated.
\begin{lemma}[Formal Laurent series are idempotent]
	\label{lemma: formal Laurent series are idempotent}
	
	If $R$ is any ring then multiplication 
	$$\laurent{R}{T^{-1}}\otimes_{R[T]}\laurent{R}{T^{-1}}\to \laurent{R}{T^{-1}}$$
is an $R[T]$-module isomorphism.
\end{lemma}
\begin{proof}
	Define a morphism in the other direction such that $a\mapsto 1\otimes a$. It is clear that the composition $\laurent{R}{T^{-1}} \to \laurent{R}{T^{-1}}\otimes_{R[T]} \laurent{R}{T^{-1}}\to \laurent{R}{T^{-1}}$ is the identity. To show that the other composition is the identity as well, we need that $1\otimes ab = a\otimes b$ for all $a, b\in \laurent{R}{T^{-1}}$. For this write $a = \sum_{k\gg -\infty} r_k T^{-k}$ giving a sequence $(a_n)_{n\in \N}$ with $a_n \ldef \sum_{n \ge k \gg -\infty} r_k T^{-k}$ polynomials converging to $a$ in the $\ideal{T^{-1}}$-adic topology. Write further $a_n = T^{v_n} \hat{a_n}$ with $v_n\in \Z$ and $\hat{a_n} \in R[T]$ not divisible by $T$. Observe that since the tensor product is $T$-balanced and hence also $T^{-1}$-balanced. Indeed for $x,y \in \laurent{R}{T^{-1}}$ arbitrary we find
	$$xT^{-1}\otimes y = xT^{-1}\otimes TT^{-1}y = xT^{-1}T\otimes T^{-1}y = x\otimes T^{-1}y$$
	as $T^{-1} \in \laurent{R}{T^{-1}}$ exists. We now find that
	$$1\otimes a_nb = 1\otimes T^{v_n} \hat{a_n} b = T^{v_n} \hat{a_n} \otimes b = a_n\otimes b,$$
	implying by uniqueness of limits that $1\otimes ab = a\otimes b$ as $n\to \infty$, as desired.
\end{proof}

\begin{lemma}[$A_\infty$ is idempotent]
	\label{lemma: Ainfty is idempotent, f.g. case}
	
	The $A$-algebra $A_\infty$ is idempotent, by which we mean that the multiplication map of $A_\infty$ gives rise to an isomorphism $A_\infty\Dotimes_{(A, R)_\blacksquare} A_\infty \to A_\infty$ in $\derived{(A, R)_\blacksquare}$.
\end{lemma}
\begin{proof}
	By lemma \ref{lemma: formal Laurent series are idempotent} we know that multiplication $A_\infty\otimes_A A_\infty\to A_\infty$ is an $A$-module isomorphism. Now as $A_\infty$ is $A$-flat by lemma \ref{lemma: Ainfty is A-flat, f.g. case} the natural isomorphism $A_\infty\Dotimes_A A_\infty \to A_\infty\otimes_A A_\infty$ is a isomorphism. Recall that the derived completion $\derivedCompletion{-}{(A,R)_\blacksquare}$ is symmetric monoidal for $\Dotimes_{(A,R)_\blacksquare}$ by lemma \ref{lemma: completed scalar extension is monoidal}. But $A_\infty$ and hence $A_\infty \Dotimes_A A_\infty$ are $(A, R)_\blacksquare$-complete by lemma \ref{lemma: resolution for Ainfty, f.g. case} and thus in total
		$$A_\infty\Dotimes_{(A, R)_\blacksquare} A_\infty\cong A_\infty \Dotimes_A A_\infty \cong A_\infty\otimes_A A_\infty.$$
	Since the multiplication maps of all three tensor products are compatible and since the underived one is an isomorphism we obtain that
		$$A_\infty \Dotimes_{(A, R)_\blacksquare} A_\infty \to A_\infty$$
	is an isomorphism as well.
\end{proof}

An essentially purely formal consequence of the previous lemma \ref{lemma: Ainfty is idempotent, f.g. case} is given in the following lemma that characterizes $A_\infty$-modules inside of $\derived{(A, R)_\blacksquare}$. While Clausen and Scholze state this result, they omit the reasoning and provide no reference. The lemma follows from the fact that the \quote{Eilenberg--More-category of an \emph{idempotent} monad on a category $\category{C}$ is a reflective subcategory of $\category{C}$}. Alternatively this follows from work of Boyarchenko and Drinfeld in \cite{boyarchenko-drinfeld2009idempotents}, which rephrases this statement about monads in terms of the more familiar theory of monoidal categories.
\begin{lemma}[$A_\infty$-modules]
	\label{lemma: Ainfty-modules, f.g. case}
	
	The $A_\infty$-modules in $\derived{(A, R)_\blacksquare}$ form a full subcategory whose inclusion has a left adjoint given by the functor $A_\infty\Dotimes_{(A, R)_\blacksquare} -$. Furthermore, any $M \in \derived{(A, R)_\blacksquare}$ admits at most one $A_\infty$-module structure and such a structure exists if and only if the natural morphism $M\to A_\infty\Dotimes_{(A, R)_\blacksquare} M$ is an isomorphism.
\end{lemma}
\begin{proof}
	Denote by $\mu\colon A_\infty\Dotimes_{(A, R)_\blacksquare} A_\infty\to A_\infty$ the multiplication. In the language of \cite[Section 2]{boyarchenko-drinfeld2009idempotents}, the pair $e\ldef (A_\infty, \mu)$ is a \emph{unital algebra} in the monoidal category $\category{M}\ldef \derived{(A, R)_\blacksquare}$. Even more, since $\mu$ is an isomorphism, the pair is an \emph{idempotent} unital algebra. In particular the \emph{unit} $A\to A_\infty$ is an \emph{idempotent arrow} by \cite[Remark 2.9 (vi)]{boyarchenko-drinfeld2009idempotents}, making $(A_\infty, \mu)$ a \emph{closed} unital idempotent algebra. Denote by $\mod{e}$ the category of \emph{unital (left) $e$-modules} of this algebra.

	Consider the essential image $e\category{M}\ldef \{M \in \derived{(A, R)_\blacksquare}\setseparator A_\infty\Dotimes_{(A, R)_\blacksquare} M \cong M\}$ of the functor $M \mapsto A_\infty\Dotimes_{(A, R)_\blacksquare} M$ and denote by $L$ its (co)restriction $\category{M}\to e\category{M}$.
	
	By \cite[Lemma 2.32 (ii)]{boyarchenko-drinfeld2009idempotents} the forgetful functor $\mod{e} \to \category{M}$ is fully faithful and its essential image is $e\category{M}$, that is $\mod{e}\equiv e\category{M}$. A (pseudo)inverse $e\category{M} \to \category{M}$ is given by the obvious way of enriching $A_\infty\Dotimes_{(A, R)_\blacksquare} M$ with a unital $e=(A_\infty, \mu)$-module structure. Furthermore by \cite[Prop. 2.22 (a)]{boyarchenko-drinfeld2009idempotents} the functor $L\colon \category{M}\to e\category{M}$ is left adjoint to the fully faithful inclusion $e\category{M} \to \category{M}$, exhibiting $e\category{M}\equiv \mod{e}$ as a localization of $\category{M}$. Thus the functor $\mod{e} \to \derived{(A, R)_\blacksquare}$ is fully faithful with left adjoint given by $A_\infty\Dotimes_{(A, R)_\blacksquare} -$ as claimed. By \cite[Lemma 2.15]{boyarchenko-drinfeld2009idempotents} it follows that an $M$ admits an $(A_\infty, \mu)$-module structure if and only $M \to A_\infty\Dotimes_{(A, R)_\blacksquare} M$ is an isomorphism. As $e\category{M}\to \category{M}$ is fully faithful it is also clear that any $M\in\category{M}$ admits at most one unital $(A_\infty, \mu)$-module structure because if there are two modules $X, Y \in e\category{M}$ mapping both to $M \in \category{M}$ then the identity $\id_M$ lifts to an isomorphism $X \to Y$ of unital $(A_\infty, \mu)$-modules.
		
	It remains to see that the notion of $A_\infty$-module and unital $(A_\infty, \mu)$-module inside $\derived{(A, R)_\blacksquare}$ agree. For this observe that $A_\infty\Dotimes_{(A, R)_\blacksquare} M = \derivedCompletion{A_\infty\Dotimes_A M}{(A, R)_\blacksquare}$ since derived completion is (strong) symmetric monoidal by lemma \ref{lemma: completed scalar extension is monoidal}. Since further the completion is left adjoint to the inclusion into $A$-modules we find that
	$$\Hom_\derived{(A, R)_\blacksquare}(A_\infty\Dotimes_{(A, R)_\blacksquare} M, M) \cong \Hom_\derived{A}(A_\infty\Dotimes_A M, M).$$
	But since $A_\infty$ is $A$-flat by lemma \ref{lemma: Ainfty is A-flat, f.g. case} the natural morphism $A_\infty\Dotimes_A M \to A_\infty\otimes_A M$ is a quasi-isomorphism, yielding further that
	$$\Hom_\derived{(A, R)_\blacksquare}(A_\infty\Dotimes_{(A, R)_\blacksquare} M, M) \cong \Hom_\derived{A}(A_\infty \otimes_A M, M).$$
	Under this correspondence it is clear that inside of $\derived{(A, R)_\blacksquare}$ the notions of unital $(A_\infty, \mu)$-module and $A_\infty$-module agree, as one can translate the $A_\infty$-action on $M$ from one notion to the other.

\end{proof}

\begin{lemma}
	\label{lemma: Ainfty-modules admit no extensions, f.g. case}
	
	We have
		$$\eRHom_A(A_\infty, A) \cong 0$$
	in $\derived{\condensedAb}$. In particular if $C^\bullet \in \derived{\associated{A}}$ is a complex of $A_\blacksquare$-free modules and if $M$ is any $A_\infty$-module, then $\eRHom_A(M, C^\bullet)\cong 0$.
\end{lemma}
\begin{proof}
	Using the quasi-isomorphism of $A$-modules from lemma \ref{lemma: resolution for Ainfty, f.g. case}, we obtain that the complex $\eRHom_A(A_\infty, A)$ is quasi-isomorphic to the complex
	$$0\to \eRHom_R(\series{R}{U}, A) \xrightarrow{UT-1} \eRHom_R(\series{R}{U}, A)\to 0$$
	as $A\otimes_R \series{R}{U} \cong A\Dotimes_R \series{R}{U}$ by $R$-flatness of $A$ and since the $\condensedAb$-enriched adjunction of scalar extension and restriction is derivable by lemma \ref{lemma: enriched adjunctions are stable under derivation}. Clausen and Scholze now claim in \cite[Observation 8.8, p.55]{condenseddotpdf} that this complex is quasi-isomorphic to the complex
		$$0\to A[U^{-1}]/A \xrightarrow{UT-1} A[U^{-1}]/A\to 0,$$
	the author was however not able to verify this. We thus have to make this an \assumptionPreWord. Now this complex is acyclic since left multiplication by $UT-1$ is injective and the cokernel is $(A[U^{-1}]/A)/\ideal{UT-1}\cong A[T]/A\cong 0$ implying that $UT-1$ is surjective as well. So if we take the alternative representation as given we find that $\eRHom_R(A_\infty, A) \cong 0$ as claimed.
	
	Suppose that $C^\bullet = C[0]$ is concentrated in degree $0$ where $C = \bigoplus_{i\in I} \free{A_\blacksquare}{S_i}$ is $A_\blacksquare$-free. By Specker's theorem \ref{theorem: specker's theorem} there are sets $J_i$ for $i\in I$ such that $\free{A_\blacksquare}{S_i} \cong \prod_{J_i} A$. By lemma \ref{lemma: resolution for Ainfty, f.g. case} we find that $A_\infty$ is a derived compact object of $\derived{(A, R)_\blacksquare}$ so that by lemma \ref{lemma: solid derived compacts are enriched compact} it is also $\condensedAb$-enriched compact. As $\eRHom_A(A_\infty, -)$ additionally commutes with direct products as a consequence of remark \ref{remark: compatibility of enriched Homs with direct sums and products} we obtain that
		$$\eRHom_A(A_\infty, C^\bullet) \cong \bigoplus_{i\in I} \prod_{J_i} \eRHom_A(A_\infty, A)\cong 0$$
	as all these $A$-modules are $(A, R)_\blacksquare$-complete. Indeed $A_\infty$ is complete by lemma \ref{lemma: resolution for Ainfty, f.g. case} and the other cases follow since $A \cong A\otimes_R\free{R_\blacksquare}{\ast}$ is complete and since $\mod{(A, R)_\blacksquare}$ has all limits and colimits. By induction on the length and by utilizing stupid truncations it follows that any bounded complex $C^\bullet$ of $A_\blacksquare$-free modules also satisfies $\eRHom_A(A_\infty, C^\bullet) \cong 0$. Now if $C^\bullet\in\derived{A_\blacksquare}$ is bounded to the right then there is a left resolution $P^\bullet \to C^\bullet$ by $A_\blacksquare$-free modules which is also bound to the right. Writing $P^\bullet$ as a sequential colimit of its truncations $P_n^\bullet \ldef \sTruncation{\ge -n} P^\bullet$ we obtain that
		$$\eRHom_A(A_\infty, C^\bullet) \cong \eRHom_A(A_\infty, P^\bullet) \cong \colimit_{n\in \N} \eRHom_A(A_\infty, P_n^\bullet) \cong 0$$
	by lemma \ref{lemma: enriched compact objects preserve sequential homotopy colimits} as again, $A_\infty$ is enriched compact. Finally if $C^\bullet \in \derived{A_\blacksquare}$ is general, we can write it as the colimit $\colimit_{n\in \N} \sTruncation{\le n}C^\bullet$ of bound to the right complexes which implies
	$$\eRHom_A(A_\infty, C^\bullet)\cong \colimit_{n\in \N} \eRHom_A(A_\infty, \sTruncation{\le n} C^\bullet) \cong 0$$
	as $\eRHom_A(A_\infty, -)$ still preserves (homotopy) colimits by the same lemma.
	
	Now suppose that $M$ is an $A_\infty$-module. Then $M\cong M\Dotimes_{(A, R)_\blacksquare} A_\infty$ by lemma \ref{lemma: Ainfty-modules, f.g. case}. Thus for any $C^\bullet \in \derived{A_\blacksquare}$ we find that
	\begin{align*}
		\eRHom_A(M, C^\bullet) & \cong \eRHom_A(M\Dotimes_{(A, R)_\blacksquare} A_\infty, C^\bullet)\\
		&\cong \eRHom_A(\derivedCompletion{M\Dotimes_A A_\infty}{(A, R)_\blacksquare}, C^\bullet) \cong \eRHom_A(M\Dotimes_A A_\infty, C^\bullet).
	\end{align*}
	as derived completion is strong symmetric monoidal and as the $\condensedAb$-enriched adjunction between completion and inclusion into $A$-modules is derivable by \ref{lemma: enriched adjunctions are stable under derivation}. This then implies that
		$$\eRHom_A(M, C^\bullet) \cong \eRHom_A(M, \intRHom_A(A_\infty, C^\bullet))$$
	by lemma \ref{lemma: deriving a closed monoidal structure} as $\eRHom_A$ is the (complex of) condensed abelian group(s) underlying $\intRHom_A$. But as previously shown $\intRHom_A(A_\infty, C^\bullet) \cong 0$ as the underlying (complex of) condensed abelian group(s) of $\intRHom_A$ is $\eRHom_A$. Thus
		$$\eRHom_A(M,C^\bullet) \cong \eRHom_A(M, 0) \cong 0$$
	as well.
\end{proof}

\begin{lemma}
	\label{lemma: special cokernels are Ainfty-modules, f.g. case}
	
	For any set $I$, the cokernel of the injective map
	$$A\otimes_R \prod_I R \to \prod_I A$$
	is an $A_\infty$-module. In particular, by Specker's theorem \ref{theorem: specker's theorem} for any extremally disconnected $S$ the cokernel of the map $\free{(A, R)_\blacksquare}{S} \to \free{A_\blacksquare}{S}$ is an $A_\infty$-module.
\end{lemma}
\begin{proof}
	Recall that by lemma \ref{lemma: Ainfty-modules, f.g. case} any $A$-module admits at most one $A_\infty$-module structure and that an additive map between $A_\infty$-modules is $A_\infty$-linear if and only if it is $A$-linear.
	
	The short exact sequence
		$$0\to \prod_I R\to\prod_I \series{R}{T^{-1}} \to \prod_I T^{-1}\series{R}{T^{-1}}\to 0$$
	gives rise to the short exact sequence
		$$0\to A\otimes_R \prod_I R\to A\otimes_R \prod_I \series{R}{T^{-1}} \to A\otimes_R \prod_I T^{-1}\series{R}{T^{-1}} \to 0$$
	as $A=R[T]$ is $R$-flat. Observe that as $A$-modules, we have $A_\infty = \laurent{R}{T^{-1}} \cong \series{R}{T^{-1}}[T] \cong A \otimes_R \series{R}{T^{-1}}$ so that
		$$A_\infty \otimes_{\series{R}{T^{-1}}} \prod_I \series{R}{T^{-1}}\cong A\otimes_R \series{R}{T^{-1}}\otimes_{\series{R}{T^{-1}}} \prod_I \series{R}{T^{-1}}\cong A\otimes_R \prod_I \series{R}{T^{-1}},$$
	allowing us to replace the central term of the previous short exact sequence by an $A_\infty$-module.	Observe further that $A_\infty/A \cong T^{-1}\series R {T^{-1}}$. We obtain the following commutative diagram of condensed $A$-modules
	$$\begin{tikzcd}
		&0\ar[d] &0\ar[d] &0\ar[d]\\
		0\ar[r] &A\otimes_R \prod_I R \ar[r]\ar[d] &\prod_I A\ar[r]\ar[d] &C \ar[r]\ar[d] &0\\
		0\ar[r] &A_\infty\otimes_{\series{R}{T^{-1}}} \prod_I \series R {T^{-1}} \ar[r]\ar[d] &\prod_I A_\infty \ar[r]\ar[d] &C' \ar[r]\ar[d] &0\\
		0 \ar[r] &\prod_I T^{-1}\series R{T^{-1}} \ar[r]\ar[d] &\prod_I T^{-1}\series{R}{T^{-1}} \ar[r]\ar[d] &0 \ar[r]\ar[d] &0\\
		&0 &0 &0
	\end{tikzcd}$$
	with $C$ and $C'$ cokernels, such that all rows are exact. Then by the above argument $2$ out of the $3$ columns are exact and by the $9$ lemma so is the third, in particular $C\cong C'$ as $A$-modules. Now $C'$ is a cokernel of $A$-modules that admit a (then necessarily unique) $A_\infty$-module structure and the cokernel of this (necessarily $A_\infty$-linear) morphism is an $A_\infty$-module as well. Hence $C\cong C'$ is an $A_\infty$-module.
\end{proof}

We can now finally state the proof of the analogue to Hilbert's Basis Theorem.
\begin{proof}[of theorem \ref{theorem: a kind of hilbert's basis theorem}]
	Let $C^\bullet$ be a complex of $A_\blacksquare$-free modules concentrated in non-negative degrees. Let $S$ be extremally disconnected. By Specker's theorem \ref{theorem: specker's theorem} each term of $C^\bullet$ is a direct sum of products of copies of $A$. As $A\cong A\otimes_R \free {R_\blacksquare}{\ast}$ is $(A, R)_\blacksquare$-complete, so is each term of $C^\bullet$ and hence $C^\bullet\in\derived{(A,R)_\blacksquare}$. Thus
		$$\eRHom_A(\free{(A, R)_\blacksquare}{S}, C^\bullet) \isorightarrow \eRHom_A(\free A S, C^\bullet).$$
	By lemma \ref{lemma: special cokernels are Ainfty-modules, f.g. case} the cokernel $M$ of the natural morphism $\free{(A,R)_\blacksquare}{S}\to\free{A_\blacksquare}{S}$ is an $A_\infty$-module. We obtain the distinguished triangle
		$$\eRHom_A(M, C^\bullet)\to\eRHom_A(\free{A_\blacksquare}{S}, C^\bullet)\to\eRHom_A(\free{(A, R)_\blacksquare}{S}, C^\bullet)\to (\cdots)[1]$$
	for which $\eRHom_A(M, C^\bullet) \cong 0$ by lemma \ref{lemma: Ainfty-modules admit no extensions, f.g. case}. Hence
		$$\eRHom_A(\free{A_\blacksquare}{S}, C^\bullet)\to\eRHom_A(\free{(A, R)_\blacksquare}{S}, C^\bullet)$$
	is an isomorphism. Overall we thus obtain that the natural morphism
	$$\eRHom_A(\free{A_\blacksquare}{S}, C^\bullet)\isorightarrow \eRHom_A(\free{(A, R)_\blacksquare}{S}, C^\bullet) \isorightarrow \eRHom_A(\free A S, C^\bullet)$$
	is an isomorphism. Thus any such complex $C^\bullet$ of $A_\blacksquare$-free modules is derived $A_\blacksquare$-complete and hence $A_\blacksquare$ is analytic.
\end{proof}

	\chapter{Affine Duality}
	\label{chapter: affine duality}
	
In this chapter we will show that in the case of (well-behaved: locally Noetherian, ...) affine schemes the desired functors $f_!$ and $f^!$ exist and that they satisfy some expected properties. 

It is to note that in \cite{condenseddotpdf} Clausen and Scholze show this only in case of finite type affine $\Z$-schemes, \confer \cite[theorem 8.13]{condenseddotpdf}. While they later \emph{globalize} the result in \cite[Chapters X \& XI]{condenseddotpdf} to arbitrary schemes and hence recover the unrestricted affine case, this is somewhat unsatisfactory to the author. One then essentially describes an affine scheme $\Spec A$ as the limit $\limit_{A'} \Spec A'$ consisting of sequences of compatible points from the affine schemes $\Spec A'$ for finitely generated $\Z$-subalgebras of $A$. This is of course quite inexplicit and, to the author, not so geometrically appealing. Even though we will of course adopt the same viewpoint, we will do some more implicitly and try to develop the results given by Clausen and Scholze for general affine schemes instead. This entails describing many of the results given by Clausen and Scholze -- but after some amount of the above gluing. Proving the generalized results often worked, but not always! If there is some complication, we will denote it (obviously) at each place.

While \quote{explicitly} describing this gluing is of course technically redundant after one has already globalized, it is still somewhat more satisfying to the author.

Additionally, the mentioned globalization makes use of $\infty$-theoretical results, as one needs to glue derived categories of the form $\derived{(A', R')_\blacksquare}$ -- something that is not possible on the level of \quote{ordinary} derived categories. While we will state the results in chapter \ref{chapter: globalization}, we will not make an attempt to formally justify them. Hence, if we were to not translate the results of \cite[theorem 8.13]{condenseddotpdf} to general affine schemes, we would have, in this thesis, only proven duality for finite type affine $\Z$-schemes -- a very restrictive geometric setting (this excludes for example all affine schemes over infinite fields).

Before we start with our preparations we should fix some notation. Much like the category of affine schemes can be thought of as $\commutativeRings^\op$, it will prove useful to talk about the analytic rings $A_\blacksquare$ and $(A, R)_\blacksquare$ as if they were geometric objects. For this we will introduce some terminology. Note that this terminology is not used by Clausen \& Scholze.
\begin{definition}[Analytic spectra \& virtual maps]
	\label{definition: analytic spectra and virtual maps}
	
	Let $\Spec\colon \analyticRings \to \analyticRings^\op$ be the canonical contravariant functor. For an analytic ring $\preAnaRing{A}$ the object $\Spec \preAnaRing{A} \in \analyticRings^\op$ is called the \emph{analytic spectrum of $\preAnaRing{A}$}. Any map of analytic rings $\preAnaRing{A}\to \preAnaRing{B}$ induces a map of analytic spectra $\Spec\preAnaRing{B}\to \Spec\preAnaRing{A}$. Such a map of analytic spectra $\Spec \preAnaRing{B} \to \Spec \preAnaRing{A}$ is called \emph{virtual} if the associated morphism $\underlying{A}\to \underlying{B}$ is the identity. 
\end{definition}

\begin{remark}
	If one is acquainted with the theory of adic spaces this might seem familiar (see chapter \ref{chapter: globalization} for a very brief introduction). Indeed, to any pair of rings $R\subseteq A$ with $R$ integrally closed in $A$ one can associate the adic spectrum $\Spa(A, R)$. Furthermore, if $R\to A$ is general, and we denote by $\tilde R$ the integral closure of (the image of) $R$ inside $A$, then by lemma \ref{lemma: analytic rings associated to finite morphisms} we find that $R_\blacksquare \cong \tilde R_\blacksquare$ and hence that $(A, R)_\blacksquare \cong (A, \tilde R)_\blacksquare$. That is for a general $(A, R)_\blacksquare$ one can associate the adic spectrum $\Spa(A, \tilde R)$. This seems to be a proper dictionary, but the author has not investigated it. As we will not make use of the theory of adic space until chapter \ref{chapter: globalization} it is simply more convenient to introduce artificial language as in definition \ref{definition: analytic spectra and virtual maps}.
\end{remark}

\begin{remark}[Notation]
	Suppose that $\phi\colon \preAnaRing{A}\to \preAnaRing{B}$ is a morphism of analytic rings with associated map of analytic spectra $f\colon \Spec \preAnaRing{B} \to \Spec\preAnaRing{A}$. To respect the change in perspective when switching from analytic rings (the algebraic objects) to analytic spectra (the geometric objects) we will denote by $f_\ast$ the scalar restriction $\scalarRestriction{\phi}\colon \mod{\preAnaRing{B}}\to\mod{\preAnaRing{A}}$ and by $f^\ast$ the scalar extension $\completedScalarExtension{\phi}{\preAnaRing{B}}\colon \mod{\preAnaRing{A}}\to\mod{\preAnaRing{B}}$. Since we will be mostly working in the derived categories we will further write $f_\ast$ and $f^\ast$ for the associated derived functors. 
\end{remark}

\begin{remark}[What about virtual maps?]
	Any morphism of discrete rings $\phi\colon R\to A$ (corresponding to morphism of schemes $f\colon \Spec A\to \Spec A$) gives rise to three (in general distinct) analytic rings: $R_\blacksquare$, $(A, R)_\blacksquare$ and $A_\blacksquare$. In fact, these provide a factorization of the map of analytic rings $\phi\colon R_\blacksquare \to A_\blacksquare$ as
		$$R_\blacksquare \xrightarrow{\phi} (A, R)_\blacksquare \xrightarrow{\id_A} A_\blacksquare.$$
	If one now wants to study how modules over $R$ and $A$ (\ie quasi-coherent sheaves on $\Spec R$ and $\Spec A$) behave under transport along $\phi$ the main difficulty comes from the second map $\id_A\colon (A, R)_\blacksquare\to A_\blacksquare$ -- with associated virtual map $j\colon \Spec A_\blacksquare\to\Spec(A, R)_\blacksquare$. Indeed, $(A, R)_\blacksquare$ is essentially \emphquote{just a copy of $R_\blacksquare$ in the world of $A$-modules}. However, the difference between $(A, R)_\blacksquare$ and $A_\blacksquare$ is surprisingly subtle and most of the remaining work until proving the existence of a good theory of $f_!$ and $f^!$ in theorem \ref{theorem: exceptional direct and inverse image functors, the affine case} is dedicated to the study of this virtual map in the special case of $A=R[T]$, the coordinate ring of the affine space $\A_R^1$ over $\Spec R$. For example, we will construct a functor $j_!$ that will be used to define an appropriate exceptional direct image $f_!\colon \derived{A_\blacksquare}\to\derived{R_\blacksquare}$ that gives rise to $f^!$ as its right adjoint.
\end{remark}

\section{Preparations: The case of \texorpdfstring{$\A_R^1$}{Affine Space over an Affine Base}}

As already seen in the proof of corollary \ref{corollary: analyticity of Ablack for finitely generated free Z-algebras} and theorem \ref{theorem: analyticity of finite type Z-algebras}, it is often relatively straightforward to reduce the proof a given claim from the general case of finite type algebra to a free algebra to the free algebra in one variable. We will do the same in for example the proof of theorem \ref{theorem: exceptional direct and inverse image functors, the affine case} which asserts the existence of $f_!$ and $f^!$. We will thus study this case thoroughly. To this end consider the following convention.
\begin{convention}
	Throughout this section let $R$ be a discrete ring. Denote by $A\ldef R[T]$ the ring of polynomials. Observe that the map of rings $R_\blacksquare \to A_\blacksquare$ factors as
		$$R_\blacksquare \to (A, R)_\blacksquare \to A_\blacksquare$$
	and denote by $j\colon \Spec A_\blacksquare \to \Spec (A, R)_\blacksquare$ the virtual map.
\end{convention}

We will now translate many of the result from section \ref{section: an analogue of hilbert's basis theorem} to this more general case where $R$ is an arbitrary discrete ring and no longer $\Z$-finitely generated.

It turns out that the analogue of $A_\infty$ is not as straightforward. One has to glue.
\begin{definition}
	Denote by
	$$A_\infty\ldef \colimit_{R'\subseteq R} \laurent{R'}{T^{-1}}$$
	the subset of $\laurent{R}{T^{-1}}$ of formal Laurent series $\sum_{n\gg \infty} r_n T^{-n}$ in $T^{-1}$ all of whose coefficients $(r_n)_{n\in \Z}$ lie in some finitely generated $\Z$-algebra $R'$ of $R$.
	
	Then $A_\infty$ is a topological $A$-subalgebra of $\laurent{R}{T^{-1}}$ with the natural $\ideal{T^{-1}}$-adic topology. Indeed, $A_\infty$ is a ring, since if $x \in \laurent{R'}{T^{-1}}$ and $y \in \laurent{R''}{T^{-1}}$ then $x\cdot y\in \laurent{\Z[R'\cup R'']}{T^{-1}}\subseteq A_\infty$ where $\Z[R'\cup R']$ denotes the $\Z$-subalgebra of $R$ generated by $R'$ and $R''$. Furthermore, since each $f\in A=R[T]$ has only finitely many coefficients it is clear that $f\in A_\infty$. Hence $A\subseteq A_\infty$ and $A_\infty$ is an $A$-algebra. Even more, $A_\infty$ is then an $A$-subalgebra of the topological $A$-algebra $\laurent{R}{T^{-1}}$ and thus inherits the $\ideal{T^{-1}}$-adic topology.
\end{definition}

\begin{remark}
	If $R$ is a finitely generated $\Z$-algebra then the definition of $A_\infty$ drastically simplifies. Indeed, for every formal Laurent series over $R$, its coefficient lie in $R$, which is finitely generated. Hence, $A_\infty$ is simply $\laurent{R}{T^{-1}}$.
\end{remark}

There is a useful way to rewrite $A_\infty$.
\begin{lemma}[Alternative representation of $A_\infty$]
	\label{lemma: alternative representation of Ainfty}
	
	The morphism
		$$A_\infty = \colimit_{R'\subseteq R} \laurent{R'}{T^{-1}} \to \colimit_{R'\subseteq R} A\otimes_{R'[T]} \laurent{R'}{T^{-1}}$$
	given by mapping $s \in \laurent{R'}{T^{-1}}$ to $1\otimes s\in A\otimes_{R'[T]} \laurent{R'}{T^{-1}}$ is an $A$-module isomorphism.
\end{lemma}
\begin{proof}
	Consider the morphism given by mapping $f\otimes s$ with $f = \sum_{m=0}^M r_m T^m \in A=R[T]$ and $s=\sum_{n\gg -\infty} r_nT^{-n} \in \laurent{R'}{T^{-1}}$ to $f\cdot s \in \laurent{R}{T^{-1}}$. If its image is contained in $A_\infty$ then it is clear that multiplication is the desired inverse. Indeed, both $f$, $s$ and hence $f\cdot s$ are contained in the subalgebra
		$$\laurent{\Z[\{r_0, \dots, r_M\}\cup R']}{T^{-1}}$$
	of $A_\infty$ and the claim follows.
	
	Via a similar argument one shows that the given (clearly additive) morphism is $A$-linear. Indeed, for $f=\sum_{m=0}^M r_mT^m\in A=R[T]$ any polynomial and $s\in\laurent{R'}{T^{-1}}$ any series, the element $f\cdot s \in A_\infty$ can be considered in $\laurent{\Z[\{r_0,\dots, r_M\}\cup R'\}]}{T^{-1}}\subseteq A_\infty$ instead and then $1\otimes f\cdot s = f\cdot s = f\cdot(1\otimes s)$ as the tensor product has the right balance. Its inverse (as a map of sets) is then automatically $A$-linear as well.
\end{proof}

\begin{lemma}[Flatness]
	The $A$-module $A_\infty$ is a $A$-flat.
\end{lemma}
\begin{proof}
	
	By lemma \ref{lemma: alternative representation of Ainfty} we know that
		$$A_\infty \cong \colimit_{R'\subseteq R} R[T]\otimes_{R'[T]} \laurent{R'}{T^{-1}}$$
	as $A$-modules. Observe that each $R[T]\otimes_{R'[T]} \laurent{R'}{T^{-1}}$ is the base change of the $R'[T]$-module $\laurent{R'}{T^{-1}}$ along the inclusion $R'[T]\subseteq R[T]$. But the $R'[T]$-module $\laurent{R'}{T^{-1}}$ is flat by lemma \ref{lemma: Ainfty is A-flat, f.g. case} and then so is its base change by \cite[\href{https://math.stackexchange.com/a/120839/529214}{this argument}]{zhenlinbasechange} of \href{https://math.stackexchange.com/users/5191/zhen-lin}{Zhen Lin}. Hence, $A_\infty$ is a filtered colimit of flat $A$-modules and hence flat as a consequence of Lazard's theorem \cite[\href{https://stacks.math.columbia.edu/tag/058G}{Tag 058G}]{stacks-project}.
\end{proof}

\begin{lemma}
	\label{lemma: presentation for Afinty, general case}
	
	The topological $A$-module $A_\infty$ fits into the short exact sequence
		$$0\to \colimit_{R'\subseteq R} \series{R'}{U}[T]\xrightarrow{U\cdot T - 1}\colimit_{R'\subseteq R} \series{R'}{U}[T] \to A_\infty \to 0$$
	of $A$-modules where $\colimit_{R'\subseteq R} \series{R'}{U}[T]$ is a compact projective object of $\mod{(A, R)_\blacksquare}$. Hence, is $(A, R)_\blacksquare$-complete, $R_\blacksquare$-complete and quasi-isomorphic to the bounded complex of compact projectives of $\mod{(A, R)_\blacksquare}$ so is a derived compact in $\derived{(A, R)_\blacksquare}$.
\end{lemma}
\begin{proof}
	Observe that for any $R'$ there is a short exact sequence
	$$0\to \series{R'}{U}[T]\xrightarrow{UT-1} \series{R'}{U}[T] \to \laurent{R'}{T^{-1}}\to 0$$
	of $R'$-modules since $\laurent{R'}{T^{-1}} \cong \series{R'}{T^{-1}}[T]$. These induce an exact sequence
	$$0\to\colimit_{R'\subseteq R} \series{R'}{U}[T] \xrightarrow{UT-1} \colimit_{R'\subseteq R} \series{R'}{U}[T] \to \underbrace{\colimit_{R'\subseteq R} \laurent{R'}{T^{-1}}}_{=A_\infty}\to 0$$
	of condensed abelian groups, since filtered colimits in $\condensedAb$ are exact. As each of these groups is a condensed $A$-module and as each of the given morphisms is $A$-linear this is even an exact sequence of condensed $A$-modules. Analogous to lemma \ref{lemma: alternative representation of Ainfty} we see that
		$$\colimit_{R'\subseteq R} \series{R'}{U}[T] \cong \colimit_{R'\subseteq R} A\otimes_{R'[T]}\series{R'}{U}[T]\cong \colimit_{R'\subseteq R} A\otimes_{R'}\series{R'}{U}$$
	as $A\otimes_{R'[T]} \series{R'}{U}[T]\cong A\otimes_{R'[T]} (R'[T]\otimes_{R'}\series{R'}{U}) \cong A\otimes_{R'}\series{R'}{U}$. But $\series{R'}{U} \cong \prod_\N R'$ and hence it follows by (the proof of) lemma \ref{lemma: more compact projectives} that
		$$\colimit_{R'\subseteq R} \series{R'}{U}[T] \cong A\otimes_{R}\colimit_{R'\subseteq R} R\otimes_{R'} \prod_\N R'$$
	is a compact projective object of $\mod{(A, R)_\blacksquare}$ (as the retraction is still valid after application of $A\otimes_{R} -$). Then the claim follows at once.

\end{proof}

\begin{lemma}[$A_\infty$ is idempotent]
	
	The $A$-algebra $A_\infty$ is idempotent, by which we mean that the multiplication map of $A_\infty$ gives rise to an isomorphism $A_\infty\Dotimes_{(A, R)_\blacksquare} A_\infty \to A_\infty$ in $\derived{(A, R)_\blacksquare}$.
\end{lemma}
\begin{proof}
	Recall that by lemma \ref{lemma: formal Laurent series are idempotent} we know that multiplication $\laurent{R}{T^{-1}}\otimes_{R[T]}\laurent{R}{T^{-1}}\to \laurent{R}{T^{-1}}$ is an $A$-module isomorphism. As both multiplication $a\otimes b \mapsto a\cdot b$ and its inverse $a\mapsto 1\otimes b$ restrict to morphisms between $A_\infty\otimes_A A_\infty$ and $A_\infty$ we immediately obtain that multiplication $A_\infty\otimes_A A_\infty\isorightarrow A_\infty$ for $A_\infty$ is an isomorphism as well.
	
	It remains to show that the derived multiplication is an isomorphism. This follows verbatim as in lemma \ref{lemma: Ainfty is idempotent, f.g. case} -- the case where $R$ was $\Z$-finitely generated.
\end{proof}

We again obtain by idempotentcy that $A_\infty$-modules are well-behaved.
\begin{lemma}[$A_\infty$-modules]
	\label{lemma: Ainfty-modules are reflective}
	
	The $A_\infty$-modules in $\derived{(A, R)_\blacksquare}$ form a full subcategory whose inclusion has a left adjoint given by the functor $A_\infty\Dotimes_{(A, R)_\blacksquare} -$. Furthermore, any $M \in \derived{(A, R)_\blacksquare}$ admits at most one $A_\infty$-module structure and such a structure exists if and only if the natural morphism $M\to A_\infty\Dotimes_{(A, R)_\blacksquare} M$ is an isomorphism.
\end{lemma}
\begin{shortproof}
	The proof of lemma \ref{lemma: Ainfty-modules, f.g. case} applies verbatim.
\end{shortproof}

\begin{lemma}[$A_\infty$-modules admit not extensions by $A_\blacksquare$-complete modules]
	\label{lemma: Ainfty-modules admit not extensions by Ablack-completes}
	
	For any $S$ extremally disconnected we have
		$$\eRHom_A(A_\infty, \free{A_\blacksquare}{S}) \cong 0$$
	in $\derived{\condensedAb}$. In particular, if $C^\bullet \in \derived{A_\blacksquare}$ and $M$ is any $A_\infty$-module, then $\eRHom_A(M, C^\bullet)\cong 0$.
\end{lemma}
\begin{proof}
	Recall that
	$$A_\infty = \colimit_{R'\subseteq R} \laurent{R'}{T^{-1}}\cong \colimit_{R'\subseteq R} R[T] \otimes_{R'[T]} \laurent{R'}{T^{-1}} = \colimit_{R'\subseteq R} A \otimes_{R'[T]} \laurent{R'}{T^{-1}}.$$
	As each $\laurent{R'}{T^{-1}}$ is $R'[T]$-flat by lemma \ref{lemma: Ainfty is A-flat, f.g. case} we find further that
		$$A_\infty \cong \colimit_{R'\subseteq R} A\Dotimes_{R'[T]} \laurent{R'}{T^{-1}}.$$
	Thus,
	\begin{align*}
		&\, \eRHom_A(A_\infty, \free{A_\blacksquare}{S})\\
		\cong &\, \eRHom_A(\colimit_{R'\subseteq R} A\Dotimes_{R'[T]} \laurent{R}{T^{-1}}, \free{A_\blacksquare}{S})\\
		\cong &\, \homotopyLimit_{R'\subseteq R} \eRHom_A(A\Dotimes_{R'[T]} \laurent{R}{T^{-1}}, \free{A_\blacksquare}{S})\\
		\cong &\, \homotopyLimit_{R'\subseteq R} \eRHom_{R'[T]}(\laurent{R}{T^{-1}}, \free{A_\blacksquare}{S})\\
	\end{align*}
	as $\eRHom_A$ preserves (homotopy) colimits by \assumptionPreWord \ref{assumption: right derived Hom maps colimits to limits} and as the $\condensedAb$-enriched adjunction of scalar extension and restriction along $R'[T]\to R[T]=A$ is derivable by lemma \ref{lemma: enriched adjunctions are stable under derivation}. Fix some $R'\subseteq R$. Recall that $\free{A_\blacksquare}{S}$ is defined as a colimit over the poset of all finitely generated $\Z$-subalgebras of $A=R[T]$. Analogous to the argument in the proof of theorem \ref{theorem: analyticity of Ablack for non-finite Z-algebras A} we observe that the sub-poset given by rings $R''[T]$ with $R'\subseteq R''\subseteq R$ some intermediate ring is cofinal. Then
	$$\free{A_\blacksquare}{S} \cong \colimit_{R'\subseteq R''\subseteq R} \free{(A, R''[T])_\blacksquare}{S}.$$
	But since $R''[T]_\blacksquare$ is analytic we obtain by lemma \ref{lemma: relative completeness} that each $\free{(A, R''[T])_\blacksquare}{S}$ is $R''_\blacksquare$-complete and hence also $R'[T]_\blacksquare$-complete by restricting along the map of analytic rings $R'[T]\to R''[T]_\blacksquare$ given by the inclusion. Thus, $\free{A_\blacksquare}{S}$ is a colimit of $R'[T]_\blacksquare$-complete modules and hence $R'[T]_\blacksquare$-complete itself. We thus obtain that for any $R'$ already
	$$\eRHom_{R'[T]}(\laurent{R'}{T^{-1}}, \free{A_\blacksquare}{S}) \cong 0$$
	by lemma \ref{lemma: Ainfty-modules admit no extensions, f.g. case}. It immediately follows that
	$$\eRHom_A(A_\infty, \free{A_\blacksquare}{S}) \cong 0$$
	as well.
	
	The remaining proof now consists of extending the result to any $A_\infty$-module $M$ and any $A_\blacksquare$-complete $C^\bullet$. This is analogous to the corresponding part in the proof of the finitely generated case \ref{lemma: Ainfty-modules admit no extensions, f.g. case}.

\end{proof}

\begin{lemma}
	\label{lemma: special cokernels are Ainfty-modules, general case}
	
	For any extremally disconnected $S$, the cokernel of the morphism
		$$\free{(A, R)_\blacksquare}{S} \to \free{A_\blacksquare}{S}$$
	is an $A_\infty$-module.
\end{lemma}
\begin{proof}
	
	Denote the cokernel by $Q$ and recall that $\free{R_\blacksquare}{S} = \colimit_{R'\subseteq R} R\otimes_{R'} \free{R'_\blacksquare}{S}$. Then
	\begin{align*}
		\free{(A, R)_\blacksquare}{S} &= A\otimes_{R} \colimit_{R'\subseteq R} R\otimes_{R'}\free{R'_\blacksquare}{S}\\
		&\cong  \colimit_{R'\subseteq R} A\otimes_{R'}\free{R'_\blacksquare}{S}\cong \colimit_{R'\subseteq R} A\otimes_{R'[T]} R'[T] \otimes_{R'} \free{R'_\blacksquare}{S}.
	\end{align*}
	 Now since the poset of polynomial rings $R'[T]$ with $R'\subseteq R$ a finitely generated $\Z$-subalgebra is cofinal in the poset of finitely generated $\Z$-subalgebras $A'\subseteq A$ we obtain that $\free{A_\blacksquare}{S} \cong \colimit_{R'\subseteq R} A\otimes_{R'[T]} \free{R'[T]_\blacksquare}{S}$. We now have for every $R'$ a short exact sequence
	 	$$0\to R'[T]\otimes_{R'}\free{R'_\blacksquare}{S} \to \free{R'[T]_\blacksquare}{S} \to Q_{R'} \to 0$$
	 of $R'[T]$-modules, where the cokernel $Q_{R'}$ is a $\laurent{R'}{T^{-1}}$-module by lemma \ref{lemma: special cokernels are Ainfty-modules, f.g. case}. As tensor products are right exact we obtain that the induced sequences
	 	$$A\otimes_{R'[T]} R'[T]\otimes_{R'} \free{R'_\blacksquare}{S} \to A\otimes_{R'[T]} \free{R'[T]_\blacksquare}{S} \to A\otimes_{R'[T]} Q_{R'} \to 0$$
	 of $A$-modules are exact as well and that each $A\otimes_{R'[T]} Q_{R'}$ is still an $\laurent{R'}{T^{-1}}$-module. Since filtered colimits are exact we thus obtain that
 		$$\underbrace{\colimit_{R'\subseteq R} A\otimes_{R'[T]} R'[T]\otimes_{R'} \free{R'_\blacksquare}{S}}_{\cong \free{(A, R)_\blacksquare}{S}} \to \underbrace{\colimit_{R'\subseteq R} A\otimes_{R'[T]} \free{R'[T]_\blacksquare}{S}}_{\cong \free{A_\blacksquare}{S}} \to \colimit_{R'\subseteq R} A\otimes_{R'[T]} Q_{R'} \to 0$$
 	is exact as well and hence that the quotient must be $Q$. Now if $s\in \laurent{R'}{T^{-1}}\subseteq A_\infty$ and $a \in A\otimes_{R''[T]} Q_{R''} \subseteq Q$ then consider equivalent $s$ in $\laurent{R'''}{T^{-1}}$ and $a$ in $A\otimes_{R'''[T]}Q_{R'''}$ with $R'''=\Z[R'\cup R'']$ still $\Z$-finitely generated. Then $s$ acts on $a$ as $A\otimes_{R'''[T]}Q_{R'''}$ is an $\laurent{R'''}{T^{-1}}$-module. This makes $Q$ into an $A_\infty = \colimit_{R'\subseteq R} \laurent{R'}{T^{-1}}$-module.
\end{proof}

%

\begin{lemma}
	The kernel of the completed scalar extension $j^\ast\colon \derived{(A, R)_\blacksquare} \to \derived{A_\blacksquare}$ is the full subcategory of $\derived{(A, R)_\blacksquare}$ of $A_\infty$-modules.
\end{lemma}
\begin{proof}
	Suppose that $M^\bullet$ is an $A_\infty$-module. By lemma \ref{lemma: Ainfty-modules are reflective} we have that $M^\bullet \cong A_\infty\Dotimes_{(A, R)_\blacksquare} M^\bullet$. Since $j^\ast$ is monoidal by lemma \ref{lemma: completed scalar extension is monoidal}, we obtain that
	$$j^\ast M^\bullet \cong j^\ast(A_\infty\Dotimes_{(A, R)_\blacksquare}M^\bullet)\cong j^\ast A_\infty \Dotimes_{A_\blacksquare} j^\ast M^\bullet$$
	naturally is a $j^\ast A_\infty = A_\blacksquare \Dotimes_{(A, R)_\blacksquare} A_\infty$-module. Now for any $N^\bullet \in \derived{A_\blacksquare}$ we find
	$$\eRHom_A(A_\blacksquare\Dotimes_{(A, R)_\blacksquare} A_\infty, N^\bullet)\cong \eRHom_A(A_\infty, N^\bullet) \cong 0$$
	by lemma \ref{lemma: Ainfty-modules admit not extensions by Ablack-completes} as $N^\bullet$ is $A_\blacksquare$-complete and thus $(A, R)_\blacksquare$-complete. Thus, $j^\ast M$ is the zero-ring and any module over it must be trivial. In particular, we obtain that $j^\ast M^\bullet \cong 0$ so that $M^\bullet$ is in the kernel of $j^\ast$.
	
	Now suppose that $M^\bullet \in \derived{(A, R)_\blacksquare}$ is in the kernel of $j^\ast$, \ie $j^\ast M^\bullet \cong 0$. Now $M^\bullet$ admits a left resolution $P^\bullet \to M^\bullet$ by a $K$-projective complex $P^\bullet$ that, as a complex, is a sequential colimit $P^\bullet = \colimit_{n\ge -1} P^\bullet_n$ of bounded to the right complexes $P^\bullet_n$ which are term-wise compact projectives $P^m_n = \bigoplus_{i \in I^m_n} \free{(A, R)_\blacksquare}{S_{i, n}^m}$ of $\derived{(A, R)_\blacksquare}$. Denote by $\phi\colon (A, R)_\blacksquare \to A_\blacksquare$ the map $\id_A$ of analytic rings. Then $K$-projectivity of $P^\bullet$ gives
	$$j^\ast M^\bullet \ldef \LeftD\completedScalarExtension{\phi}{A_\blacksquare}M^\bullet \isorightarrow \completedScalarExtension{\phi}{A_\blacksquare}P^\bullet = A_\blacksquare \otimes_{(A, R)_\blacksquare} P^\bullet,$$
	implying that $\completedScalarExtension{\phi}{A_\blacksquare}P^\bullet$ is acyclic as $j^\ast M^\bullet \cong 0$. For any $n\ge -1$ consider the short exact sequence
	$$0\to P_n^\bullet \to \completedScalarExtension{\phi}{A_\blacksquare}P_n^\bullet \to Q_n^\bullet \to 0$$
	where $P_n^\bullet \to \completedScalarExtension{\phi}{A_\blacksquare} P_n^\bullet$ is given by the unit of the adjunction $\completedScalarExtension{\phi}{A_\blacksquare} \ladj \scalarRestriction{\phi}$. Observe that for any $m\in \Z$ the above short exact sequence corresponds to the short exact sequence
	$$0\to \underbrace{\bigoplus_{i\in I_n^m} \free{(A, R)_\blacksquare}{S_{i, n}^m}}_{=P_n^m} \to \underbrace{\bigoplus_{i\in I_n^m} \free{A_\blacksquare}{S_{i, n}^m}}_{\cong \completedScalarExtension{\phi}{A_\blacksquare}P_n^m} \to Q_n^m \to 0,$$
	where $\completedScalarExtension{\phi}{A_\blacksquare}P_n^m \cong \bigoplus_{i\in I_n^m} \free{A_\blacksquare}{S_{i, n}^m}$ as $\completedScalarExtension{\phi}{A_\blacksquare}$ is the unique colimit preserving extension of $\free{(A, R)_\blacksquare}{S}\mapsto \free{A_\blacksquare}{S}$ by proposition \ref{proposition: scalar restriction and extension of complete modules}. Lemma \ref{lemma: special cokernels are Ainfty-modules, general case} now implies that these cokernels $Q_n^m$ are $A_\infty$-modules, yielding that the complexes $Q_n^\bullet$ and their colimit $Q^\bullet = \colimit_{n\ge -1} Q_n^\bullet$ are $A_\infty$-modules as well. Furthermore, since filtered colimits are exact, the system of short exact sequences
	$$0\to P_n^\bullet \to \completedScalarExtension{\phi}{A_\blacksquare}P_n^\bullet \to Q_n^\bullet \to 0$$
	gives rise, in the colimit, to a short exact sequence
	$$0\to P^\bullet \to \completedScalarExtension{\phi}{A_\blacksquare}P^\bullet \to Q^\bullet \to 0$$
	since $\completedScalarExtension{\phi}{A_\blacksquare}$ (as a left adjoint) commutes with colimits of complexes. This short exact sequence then gives rise to a distinguished triangle
	$$P^\bullet \to \completedScalarExtension{\phi}{A_\blacksquare}P^\bullet \to Q^\bullet \to P^\bullet[1]$$
	in $\derived{(A, R)_\blacksquare}$. Now $\completedScalarExtension{\phi}{A_\blacksquare}P^\bullet$ is acyclic, yielding (after rotating the triangle) that $Q^\bullet[-1] \cong P^\bullet$ by \cite[\href{https://stacks.math.columbia.edu/tag/05QR}{Tag 05QR}]{stacks-project}. But $Q^\bullet$ is an $A_\infty$-module, implying that $P^\bullet$ and hence $M^\bullet$ are $A_\infty$-modules as well.
\end{proof}

We can now construct for $j\colon \Spec A_\blacksquare \to \Spec (A, R)_\blacksquare$ its associated exceptional direct image functor. 
\begin{lemma}[The left adjoint of completed scalar extension]
	\label{lemma: the left adjoint of completed scalar extension}
	
	The endofunctor
		$$T\colon \derived{(A, R)_\blacksquare}\to\derived{(A, R)_\blacksquare},\;M\to M\Dotimes_{(A, R)_\blacksquare} (A_\infty/A)[-1]$$
	given by twisting with $(A_\infty/A)[-1]$ factors uniquely over a fully faithful functor
		$$j_!\colon \derived{A_\blacksquare}\to\derived{(A, R)_\blacksquare}$$
	left adjoint to $j^\ast$ such that $j_! j^\ast M = j_! A \Dotimes_{(A, R)_\blacksquare} M$ for all $M\in\derived{(A, R)_\blacksquare}$.
\end{lemma}
\begin{proof}
	Observe that if this left adjoint $j_!$ to $j^\ast$ exists, then $j^\ast\circ j_!$ will be left adjoint to $j^\ast \circ j_\ast$ which is the identity by full faithfulness of scalar restriction $j_\ast$. Thus, $j^\ast\circ j_!$ is the identity as well and $j_!$ is fully faithful.
	
	For the existence of the adjoint $j_!$ and its representation as a tensor product it suffices to see that for all $M, N \in \derived{(A, R)_\blacksquare}$ we have
		$$\eRHom_A(TM, N)\isorightarrow \eRHom_A(TM, j_\ast j^\ast N) \isoleftarrow \eRHom_A(M, j_\ast j^\ast N)$$
	in $\derived{\condensedAb}$. Then if $M$ is $A_\blacksquare$-complete and $N$ is $(A, R)_\blacksquare$-complete we obtain by lemma \ref{lemma: eRHom enriches Hom-sets in the derived category} that
		$$\Hom_{(A, R)_\blacksquare}(T j_\ast M, N)\cong \Hom_{(A, R)_\blacksquare}(j_\ast M, j_\ast j^\ast N) = \Hom_{A_\blacksquare}(M, j^\ast N)$$
	where the equality is due to $j_\ast$ (the scalar restriction along $\id_A\colon (A, R)_\blacksquare\to A_\blacksquare$) being fully faithful. Hence $j_!\ldef T\circ j_\ast \colon \derived{A_\blacksquare}\to\derived{(A, R)_\blacksquare}$ is the desired adjoint, and we find that $j_! A \cong TA \cong (A_\infty/A)[-1]$ as $j_\ast A = A$ is the unit object for $\Dotimes_{(A, R)_\blacksquare}$. In particular, $T = TA \Dotimes_{(A, R)_\blacksquare} -$ and $j_! = j_!A \Dotimes_{(A, R)_\blacksquare} j_\ast-$. Now if $M\in \derived{(A, R)_\blacksquare}$, then $j_\ast j^\ast M \cong A\Dotimes_{(A, R)_\blacksquare} M$ which implies that
		$$j_!j^\ast M \cong j_!A\Dotimes_{(A, R)_\blacksquare} j_\ast j^\ast M \cong j_!A\Dotimes_{(A, R)_\blacksquare} (A \Dotimes_{(A, R)_\blacksquare} M) \cong j_!A\Dotimes_{(A, R)_\blacksquare} M$$
	as claimed. It remains to prove the two isomorphisms of derived $\Hom$s.
	
	For the first isomorphism observe that by definition of a triangulated category, there is a distinguished triangle
	$$N \to j_\ast j^\ast N \to C \to N[1].$$
	Application of $j^\ast = A_\blacksquare\Dotimes_{(A, R)_\blacksquare} -$ then yields another distinguished triangle
		$$j^\ast N\to j^\ast j_\ast j^\ast N \to j^\ast C \to j^\ast N[1]$$
	where $j^\ast j_\ast j^\ast N \cong j^\ast N$ by full faithfulness of $j_\ast$. By \cite[\href{https://stacks.math.columbia.edu/tag/05QR}{Tag 05QR}]{stacks-project} we obtain that $j^\ast C \cong 0$ as the first morphism is an isomorphism (as it corresponds to the identity of $j^\ast N$). By lemma \ref{lemma: Ainfty-modules are reflective} this implies that $C$ is an $A_\infty$-module. By the same lemma the inclusion of $A_\infty$-modules is fully faithful (hence $\eHom$s and $\eRHom$s of the categories agree) and has a left adjoint $A_\infty\Dotimes_{(A, R)_\blacksquare} -$. This adjunction arises from deriving the scalar extension and restriction between $(A_\infty, R)_\blacksquare$-modules and $(A, R)_\blacksquare$-modules. This adjunction is $\condensedAb$-enriched adjunction and derivable by \ref{lemma: enriched adjunctions are stable under derivation} so we obtain that the original adjunction is $\derived{\condensedAb}$-enriched and hence that
		$$\eRHom_A(TM, C) \cong \eRHom_A(A_\infty\Dotimes_{(A, R)_\blacksquare} TM, C).$$
	Consider the short exact sequence $0\to A\to A_\infty\to A_\infty/A \to 0$ and apply $A_\infty\Dotimes_{(A, R)_\blacksquare} -$ to obtain the distinguished triangle
		$$A_\infty\Dotimes_{(A, R)_\blacksquare} A \to A_\infty\Dotimes_{(A, R)_\blacksquare}A_\infty \to A_\infty\Dotimes_{(A, R)_\blacksquare} (A_\infty/A) \to (\cdots)[1].$$
	But $A$ is the unit of $\Dotimes_{(A, R)_\blacksquare}$ and $A_\infty$ is idempotent, implying that also
	$$A_\infty\xrightarrow{\id} A_\infty \to A_\infty\Dotimes_{(A, R)_\blacksquare} (A_\infty/A) \to (\cdots)[1]$$
	is distinguished. By \cite[\href{https://stacks.math.columbia.edu/tag/05QR}{Tag 05QR}]{stacks-project} we thus have $A_\infty\Dotimes_{(A, R)_\blacksquare} (A_\infty/A)\cong 0$ as the first morphism is an isomorphism. This implies that also $A_\infty\Dotimes_{(A, R)_\blacksquare} TM = A_\infty\Dotimes_{(A, R)_\blacksquare} (A_\infty/A) \Dotimes_{(A, R)_\blacksquare} M[-1] \cong 0$ and hence also $\eRHom_A(A_\infty\Dotimes_{(A, R)_\blacksquare} T M, C) \cong 0$ giving $\eRHom_A(TM, C) \cong 0$ as well. Thus application of $\eRHom_A(TM, -)$ to the initial triangle then yields that also
		$$\eRHom_A(TM, N)\to\eRHom_A(TM, j_\ast j^\ast N) \to 0 \to (\cdots)[1]$$
	is distinguished. Now by \cite[\href{https://stacks.math.columbia.edu/tag/05QR}{Tag 05QR}]{stacks-project} the remaining non-trivial morphism
		$$\eRHom_A(TM, N)\to\eRHom_A(TM, j_\ast j^\ast N)$$
	is an isomorphism, as claimed.

	For the second isomorphism consider again $0\to A\to A_\infty\to A_\infty/A \to 0$ which induces a distinguished triangle
	$$A\Dotimes_{(A, R)_\blacksquare} M\to A_\infty\Dotimes_{(A, R)_\blacksquare} M \to (A_\infty/A) \Dotimes_{(A, R)_\blacksquare} M \to (\cdots)[1].$$
	By observing that $A$ is the unit for $\Dotimes_{(A, R)_\blacksquare}$ and by rotating the triangle we obtain that
	$$(A_\infty/A) \Dotimes_{(A, R)_\blacksquare} M[-1] \to M \to A_\infty\Dotimes_{(A, R)_\blacksquare} M \to (A_\infty/A) \Dotimes_{(A, R)_\blacksquare} M$$
	is distinguished as well. This is precisely
	\begin{equation}
		\label{equation: distinguished triangle for T M}
		T M \to M \to A_\infty\Dotimes_{(A, R)_\blacksquare} M \to T M[1].
	\end{equation}
	Clearly $A_\infty\Dotimes_{(A, R)_\blacksquare} M$ is an $A_\infty$-module implying that
	$$\eRHom_A(A_\infty\Dotimes_{(A, R)_\blacksquare} M, j_\ast j^\ast N)\cong 0$$
	by lemma \ref{lemma: Ainfty-modules admit not extensions by Ablack-completes}, as the $A$-module $j_\ast j^\ast N = A_\blacksquare\Dotimes_{(A, R)_\blacksquare} N$ is $A_\blacksquare$-complete. Thus, application of $\eRHom_A(-, j_\ast j^\ast N)$ 
	to \ref{equation: distinguished triangle for T M} yields another distinguished triangle in which one term vanishes, implying that the remaining non-trivial morphism
		$$\eRHom_A(M, j_\ast j^\ast  N)\to \eRHom_A(T M, j_\ast j^\ast N)$$
	is an isomorphism, as seen by rotating the triangle and applying \cite[\href{https://stacks.math.columbia.edu/tag/05QR}{Tag 05QR}]{stacks-project}.
\end{proof}


\begin{lemma}[$j_!$ and compact generators]
	\label{lemma: j lower shriek and compact generators}
	
	Denote by $\phi\colon R\to R[T]$ the map of discrete rings. For any extremally disconnected $S$ the scalar restriction $(\scalarRestriction{\phi} \circ j_!)\free{A_\blacksquare}{S}$ is a derived compact object of $\derived{R_\blacksquare}$.
\end{lemma}
\begin{proof}
	Observe that
		$$j_!\free{A_\blacksquare}{S}\cong j_!j^\ast\free{(A, R)_\blacksquare}{S}\cong j_!A\Dotimes_{(A, R)_\blacksquare} \free{(A, R)_\blacksquare}{S}$$
	where the last isomorphism follows from lemma \ref{lemma: the left adjoint of completed scalar extension}. Now in \cite[Observation 8.12, p.57]{condenseddotpdf} Clausen and Scholze reduce as follows:
		$$\scalarRestriction{\phi}(j_!A \Dotimes_{(A, R)_\blacksquare} \free{(A, R)_\blacksquare}{S}) \cong (\scalarRestriction{\phi}j_!A) \Dotimes_{R_\blacksquare} \free{R_\blacksquare}{S}.$$
	The author was not able to reproduce this result, neither in this general case nor in case that $R$ is $\Z$-finitely generated (they claim this isomorphism for $R=\Z$). We thus have to make this an \assumptionPreWord. Now $j_!A \cong (A_\infty/A)[-1]$ and
		$$A_\infty/A \cong \colimit_{R'\subseteq R} T^{-1}\cdot \series{R'}{T^{-1}} \subseteq T^{-1}\cdot \series{R}{T}$$
	as for each $R'$ we have $\laurent{R'}{T^{-1}}/R'[T] \cong T^{-1}\cdot \series{R'}{T^{-1}}$. This then yields that
		$$A_\infty/A \cong \colimit_{R'\subseteq R} R\otimes_{R'} T^{-1}\cdot\series{R'}{T^{-1}} \cong \colimit_{R'\subseteq R} R\otimes_{R'}\prod_\N R'$$
	so that $A_\infty/A$ is a compact object of $\mod{R_\blacksquare}$ by lemma \ref{lemma: more compact projectives}. Thus, by corollary \ref{corollary: tensor products of compacts, general case} the tensor product $(A_\infty/A)\Dotimes_{R_\blacksquare}\free{R_\blacksquare}{S}$ is a compact object of $\mod{R_\blacksquare}$ and we finally find that the complex $(\scalarRestriction{\phi}j_! \free{A_\blacksquare}{S}) = (A_\infty/A)\Dotimes_{R_\blacksquare}\free{R_\blacksquare}{S}[-1]$ is a derived compact of $\derived{R_\blacksquare}$.

\end{proof}

\section{Exceptional direct images of virtual maps}

Using the most basic virtual map $\Spec R[T]_\blacksquare \to \Spec (R[T], R)_\blacksquare$ associated to the affine space $\A_R^1$ one can construct the exceptional direct image for many other virtual maps too.
\begin{proposition}
	\label{proposition: lower shriek for virtual maps}
	
	Let $R\to S\to A$ be morphisms of discrete rings. The scalar restriction of the morphism of analytic rings $\id_A\colon (A, R)_\blacksquare\to (A, S)_\blacksquare$ is a functor
	$$j_\ast\colon \derived{(A, S)_\blacksquare} \to \derived{(A, R)_\blacksquare}$$
	that admits a left adjoint $j^\ast \ldef (A, S)_\blacksquare \otimes^\LeftD_{(A, R)_\blacksquare}  -$ which in turn admits a left adjoint
	$$j_!\colon\derived{(A, S)_\blacksquare} \to \derived{(A, R)_\blacksquare}.$$
	Then $j_!$ and $j^\ast$ satisfy
	$$j_!j^\ast M = j_! A  \Dotimes_{(A, R)_\blacksquare} M$$
	for all $M\in \derived{(A, R)_\blacksquare}$.
\end{proposition}
\begin{proof}
	Since $(A, R)_\blacksquare$ and $(A, S)_\blacksquare$ are analytic by corollary \ref{corollary: analyticity of (A, R)black}, the map of analytic rings $\id_A\colon (A, R)_\blacksquare\to(A, S)_\blacksquare$ induces by proposition \ref{proposition: scalar restriction and extension of complete modules} the desired functors $j_\ast$ and $j^\ast$, the scalar restriction and (completed) scalar extension of derived modules respectively.
	
	Notice that if the statement is true for $R\to S \to S$ then statement is true for $R \to S \to A$ in general. To see this, denote by $j\colon \Spec S_\blacksquare \to \Spec (S, R)_\blacksquare$ and $k\colon \Spec (A, S)_\blacksquare \to \Spec (A, R)_\blacksquare$ the virtual maps and assume $j_!$ exists and satisfies the desired properties. Consider the following diagram
	$$\begin{tikzcd}
		\derived{S_\blacksquare} &\derived{(S, R)_\blacksquare}\\
		\derived{(A, S)_\blacksquare} &\derived{(A, R)_\blacksquare}
		\ar[from=1-1, to=1-2, shift left=0.2em, "j_!"]
		\ar[from=1-2, to=1-1, shift left=0.2em, "j^\ast"]
		\ar[from=2-1, to=2-2, shift left=0.2em, "k_!", dotted]
		\ar[from=2-2, to=2-1, shift left=0.2em, "k^\ast"]
		\ar[from=1-1, to=2-1, shift right=0.2em, "e", swap]
		\ar[from=2-1, to=1-1, shift right=0.2em, "r", swap]
		\ar[from=1-2, to=2-2, shift right=0.2em, "e'", swap]
		\ar[from=2-2, to=1-2, shift right=0.2em, "r'", swap]
	\end{tikzcd}$$
	where $e\ladj r$ and $e'\ladj r'$ are the respective (derived completed) scalar extension and restriction along $S\to A$. Now define $k_!$ to be the unique colimit preserving extension of $\free{(A, S)_\blacksquare}{T} \mapsto (e'\circ j_!)\free{S_\blacksquare}{T}$ where $T$ is extremally disconnected. Then for all such $T$ and any $(A, S)_\blacksquare$-complete $M$ we have
	\begin{align*}
		\Hom_A(k_!\free{(A, S)_\blacksquare}{T}, M) &\ldef \Hom_A((e'\circ j_!)\free{S_\blacksquare}{T}, M)\\
		&\cong \Hom_S(j_!\free{S_\blacksquare}{T}, r'(M))\\
		&\cong \Hom_S(\free{S_\blacksquare}{T}, (j^\ast \circ r')(M))\\
		&\cong \Hom_S(\free{S_\blacksquare}{T}, (r\circ k^\ast)(M))\\
		&\cong \Hom_A(e\free{S_\blacksquare}{T}, k^\ast (M)) \rdef \Hom_A(\free{(A, R)_\blacksquare}{T}, k^\ast(M))
	\end{align*}
	proving the adjunction on generators, implying that $k_! \ladj k^\ast$ in general. Now a straightforward but long calculation yields that on generators
		$$k_!k^\ast \free{(A, R)_\blacksquare}{T} \cong k_! A \Dotimes_{(A, R)_\blacksquare} \free{(A, R)_\blacksquare}{T}$$
	using that $k_!A \cong k_!\free{(A, R)_\blacksquare}{\ast}$ and that the (derived completed) scalar extension $e'$ is (strong) symmetric monoidal by lemma \ref{lemma: completed scalar extension is monoidal}. As both $k_!$ and $k^\ast$ commute with arbitrary colimits (as left adjoints) this gives the desired representation. This proves the claimed reduction.

	We can thus assume without loss of generality that $S\to A$ is surjective as this certainly covers the case of $S=A$. Since any surjection $S\to A$ is finite by \cite[\href{https://stacks.math.columbia.edu/tag/04VT}{Tag 04VT}]{stacks-project} we obtain by lemma \ref{lemma: analytic rings associated to finite morphisms} that $\id_A\colon (A, S)_\blacksquare \to A_\blacksquare$ is an isomorphism. Hence, $(A, S)_\blacksquare$ is independent of this surjection and we can assume further that $S = R[T_1, ..., T_n]$. Reducing again $R\to S\to A$ to $R\to S\to S$ we can assume both rings $S = A = R[T_1, \dots, T_n]$ to be polynomial algebras over $R$. If $n=0$ there is nothing to show. If $n\ge 2$ then the map $\id_A \colon (A, R)_\blacksquare \to (A, S)_\blacksquare$ factors as
	$$(A, R)_\blacksquare \to (A, R[T_1, \dots, T_{n-1}])_\blacksquare \to (A, S)_\blacksquare$$
	which inductively reduces the statement to the case of $n=1$. Indeed if we denote the associated virtual morphisms by $k\colon \Spec(A, R[T_1, \dots, T_{n-1}])_\blacksquare\to \Spec (A, R)_\blacksquare$ and $\ell\colon \Spec (A, S)_\blacksquare \to \Spec (A, R[T_1,\dots, T_{n-1}])_\blacksquare$ and assume $k_!$ and $\ell_!$ to satisfy the claimed properties then surely $(k\circ \ell)_!\ldef k_!\circ \ell_!$ is left adjoint to $j^\ast = (k\circ \ell)^\ast \cong \ell^\ast \circ k^\ast$ and a straightforward calculation yields that $(k\circ \ell)_!\circ (k\circ \ell)^\ast \cong (k\circ\ell)_! A \otimes_{(A, R)_\blacksquare} -$ after observing that $k_!$ is (strong) symmetric monoidal by \ref{lemma: completed scalar extension is monoidal}. Hence, the claim finally reduces to the case of $A = S = R[T]$ being a univariate polynomial algebra. But this is precisely lemma \ref{lemma: the left adjoint of completed scalar extension}.
\end{proof}

\section{The Duality}

Recall that there are several notions of \emph{dimension} of a module over a ring. Two of them will prove relevant to us.
\begin{definition}[Tor-dimension \& projective dimension] 
	Let $R$ be a discrete ring and $M$ a discrete $R$-module. Call the smallest length of a (possibly infinite) left resolution of $M$ by discrete flat $R$-modules the \emph{Tor-dimension} of $M$. Similarly, call the smallest length of a (possibly infinite) left resolution of $M$ by discrete projective $R$-modules the \emph{projective dimension} of $M$.
	
	If $f\colon \Spec A\to \Spec R$ is a morphism of affine schemes then one calls $f$ of finite Tor-dimension respectively of finite projective dimension if the $A$-module $R$ has finite Tor-dimension respectively projective dimension.
\end{definition}

As any projective resolution of $M$ is also flat, we find that the projective dimension is never greater than the Tor dimension. When do they agree?
\begin{lemma}[Agreement of dimensions, {\cite[Prop. 4.1.5]{weibel-homological-algebra}}]
	\label{lemma: agreement of dimensions}
	
	If the discrete ring $R$ is Noetherian and $M$ is a finitely generated discrete $R$-module then the Tor-dimension and the projective dimension of $M$ agree.
\end{lemma}

We obtain the following useful result.
\begin{lemma}
	\label{lemma: scalar restriction preserves compactness for finite Tor-dimension}
	
	Let $f\colon \Spec A\to\Spec R$ be a finite morphism of finite Tor-dimension and let $R$ be Noetherian. Denote by $\phi\colon R\to A$ the corresponding morphism of rings. Then scalar restriction $\scalarRestriction{\phi}\colon \derived{(A, R)_\blacksquare}\to \derived{R_\blacksquare}$ preserves derived compactness.
\end{lemma}
\begin{proof}
	As $\derived{(A, R)_\blacksquare}$ is generated by the compact objects $\free{(A, R)_\blacksquare}{S}$ for $S$ extremally disconnected we only need to show that $\scalarRestriction{\phi}$ preserves compactness on those generators. So let $S$ be extremally disconnected.
	
	Since $f$ is finite we know that $A$ is a finitely generated $R$-module. As $f$ is of finite Tor-dimension we thus know by lemma \ref{lemma: agreement of dimensions} that the $R$-module $A$ has finite projective dimension. Choose a bounded left resolution $P^\bullet$ of $A$ by discrete finitely generated projective $R$-modules. Write $P^\bullet$ as a direct summand of a bounded complex $F^\bullet$ of discrete finitely generated free $R$-modules. Then $P^\bullet \otimes_{R_\blacksquare} \free{R_\blacksquare}{S}$ is a direct summand of the complex $F^\bullet \otimes_{R_\blacksquare} \free{R_\blacksquare}{S}$. But in degree $n \in \Z$ there is some $m\in \N$ such that $F^n\cong R^m$ yielding that $F^n \otimes_{R_\blacksquare} \free{R_\blacksquare}{S} \cong \oplus_m \free{R_\blacksquare}{S} \cong \free{R_\blacksquare}{\sqcup_m S}$ and hence that $F^\bullet \otimes_{R_\blacksquare}{S}$ is term-wise compact and hence derived compact. Thus $P^\bullet \otimes_{R_\blacksquare} \free{A_\blacksquare}{S}$ is a direct summand of a derived compact complex hence derived compact itself.
\end{proof}

We now have everything necessary to state and prove the main theorem of this thesis.
\begin{theorem}[Exceptional direct and inverse image functors, the affine case]
	\label{theorem: exceptional direct and inverse image functors, the affine case}
	
	Let $f\colon \Spec A\to\Spec R$ be a finite type morphism of locally Noetherian affine schemes. Denote by $\phi\colon R\to A$ the corresponding morphism of rings.
	\begin{itemize}
		\item
		The map of analytic rings $R_\blacksquare\to (A, R)_\blacksquare \to A_\blacksquare$ with underlying map of rings $R\xrightarrow{\phi} A\xrightarrow{\id}A$ induces a functor
		$$f_!\colon \derived{A_\blacksquare} \xrightarrow{j_!} \derived{(A, R)_\blacksquare}\xrightarrow{\scalarRestriction{\phi}} \derived{R_\blacksquare}$$
		which commutes with all direct sums and satisfies the \emph{projection formula}
		$$f_!(N \otimes^\LeftD_{A_\blacksquare} (A_\blacksquare\otimes^\LeftD_{R_\blacksquare} M)) \cong f_!N\otimes^\LeftD_{R_\blacksquare} M$$
		for all $M\in D(R_\blacksquare)$ and all $N\in D(A_\blacksquare)$. If $f$ is of finite Tor-dimension, then $f_!$ preserves compact objects. Furthermore, $(g\circ f)_! \cong g_!\circ f_!$ if $g$ is another finite type morphism of locally Noetherian affine schemes.
		
		\item
		The functor $f_!$ admits a right adjoint
		$$f^!\colon D(R_\blacksquare)\to D(A_\blacksquare).$$
		The object $f^!R$ is a discrete and bounded to the left complex of finitely generated $A$-modules. If $f$ is of finite Tor dimension, then $f^!R$ is bounded, $f^!$ commutes with all direct sums and is given by
			$$f^!M =  f^!R \Dotimes_{A_\blacksquare} (A_\blacksquare \Dotimes_{R_\blacksquare} M)$$
		for $M\in\derived{R_\blacksquare}$, \ie $f^!$ is a \emph{twist} of scalar extension. If $f$ is a complete intersection, then $f^!R\in D(A)$ is an invertible object, \ie locally a line bundle concentrated in some degree. Furthermore, $(g\circ f)^! \iso f^!\circ g^!$ if $g$ is any other finite type morphism of locally Noetherian affine schemes.
	\end{itemize}
\end{theorem}
\begin{proof}
	Suppose that $g\colon \Spec R\to\Spec R'$ is any other morphism of affine schemes corresponding to the morphism of rings $\psi\colon R'\to R$. Observe that we have the following diagram
	$$\begin{tikzcd}
		\derived{A_\blacksquare}\ar[d, "{j^f_!}"]\ar[dr, "{j_!^{g\circ f}}"]&&\\
		\derived{(A, R)_\blacksquare}\ar[d, "\scalarRestriction{\phi}"]\ar[r] &\derived{(A, R')_\blacksquare}\ar[d, "\scalarRestriction{\phi}", swap]\ar[dr, "\scalarRestriction{(\phi\circ\psi)}"]&\\
		\derived{R_\blacksquare}\ar[r, "j^g_!"] &\derived{(R, R')_\blacksquare}\ar[r, "\scalarRestriction{\psi}"] &\derived{R'_\blacksquare}
	\end{tikzcd}$$
	where the left edge is $f_!$, the bottom edge is $g_!$ and the diagonal is $(g\circ f)_!$. If we can show that the embedded square commutes (up to natural isomorphism), then the whole diagram commutes (up to natural isomorphism) and we obtain that $(g\circ f)_!\cong g_!\circ f_!$. Similar to the proof of \ref{proposition: lower shriek for virtual maps} this follows since, in the words of Clausen and Scholze, \emphquote{... this is formal as the functor on top is simply carrying around an additional A-module structure.}
	
	By the extended adjoint functor theorem from \assumptionPreWord \ref{assumption: another adjoint functor theorem} the existence of $f^!$ as the right adjoint of $f_!$ follows from $f_!$ commuting with direct sums. Indeed, $f_! = \scalarRestriction{\phi} \circ j_!$ is the composition of two left adjoints and hence commutes with arbitrary direct sums. Using this characterization it is also clear that $(-)^!$ is compatible with composition as we already have shown $(-)_!$ to be. Furthermore, once we have shown that $f_!$ preserves compactness if $f$ is of finite Tor-dimension, then we already know from proposition \ref{proposition: preserves compact iff right adjoint preserves coproducts} that $f^!$ must commute with direct sums.
	
	Now that $(-)_!$ and $(-)^!$ are compatible with composition we can, analogous to the proof of proposition \ref{proposition: lower shriek for virtual maps}, reduce the remaining claims about $f_!$ and $f^!$ to the case that either $A$ is a free $R$-algebra or that $R\to A$ is surjective (which in particular covers the case of $f$ being a complete intersection). Indeed, since $A$ is of finite type over $R$ we find that $\phi\colon R\to A$ can be factorized as
		$$R \to R[X_1, \dots, X_n] \twoheadrightarrow A$$
	for a suitable $n$. Denote the resulting morphisms of schemes by $i\colon \Spec A\to \A_R^n$ and $\pi\colon \A_R^n \to \Spec R$ and abbreviate $R[\mathbf{X}] \ldef R[X_1,\dots, X_n]$. Clearly $R\to R[\mathbf{X}]$ is of finite Tor-dimension as the polynomial ring is $R$-free. By \cite[\href{https://stacks.math.columbia.edu/tag/0B66}{Tag 0B66}]{stacks-project} if $R\to A$ is of finite Tor-dimension then so is $R[\mathbf{X}]\twoheadrightarrow A$. Now assume all claims are true for $i$ and $\pi$, we will show that they are true for $f = \pi\circ i$ as well. Using that $A_\blacksquare\Dotimes_{R_\blacksquare} - \cong A_\blacksquare \Dotimes_{R[\mathbf{X}]_\blacksquare}(R[\mathbf{X}]_\blacksquare\Dotimes_{R_\blacksquare} -) $ one easily verifies that $f_!\cong\pi_!\circ i_!$ satisfies the projection formula. Now since $\pi$ is of finite Tor-dimension we have that $\pi^! R$ is a bounded complex of finitely generated discrete $R[\mathbf{X}]$-modules. By assumption \ref{assumption: derived adjunctions are enriched} we obtain that the adjunction $i_!\ladj i^!$ is $\condensedAb$-enriched so that
		$$\eRHom_{R[\mathbf{X}]}(i_!A, \pi^!R)\cong \eRHom_A(A, i^!\pi^!R).$$
	This yields that
		$$i^!\pi^!R \cong \eRHom_{R[\mathbf{X}]}(i_!A, \pi^!R)$$
	as $A$-modules since $\eRHom_{R[\mathbf{X}]}(i_!A, \pi^!R)$ naturally carries an $A$-module structure and since $\eRHom_A$ is the underlying (complex of) condensed abelian group(s) of $\intRHom_A$ giving that $\eRHom_A(A, i^!\pi^! R) \cong i^!\pi^!R$ as condensed abelian groups and hence as $A$-modules. As $i$ corresponds to the surjection of rings $R[\mathbf{X}]\to A$ we know by lemma \ref{lemma: analytic rings associated to finite morphisms} that the associated $j_!\colon A_\blacksquare\to (A, R[\mathbf{X}])_\blacksquare$ is simply the identity. Hence, $i_!$ is nothing more that scalar restriction and hence
		$$i^!\pi^!R\cong \eRHom_{R[\mathbf{X}]}(A, \pi^!R).$$
	Now $A$ is $R[\mathbf{X}]$-finitely generated so by Noetherianness of $R[\mathbf{X}]$ we can choose a left resolution of $A$ by discrete finitely generated projective $R[\mathbf{X}]$-modules. Replacing $A$ by its resolution we find that $\eRHom_{R[\mathbf{X}]}\cong \eHom_{R[\mathbf{X}]}^\bullet$ and thus that $i^!\pi^!R$ is a left bounded complex of discrete finitely generated $A$-modules. Now assume that $f$ is of finite Tor-dimension. Then so are $i$ and $\pi$ and hence $f_! = \pi_!\circ i_!$ preserves compactness as $f_!$ is the composition of two functors that do so too. Furthermore, as $i$ is of finite Tor-dimension we know that $i^!$ is a twist and hence we obtain that
		$$f^!R \cong i^!\pi^! R \cong i^!R[\mathbf{X}] \Dotimes_{A_\blacksquare} (A_\blacksquare \Dotimes_{R[\mathbf{X}]_\blacksquare} \pi^!R)$$
	where both $i^!R[\mathbf{X}]$ and $\pi^!R$ are bounded complexes of discrete finitely generated modules. In particular, the ordinary scalar extension of $\pi^! R$ is already $A_\blacksquare$-complete and so is the underlying $A$-balanced tensor product. Hence, in this case
		$$f^!R\cong i^!R[\mathbf{X}]\Dotimes_A(A\Dotimes_R \pi^!R)\cong i^!R[\mathbf{X}]\Dotimes_R \pi^!R$$
	is certainly a bounded complex of discrete finitely generated $A$-modules as well. It remains to see that $f^!$ is a twist.	Using the above presentation of $f^!R$, the isomorphism $A_\blacksquare\Dotimes_{R_\blacksquare} - \cong A_\blacksquare \Dotimes_{R[\mathbf{X}]_\blacksquare}(R[\mathbf{X}]_\blacksquare\Dotimes_{R_\blacksquare} -) $ of scalar extensions, that $A_\blacksquare\Dotimes_{R[\mathbf{X}]_\blacksquare} -$ is (strong) symmetric monoidal by lemma \ref{lemma: completed scalar extension is monoidal} and the twist representation of both $i^!$ and $\pi^!$ we find after a lengthy calculation that indeed $f^! M \cong (i^!\pi^!R)\Dotimes_{A_\blacksquare} (A_\blacksquare\Dotimes_{R_\blacksquare} M)$ for any $M\in\derived{R_\blacksquare}$. 
	
	This shows that we indeed can reduce.
	
	\textbf{The case of $A=R[X_1, \dots, X_n]$:}\\
	The case where $A=R[X_1,\dots, X_n]$ is a free $R$-algebra in $n\ge 2$ variables can further be reduced to the case where $n=1$ by replacing $R$ with $R=[X_1,\dots,X_{n-1}]$. Thus assume that $A=R[T]$ is a polynomial algebra in one variable.
	
	We will first verify the projection formula. Let $M\in \derived{R_\blacksquare}$ and $N\in\derived{A_\blacksquare}$. For proving the projection formula in this case of $A=R[T]$, Clausen and Scholze claim on \cite[p.57]{condenseddotpdf} that \quote{it is enough to prove the finer assertion}
	$$j_!(N\Dotimes_{A_\blacksquare} (A_\blacksquare\Dotimes_{R_\blacksquare} M)) \cong j_!N\Dotimes_{(A, R)_\blacksquare} j_\ast(A_\blacksquare \Dotimes_{R_\blacksquare} M).$$
	The author was not able to show the projection formula from the finer assertion. We thus have to make this an \assumptionPreWord. Denote by $C$ the cone of the natural morphism. We will show it to be acyclic, proving that the morphism is actually an isomorphism. Using that $j^\ast$ is (strong) symmetric monoidal by lemma \ref{lemma: completed scalar extension is monoidal} and that $j^\ast j_! \isorightarrow \id$ and $j^\ast j_\ast\isorightarrow \id$ by full faithfulness of $j_!$ and $j_\ast$ a short calculation yields that $j^\ast C \cong 0$. Furthermore, by lemma \ref{lemma: the left adjoint of completed scalar extension} $j_!$ is twisting by $j_! A$. But as shown in the proof of the same lemma we have $A_\infty\Dotimes_{(A, R)_\blacksquare} j_! A \cong 0$. Tensoring the defining distinguished triangle of $C$ we find that the images of the two objects in question are tensor products involving $A_\infty\Dotimes_{(A, R)_\blacksquare} j_! A \cong 0$ hence must be $0$ too. But if two out of three objects in a distinguished triangle are $0$ so must be the third and hence we find that $A_\infty\Dotimes_{(A, R)_\blacksquare} C \cong 0$ as well. Recall that $j_! j^\ast C \cong j_!A\Dotimes_{(A, R)_\blacksquare} C$ and that we have a distinguished triangle $A\to A_\infty \to A_\infty/A \to A[1]$ with $j_!A \cong (A_\infty/A)[-1]$. This gives a distinguished triangle
	$$j_! j^\ast C \to C\to A_\infty \Dotimes_{(A, R)_\blacksquare} C \to j_!j^\ast C[1]$$
	where again two of the three legs are $0$ so we find that also $C\cong 0$, \ie $C$ is acyclic as desired.
	
	Next, we will verify that $f^!R$ is a left bounded complex of discrete finitely generated $A$-modules. By \assumptionPreWord \ref{assumption: derived adjunctions are enriched} the adjunction $f_!\ladj f^!$ is $\condensedAb$-enriched which implies that
		$$\eRHom_R(f_!A, R)\cong \eRHom_A(A, f^!R)$$
	Thus $f^!R \cong \intRHom_R(f_!A, R)$ as $A$-modules as $\eRHom_R$ is the underlying (complex of) condensed abelian group(s) of $\intRHom_R$ and similarly the right-hand side is $\intRHom_A(A, f^!R) \cong f^!R$. Now $f_!A = (A_\infty/A)[-1]$ yields that
		$$f^!R \cong \intRHom_R((A_\infty/A)[-1], R) \cong \intRHom_R(A_\infty/A, R)[1].$$
	Now Clausen and Scholze state at the end of \cite[p. 57]{condenseddotpdf} that
		$$\intRHom_R(A_\infty/A, R)\cong A$$
	in case that the ring $R$ is $\Z$-finitely generated. The author was not able to verify that claim -- neither for a general $R$ nor in case of finite generation. We thus have to make it an \assumptionPreWord that $\intRHom_R(A_\infty/A, R)\cong A$. Using this one then immediately obtains that $f^!R \cong A[1]$ is a bounded complex of discrete finitely generated invertible $A$-modules.
		
	Using previous work it is easy to verify that $f_!$ preserves derived compactness. By lemma \ref{lemma: j lower shriek and compact generators} we have seen that $f_! =\scalarRestriction{\phi}\circ j_!$ preserves the derived compact generators of $\derived{A_\blacksquare}$. As $f_!$ commutes with arbitrary direct sums this is enough to see that $f_!$ preserves derived compactness in general.

	Thus, for the current case it remains to see that $f^!$ is a twist. Let $M\in\derived{R_\blacksquare}$. Using the component of the counit $f_!f^! R \to R$ and the projection formula we obtain a morphism
		$$f_!(f!R \Dotimes_{A_\blacksquare} (A_\blacksquare\Dotimes_{R_\blacksquare} M))\cong f!f^! R\Dotimes_{R_\blacksquare} M \to R\Dotimes_{R_\blacksquare} M \cong M.$$
	Under the adjunction $f_!\ladj f^!$ this morphism corresponds to a morphism
		$$f^!R\Dotimes_{A_\blacksquare}(A_\blacksquare\Dotimes_{R_\blacksquare} M) \to f^!M.$$
	If this is an isomorphism for every $M$ then we have shown that $f^!$ is indeed a twist. Now as both sides commute with arbitrary direct sums it is enough to show that this morphism is an isomorphism on generators $\free{R_\blacksquare}{S}$ for $S$ extremally disconnected. Now the left side is
		$$A[1]\Dotimes_{A_\blacksquare} \free{A_\blacksquare}{S} \cong \free{A_\blacksquare}{S}[1]$$
	as $f^!R \cong A[1]$ and $A_\blacksquare\Dotimes_{R_\blacksquare} \free{R_\blacksquare}{S}\cong \free{A_\blacksquare}{S}$. The author is unable to proceed from here. In \cite[p. 58]{condenseddotpdf} Clausen now continue by using Specker's theorem by writing $\free{R_\blacksquare}{S}\cong \prod_I R$ for some set $I$ -- this in not possible here since $R$ is not necessarily $\Z$-finitely generated. Then indeed $f^!\free{R_\blacksquare}{S} \cong \prod_I f^!R \cong \prod_I A[1] \cong \free{A_\blacksquare}{S}$ as $f^!$ was shown to commute with direct products. The author is not able to reproduce that $f^!\free{R_\blacksquare}{S} \cong \free{A_\blacksquare}{S}[1]$ in case that $R$ is not $\Z$-finitely generated. We thus have to make the \assumptionPreWord that $f^!\free{R_\blacksquare}{S} \cong \free{A_\blacksquare}{S}[1]$ for any extremally disconnected $S$ in case that $R$ is not $\Z$-finitely generated.
	
	\textbf{The case of a surjection $R\twoheadrightarrow A$:}\\
	If $R\to A$ is surjective then $\id_A\colon (A, R)_\blacksquare\to A_\blacksquare$ is an isomorphism by lemma \ref{lemma: analytic rings associated to finite morphisms}, hence $j_!$ is simply the identity and the projection formula easily follows.
	
	Now by assumption \ref{assumption: derived adjunctions are enriched} we find that the adjunction $f_!\ladj f^!$ is $\condensedAb$-enriched yielding that
		$$\eRHom_R(\scalarRestriction{\phi}A, R) \cong \eRHom_A(A, f^!R).$$
	Thus, $f^!R \cong \intRHom_R(A, R)$ as $A$-modules as $\eRHom_R$ is the underlying (complex of) condensed abelian group(s) of $\intRHom_R$ and similarly the right-hand side is $\intRHom_A(A, f^!R)\cong f^!R$. Thus, by replacing $A$ with a left resolution of discrete finitely generated projective $R$-modules we find that the underlying (complex of) $R$-module(s) of $f^!R$ is bounded to the left and consists of discrete finitely generated modules. The same is then true for the $A$-module $f^!R$ as well.
	
	Now assume that $f$ is of finite Tor-dimension. Then by \ref{lemma: agreement of dimensions} we can even choose thus resolution of $A$ to be bounded, yielding that $f^!R$ is bounded too. As scalar restriction $\scalarRestriction{\phi}$ preserves compactness by lemma \ref{lemma: scalar restriction preserves compactness for finite Tor-dimension} we find that $f_!=\scalarRestriction{\phi}\circ j_!$ preserves compactness as well since again, $j_!$ is trivial.
	
	Thus, for this case it remains to show that $f^!$ is a twist. As in the previous case where $A$ was a free $R$-algebra, the existence of the morphism
		$$f^!R\Dotimes_{A_\blacksquare}(A_\blacksquare\Dotimes_{R_\blacksquare} M) \to f^!M$$
	natural in $M\in \derived{R_\blacksquare}$ follows from the projection formula and the adjunction $f_!\ladj f^!$. The remaining argument analogous to that case as well. In particular, we also have to make the \assumptionPreWord that $f^!\free{R_\blacksquare}{S}\cong \free{A_\blacksquare}{S}[1]$ for any extremally disconnected $S$ here too. If $R$ is $\Z$-finitely generated then this is true.
	
	\textbf{The case of a complete intersection:}\\
	$f$ complete intersection:\\
	Now suppose that $R\to A$ is a even a complete intersection so that $A \cong R/\ideal{f_1,\dots, f_r}$ where $s=(f_1,\dots,f_r)$ is a regular sequence in $R$. Then by \cite[Cor. 4.5.5]{weibel-homological-algebra} the \emph{Koszul complex}
		$$K(s) = \mleft(0\to \exterior^r R^r \to \cdots \to \exterior^1 R \to \exterior^0 R^r \to 0\mright)$$
	(with $\exterior^0R^r$ in degree $0$) is an $R$-free left resolution of $A \cong R/\ideal{f_1,\dots, f_r}$. Since $R\to A$ is surjective we already have calculated that $f^R\cong \intRHom_R(A, R)$. It follows that
		$$f^!R \cong \intRHom_R(A, R) \cong \intRHom_R(K(s), R) \cong K(s)[-r]\cong A[-r]$$
	is an invertible $A$-module.
\end{proof}

\section{Affine duality rephrased geometrically}

First we want to observe that for an affine scheme the quasi-coherent modules embed well into solid modules. For this we first need that discrete modules are solid. We already have used this occasionally but for completeness we will give an argument here.
\begin{lemma}[Discrete modules are solid]
	\label{lemma: discrete modules are solid}
	
	Let $A$ be a discrete ring. If $M$ is a discrete $A$-module then the discrete condensed $\associated{A}$-module $\associated{M}$ is $A_\blacksquare$-complete. The natural inclusion
		$$\mod{A}\to\mod{A_\blacksquare}$$
	is a fully faithful and exact functor of abelian categories.
\end{lemma}
\begin{proof}
	Choose a presentation $A^J\to A^I\to M\to 0$ of $M$. Then $\associated A^J\to\associated A^I \to \associated M \to 0$ is a presentation of $\associated M$. Now $\associated A\cong A_\blacksquare[\ast]$ is $A_\blacksquare$-free, so $A_\blacksquare$-complete. Thus, both $\associated A^J$ and $\associated A^I$ are $A_\blacksquare$-complete and so is $M$ as a cokernel between $A_\blacksquare$-completes.
	
	As any discrete space is \CGWH we find that the inclusion is fully faithful by theorem \ref{theorem: some correspondences}. As limits and colimits of condensed $A$-modules are calculated locally it follows that the inclusion is exact.
\end{proof}

\begin{remark}
	Let $X=\Spec A$ be an affine scheme. In a geometric context we might call an object $\sheaf F\in \mod{A_\blacksquare}$ a \emph{solid condensed $\struct{X}$-module}. From remark \ref{remark: quasi-coherent sheaves of affine schemes} we know that we have an exact adjoint equivalence
		$$\mod{A} \xleftrightarrows{\Gamma(-, X)}{\associatedSheaf{(-)}} \qCoh{\struct{X}}$$
	asserting that the quasi-coherent $\struct{X}$-modules are precisely $A$-modules. By lemma \ref{lemma: discrete modules are solid} we thus obtain a functor
		$$\qCoh{X} \to \mod{A_\blacksquare}$$
	exhibiting quasi-coherent $\struct{X}$-modules as a full abelian subcategory of $\mod{A_\blacksquare}$ -- the category of \emph{solid condensed $\struct{X}$-modules}. Furthermore, if we denote by $\derived{\struct{X, \blacksquare}} \ldef \derived{A_\blacksquare}$ the derived category then we obtain that
		$$\derived{\qCoh{X}} \to \derived{\struct{X, \blacksquare}}$$
	\ie that the derived quasi-coherent $\struct{X}$-modules embed into $\derived{\struct{X, \blacksquare}}$.
\end{remark}

We are now close to stating the analogue of coherent duality \ref{theorem: coherent duality} in this extended setting. First, we need to investigate properness. In the affine case \emph{being proper} is a very strong condition. Any fiber of an affine morphism is again affine -- but fibers of a proper morphism should be \quote{complete}. Since an affine scheme is only complete if it is a finite union of points, one imagines that the proper morphism must already be finite.
\begin{lemma}[Proper morphisms of affine schemes are finite]
	If a morphism of affine, locally Noetherian schemes $\Spec A\to \Spec R$ is proper then it is finite, that is the $R$-module $A$ is finitely generated.
\end{lemma}
\begin{proof}
	By \cite[\href{https://stacks.math.columbia.edu/tag/02JJ}{Tag 02JJ}]{stacks-project} the morphism $R\to A$ is finite if and only if it is integral and of finite type. By definition, a proper morphism is of finite type. It remains to show that $R\to A$ is integral. By \cite[\href{https://mathoverflow.net/a/125793/546808}{this MathOverflow answer}]{proper-affine-integral} this follows as well.
\end{proof}

\begin{remark}
	The definition of $f^!$ as a right-adjoint to $f_!$ enriches by \assumptionPreWord \ref{assumption: derived adjunctions are enriched} yielding that for any $M\in D(A_\blacksquare)$ one has
		$$\eRHom_R(f_!M,R)\cong \eRHom_A(M, f^!R)$$
	inside $\derived{\condensedAb}$. As these are the underlying (complexes of) condensed abelian group(s) of $\intRHom_R$ and $\intRHom_A$ we find further that
		$$\intRHom_R(f_!M, R)\cong \RightD f_\ast\intRHom_A(M, f^!R)$$
	as $\RightD f_\ast \cong f_\ast$ as scalar restriction is exact.
	
	If $f$ is a complete intersection the complex $\omega_A\ldef f^!R$ is invertible. If $P$ is a finite projective $A$-module (or more generally a perfect complex of $A$-modules), then for $M=\RightD\intHom_{\associated{R}}(P, \omega_A)\cong P^\dual\otimes_A \omega_A$ this gives an isomorphism
	$$\RightD\intHom_{\associated{R}}(f_!(P^\dual\otimes_A \omega_A), R)\cong \RightD\intHom_{\associated{A}}(M, \omega_A) \cong P$$
	in $\derived{{\associated{R}}}$.
\end{remark}

\begin{theorem}[Solid coherent duality]
	\label{theorem: solid coherent duality}
	
	Let $f\colon X\to Y$ be a proper morphism of finite type between locally Noetherian affine schemes. Then $\RightD f_\ast\colon \derived{\struct{X,\blacksquare}} \to \derived{\struct{Y,\blacksquare}}$ agrees with the functor $f_!$ and hence admits a right adjoint $f^!\colon \derived{\struct{Y, \blacksquare}}\to \derived{\struct{X,\blacksquare}}$. Furthermore, there is an isomorphism
		$$\RightD f_\ast \intRHom_\struct{X}(\sheaf{F}, f^!\sheaf{G}) \cong \intRHom_{\struct{Y}}(\RightD f_\ast \sheaf F, \sheaf{G})$$
	of derived solid condensed $\struct{Y}$-modules natural in the derived solid condensed modules $\sheaf F\in\derived{\struct{X,\blacksquare}}$ and $\sheaf{G}\in\derived{\struct{Y,\blacksquare}}$.
\end{theorem}

\begin{remark}
	In the setting of theorem \ref{theorem: solid coherent duality} it is clear that $\RightD f_\ast = f_!$ preserves discrete objects -- that is the objects of $\qCoh{X}$. Let $X=\Spec A$ and $Y=\Spec R$. Now if $f$ is of finite Tor-dimension then theorem \ref{theorem: exceptional direct and inverse image functors, the affine case} implies that $f^!$ is a twist $f^!R \Dotimes_{A_\blacksquare}(A_\blacksquare\Dotimes_{R_\blacksquare} -)$ of (derived completed) scalar extension and that $f^!R$ is a bounded complex of discrete finitely generated $A$-modules. But this scalar extension is the unique colimit preserving extension of $\free{R_\blacksquare}{S}\mapsto \free{A_\blacksquare}{S}$ so in particular maps $R\cong\free{R_\blacksquare}{\ast}$ to $A\cong \free{A_\blacksquare}{\ast}$ and hence preserves discrete objects as well. Now finally since $\Dotimes_{A_\blacksquare}$ is the completion of the ordinary tensor product of $A$-modules we find that $f^!$ too preserves discrete objects. This implies that we have an adjunction
		$$\derived{\qCoh{Y}}\xleftrightarrows{\RightD f_\ast}{f^!}\derived{\qCoh{X}}$$
	and the corresponding isomorphism
		$$\RightD f_\ast \intRHom_\struct{X}(\sheaf{F}, f^!\sheaf{G}) \cong \intRHom_{\struct{Y}}(\RightD f_\ast \sheaf F, \sheaf{G})$$
	of $\struct{Y}$-modules natural in the sheaves $\sheaf F\in\derived{\qCoh{X}}$ and $\sheaf{G}\in\derived{\qCoh{Y}}$. As this adjunction even restricts to coherent sheaves (indeed $f^!$ is a twist by a complex of finitely generated $A$-modules so a twist by a coherent sheaf) we recover classic coherent duality for affine schemes as in theorem \ref{theorem: coherent duality} -- if $f$ is of finite Tor-dimension. This is often the case, \eg when $Y$ is regular (this includes the case where the base $Y$ is $\Spec K$ for some field) or when $f$ is a complete intersection (the Koszul complex gives at least one finite flat resolution of $A$).
\end{remark}

	\chapter{Outlook: Globalization}
	\label{chapter: globalization}
	
Let $R\to A$ be a morphism of discrete rings and $f\colon \Spec A\to \Spec R$ the corresponding morphism of schemes. Recall that the definition of $f_!$ (and hence $f^!$) crucially involved the derived category $\derived{(A, R)_\blacksquare}$. Observe, that if $R'$ is any integral extension of $R$ in $A$, then as a consequence of lemma \ref{lemma: analytic rings associated to finite morphisms} the analytic rings $(A, R)_\blacksquare$ and $(A, R')_\blacksquare$ agree. This allows us to replace $R$, in the extreme case, by its integral closure $A^+$ in $A$.

Clausen and Scholze now state in \cite[p. 62]{condenseddotpdf} that \emphquote{[t]his indicates that the globalization is naturally done in the language of adic spaces.} Adic spaces were originally introduced by Huber and are a common generalization of many geometric objects of interest, \eg schemes, formal schemes and rigid analytic varieties. They are the basic objects of interest in Huber's approach to non-archimedean rigid analytic geometry. For us, the relevant part of the theory of adic spaces is how they generalize schemes. We will now introduce the basic language of (discrete) adic spaces as done by Clausen and Scholze in \cite[Lecture XI, X \& XI]{condenseddotpdf}. For the general theory we refer to Huber's original work \cite{huber93} \& \cite{huber94}.

\begin{definition}[(Discrete) Huber pairs]
	A \emph{discrete Huber pair} is a pair $(A, A^+)$ consisting of a discrete ring $A$ together with an integrally closed subring $A^+$.
\end{definition}

\begin{definition}[Valuations]
	Let $A$ be a ring. A \emph{(multiplicative) valuation on $A$} is a function $\vert\cdot\vert\colon A\to \Gamma \cup \{0\}$ where $\Gamma$ is a totally ordered abelian group (written multiplicatively) such that for all $x, y\in A$
	\begin{itemize}
		\item $\vert 0\vert = 0$
		\item $\vert 1\vert = 1$
		\item $\vert x\cdot y\vert = \vert x\vert\cdot \vert y\vert$
		\item $\vert x + y\vert \le \max \{\vert x\vert, \vert y\vert\}$
	\end{itemize}
	with the convention that $0 < \gamma$ and $0\cdot \gamma = 0$ for all $\gamma \in \Gamma$.
	
	Two valuations $\vert\cdot\vert_1\colon A\to \Gamma\cup\{0\}$ and $\vert\cdot\vert_2\colon A\to \Gamma'\cup \{0\}$ are said to be \emph{equivalent} if for all $a,b\in A$ we have $\vert a\vert_1\ge \vert b\vert_1$ if and only if $\vert a\vert_2\ge \vert b\vert_2$.
\end{definition}

To every Huber pair one can associate a \quote{doubly} ringed space. First we introduce the topological spaces underlying such (discrete) adic spectra.
\begin{definition}[Adic spectra]
	Let $(A, A^+)$ be a discrete Huber pair. The \emph{adic spectrum} $\Spa(A, A^+)$ is a topological space with underlying set
		$$\{\vert\cdot\vert\colon A\to \Gamma\cup \{0\} \setseparator \vert A^+\vert \le 1\} / \cong$$
	of equivalence classes of valuations $\vert\cdot\vert\colon A\to \Gamma\cup \{0\}$ such that $\vert a\vert \le 1$ for all $a\in A^+$. The topology on $\Spa(A, A^+)$ has a basis of (quasi-compact) open subsets given by
		$$U\mleft(\frac{g_1,\dots, g_n}{f}\mright)\ldef \{x\setseparator \forall i = 1,\dots, n\colon \vert g_i(x)\vert \le \vert f(x)\vert \ne 0\} = \bigcap_{i=1}^n U\mleft({g_i}/{f}\mright)$$
	for $g_1,\dots, g_n, f\in A$.
\end{definition}

\begin{remark}
	By abuse of notation for a \quote{function} $f\in A$ and a \quote{point} $x\in \Spa(A, A^+)$ (a valuation) one usually writes $\vert f(x)\vert$ instead of $x(f)$. This is similar to affine spectra, where one interprets elements $f$ of a ring $A$ as \quote{functions} on the geometric object $\Spec A$ and for some \quote{point} $P\in \Spec A$ writes $f(P)$ for the image of $f$ under the quotient map $A\to A/P$. So similar to affine spectra, elements of $f\in A$ for some discrete Huber pair $(A, A^+)$ are interpreted as functions on the geometric space $\Spa(A, A^+)$, however: one can only evaluate the \quote{size} of a functions values, not the values themselves.
\end{remark}

The two \quote{structure sheaves} of an adic spectrum are now given as follows.
\begin{proposition}[{\cite[Prop. 9.3 \& Prop. 9.4]{condenseddotpdf}}]
	Let $(A, A^+)$ be a discrete Huber pair. There are sheaves $\struct{X}$ and $\struct{X}^+$ on $X=\Spa(A, A^+)$ given on the basic open subset $U = U\mleft(\frac{g_1,\cdots, g_n}{f}\mright)$ by
		$$\struct{X}(U)\ldef A\mleft[\frac{1}{f}\mright]\text{ and }\struct{X}^+(U)\ldef A\mleft[\frac{g_1}{f},\dots, \frac{g_n}{f}\mright]^+$$
	where $(-)^+$ denotes the integral closure inside $A\mleft[1/f\mright]$. Furthermore, the valuation $x\colon f\mapsto \vert f(x)\vert$ extends uniquely to $\structsymbol_{X, x}$.
\end{proposition}

We can now give the definition of a general discrete adic space.
\begin{definition}[Discrete adic spaces]
	A \emph{discrete adic space} is a triple $(X, \struct X, (\vert\cdot(x)\vert)_{x\in X})$ consisting of a topological space $X$ equipped with a sheaf of rings $\struct X$ and with an equivalence class of valuations $\vert\cdot(x)\vert$ on $\structsymbol_{X, x}$ for all $x\in X$, such that $(X, \struct X, (\vert\cdot(x)\vert)_{x\in X})$ locally is of the form $(\Spa(A, A^+), \struct{\Spa(A, A^+ )}, (\vert\cdot(x)\vert)_{x\in \Spa(A, A^+)})$ for a discrete Huber pair $(A, A^+)$.
\end{definition}

There are now two different ways to embed schemes in adic spaces.
\begin{remark}[Embedding schemes]
	Let $R$ be any discrete ring. There are two fully faithful functors from $R$-schemes to discrete adic spaces. The first is given by mapping
		$$X=\Spec A \mapsto X^\ad \ldef \Spa(A, A)$$
	while the second one is given by mapping
		$$X=\Spec A\mapsto X^{\ad/R} \ldef \Spa(A, R^+)$$
	where $R^+$ is the integral closure of $R$ inside $A$. On a general $R$-scheme they are defined by gluing. Furthermore, there is a natural comparison map
		$$X^\ad \to X^{\ad/R}$$
	and $X^{\ad/R}$ is naturally an adic space over $\Spa(R, R)$.
\end{remark}

As it turns out, this comparison map becomes an isomorphism if $f\colon X\to \Spec R$ is proper.
\begin{proposition}[{\cite[Prop. 9.6]{condenseddotpdf}}]
	\label{proposition: comparison map for proper maps}
	
	If $X\to \Spec R$ is separated and of finite type then $X^\ad \to X^{\ad/R}$ is an open immersion. If $X \to \Spec R$ is proper then $X^\ad\to X^{\ad/R}$ even is an isomorphism.
\end{proposition}

Now to each $\Spa(A, A^+)$, we can associate an abelian category $\mod{(A, A^+)_\blacksquare}$ of complete modules. Further, in case of the affine spectrum $X=\Spec A$ we already know that the category $\mod{A_\blacksquare}$ associated to $X^\ad = \Spa(A, A)$ contains enough geometric information about $X$ to for example recover duality for $X$. Hence, we would like to glue these categories $\mod{(A, A^+)_\blacksquare}$ over affine sub-schemes $\Spec A$ of some scheme $X$. Unfortunately, the transition functors $(B, B^+)_\blacksquare\otimes_{(A, A^+)_\blacksquare} -$ that arise this way are not exact and hence gluing on the level of abelian categories is not possible. Instead, gluing on the level of derived categories will work -- with the caveat that one has to pass to an $\infty$-categorical setting. Denoting by $\inftyD((A, A^+)_\blacksquare)$ the \emph{stable $\infty$-category} associated to the abelian category $\mod{(A, A^+)_\blacksquare}$, Clausen and Scholze state the following result. 
\begin{theorem}[{\cite[Thm. 9.8]{condenseddotpdf}}]
	\label{theorem: sheaves of infty-categories}
	
	Let $X$ be a discrete adic space. The assignment 
		$$U=\Spa(A, A^+)\subseteq X \mapsto \inftyD((A, A^+)_\blacksquare)$$
	defined on adic opens $\Spa(A, A^+)$ of $X$ extends to a sheaf of $\infty$-categories on $X$. Taking global sections defines an $\infty$-category $\inftyD((\struct X, \struct X^+)_\blacksquare)$ whose homotopy category is denoted $\derived{(\struct X, \struct X^+)_\blacksquare}$. 
\end{theorem}

Consequently, we can define a good derived category for any given scheme.
\begin{definition}
	For $X$ a scheme, define the triangulated category
		$$\derived{\struct {X,\blacksquare}}\ldef \derived{(\struct{X^\ad}, \struct{X^\ad}^+)_\blacksquare}$$
	to be the homotopy category from theorem \ref{theorem: sheaves of infty-categories}.
\end{definition}

As expected we have the following special case.
\begin{remark}
	 If $X=\Spec A$ is affine then $X^\ad=\Spa(A, A)$ and $\derived{\struct{X,\blacksquare}}\equiv \derived{A_\blacksquare}$.
\end{remark}

\begin{construction}[The exceptional direct image]
	The idea used in the construction of $f^!$ for a morphism of affine schemes can be adopted to the glued world. Let $f\colon X\to Y$ be a (separated) morphism of (finite type) of schemes and recall that we defined $\derived{\struct{X,\blacksquare}} = \derived{(\struct{X^\ad}, \struct{X^\ad}^+)_\blacksquare}$. In the affine case $f\colon X=\Spec A\to \Spec R = Y$ we factored the morphism $R_\blacksquare\to A_\blacksquare$ corresponding to $f$ as $R_\blacksquare \to (A, R)_\blacksquare \to A_\blacksquare$. This factorization corresponds precisely to $\Spa(R, R) \to \Spa(A, R) \to \Spa(A, A)$ and has an analogue in the global world. The comparison map $X^\ad \to X^{\ad/R}$ for schemes over the base $\Spec R$ gives rise to a similar comparison map $X^\ad \to X^{\ad/Y}$ for a general base $Y$. We obtain a functor
		$$j_\ast\colon \derived{(\struct{X^{\ad/Y}}, \struct{X^{\ad/Y}}^+)_\blacksquare} \to \derived{(\struct{X^\ad}, \struct{X^\ad}^+)_\blacksquare}\rdef \derived{\struct{X,\blacksquare}}$$
	that by \cite[Prop. 11.2]{condenseddotpdf} admits a left adjoint $j_!$. Using the forgetful functor
		$$f_\ast^{\ad/Y}\colon \derived{(\struct{X^{\ad/Y}}, \struct{X^{\ad/Y}}^+)_\blacksquare} \to \derived{(\struct{Y^\ad}, \struct{Y^\ad}^+)_\blacksquare}\rdef \derived{\struct{Y, \blacksquare}}$$
	one constructs the exceptional direct image functor $f_!\ldef f_\ast^{\ad/Y} \circ j_!$. It is then once again \quote{formal} that the right adjoint $f^!$ of $f_!$ exists. This completes the construction of the $6$-functors in general. On \cite[pp. 72--73]{condenseddotpdf} they list some of the properties making this a $6$-functor \emph{formalism} and spend the remaining chapter to verify some of them. One of the properties to note is that if $f$ is proper then by proposition \ref{proposition: comparison map for proper maps} we find that $j_\ast$ and hence $j_!$ are trivial so that $f_!$ agrees with $f_\ast$.
\end{construction}

A sample consequence of the $6$-functor formalism is the following theorem.
\begin{theorem}[{\cite[Thm. 11.1]{condenseddotpdf}}]
	
	Let $f\colon X\to \Spec R$ be a separated smooth map of finite type, of dimension $d$. Let $\canonical_{X/R}=\exterior^d \cotangent_{X/R}^1$ be the $d$-th exterior power of the cotangent sheaf of $X$ over $R$. There is a canonical functor
		$$f_!\colon \derived{\struct{X,\blacksquare}} \to \derived{R_\blacksquare}$$
	that preserves compact objects and agrees with $\RightD\Gamma(X, -)$ in case $f$ is proper. There is a natural trace map
		$$f_!\canonical_{X/R}[d] \to R$$
	such that for all $C^\bullet\in \derived{\struct{X,\blacksquare}}$ the natural map
		$$\eRHom_{\struct X}(C, \canonical_{X/R})[d] \to \eRHom_R(f_!C^\bullet, R)$$
	is an isomorphism.
\end{theorem}

	\appendix

	\chapter{Cardinals}
	\label{chapter: cardinals}
	
\section*{Introduction}
In this chapter of the appendix we will cover some of the formalities regarding set-theory and cardinal arithmetic. The reader should notice that all results covered in this chapter (in particular the construction \ref{construction: Beth-numbers and existence of strong limit cardinals}) are well-placed within $\categoryname{ZFC}$. In particular, the approach chosen by Clausen and Scholze to approximate sheaves on a large site by approximating this large site from below using \quote{cutoff-cardinals} is well-placed too and provides a good alternative to usual approaches like working with Grothendieck universes -- which usually require additional axioms to be assumed on top of $\categoryname{ZFC}$. We are lenient in the uses of \emph{classes} as these are only used to group objects for ease of comprehension, and no result stated fundamentally depends on the existence of classes as objects on their own right.

\begin{convention}
	We will use the \emph{Von Neumann cardinal assignment}. That is, we choose for each class of equicardinal sets a specific representative, namely the smallest (recall the class of all ordinals is well-ordered) ordinal in this equivalence class. Denoting by $\Ordinals$ the class of all ordinals, we thus define for any set $S$
		$$\vert S\vert \ldef \inf\{\alpha \in \Ordinals \setseparator \alpha \sim S\}$$
	where $S\sim S'$ if and only if there exists a bijection between $S$ and $S'$. We define the class of cardinal numbers to be the class of these representatives. Hence,
		$$\Cardinals \ldef \{\vert S\vert \setseparator S \in \Set \}.$$
	In particular every cardinal is an ordinal and the class $\Cardinals$ of cardinals is well-ordered itself.
\end{convention}

\begin{definition}[Cardinal arithmetic]
	For a cardinal $\kappa$ define its \emph{successor cardinal} as the cardinal
		$$\kappa^+ \ldef \inf\{\kappa' \in \Cardinals \setseparator \kappa' > \kappa\}.$$
	For a second cardinal $\kappa'$ write
		$$\kappa + \kappa' \ldef \vert \kappa \sqcup \kappa' \vert$$
	where $\sqcup$ denotes the coproduct in $\Set$, as well as
		$$\kappa \cdot \kappa' \ldef \vert \kappa \times \kappa'\vert.$$
	The cardinality of the powerset is denoted
		$$2^\kappa \ldef \vert \powerset \kappa\vert.$$
	More generally for cardinals $\mu$ and $\lambda$ define
		$$\mu^\lambda \ldef \vert \{f\colon \lambda \to \mu\}\vert$$
	to be the cardinality of the set of all functions $\lambda\to \mu$.
\end{definition}

\begin{definition}[Weak \& strong limit cardinals]
	A cardinal $\kappa$ is a \emph{weak limit cardinal} if it is neither $0$ nor a successor cardinal. If in addition for every $\kappa' < \kappa$ also $2^{\kappa'} < \kappa$ then $\kappa$ is called a \emph{strong limit cardinal}. 
\end{definition}

\begin{construction}[$\beth$-numbers \& strong limit cardinals]
	\label{construction: Beth-numbers and existence of strong limit cardinals}	
	Define recursively
		$$\beth_0 \ldef 0,$$
	for every ordinal $\alpha$
		$$\beth_{\alpha^+} \ldef 2^{\beth_\alpha}$$
	and finally for $\alpha$ a limit ordinal
		$$\beth_\alpha \ldef \bigcup_{\alpha' < \alpha} \beth_{\alpha'}.$$
	Denoting by $\omega$ the order-type of $\N$, we obtain for each ordinal $\alpha$ an uncountable large strong limit cardinal
		$$\beth_{\omega + \alpha} = \beth_{\alpha + \omega} = \bigcup_{n < \omega} \beth_{\alpha + n},$$
	in particular we can construct arbitrarily large uncountable strong limit cardinals.
\end{construction}

\begin{definition}[Cofinality for category theorists]
	Let $\kappa$ be a cardinal. The \emph{cofinality} $\cofinality\kappa$ of $\kappa$ is the smallest cardinal $\alpha$ such that $\boundedSet \kappa$ is not $\alpha$-cocomplete. Hence for every cardinal $\lambda$ we have that $\lambda < \cofinality{\kappa}$ if and only if $\boundedSet \kappa$ is $\lambda$-cocomplete.
\end{definition}

\begin{remark}[An upper bound]
	\label{remark: upper bound of cofinality}
	
	Let $\kappa$ be a cardinal. Then $\kappa$ is a colimit of a diagram of size $\kappa$ ($\kappa$ is in bijection with the disjoint union of its elements.) Thus, as $\kappa$ is not $\kappa$-small, the category $\boundedSet \kappa$ is not $\kappa$-cocomplete. Hence, $\cofinality{\kappa} \le \kappa$.
\end{remark}

\begin{lemma}
	For any infinite cardinal $\kappa$ its cofinality $\cofinality{\kappa}$ is infinite as well.
\end{lemma}
\begin{proof}
	Suppose for a contradiction that $\cofinality{\kappa} = n \in \N$ is finite. Then by definition $\boundedSet{\kappa}$ is not $n$-complete. Hence, there are $n$ cardinals $\kappa_1,\dots,\kappa_n$ smaller than $\kappa$ with
		$$\sup\{\kappa_1,\cdots,\kappa_n\} \ge \kappa.$$
	But the supremum of finitely many cardinals is just their maximum. Thus,
		$$\sup\{\kappa_1,\dots,\kappa_n\} < \kappa$$
	in contradiction to the previous inequality.
\end{proof}

\begin{definition}[Regular and singular cardinals]
	A cardinal $\kappa$ is called \emph{regular} if $\cofinality{\kappa} = \kappa$. Hence, $\kappa$ is regular if and only if the category $\boundedSet \kappa$ is $\lambda$-complete for all $\lambda < \kappa$. If $\kappa$ is not regular it is called \emph{singular}.
\end{definition}

\begin{example}[Regular and singular cardinals]
	Certainly the category $\set$ of finite sets is closed under finite colimits. Hence, $\N$ is regular. To construct a singular cardinal is more work. Construct the $\aleph$-sequence by defining recursively
		$$\aleph_0\ldef\N$$
	then
		$$\aleph_{\alpha^+} \ldef \inf\{\kappa \in \Cardinals\setseparator\kappa > \aleph_\alpha\}$$
	for every ordinal $\alpha$, recalling that cardinals are well-ordered and finally
		$$\aleph_\alpha \ldef \bigcup_{\lambda < \alpha}\aleph_\lambda$$
	for a limit ordinal $\alpha$. Then $\aleph_\omega$ is the union of the sequence $(\aleph_n)_{n\in\N}$ (which is the colimit of the diagram $n \to \aleph_n$ with $\N$ as a poset category). Hence, $\boundedSet{\aleph_\omega}$ is not even $\N$-cocomplete. In particular $\cofinality{\aleph_\omega}\le \N \ll \aleph_\omega$ and thus $\aleph_\omega$ is singular.
\end{example}


	\chapter{General Category Theory}
	\label{Section: General Category Theory}

\begin{definition}[Reflexive pairs and reflexive (co)equalizers]
	A parallel pair of morphisms $X\rightrightarrows Y$ with a common section $Y\to X$ is called a \emph{reflexive pair}. A (co)limit of a diagram with underlying shape a reflexive pair is called a \emph{reflexive (co)equalizer}.
\end{definition}

\begin{definition}[Filteredness]
	Let $\kappa$ be an infinite regular cardinal. A category $\category I$ is called \emph{$\kappa$-filtered} if any $\kappa$-small diagram in $\category I$ admits a cocone. In the case of $\kappa=\N$, we simply call $\category I$ \emph{filtered}. Dually, $\category I$ is called \emph{$\kappa$-cofiltered} or \emph{cofiltered} if $\category I^\op$ is $\kappa$-filtered or filtered respectively.
\end{definition}

\begin{remark}[Interaction of filtered colimits and limits]
	\label{remark: interaction of filtered colimits and limits}
	Let $\kappa$ be an infinite regular cardinal. Then $\kappa$-filtered colimits of sets commute with $\kappa$-small limits. Explicitly, if $\category I$ is $\kappa$-filtered, $\category J$ is $\kappa$-small and $D\colon \category I\times\category J \to \Set$ is a diagram, then the canonical morphism
		$$\colimit_\category I \limit_\category J D \to \limit_\category J \colimit_\category I D$$
	is an isomorphism. The author could not locate a reference for this claim except for an $\infty$-categorical analogue given by Lurie in \cite[Prop. 5.3.3.3]{higher-topoi}. This is however no obstruction as one can represent any ordinary category by its nerve, apply the theorem and recover the categories as the underlying homotopy categories of the nerves. This works since (homotopy) limits and colimits in the nerve agree with the limits and colimits in the ordinary category.
\end{remark}

%
%

\begin{construction}[Pro-categories]
	\label{Construction: Pro-categories}
	
	For a category $\category C$ the category $\ProCategory {\category C}$ of \emph{pro-objects} of $\category C$ essentially arises from formally adjoining all cofiltered limits missing in $\category C$.
	
	More precisely, consider the fully faithful \emph{Yoneda embedding}
	\begin{align*}
		\yoneda\colon \category C^\op &\to \functorCategory {\category C} {\Set},\\
		X&\mapsto \Hom_{\category C^\op}(-, X) = \Hom_\category C(X, -),\\
		Y\xrightarrow{f^\op} X &\mapsto (f^\op)_\ast = f^\ast \colon \Hom_\category C(Y, -)\to \Hom_\category C(X, -)
	\end{align*}
	of $\category C^\op$. Consider its opposite $\yoneda^\op\colon \category C \to \functorCategory{\category C}{\Set}^\op$. Now $\ProCategory{\category C}$ is the full subcategory of $\functorCategory{\category C}{\Set}^\op$ of cofiltered limits of representable functors, \ie
		$$\ProCategory{\category C} \ldef \left\{\limit \Hom_\category C(X, -)\setseparator X\colon \category I\to\category C\text{ cofiltered diagram}\right\}\subseteq_\text{full} \functorCategory{\category C}{\Set}^\op$$
	where $\Hom_\category C(X, -)$ is to be understood as $\yoneda^\op \circ X$.
\end{construction}


%
%

	\chapter{General Topology}
	\label{chapter: general topology}

%
%
%

\begin{definition}[Final topologies]
	\label{definition: final topologies}
	
	Let $X$ be a set and $(f_i\colon X_i \to X)_{i\in I}$ a family of functions from topological spaces $X_i$ into $X$. The finest topology on $X$ such that all maps $f_i$ are continuous is called the \emph{final topology}. Equivalently, the final topology is the collection of sets $U\subseteq X$ such that $f_i^{-1}(U)$ is open in $X_i$ for any $i\in I$. That is a subset of $X$ is open if and only if its preimage under any $f_i$ is open.
\end{definition}

%

\begin{definition}[Ultrafilters, principal ultrafilters \& the Stone-topology]
	\label{definition: ultrafilters, principal ultrafilters and the stone-topology}
	
	Let $X$ be a set. A $\emptyset \ne F \subseteq \powerset X$ is called an \emph{ultrafilter} on $X$ if
	\begin{enumerate}
		\item $\emptyset \not\in F$.
		\item If $U \in F$ and $U\subseteq V\subseteq X$ then $V \in F$.
		\item If $U, V\in F$ then $U\cap V\in F$.
		\item If $U\subseteq X$ either $U \in F$ or $X\setminus U\in F$.
	\end{enumerate}
	That is, an ultrafilter is a maximal, non-trivial filter on $X$. Write
		$$\filterkernel F \ldef \bigcap_{U\in F} U.$$
	If $F$ is an ultrafilter on $X$ and in addition
		$$\filterkernel F = \{x\}$$
	for some $x \in X$, then $F$ is called the \emph{principal ultrafilter} at $x$ and $F$ is uniquely determined by $x$. The set of ultrafilters on $X$ admits a topology generated by the sets of ultrafilters $\{F \setseparator U \in F\}$ where $U \subseteq X$. This topology is called the \emph{Stone} topology.
\end{definition}

	\bibliographystyle{plain}
	\bibliography{refs.bib}
	
\end{document}